\documentclass[11pt,english]{book}

\usepackage{tikz}
\usepackage[all]{xy}
\usepackage[T1]{fontenc}
\usepackage[latin1]{inputenc}
\usepackage{babel}

\usepackage{setspace}
\usepackage{indentfirst}

\makeatletter

\usepackage{verbatim}

\makeatother

\usepackage{enumerate}
\usepackage{amsmath,amscd}
\usepackage{amssymb}
\usepackage[all]{xy}

\usepackage{amsmath,amscd,amsthm}

           {\nolinebreak $\Box$ \end{trivlist}}

\newtheorem{prop}{Proposition}[section]
\newtheorem{lem}[prop]{Lemma}

\newtheorem{cor}[prop]{Corollary}
\newtheorem{them}[prop]{Theorem}

\newtheorem{defi}[prop]{Definition}

\newtheorem{nota}[prop]{Notation}

\theoremstyle{definition}
\newtheorem{rema}[prop]{Remark}
\newtheorem{exampl}[prop]{Example}
\newenvironment{rem}{ \begin{rema}   }{$\qed $ \end{rema}}
\newenvironment{examp}{ \begin{exampl}   }{$\qed $ \end{exampl}}

\newcommand{\Cechs}{ N[{\mathcal U}] }

\newcommand{\toto}{\rightrightarrows}
\newcommand{\BBB}{B}

\newcommand{\id}{{\mathrm{ Id}}}
\newcommand{\Aut}{{\mathrm{ Aut}}}

\newcommand{\diff}{\mathrm{ d}}
\newcommand{\g}{{\mathfrak g}}
\newcommand{\G}{{\mathbf G}}
\newcommand{\gggg}{{\mathbf g}}

\newcommand{\vvvv}{{\mathbf v}}
\newcommand{\fatg}{{\mathbf g}}

\newcommand{\e}{{\mathrm{ e}}}

\newcommand{\RP}{\}_{_{\pi}}^{\mu}}
\newcommand{\GGG}{{\mathcal G}}

\newcommand{\lra}{{\longrightarrow}}
\newcommand{\HH}{{\mathbf H}}

\newcommand{\gpdproductR}{\, \bullet_{\mathcal R} \,}
\newcommand{\gpdactkernel}{\, \bullet_{{\mathcal R},{K}} \, }
\newcommand{\gpdactband}{\, \bullet_{{\mathcal R},{Band}} \,}
\newcommand{\gpdactP}{\, \bullet_{{\mathcal R},{P}} \,}

\title{\textbf{Higher structures:}\\ gerbes\\ and\\ Nijenhuis forms}
\author{Azimi} 
\begin{document}
\pagenumbering{roman}
\maketitle
\thispagestyle{empty}
\begin{center}
{\large{\textsc{Mohammad Jawad Azimi}}}\\
\vspace{7.8cm}
{\Large \textsc{}}
\end{center}
\vspace{5.4cm}
\begin{flushright}\parbox[]{7.1cm}{Disserta\c{c}\~{a}o apresentada \`{a} Faculdade de Ci\^{e}ncias
e Tecnologia da Universidade de Coimbra, para
a obten\c{c}\~{a}o do grau de Doutor em Matem\'{a}tica.}\end{flushright}
\vspace{2.5cm}
\begin{center}
Coimbra\\2013
\end{center}

\tableofcontents

\backmatter
\newpage
\chapter{}
\section*{Resumo}
\addcontentsline{toc}{section}{Resumo}
Esta tese trata de estruturas de ordem superior, designa\c c\~ao gen\'erica para todas as cole\c c\~oes de par\^entesis ou produtos $n$-uplos que, no caso de $n=2$, se reduzem aos usuais. Exemplos destas estruturas incluem os $2$-grupos e as no\c c\~oes com eles relacionadas de fibrados principais, isto \'e, gerbes n\~ao-Abelianos, e estruturas $L_\infty$. Estes dois exemplos importantes s\~ao os objetos centrais dos dois cap\'{\i}tulos desta disserta\c c\~ao.

No primeiro cap\'{\i}tulo, apresentamos uma descri\c c\~ao geral e precisa de gerbes com valores em m\'odulos cruzados arbitr\'arios e sobre stacks diferenciais arbitr\'arios. Para esta descri\c c\~ao usamos grup\'oides de Lie, ou seja, apenas geometria diferencial cl\'assica, considerando os stacks diferenciais como sendo grup\'oides de Lie, m\'odulo equival\^encia de Morita.
Provamos que a descri\c c\~ao apresentada conduz a uma no\c c\~ao que \'e equivalente \`as  j\'a existentes, comparando a nossa constru\c c\~ao com a cohomologia n\~ao-Abeliana. Mais exatamente, introduzimos a no\c c\~ao chave de extens\~ao de grup\'oide de Lie com valores num m\'odulo cruzado, relacionamo-la com $1$-cociclos n\~ao-Abelianos de Dedecker e provamos, em seguida, que a equival\^encia de Morita se traduz em cobordos, abrindo assim o caminho para uma defini\c c\~ao geral de gerbes com valores num m\'odulo cruzado sobre um stack diferencial.

No segundo cap\'{\i}tulo, desenvolvemos a teoria de formas de Nijenhuis em \'algebras $L_\infty$. Come\c camos por apresentar uma defini\c c\~ao de par\^enteses
de Richardson-Nijenhuis para formas sim\'etricas graduadas a valores vetoriais, num espa\c co vetorial graduado. Para este par\^enteses, as estruturas $L_\infty$ s\~ao simplesmente elementos de tipo Poisson. Dada uma \'algebra $L_{\infty}$, uma forma a valores vetoriais, de grau zero, que deforma um elemento de Poisson num outro elemento de Poisson, diz-se uma forma fraca de Nijenhuis. Aqui, a deforma\c c\~ao consiste em tomar o par\^enteses da forma fraca de Nijenhuis com o elemento. As formas de Nijenhuis ${\mathcal N}$ s\~ao aquelas para as quais deformar duas vezes por ${\mathcal N}$ \'e o mesmo que deformar uma vez por uma forma ${\mathcal K}$, que \'e dita o quadrado de ${\mathcal N}$. Neste contexto, obtemos uma hierarquia infinita de \'algebras $L_\infty$.

De entre os exemplos de deforma\c c\~oes de Nijenhuis, contam-se a aplica\c c\~ao de Euler numa \'algebra $L_\infty$ arbitr\'aria, bem como os elementos de Poisson e de Maurer Cartan numa \'algebra de Lie diferencial graduada.

Efetuamos a classifica\c c\~ao das formas de Nijenhuis em $2$-\'algebras de Lie com \^ancora nula. Mostramos tamb\'em que, sobre certas condi\c c\~oes, existe uma correspond\^encia biun\'{\i}voca entre as formas de Nijenhuis a valores vetoriais, na $2$-\'algebra de Lie associada a um algebr\'oide de Courant, e as aplica\c c\~oes de Nijenhuis $\mathcal{C}^{\infty}$-lineares no mesmo algebr\'oide de Courant. Apresentamos exemplos de formas de Nijenhuis a valores vetoriais nas $n$-\'algebras de Lie associadas a variedades $n$-pl\'eticas. Explicamos tamb\'em como tensores de Nijenhuis num algebr\'oide de Lie podem ser vistos como formas de Nijenhuis numa certa \'algebra de Gerstenhaber, considerada como \'algebra $L_{\infty}$. Al\'em disso, para esta \'ultima estrutura de \'algebra $L_{\infty}$, estruturas $\Omega N$ e estruturas de Poisson-Nijenhuis podem tamb\'em ser vistas como formas de Nijenhuis.


\newpage
\section*{Abstract}
\addcontentsline{toc}{section}{Abstract}

The thesis is devoted to higher structures, which is a generic name for all those collections of $n$-ary brackets or products
reducing for $n=2$ to the ordinary ones. Among examples of those are $2$-groups, and their related notions of principal bundles,
i.e. non-Abelian gerbes, and $L_\infty$-structures.
These two major examples are the central objects of the two chapters of the present work.

In the first chapter, we give a precise and general description of gerbes valued in  arbitrary crossed module and over an arbitrary differential stack.
We do it using only Lie groupoids, hence ordinary differential geometry, by considering differential stacks as being Lie groupoids up to Morita equivalence.
We prove the coincidence with the existing notions by comparing our construction with non-Abelian cohomology. More precisely, we introduce the key notion of extension of Lie groupoids valued in a crossed-module.
 We relate it with Dedecker's non-Abelian $1$-cocycles, and we then show that Morita equivalence amounts to co-boundaries, paving the way for a general definition
 of gerbes valued in a crossed-module over a differential stack.

In the second chapter, we develop the theory of Nijenhuis forms on $L_\infty$-algebras.
First, we recall a convenient notion of Richardson-Nijenhuis bracket on the graded symmetric vector valued
forms on a graded vector space, bracket for which $L_{\infty}$-algebras are simply Poisson elements. Weak Nijenhuis
 vector valued forms for a given $L_{\infty}$-algebra are defined to be forms of degree $0$
 deforming (i.e. taking bracket) that Poisson element into an other Poisson element.
 Nijenhuis forms are those forms ${\mathcal N}$ for which deforming twice by ${\mathcal N}$ is like deforming once by a form ${\mathcal K}$ called the square of ${\mathcal N }$.
 We obtain in this context an infinite hierarchy of $L_{\infty}$-algebras.

 Among examples of such Nijenhuis deformations are the Euler map on an arbitrary $L_{\infty}$-algebra or Poisson
and Maurer Cartan elements on a differential graded Lie algebra.
A classification of Nijenhuis forms on  anchor-free Lie $2$-algebras can be completed. We also show that there is, under adequate conditions, a one to one correspondence between the Nijenhuis vector valued forms
$\mathcal{N}$ with respect to the Lie $2$-algebra associated to a Courant algebroid and Nijenhuis $\mathcal{C}^{\infty}$-linear
maps on the Courant algebroid itself. We give examples of Nijenhuis vector valued forms on the Lie $n$-algebras associated
 to $n$-plectic manifolds. We also explain how Nijenhuis tensors on a Lie algebroid
 are indeed Nijenhuis forms of some Gerstenhaber algebra, considered as an $L_\infty$-algebra. For the latter $L_{\infty}$-algebra structure, moreover, $\Omega N$-structures and
 Poisson-Nijenhuis structrures can also be seen as Nijenhuis forms. 
\backmatter

\chapter{Introduction}

\section*{Generalities on higher structures}

There is no precise mathematical definition of what a "higher structure" is. The name appears in the literature as soon as binary operations, e.g. Lie algebra brackets or Lie group products, are being replaced by collections of $n$-ary operations, an $n$-ary operation on a set $S$, for us, being simply a map from $S \times \cdots \times S$ ($n$ times) to $S$. To deserve the name "higher Lie algebras" or "higher Lie groups", these operations are assumed, in general, to satisfy some quadratic relations that reduce to Jacobi identity or associativity when  the only non-trivial operations are binary.
In that manner, $L_\infty$-algebras are being defined.
As particular cases, when one only has $1$-ary and $2$-ary operations, we recover differential graded
Lie algebras and crossed-modules of groups.
When all the $k$-ary operations for $k \geq n$ are trivial,
we obtain, among others, Lie $n$-algebras and Lie $n$-groups, $n \in \mathbb{N}$.
The name higher structure is also used for gerbes, Abelian or not, because they can be seen as being
higher-principal bundles, i.e. principal bundles on higher groups, Lie $2$-groups, for our purpose.

Higher structures are objects of growing importance. For instance, $L_\infty$-algebras,
once called strongly homotopy Lie algebras \cite{Lada-Sta}, introduced first
by string theorists \cite{Zwiebach}, gained notoriety when Kontsevitch
used $L_\infty$-morphisms to prove the existence of star-products on Poisson manifolds \cite{Kont}.
In fact, the approach by Kontsevitch avoids $L_\infty$-algebras, but Voronov \cite{Voronov} using a very general construction, associates an $L_\infty$-algebra to a Poisson element and an Abelian sub-algebra of a differential graded Lie algebra.
Several authors have shown that an $L_\infty$-algebra encodes a Poisson structure in a
neighborhood of a coisotropic submanifold, provided that a linear transversal is given, see \cite{Cattaneo-Schaaetz} and \cite{Cattaneo-Felder}.
This makes $L_\infty$-algebras a central tool for studying Poisson brackets, but there are more relations.
 Roytenberg and  Weinstein \cite{RoytenbergWeinstein} have given a description of the so-called Courant algebroids,
i.e. the ambient space on which Dirac structures \cite{Courant} live, in terms of Lie $2$-algebras.
 In  the same vein, Rogers \cite{CRoger,CRogerspaper} encodes $n$-plectic manifolds by Lie $n$-algebras
and Fr\'{e}iger, Roger and Zambon \cite{FRZ} used this formalism to construct moment maps of those.

In the meantime, there has been a continuous interest into Abelian and non-Abelian gerbes.
Given a central extension of Lie groups $\tilde{K} \to K$, Abelian gerbes first appeared as being the obstruction of a $K$-principal bundle
to  come from a $\tilde{K}$-principal bundle \cite{Brylinski}, \cite{Grothendieck1} and \cite{Grothendieck2}.
Abelian gerbes where later on studied for themselves, for instance by Giraud \cite{Giraud}  Chatterjee \cite{Chatterjee}.
Again, the origin of this attention comes partly from physics \cite{GR}, \cite{Hitchin} and \cite{Witten}.
But there are also beautiful interpretations of those in terms of category, see for instance Baez and Schreiber \cite{BaezSchreiber},
or Breen and Messing \cite{BM}. With many variations (see \cite{NikolausWaldorf}), all these interpretations consist
in seeing gerbes as being principal bundles but for a Lie $2$-group (Lie $2$-groups being more or less the same thing as crossed-modules, see \cite{GinotStienon})
instead of for a Lie group. Since principal bundles admit connections and these connections admit curvature, there has been a great interest
for defining those on gerbes. For the Abelian case, Brylinski \cite{Brylinski}, Murray \cite{Murray}, Chatterjee \cite{Chatterjee} and Hitchin \cite{Hitchin} give approaches on the matter.
 The non-Abelian case was settled by Breen and Messing in the pioneering paper \cite{BM}.
 Later on, a more physicist point of view was adopted \cite{AschieriCantiniJurco},
  and a Lie groupoid interpretation of these connections and curvatures were given in \cite{LaurentStienonXu}.
  Also, categorical  interpretations were given \cite{SchreiberWaldorf}.


\subsection*{A first glance at $L_\infty$-algebras}\label{Int:1stglanceatLinfty}

As previously said, $L_\infty$-algebras are collections of $n$-ary operations, assumed to satisfy some quadratic relations that reduce to Jacobi identity, when only the binary operation is not trivial. Let us explain briefly how such objects could appear.
Let $E$ be a vector space of finite dimension. It is well-known that there is a one to one
correspondence between
\begin{enumerate}
\item[(i)] Lie algebra brackets on $E$,
\item[(ii)] derivations of degree $+1$ squaring to zero of $\wedge^\bullet E^*$
(with the understanding that elements of $\wedge^k E^*$ are of degree $k$, for all $k \geq 0$).
\end{enumerate}
The correspondence consists in associating to a Lie bracket on $E$ its Chevalley-Eilenberg differential,
i.e. the unique derivation of $\wedge^\bullet E^*$ whose restriction to $\wedge^1 E^* \simeq E^*$
is  the dual of the bilinear bracket $\wedge^2 E \to E$.
This is clearly injective, and the surjectivity comes from the fact that every derivation of degree $+1$
of $\wedge^\bullet E^*$ restricts to a map $E^* \to \wedge^2 E^*$, whose dual map is a Lie bracket
if and only if the derivation squares to zero.

Now, let $E= \oplus_{i\in \mathbb{Z}} E_i$ be a graded vector space whose components $E_i$ are of finite dimension for all $i \in \mathbb{Z}$. The morphism above can be generalized to provide a map from graded Lie algebra brackets on $E$ to derivations of degree $+1$ squaring to zero of $\wedge^\bullet E^*$
by precisely the same construction. But this map is not a one to one correspondence anymore,
because derivations of degree $+1$ of $\wedge^\bullet E^*$ do not need to restrict to maps
$\wedge^1 E^*\simeq  E^* \to \wedge^2 E^*$ anymore (here, the degree of
 $ \alpha_1 \wedge \cdots \wedge\alpha_k \in \wedge^k E^*$
 is defined to be $a_1 + \dots + a_k  + k$ for all  $ \alpha_1 \in E_{a_1}^*,  \cdots , \alpha_k \in E_{a_k}^* $).
 Hence, when we are given a graded vector space, graded Lie algebra brackets form a (in general strict) subset
 of the set of all derivations of degree $+1$ squaring to zero of $\wedge^\bullet E^*$.
 To obtain a one to one correspondence one has to consider a collection of $n$-ary brackets:
  $$ l_n : \wedge^n E \to E $$
  and to associate to it the unique derivation of $\wedge^\bullet E^*$
  whose restriction to $E^*$ is the map from $E^* $ to $\oplus_{n\geq 0} \wedge^n E^*$
  mapping $\alpha$ to $D(\alpha)=\sum_n l_n^* (\alpha)$.
  By counting the degree, one checks that $l_n$ has to map $  E_{a_1} \times \dots \times E_{a_n}$ to $E_{a_1+\dots +a_n+n-2} $
  so ensure the degree of $D$ to be $+1$. Such a collection of brackets shall be called an $L_\infty$-structure when the corresponding derivation $D$ squares to zero, which amounts to the relation:
   \begin{equation}\label{L_infty-def-skew}
\sum_{i+j=n+1}\sum_{\sigma \in Sh(i,j-1)}(-1)^{i(j-1)}\chi(\sigma)l_j(l_i(X_{\sigma(1)},\cdots,X_{\sigma(i)}),\cdots,X_{\sigma(n)})=0,
\end{equation}
for all $n\geq 1$ and all $X_1,\cdots,X_n \in E$. Here, $Sh(i,j-1)$ stands for the set of $(i,j-1)$-unshuffles
 and $\chi(\sigma)=\epsilon(\sigma).sign(\sigma)$, with $sign(\sigma)$ being the sign of the permutation while $\epsilon(\sigma)$ is the Koszul sign. Given a permutation $\sigma$ in the group of permutations of $n$ elements, the Koszul sign $\epsilon(\sigma)$ is defined as follows
   \begin{equation}
   X_{\sigma(1)}\otimes\cdots \otimes X_{\sigma(n)}=\epsilon(\sigma)X_{1}\otimes\cdots \otimes X_{n},
   \end{equation}
for all $X_{1}, \cdots ,X_{n} \in E,$
 see \cite{Lada-Sta}. The family of skew-symmetric vector valued forms $(l_{i})_{i\geq 1}$ is called an $L_{\infty}$-structure on the graded vector space $E$.

 In this thesis, it shall be more convenient to work with symmetric brackets on $E[1]$, with $E[1]$ being the graded vector space whose component of degree
 $i-1$ is formed by elements of $E_i$, so that the previous signs do not match our future conventions,
  but are equivalent to them.
  This may come a surprise, since skew-symmetric and symmetric objects are in general considered as being totally different in nature.
  Symplectic geometry, for instance, does not follow the same pattern as Riemanian geometry.
  However, in the world of graded mathematics, they are in some sense not so different.
  For instance, a graded Lie algebra structure is in general defined to be a graded skew-symmetric bilinear assignment
     $ [.,.] : E \times E \to E$, satisfying the graded Jacobi identity.
  Here, by graded skew-symmetric bilinear assignment, we mean
   $$ [ X,Y ]=-(-1)^{|X||Y|}[Y,X],$$
   where $|X|$ denotes the degree of $X\in E$.
  But we can also see this bracket as a graded symmetric map. For this purpose one has to consider $E[1]$, and to define
   \begin{equation} \label{eq:decalageFor_n=2}[X,Y]' =(-1)^{|X|} [X,Y].\end{equation}
  A direct computation now gives
   $$ [X,Y]' = (-1)^{(|X|+1)(|Y|+1)}[Y,X]',$$
   that is to say, $ [.,.]'$ is a graded symmetric bilinear assignment on $E[1]$.

   This discussion, indeed, can be extended for $L_\infty$-algebras, and it can be shown that there is a one to one correspondence between
   $L_\infty$-algebras structure as defined above and $L_\infty $-algebras structures as defined next.
   \begin{defi}\label{def:SymLinfty}
An $L_{\infty}$-algebra is a graded vector space $E$ together with a family of symmetric vector valued forms (symmetric multi-linear maps) $(l_i)_{i\geq 1}$ of degree $+1$, with $l_i: \otimes^i E\to E$ satisfying the following relation:
\begin{equation}\label{L_infty-def}
\sum_{i+j=n+1}\sum_{\sigma \in Sh(i,j-1)}\epsilon(\sigma)l_j(l_i(X_{\sigma(1)},\cdots,X_{\sigma(i)}),\cdots,X_{\sigma(n)})=0
\end{equation}
for all $n\geq 1$ and all homogeneous  $X_1,\cdots,X_n \in E$, where $\epsilon(\sigma)$ is the Koszul sign. The family of symmetric vector valued forms $(l_{i})_{i\geq 1}$ (or sometimes the vector valued form $\mu:=\sum_{i\geq 1} l_i$) is called an $L_{\infty}$-structure on the graded vector space $E$.
\end{defi}

   The correspondence between both definitions of $L_\infty$-structure is through the so-called d\'ecalage map, which
   is a map from the space of  multi-linear maps of degree $2-k$ on a graded vector space $E$ to the space of multi-linear maps of degree $+1$
   on $E[1]$ given by
   \begin{equation}\label{def:L_inftysymmetric}
   l_n(X_1,\cdots,X_n)\mapsto (-1)^{(n-1)|X_1|+(n-2)|X_2|+\cdots+|X_{n-1}|}l'_n(X_1,\cdots,X_n).
   \end{equation}
  A direct computation shows that $l_n$ is graded skew-symmetric on $E$ if and only if $l'_n$ is graded symmetric on $E[1]$.
   Also, when $n=2$, we recover (\ref{eq:decalageFor_n=2}). Last, having an $L_{\infty}$-algebra in terms of graded skew-symmetric brackets $(l_i)_{i\geq 1}$ and applying the d\'ecalage map to the brackets we get an $L_{\infty}$-algebra in terms of graded  symmetric brackets $(l'_i)_{i\geq 1}$, that is, equation (\ref{L_infty-def-skew}) holds for the graded skew-symmetric brackets $(l_i)_{i\geq 1}$ if and only if equation (\ref{L_infty-def}) holds for graded symmetric brackets $(l'_i)_{i\geq 1}$.

\subsection*{A first glance at gerbes}\label{Int:1stglancegerbs}

Groups and groupoids, like their counterpart Lie algebras and Lie algebroids, are binary products.
In a groupoid for instance, the product of two elements may be defined or not, but, when it is defined,
it is a binary product.
  There exists a notion of $n$-groups with $n$-ary products.
  The exact definition of those is beyond the scope of this thesis, but we can give an intuition
  on those by describing some examples.

  Let $S$ be  a set.
  By an \emph{oriented triangle} on $S$, we mean an oriented triplet of elements in $S$,
  i.e. an element in ${\bf T}_S := S^3/\sim $ where $ \sim $ identifies elements that differ by cyclic permutation:
   $$ (i,j,k)\sim (k,i,j) \sim (j,k,i) .$$
  We invite the reader to see them as being oriented triangles.
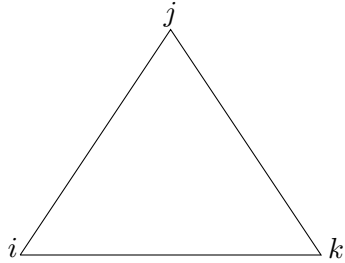
\begin{center}
   \begin{tikzpicture}
\draw (-2,0) -- (2,0) -- (0,3)-- (-2,0);
\node at (-2.1,0.1){$i$};
\node at (2.2,0.1){$k$};
\node at (0,3.2){$j$};
\end{tikzpicture}
\end{center}
%
%
We define a $3$-ary product $(.,.,.)_3$ on ${\bf T}_3$ as follows. It is not always defined: in fact it is only defined when
  the three oriented triangles are the faces of a tetrahedron, and have compatible orientations.
  The product then consists in applying to those the fourth face of the tetrahedron, equipped with a compatible orientation. In equation
   $$((ijk)(ikl)(ilj))_3 \mapsto (kjl)$$
   which can be geometrically understood as follows:
     
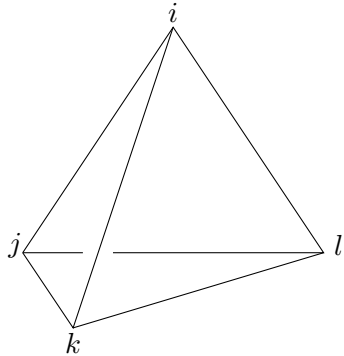
\begin{center}
   \begin{tikzpicture}
\draw (-2,0) -- (-1.2,0);
\draw (-0.8,0)-- (2,0) -- (0,3)-- (-2,0);
\draw (-2,0)--(-1.33,-1)--(0,3);
\draw (-1.33,-1)--(2,0);
\node at (-2.1,0.1){$j$};
\node at (2.2,0.1){$l$};
\node at (0,3.2){$i$};
\node at (-1.33,-1.2){$k$};
\end{tikzpicture}
\end{center}
 Another, and more complicated $3$-ary operation can be constructed out of a crossed-module
   $\G \stackrel{\rho}{\to} \HH$.
   By an \emph{oriented triangle of a crossed-module $\G \stackrel{\rho}{\to} \HH$}, we mean a representative of the quotient of the
      $ \{ (h_1,h_2,h_3,g) \in \HH^3 \times \G \hbox{ s.t. } h_1h_2h_3 =\rho (g) \}  $
      under the action of the cyclic group with three elements acting by:
        $$(h_1,h_2,h_3,g) \sim (h_2,h_3,h_1,h_1(g))  \sim (h_3,h_1,h_2,h_1(h_2(g)).$$
   We call $T_{\G\to \HH}$ this quotient set.  Again, we invite the reader to see those elements as oriented triangles:
      
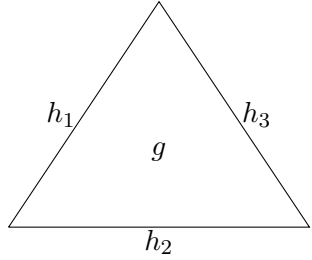
\begin{center}
   \begin{tikzpicture}
\draw (-2,0) -- (2,0) -- (0,3)-- (-2,0);
\node at (0,-0.2){$h_2$};
\node at (1.3,1.5){$h_3$};
\node at (-1.3,1.5){$h_1$};
\node at (0,1){$g$};
\end{tikzpicture}
\end{center}

     We define a $3$-ary product $(.,.,.)_3$ on ${\bf T}_3$ as follows. It is not always defined: in fact it is only defined when
  the three oriented triangle are the faces of a tetrahedron, and have compatible orientations.
  
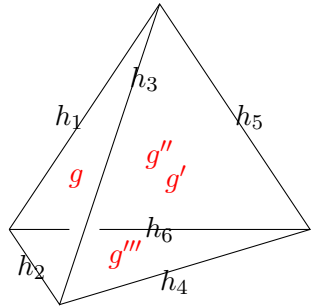
\begin{center}
   \begin{tikzpicture}
\draw (-2,0) -- (-1.2,0);
\draw (-0.8,0)-- (2,0) -- (0,3)-- (-2,0);
\draw (-2,0)--(-1.33,-1)--(0,3);
\draw (-1.33,-1)--(2,0);
\node at (0.22,0.67){${\color{red}g'}$};
\node at (-1.11,0.68){${\color{red}g}$};
\node at (0,1){${\color{red}g''}$};
\node at (-0.45,-0.33){${\color{red}g'''}$};
\node at (-1.7,-0.5){$h_2$};
\node at (0.2,-0.7){$h_4$};
\node at (-1.2,1.5){$h_1$};
\node at (-0.2,2){$h_3$};
\node at (1.2,1.5){$h_5$};
\node at (0,0){$h_6$};
\end{tikzpicture}
\end{center}
  Explicitly, the product is given by:
   $$((h_1,h_2,h_3,g)(h_3^{-1},h_4,h_5,g')(h_5^{-1},h_6,h_1^{-1},g'')_3 \mapsto (h_2,h_4,h_6,g'''),$$
   where $g'''  $ is the unique element which can naturally construct so that
    $$h_2h_4h_6=\rho(g''').$$
 Let us find this element. If we know that  $$h_1h_2h_3=\rho(g) \hbox{ and } h_3^{-1} h_4 h_5=\rho(g') \hbox{ and } h_5^{-1}h_6h_1^{-1}=\rho(g'') .$$
then we have:
   \begin{equation*}
   \begin{array}{rcl}
    \rho(h_1(g'''))&=& h_1(h_2h_4h_6)h_1^{-1}\\
                  &=&h_1h_2h_3h_3^{-1}h_4h_5h_5{-1}h_6)h_1^{-1}\\
                  &=&\rho(g)\rho(g')\rho(g'')\\
                  &=&\rho(gg'g'')
   \end{array}
   \end{equation*}
  so that a natural choice for $g'''$ is  $  g'''=h_1^{-1} (gg'g'')$, or, equivalently:
   \begin{equation}\label{eq:product3ary}
        g^{-1}  h_1(g''')= g'g''.
   \end{equation}

Non-Abelian $1$-cocycle, associated to a crossed-module, which are defined below, can be interpreted as being morphisms
for these $3$-ary products.
Let $(U_i)_{i \in I}$ be an open cover on a manifold $M$.
To every points $m\in M$, one can consider all the oriented triangles $ T_{S_m}$ with $S_m$ being the subset of $I$ made of al indices such that
$m $ belongs to $U_i$. In some sense, this is a bundle  of set equipped with $3$-ary product over $M$.

Assume that we are given, for all $i,j \in I$
maps $h_{ij} : U_i \cap U_j \to H$
and maps $ g_{ijk} : U_{i} \cap U_{j} \cap U_{k} \to G$.
Let us see what it takes to impose that the map defined by
  $$ (ijk) \to (h_{ij},h_{jk},h_{kl},g_{ijk}) $$
be a morphism for the $3$-ary bracket.
 $$ \Phi_m : T_{S_m} \to  T_{G \to H} 
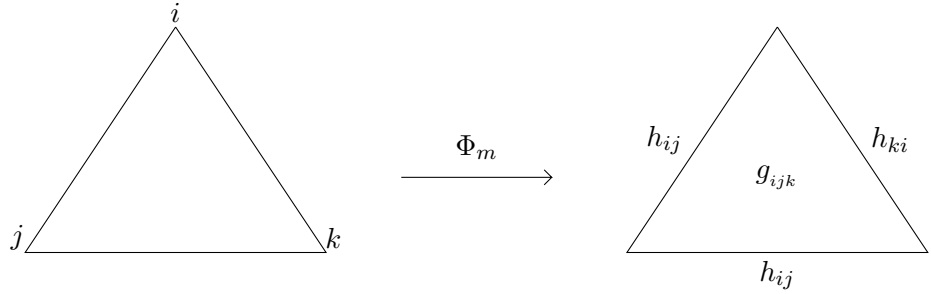
$$
 \begin{center}
  \begin{tikzpicture}
  \draw (-4,3)--(-6,0)--(-2,0)--(-4,3);
  \node at (-4,3.2){$i$};
  \node at (-6.1,0.2){$j$};
  \node at (-1.9,0.2){$k$};
  \draw (-1,1)--(1,1)--(0.9,1.1);
  \draw (1,1)--(0.9,0.9);
  \draw (4,3)--(6,0)--(2,0)--(4,3);
  \node at (2.5,1.5){$h_{ij}$};
  \node at (5.5,1.5){$h_{ki}$};
  \node at (4,-0.3){$h_{ij}$};
  \node at (4,1){$g_{_{ijk}}$};
  \node at (0,1.4){$\Phi_m$};
  \end{tikzpicture}
 \end{center}

First, it has to takes values in $T_{G \to H}$
which amounts to
 \begin{equation}\label{crossedmoduleequation1}
  h_{ij}h_{jk}h_{ki} = \rho(g_{ijk}) .
  \end{equation}
But, second, we also need to have, for every tetrahedron in $T_{S_m} $, i.e. for any quadruple of indices $i,j,kl$ such that $m \in U_{ijkl}$,
the following relation;
 $$  (\Phi_m(ijk),\Phi_m(ikl),\Phi_m(ilj))_3 = \Phi_m (kjl) .$$
 
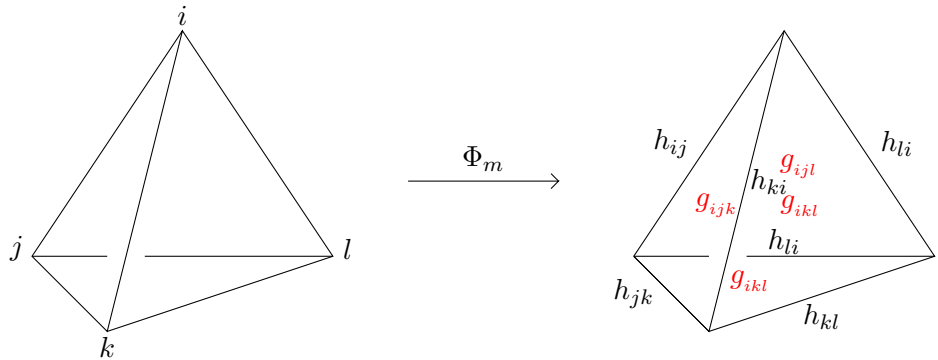
\begin{center}
 \begin{tikzpicture}
\draw (2,0) -- (3,0);
\draw (3.5,0)-- (6,0) -- (4,3)-- (2,0);
\draw (2,0)--(3,-1)--(4,3);
\draw (6,0)--(3,-1)--(2,0);

\node at (2.5,1.5){$h_{ij}$};
\node at (5.5,1.5){$h_{li}$};
\node at (3.8,1){$h_{ki}$};
\node at (2,-0.5){$h_{jk}$};
\node at (4.5,-0.8){$h_{kl}$};
\node at (4,0.2){$h_{li}$};

\node at (4.22,0.67){${\color{red}g_{_{ikl}}}$};
\node at (3.11,0.68){${\color{red}g_{_{ijk}}}$};
\node at (4.2,1.2){${\color{red}g_{_{ijl}}}$};
\node at (3.55,-0.33){${\color{red}g_{_{ikl}}}$};
\draw (-1,1)--(1,1)--(0.9,0.9);
\draw (1,1)--(0.9,1.1);
\draw (-4,3) -- (-6,0)--(-5,-1)--(-4,3)--(-2,0)--(-5,-1);
\draw (-5,0)-- (-6,0);
\draw (-4.5,0)--(-2,0);
\node at (0,1.3){$\Phi_m$};
\node at (-5,-1.2){$k$};
\node at (-6.2,0.1){$j$};
\node at (-4,3.2){$i$};
\node at (-1.8,0.1){$l$};
\end{tikzpicture}
 \end{center}
which is equivalent, in view of (\ref{eq:product3ary}), to
\begin{equation}\label{crossedmoduleequation2}
g_{ijk}g_{ikl}=h_{ij}(g_{jkl})g_{ijl}.
\end{equation}
 The relations (\ref{crossedmoduleequation1}) and (\ref{crossedmoduleequation2}) precisely define non-Abelian 1-cocycles.

 \section*{Summary of the thesis}

The first part of the thesis can be considered as an attempt to "stay classical", i.e. to use simply Lie groupoids in order to investigate non-Abelian gerbes. Our claim is that there is a definition of non-Abelian gerbes that uses only the notions of Lie groupoid and crossed-modules of Lie groups, and that this definition
make it possible to clarify the notion of non-Abelian gerbe over a differential stack.
The second part of the thesis is, on the contrary, an attempt to think purely in terms of higher structures
and to define Nijenhuis forms on $L_\infty$-algebras. Here, by Nijenhuis forms, we mean a generalization of
the notion of Nijenhuis $(1-1)$-tensors on manifolds, i.e. $(1-1)$-tensors whose Nijenhuis torsion vanish.
On manifolds, Nijenhuis tensors are 1-ary operations on  the Lie algebra of vector fields.
Since, when dealing with $L_\infty$-algebras, one has to replace Lie algebra brackets by collections of $n$-ary brackets for all integers $n \geq1$,
we also want to define Nijenhuis forms that are collections of $n$-ary operations for all integers $n \geq1$.

\subsection*{Non-Abelian gerbes with groupoids}

Differential gerbes appeared from the very beginning  as being classes in some "higher" cohomology \cite{Giraud}, e.g. non-Abelian gerbes correspond to non-Abelian $1$-cohomology in the sense of Dedecker \cite{BaezSchreiber,BM,Dedecker}, and it is the form under which it appears in theoretical physics \cite{GR,Witten}.
 But differential gerbes can also be thought of as being a certain class of bundles over a differential stack, and, to quote \cite{BehrendXu}, "there is a dictionnary between differential stacks and Lie groupoids".
The purpose of first chapter thesis is to add one entry to that dictionary, namely to define with great care in terms of Lie groupoids and for all
crossed modules $\G \to \HH $ the notion of $\G \to\HH $-gerbes  and to justify that definition by showing the coincidence of the notion introduced with non-Abelian $1$-cohomology.

There are, of course, several other manner to define non-Abelian gerbes, and to give properties of those.
In a recent work \cite{NikolausWaldorf} these numerous definitions have been carefully enumerated and shown, in a rigorous manner, to coincide. More precisely, the authors of \cite{{NikolausWaldorf}} have merged four definitions of smooth $\Gamma$-gerbes, with $\Gamma$ a strict 2-groups (notice that strict $2$-groups are indeed in one to one correspondence with crossed modules):
\begin{enumerate}
\item smooth $\Gamma$-valued $1$-cocycles (for which they refer to \cite{BM}, but which matches by construction the definition in terms of non-Abelian cohomology in the sense of \cite{Dedecker} just mentionned), see also \cite{AJC}.
\item classifying maps valued in the realization $B\Gamma$ of the simplicial tower of $\Gamma $,
\item bundle gerbes in the sense of \cite{Waldorf},
\item principal $\Gamma$-bundles in the sense of Bartels \cite{Bartels}, the idea being to generalize the notion of principal bundle from Lie groups to Lie (strict) 2-groups.
\end{enumerate}
 But, as mentioned in \cite{NikolausWaldorf}, Example 3.8, there is in particular case of gerbes over manifolds and of the crossed-module $\G \to \HH$ a fifth equivalent definition which is in terms of Lie groupoid extensions, a definition that avoids totally categorical language. We can restate our purpose by saying that it consists in giving this fifth description in the general setting of arbitrary crossed-module (and not only $\G \to \Aut (\G)$). Also, we give a definition that makes sense  when the base space is not a manifold but an arbitrary differential stack.

  Although the four approaches just mentioned can be remarkably effective in the sense that the objects have short, simple and workable definitions, it always requires a deep familiarity with category theory (or even toposes and higher categories) making them hardly accessible for a mathematician not used  to these technics. Our manner is maybe more difficult in the sense that the objects are always defined as classes of -oids up to Morita equivalences, which sometimes yield long definition, and forces us to check that  properties are Morita invariant, but it is certainly simpler in the sense that it uses
the ordinary language of differential geometry (manifolds, -oids, maybe \v{C}ech cocycles) from the beginning to the end.

  The present work is also in the continuation of \cite{BehrendXu} (where $S^1 $-gerbes over a differential stack are extensively studied using this Lie groupoid point of view),of \cite{LaurentStienonXu} (where the case of non-Abelian$\G\to\Aut(\G)$-gerbs over Lie groupoids is investigated, but the correspondence with non-Abelian  $1$-cocycles is not dealt with very precisely.
 and of \cite{BL}, (where the previous construction is investigated in detail for $\G$-gerbs and extended to connections). Our work is definitively in the same line of those, but there are important differences that we now outline. Abelian gerbes in the sense of \cite{BehrendXu} (resp. $\G$-gerbes in the sense of \cite{LaurentStienonXu}) corresponds to the case where the crossed-modules in which the gerbe takes values is $S^1 \to pt$ (resp. $\G \to \Aut (\G)$), so that our work generalizes both. Second, we made more precise the notion of gerbes over an object (manifold, Lie groupoid, or differential stack). This means that , unlike \cite{LaurentStienonXu}, we do not simply define gerbes as being $\G$-extensions up to Morita equivalence, and this is for two reasons:
\begin{enumerate}
\item first, as already stated, we wish to make precise over what object our gerbe is, which
means that we only allow ourself to take Lie groupoid extensions $X \to Y $ where
the "small" Lie groupoid $Y$ is itself "over" a given object $B$ (manifold or Lie groupoid or differential stack).
 By "over", we mean that "Y" is obtained by taking a pull-back of $B$. Also, Morita equivalence should be taken in such a way that the base manifold or groupoid is not "changed". This last issue is easily understandable, and always appear in differential geometry: the space of principal bundles over a manifold $M$, in a similar fashion, is not obtained by considering all possible principal bundles $P \to M$ modulo principle bundles isomorphisms, but modulo principal bundles isomorphisms over the identity of $M$.
\item second, when taking an arbitrary crossed-module $\G \to \HH $, Lie groupoid extensions are not enough.
By spelling out the manifold case, and knowing that we wish to have a correspondence with crossed-module valued non-Abelian $1$-cocyle, we arrived at the conclusion that we need to consider a Lie groupoid $\G$-extension together with a $\HH$-principal bundle. These two structures are not independent, and, having in mind the manifold case again, one sees that we need this principle bundle to be equipped with a principal bundle morphism taking values in the band of the Lie groupoid extension, map on which still two constraints have to be imposed.
\end{enumerate}

Chapter $I$ is organized as follows.
In section \ref{sec:GpdExt}, we recall from \cite{LaurentStienonXu} the notion of $\G$-extensions of Lie groupoids, i.e. a surjective submersion morphism
 of Lie groupoids over the same base ${\mathcal R} \stackrel{\phi}{\to} \GGG $, for which the kernel is a locally trivial bundle of groups with typical fiber $ \G$. We then recall, following \cite{LaurentStienonXu}, the notion of the band of the $\G$-extension, which is some principal bundle over the Lie groupoids
${\mathcal R}$. We then define $ \G \to \HH$-extensions, namely $\G$-extensions ${\mathcal R} \stackrel{\phi}{\to} \GGG $ endowed with some
$ \HH$-principal bundle which admits the band as a quotient, see definition \ref{def: Azimi-extension} for a more precise description.

We then recall the definition of Dedecker's non-Abelian $1$-cocyle (resp. non-Abelian $1$-coboundaries, non-Abelian $1$-cohomology)
on an open cover of a given manifold $N$ and describe a dictionary between these objects and $ \G \to \HH$-extensions.
More precisely, we define, given an open cover of a manifold, a subclass of $ \G \to \HH$-extensions called
adapted $ \G \to \HH$-extensions of the \v{C}ech groupoid, and we show the following points, given an open cover ${\mathcal U} $ on the manifold $ N$:
\begin{itemize}
 \item \emph{Proposition \ref{prop:cocycles}}, There is a one to one correspondence between:
         \begin{enumerate}
          \item[(i)] $ \G \to \HH$-valued non-Abelian $1$-cocycles w.r.t. ${\mathcal U} $
          \item[(ii)] adapted $ \G \to \HH$-extensions of the \v{C}ech groupoid $ N[{\mathcal U}]$
         \end{enumerate}
 \item \emph{Proposition \ref{prop:coboundaries}}, There is a one to one correspondence between:

         \begin{enumerate}
          \item[(i)] $ \G \to \HH$-valued non-Abelian $1$-coboundaries w.r.t. ${\mathcal U} $
          \item[(ii)] isomorphisms of adapted $ \G \to \HH$-extensions of the \v{C}ech groupoid $ N[{\mathcal U}]$
         \end{enumerate}
 \item \emph{Theorem \ref{cor:1-1Cohomology-adapted}}, There is a one to one correspondence between:
    \begin{enumerate}
          \item[(i)] $ \G \to \HH$-valued non-Abelian $1$-cohomology w.r.t. ${\mathcal U} $
          \item[(ii)] isomorphism classes of adapted $ \G \to \HH$-extensions of the \v{C}ech groupoid $ N[{\mathcal U}]$ up to Morita
 equivalence over the identity.
\item[(iii)] (assuming the covering to be a good one) isomorphism classes of $ \G \to \HH$-extensions of the \v{C}ech
groupoid $ N[{\mathcal U}]$ up to isomorphisms over the identity of $N[{\mathcal U}]$.
         \end{enumerate}
\end{itemize}

The first purpose of section \ref{sec:MoritaExt} is to overcome of the choice of an open cover and to reach therefore $ \G \to \HH$-valued non-Abelian $1$-cohomology
in its full generality. This requires to define the notion of Morita equivalence of  $ \G \to \HH$-extensions, which, in turn, allows to complete the previous isomorphisms to eventually obtain the one we are really interested in:
\begin{itemize}
\item \emph{Theorem \ref{th:coho=gerbes}}, There is a one to one correspondence between:
    \begin{enumerate}
          \item[(i)] $ \G \to \HH$-valued non-Abelian $1$-cohomology,
\item[(ii)]    $ \G \to \HH$-extensions of a pull-back of the
groupoid $ N \toto N $ up to Morita equivalence over the identity of $N$.
         \end{enumerate}
\end{itemize}
The point of this last theorem gives a clear hint of what a $ \G \to \HH$-gerbe over a given Lie-groupoid $\BBB$ should be, namely the
$ \G \to \HH$-extensions of a pull-back of the groupoid $ \BBB $ up to Morita equivalence over the identity of $\BBB$.
We conclude by saying that Morita equivalent Lie groupoids $\BBB$ and $\BBB'$ have the same
 $ \G \to \HH$-gerbe over them, making sense therefore of the notion of  $ \G \to \HH$-gerbes over a differential stack.

\subsection*{Deformations of $L_\infty$-algebras by Nijenhuis forms}\label{sec:deformationofLinfty}

Unlike the first part of the thesis, the second part of the thesis is an attempt to think purely in terms of higher structures
in order to define Nijenhuis forms on $L_\infty$-algebras.

 Given a Lie algebra $({\mathfrak g},[.,.])$
and a linear endomorphism $N$ of ${\mathfrak g}$,
the deformed bracket of $[.,.]$ by $N$ is, by definition, the bilinear map $[.,.]_N$ given by $[X,Y]_N = [NX,Y]+[X,NY] - N([X,Y])$.
The Nijenhuis torsion of $N$ \cite{Nijenhuisformingsets} is defined by:
 \begin{equation}\label{tortionintro}
  TN(X,Y) := [NX,NY] - N([X,Y]_N)
 \end{equation}
and $N$ is said to be Nijenhuis if the Nijenhuis torsion of $N$ vanishes, that is to say, $N$ is a Lie algebra morphism from $({\mathfrak g},[.,.]_N)$
to $({\mathfrak g},[.,.])$.
A simple computation shows that $N$ is Nijenhuis if, and only if,
\begin{equation}\label{eq:NijenhuiswithSquaresIntro}
[.,.]_{N,N} =[.,.]_{N^2},
\end{equation} so that there are two interpretations of Nijenhuis endomorphisms:
\begin{enumerate}
\item An endomorphism $N$ is Nijenhuis if $N$ itself is a morphism from the modified bracket to the initial one,
\item An endomorphism $N$ is Nijenhuis if deforming twice by $N$ the original bracket yields the original bracket deformed by $N^2$.
\end{enumerate}
An important point is that if $N$ is Nijenhuis, then the deformed bracket $[.,.]_N$ is still a Lie bracket.


 If $({\mathfrak g},[.,.])$ is the Lie algebra of smooth vector fields on a manifold $M$, then $TN(X,Y)=[NX,NY]-N([X,Y]_N) $ is a $(2-1)$-tensor called Nijenhuis torsion. This is certainly the most famous appearance of this notion: this goes back to the Newlander-Nirenberg's theorem,
which states that an almost complex structure comes from a complex structure if and only if it is a Nijenhuis $(1-1)$-tensor \cite{NewlanderNirenberg}.
In the seventies, Nijenhuis tensors have been proved to appear naturally while studying integrable systems, more precisely bi-Hamiltonian systems.
The pioneering article by Franco \cite{Magri} led to the notion of Poisson-Nijenhuis structures \cite{MMR} following notes by Magri and Morosi \cite{MM}.

 The notions of Nijenhuis tensors and Poisson-Nijenhuis manifolds were later understood in a more algebraic setting by Kosmann-Schwarzbach and Magri \cite{YKS}. These settings allowed to define those notions for Lie algebroids, see
 Kosmann-Schwarzbach \cite{YKSbialgebroid}, Grabowski and Urbanski \cite{GrabowskiUrbanski} and Clemente-Gallardo and Nunes da Costa \cite{Clemente-Nunes} for the more general case of Dirac-Nijenhuis structures. In \cite{CGM2}, Cariñena, Grabowski and Marmo extended Nijenhuis tensors to general (binary) algebraic structures, while Cariñena, Grabowski, Marmo \cite{CGM}, Kosman-Schwrzbach \cite{YKSBrazil} and Antunes, Laurent-Gengoux and Nunes da Costa \cite{CJP} studied Nijenhuis structures on Loday algebras and Courant algebroids. Through  this generalization, Poisson structures on Lie algebroids become particular cases of Nijenhuis tensors on some Courant structure \cite{AntunesCosta}, so that
Poisson-Nijenhuis are pairs of compatible Nijenhuis tensors \cite{CJP}.

Jump to $n$-ary. This is not easy: for binary operations, the vanishing of the torsion (\ref{tortionintro}) can be interpreted as meaning that $N$ is a morphism from the deformed binary operation to the initial one, as we said above. It does not seem reasonable  (at least, examples make us believe that it is not reasonable) to hope that such an interpretation can still be valid when working with collections of $n$-ary operations. By chance, the second interpretation of Nijenhuis
structures defined above can be generalized, with the help of a natural bracket that exists on
the vector space of multi-linear skew-symmetric (or symmetric) endomorphisms of a graded vector space : the Schouten-Nijenhuis bracket.

The Schouten-Nijenhuis bracket $[.,.]_{_{SN}}$, defined on the exterior algebra $\Gamma(\wedge TM)=\oplus_{i\geq 0}\\\Gamma(\wedge^i TM)$ of multi-vector fields on a manifold $M$, was introduced by Schouten and Nijenhuis who studied its properties in
 \cite{Scouten-differential-operators} and \cite{Nijenhuis-Jacobi-type}.

It becomes a remarkable tool in Poisson geometry when Lichnerowicz \cite{Lichnerowicz} proved that a bivector field $\pi$ on a manifold is Poisson if, and only if, the Schouten-Nijenhuis bracket of $\pi$ with itself is zero, $[\pi,\pi]_{_{SN}}=0$.

The Schouten-Nijenhuis bracket of two vector fields on a manifold, is just the Lie bracket of the vector fields. Moreover, assuming that $P\in\Gamma(\wedge^pTM)$ has degree $p-1$, the pair $(\Gamma(\wedge TM),\, [.,.]_{_{SN}})$ is a graded Lie algebra whhile the triple $(\Gamma(\wedge TM),\,\wedge.\, [.,.]_{_{SN}})$ is a Gerstenhaber algebra.

It turns out that, roughly speaking, the role of $\Gamma(\wedge TM)$ can be replaced by the space of multi-sections of a Lie algebroid $A$ and the extension, by derivation, of the Lie bracket on $\Gamma(A)$ determines a bracket on $\Gamma(\wedge A)$ which is also called the Schouten-Nijenhuis bracket of the Lie algebroid.

But despite this powerful tool, it is not that trivial to generalize the characterization (\ref{eq:NijenhuiswithSquaresIntro}) of Nijenhuis structures,
 because there is no clear definition of what the square of a collection of $n$-ary operations is - unlike for $1$-ary operation.
 In fact, the square for us will be basically any (in general quadratic) expression in the $n$-ary operations that compose $N$ and that commutes with $N$.

 The Richardson-Nijenhuis bracket of two symmetric vector valued forms, to be introduced in Section \ref{section:3.1}, is essential in our definition of Nijenhuis form in this more general context.
 \begin{defi}
 Given a symmetric vector valued form $\mu$ of degree $1$ on a graded vector space $E$, a vector valued form $\mathcal{N}$ of degree zero is called a Nijenhuis (vector valued) form, with respect to $\mu$ if there exists a vector valued form $\mathcal{K}$ of degree zero such that
   $$ [\mathcal{N},[\mathcal{N},\mu]_{_{RN}}]_{_{RN}}=[\mathcal{K},\mu]_{_{RN}} \hbox{ and } [\mathcal{N},\mathcal{K}]_{_{RN}}=0, $$
    where $[.,.]_{_{RN}}$ stands for the Richardson-Nijenhuis bracket.
 \end{defi}

 This definition has to be justified. A fundamental result, \cite{YKS} shows that when $N$ is Nijenhuis for a Lie algebra bracket
 $[.,.] $ then for all integers $ n,m \geq 1$, we have
 \begin{equation}\label{eq:powersofN}
  [.,.]_{N^n,N^m} = [.,.]_{N^{n+m}}.
  \end{equation}
 In particular, $[.,.]_{N^{n}}$ is a Lie bracket for all integer $n > 0 $, so that a whole hierarchy of Lie algebra brackets can be constructed.
 We have been able to mimic these results, but for these we have to translate them without using powers of $N$:
 (\ref{eq:powersofN}) means that $ [.,.]_{N,\dots,N} = [.,.]_{N^k}$
 so that it can be restated as meaning that all the brackets $ [.,.]_{N,\dots,N} $ obtained by successive deformations are Lie brackets,
 and that $ N$ is Nijenhuis for all of them. Theorem \ref{theo:Hierarchy} states that this remains true with our definition of Nijenhuis
 forms. Also, we explain in Remark \ref{rem:whynotsquares} why we could not replace $N^2$ by the pre-Lie product that helps defining
 the Richardson-Nijenhuis bracket, an idea that could seem natural at first.

 Of course, this definition has also to be proved to be non-empty, and to contain more than just ordinary Nijenhuis endomorphisms
 on ordinary Lie algebras. By making it non-empty, we also mean to find examples for which the Nijenhuis is not simply an endomorphism,
 but the sum of an endomorphism with a bilinear map with a trilinear map and so on.

 We now list the examples that we have so far. The first example is universal, in the sense that every $L_\infty$-structure admits
 it: the Euler map $S$, that multiplies an element by its degree.
 Of course, ordinary Nijenhuis tensors on ordinary graded Lie algebras are among the most trivial examples.
 Poisson elements, and more generally, Maurer-Cartan elements of differential graded Lie algebras
 are also examples, which are not purely made of vector valued $1$-forms, but which are the sum of a vector valued $1$-form with a vector valued $0$-form.
 \begin{itemize}
        \item \emph{Corollary \ref{cor1poissononGLA}}, Let $\mu=l_2$ be a symmetric graded Lie algebra structure on a graded vector space $E=\oplus_{i\in \mathbb{Z}}E_i$ and $\pi\in E_0$. Then $\pi+S$ is a Nijenhuis vector valued form, with respect to $\mu$ and with square $2\pi+S$ if, and only if, $\pi$ is a Poisson element.

  \end{itemize}

\begin{itemize}
        \item \emph{Corollaries \ref{cor2poissononDGLA} and \ref{cor3poissononDGLA}}, Let $l_1+l_2$ be a symmetric DGLA structure on a graded vector space $E=E_{-2}\oplus E_{-1}$ and $\pi\in E_0$. Then
                \begin{enumerate}
                    \item $\pi+S$ is a Nijenhuis vector valued form with respect to $\mu$, with square $2\pi+S$ if, and only if, $\pi$ is a Poisson element.
                     \item $Id_E+\pi$ is Nijenhuis vector valued form, with respect to $\mu$ and with square itself if, and only if, $\pi$ is a Maurer-Cartan element of the  DGLA $(E,\mu)$,
                \end{enumerate}
\end{itemize}

   Less trivial examples are given on the Lie $n$-algebras.
 On those, we can expect to have Nijenhuis forms which are not purely vector valued $1$-forms, but which are the sum of a family of vector valued $k$-forms.

\begin{itemize}
\item \emph{Proposition \ref{corlast}}, Let $(E=E_{-n}\oplus \cdots\oplus E_{-1}, \mu=l_1+\cdots +l_{n+1})$ be a Lie $n$-algebra. Let $N_1,\cdots,N_l$ be a family of symmetric vector valued $k_1,\cdots,k_l$-forms, respectively, of degree zero on $E$, with $\frac{n+3}{2}\leq k_1\leq \cdots\leq k_l \leq n+1$. Then $S+\sum_{i=1}^l N_i$ is a Nijenhuis vector valued form with respect to $\mu$, with square $S+2\sum_{i=1}^l N_i$.
\end{itemize}
 For Lie $n$-algebras, there is another class of examples that we have in mind:
 $n$-plectic manifolds. It is shown that $n$-plectic manifolds give a Lie $n$-algebra structure that can be
 deformed by a family of forms.
 \begin{itemize}
\item \emph{Theorem \ref{nplecticlasttheorem}}, Let $(\eta^j)_{j\geq 1}$ be a family of $n$-forms on an $n$-plectic manifold $(M,\omega)$. Let $(E=E_{-n}\oplus \cdots\oplus E_{-1}, \mu=l_1+\cdots+l_{n+1})$ be the associated Lie $n$-algebra to $(M,\omega)$. For each $2\leq i\leq n$ define the vector valued $i$-forms $\widetilde{\eta^j_i}$ as
\begin{equation*}
\widetilde{\eta^j_i}(\beta_1,\cdots,\beta_i)=\begin{cases}
                                             \iota_{\chi_{\beta_1}}\cdots\iota_{\chi_{\beta_i}}\eta^j, & \mbox{if} \,\,\,\ \beta_k \in E_{-1}\,\,\, \mbox{for all} \,\,\,1\leq k\leq i\\
                                             0, & \mbox{otherwise}
                                             \end{cases}
\end{equation*}where $ \chi_{\beta_1},\cdots,\chi_{\beta_i}$ are the unique Hamiltonian vector fields associated to the Hamiltonian forms $\beta_1,\cdots,\beta_i$ respectively. Then $\mathcal{N}:=S+\sum_{j\geq 1}\sum_{i=2}^{n}\widetilde{\eta^j_i}$ is a Nijenhuis vector valued form with respect to the Lie $n$-algebra structure $\mu=l_1+\cdots+l_{n+1}$, associated to the $n$-plectic manifold $(M,\omega)$.
\end{itemize}
The case of Lie $2$-algebra is treated separately. Here, we have Nijenhuis forms which are the sum of a vector valued $1$-form with a vector valued $2$-form.
\begin{itemize}
\item \emph{Theorem \ref{thm:alphalpha}}, Let $\mu=l_1+l_2+l_3$ be a Lie $2$-algebra structure on a graded vector space $E=E_{-2}\oplus E_{-1}$ and $\alpha$ be a symmetric vector valued $2$-form of degree $0$. Then  $S+\alpha$ is a Nijenhuis vector valued form with respect to $\mu$ with square of $S+2\alpha $ if and only if
 $$ \alpha(l_1 \alpha (X,Y),Z ) +c.p.=0,$$
 for all $ X,Y,Z \in E_{-1}$.
\end{itemize}

  Indeed, we can make sense of the deformations of such vector valued $2$-forms.
 For strict Lie algebra for instance, which is associated to a $3$-cocycle on an ordinary Lie algebra,
 adding a coboundary to this cocycle corresponds to deformation by a Nijenhuis form.
 Some classes of Lie $2$-algebras are shown in Proposition \ref{prop:strict+trivial}
 to be obtained by a Nijenhuis deformation of strict and trivial Lie $2$-algebras.
\begin{itemize}
 \item \emph{Proposition \ref{prop:strict+trivial}}, Given a Lie $2$-algebra structure $l_1+l_2+l_3$ on a graded vector space $E=E_{-2}\oplus E_{-1}$ such that $l_2$ vanishes on the elements of degree $-3$, there exists a Nijenhuis transformation
of the form $S+\alpha$ with $\alpha$ a vector valued $2$-form of degree zero, such that the deformed bracket $[S+\alpha,l_1+l_2+l_3] $
is the direct sum of a strict Lie $2$-algebra with a trivial $L_\infty$-algebra.
\end{itemize}

 The next example is Courant algebroids. There are already quite a lot of works on Nijenhuis structures on Courant algebroids,
 as already said. We have been able to show that these Nijenhuis tensors give examples of Nijenhuis tensors on the Lie $2$-algebra associated with,
 in Proposition \ref{prop:NijenhuisonCourant}.
 Then, we have stated several results showing that, in some sense, there is no hope to get more examples than those already known,
 at least if we impose the (reasonable) condition that the Nijenhuis tensors have to be ${\mathcal C}^\infty (M)$-linear.
  Indeed, Corollary \ref{thewrongcorollary} classifies (almost) entirely Nijenhuis tensors on the Lie $2$-algebra associated to a Courant structure.
 More precisely,
\begin{itemize}
\item let $ (E, \circ, \rho, \langle ,\rangle ) $ be a Courant algebroid with associated symmetric Lie $2$-algebra $\mu=l_1+l_2+l_3$, on the graded vector space $V=\mathcal{C}^{\infty}(M)\bigoplus \Gamma(E)$.
          \begin{enumerate}
               \item \emph{Proposition \ref{prop:NijenhuisonCourant}}, Let $N:\Gamma(E)\to \Gamma(E)$  be a Nijenhuis $(1-1)$-tensor satisfying
                   \begin{equation*}
                      \begin{cases}
                           N+N^*=\lambda I,\\
                           N^2+(N^2)^*=\gamma I.
                      \end{cases}
                   \end{equation*}
              with $\lambda,\gamma$ being Casimir functions.
              Define   ${\mathcal N}$ and ${\mathcal K}$ as
                 \begin{equation*}
                      {\mathcal N}|_{\Gamma(E)}= N \,\,\,\mbox{and}\,\,\,{\mathcal N}|_{\mathcal{C}^{\infty}(M)}=\lambda I_{\mathcal{C}^{\infty}(M)},
                 \end{equation*}
                 \begin{equation*}
                    {\mathcal K}|_{\Gamma(E)}= N^2 = \lambda N - \frac{\gamma - \lambda^2}{2} Id \,\,\,\mbox{and}\,\,\,{\mathcal K}|_{\mathcal{C}^{\infty}(M)}=\gamma I_{\mathcal{C}^{\infty}(M)}
                 \end{equation*}
               Then $\mathcal N$ is a Nijenhuis vector valued $1$-form with respect to $\mu$, with square $\mathcal K$.
         \item \emph{Corollary \ref{thewrongcorollary2}}, There is  one to one correspondence between:
           \begin{enumerate}
 \item[(i)] quadruples $(N,K,\lambda,\gamma)$ with $N,K$ being $(1-1)$-tensors and $\lambda,\gamma$ being Casimir functions satisfying the following conditions:
    $$ \begin{cases}
     \circ^{N,N}=\circ^K,\\
     N K - K N =0,\\
     N+N^*=\lambda Id_{\Gamma(E)} ,\\
     K+K^* = \gamma Id_{\Gamma(E)}.\\
     \end{cases}$$
 \item[(ii)] Nijenhuis vector valued forms ${\mathcal N} $ with respect to $\mu$, with square ${\mathcal K} $ such that the deformed brackets $[ {\mathcal N}, \mu]_{_{RN}} $ is a pre-Lie $2$-algebras associated to  pre-Courant structures with the same scalar product.
     \end{enumerate}
 \end{enumerate}
%
\end{itemize}

  In the last section, we investigate several cases around Lie algebroids, that we see as being graded Lie algebras
  by considering their Schouten-Nijenhuis bracket. We then show that the following objects give examples of Nijenhuis forms
  (or at least weak Nijenhuis or co-boundary Nijenhuis forms)
  on this graded Lie algebras:Nijenhuis tensors on Lie algebroids, $\Omega N$-structures, Poisson-Nijenhuis structures and $\Pi \Omega$-structures.
 \begin{itemize}
\item \emph{Proposition \ref{prop:NijenhuisAsNijenhuis}},
 For every Nijenhuis tensor field $N$ on a Lie algebroid $(A,\left[.,.\right],\rho)$,
the extension of $N$ by derivation, $\underline N $, is a Nijenhuis vector valued $1$-form with respect to the multiplicative GLA-structure $l_2^{\left[.,.\right]_N}$ on the graded vector space $\Gamma(\wedge A)[2]$, with square ${\underline{(N^2)}}$.
 \end{itemize}
  \begin{itemize}
  \item \emph{Corollary \ref{NijenhuisonOmegaNstructures}},
  Let $(A,[.,.], \rho)$ be a Lie algebroid, with the associated de Rham differential $\diff^A$ and with the associated multiplicative GLA structure $l_2^{\left[.,.\right]}$ on the graded vector space $\Gamma(\wedge A)$. Let $(N,\alpha)$ be an $\Omega N$-structure on the Lie algebroid $A$, then $\underline{N}+\underline{\alpha}$ is a Nijenhuis vector valued form, with respect $l_2^{\left[.,.\right]}$, with square $\underline{N^2}+\underline{\alpha_{_{N}}}$.
  \end{itemize}
\begin{itemize}
\item \emph{Proposition \ref{NijenhuisonPoissonNijenhuisalgeroid}},
  Let $(N,\pi)$ be a Poisson-Nijenhuis structure on a Lie algebroid $(A,[.,.],\rho)$,
then the derivation $\underline N $ is a weak-Nijenhuis tensor
for the $L_\infty$-structure $\mu=l_1^{\left[.,.\right],\pi}+l_2^{\left[.,.\right]}$ (on the graded vector space $\Gamma(\wedge^A)[2]$), associated to the Poisson structure $\pi$ which is given by
\begin{equation*}
l_1^{\left[.,.\right],\pi}(P)=\left[\pi,P\right]_{_{SN}} \quad \mbox{and} \quad l_2^{\left[.,.\right]}(P,Q)=(-1)^{p-1}\left[P,Q\right]_{_{SN}}.
\end{equation*}
\end{itemize}
\begin{itemize}
\item \emph{Proposition \ref{Poissonifandonlyif}},
  Let $(A,[.,.],\rho)$ be a Lie algebroid , $\pi\in \Gamma(\wedge^2 A)$ be a bi-vector and $N:\Gamma(A)\to \Gamma(A)$ be a $(1-1)$-tensor field such that
\begin{equation*}
N\pi^{\#}=\pi^{\#}N^*.
\end{equation*}
Then
 $\underline{N} + \pi $ is a co-boundary Nijenhuis form, with respect to the multiplicative GLA-structure $l_2^{\left[.,.\right]}$
with square ${\underline{N^2}} $, if and only if $(N,\pi)$ is a Poisson-Nijenhuis structure on the Lie algebroid $(A,[.,.],\rho)$.
\end{itemize}
\begin{itemize}
\item \emph{Proposition \ref{NijenhuisonPoissonOmega}},
  Let $(\pi,\omega)$ be a  $\Pi\Omega$-structure on a Lie algebroid $(A,\left[.,.\right],\rho)$. Then,
 ${\mathcal N}=\underline{\omega} + \pi $ is a co-boundary Nijenhuis form, with respect to the multiplicative GLA-structure $l_2^{\left[.,.\right]}$
with square $\underline{N} $, where $N = \pi^{\#} \circ \omega^b$.
\end{itemize} 
\mainmatter
\pagenumbering{arabic}

\chapter{Non-Abelian gerbes as Lie groupoid extensions}

The purpose of this section is to give a definition of non-Abelian gerbes purely with the help of Lie groupoids.
We refer to the introduction for the motivations.
We start by reviewing several notions related to Lie groupoids, culminating in the definition
of the object called differential stacks.

\section{Definitions and notations}

\emph{Notations related to open covers on manifolds.}\label{notation:Cech}
For ${\mathcal U}= (U_i)_{i \in I} $ an open cover on a manifold $N$,
we use the shorthand $U_{ij} = U_i \cap U_j$ for all $i,j \in I$, and introduce the convenient notation $$U_{i_1 \dots i_n} := U_{i_1} \cap \dots \cap  U_{i_n} $$  for all $n \in {\mathbb N}$ and all $i_1,\dots,i_n \in I$. We warn the reader that $U_{ij}$ is not equal to $U_{ji}$ for $i \neq j$, and, more generally, $U_{i_1 \dots i_n} $ is not equal $U_{i_{\sigma(1)} \dots i_{\sigma(n)}} $ (for $\sigma \in \Sigma_n$ a permutation, and $i_1,\dots,i_n$ distinct).

An extremely common notation in the literature dealing with gerbes is to denote by $x_i$ (resp. $x_{ij},x_{ijk}$) an element $x\in M$ that happens to belong to some open subset $ U_i$ (resp. $U_{ij}, U_{ijk} $), when it is seen as an element in $U_i$ (resp. $U_{ij},U_{ijk} $).
  We extend this convention for all kind of objects: for instance, for a function $ \lambda$ whose domain of definition is  $\coprod_{i_1,\dots,i_n \in I}  U_{i_1 \dots i_n}$,
we write $ \lambda_{i_1\dots i_n} $ for its restriction to~$U_{i_1 \dots i_n}$.
%





\bigskip
\emph{Lie groupoids : notations and basic facts.}
  Given $ M,N, P$ smooth manifolds and $f:M \to P $, $g: N \to P$ smooth maps, we define the \emph{fibered product} to be the closed subset of $M \times N  $ made of all pairs $(m,n)$ with $f(m)=g(n)$, we denote it by $M \times _{f,P,g} N$ in general, and sometimes by $M \times _{P} N$ when there is no risk of confusion. The following is extremely classical:
 \begin{lem} \label{lem:manifold}\cite{BergerGostiaux}
  Let $ M,N, P$ be smooth manifolds. If at least one of the smooth maps $f:M \to P $ or $g: N \to P$ is a surjective submersion, then the set $M \times _{f,P,g} N$ is a smooth manifold.
  \end{lem}
We refer to \cite{McK} for the definition of Lie groupoids, but we wish to clarify some notations. When introducing a Lie groupoid, we shall in general simply mention the names of the manifolds of objects and the manifolds of arrows, using the notation $\Gamma \toto M $. Indeed, the source, target and unit maps for all Lie groupoids $\Gamma \toto M$ shall be denoted by the same letters $s,t$ and $\epsilon$ respectively. In general, the product shall be either denoted by the fat dot $\bullet$
or simply skipped, and the inverse by the exponent $-1$. However, at some point, we shall have to consider pairs of manifolds that admit several different Lie groupoid structures, that, fortunately,
have the same source, target and unit maps. We will then introduce a  notation for the product (and inverse) that will distinguish them.
Last, our convention is that the product $x\bullet x'$ of two elements $x,x'$ in a Lie groupoid is defined when $t(x)=s(x')$.

 A left-action of Lie groupoid $\mathcal B \toto \mathcal B_0$ on a  manifold $X$ with respect to a surjective submersion $p:X \to \mathcal B_0$ is a map
      $$ \mathcal B \times_{t,\mathcal B_0,p} X\longrightarrow X,$$
 (denoted by $(b,x)\mapsto b\cdot x$) such that $p(b\cdot x)= s(b)$ and subject to the following axioms, analogous to those of group actions:
 $$b\cdot (a\cdot x)=(b  a)\cdot x \hbox{ and } \epsilon({p(x)}) \cdot x=x, $$
 for all admissible $a,b \in \mathcal B$ and $x\in X$. We shall often say action for left-action for the sake of simplicity.
   Since we may have to deal with situations where there are more than one Lie groupoid or more than one manifold involved, it will be convenient to write an action by $b\bullet_{\mathcal B, X}x$, mentioning therefore in the notation itself which groupoid acts and which manifold is acted upon.

\section{Differential stacks}

The purpose of this section is to define differential stacks with the help of Lie groupoids, in the spirit of \cite{BehrendXu}.
To ensure a self-contained exposition, we recall all the steps of this construction, following \cite{LaurentHabilitation}
as a guideline. For dealing with gerbes, we shall need to define Morita equivalence with the help of
pull-back Lie groupoids, which goes with some technical difficulties.
For this purpose, we recall successively:
\begin{enumerate}
\item the notion of Lie groupoid pull-back,
\item the notion of Lie groupoid Morita equivalence, as defined with the help of pull-back Lie groupoids.
\end{enumerate}
From this, we shall define differential stacks as being Lie groupoids up to Morita equivalence.

We start with the notion of pull-back Lie groupoid.
Notice first that, for a given manifold $B$, Lie groupoids that admit  $B$ as the unit manifold form a category, with morphisms being
Lie groupoid morphisms over the identity of $B$. Similarly, topological groupoids that admit  $B$ as the unit manifold  form a category.

\begin{defi}\label{def:pback}\cite{McK}
Let $p:M \to B  $ be a smooth map. The assignments below define a functor from the category of Lie groupoids
over $B$ to the category of topological groupoids over $M$:
\begin{enumerate}
\item {\emph (On objects)}
Let $\GGG  \toto B$ be a Lie groupoid over a manifold $B$. Then the set $ \GGG[p,M]:= M\times_{p,B,s}\GGG\times_{t,B,p} M $ (sometimes simply denoted by  $\GGG[p]$)is endowed with a topological groupoid structure over $M$ given as follows: the source and target $ s, t:\GGG[p,M] \rightarrow M$ are the projections on the first and the third components respectively, the unit map is given for all $x \in M $ by $x \mapsto (x,\varepsilon\circ p(x),x)$, where $\varepsilon$ is the unit map of the Lie groupoid $\GGG\toto B $. Last, the multiplication and the inverse are given by:
     $$(x,\gamma, y)\bullet(y,\gamma' , z)=(x, \gamma\bullet\gamma',z) \hbox{ and }
   (x,\gamma ,y)^{-1}=(y,\gamma^{-1},x). $$
   for all $x,y,z \in M$ and $\gamma , \gamma' \in \GGG.$
\item {\emph (On arrows)} Let $\phi : \GGG \to \GGG' $ be a Lie groupoid homomorphism over the identity of $B$. We set $\phi[p,M]$ to be $(n,r,n') \mapsto (n,\phi(r),n')$ for all $(n,r,n' ) \in \GGG[p,M] = M \times_{p,B,s} \GGG \times_{t,B,p} M $. By construction, $\phi[p,M] $ is a  Lie groupoid homomorphism over the identity of $M$ from $ \GGG[p,M]$ to $ \GGG'[p,M] $.
 \end{enumerate}
%
     The topological groupoid $\GGG[p,M] \toto M $ is called the \emph{pull-back of $\GGG \toto B$ with respect to $p: M \rightarrow B $}, or simply the pull-back groupoid when there is no risk of confusion.
\end{defi}

Indeed, the previous functor takes values in the category of Lie groupoids when $p$ is a surjective submersion.
More generally  \cite{CrainicMoerdijk}:

 \begin{lem}\label{lem: generalized surjective submersion}
 Let $\GGG \toto B$ be a Lie groupoid, $M$ be a manifold and $p:M \to B$ a smooth map. Then  $\GGG[p,M] $ admits a structure of Lie groupoid on the manifold $M$ if the map $\phi : M \times_{p,B,s}\GGG \to B$ given by $(m,\gamma) \mapsto t(\gamma) $, for all $(m,\gamma)\in M\times_{P,B,s}\GGG$, is a surjective submersion (in which case $p$ is called a \emph{generalized surjective submersion for the Lie groupoid $\GGG \toto B$}).
 \end{lem}
 \begin{proof}
 Lemma \ref{lem:manifold} applied to $\phi : M\times_{p,B,s}\GGG \to B$ and $p: M \to B$ implies that $(M \times_{p,B,s}\GGG) \times_{t,B,p} M$ is a manifold. It is routine to check that $(M \times_{p,B,s}\GGG) \times_{t,B,p} M$, together with the structure maps defined in definition \ref{def:pback} is a Lie groupoid.
 \end{proof}

We can now define Morita equivalence.
In fact, we define two different types  of Morita equivalence, at this point,
later on, we shall show that they coincide, in the sense that the equivalence relation they define
on Lie groupoids, coincide.

\begin{defi}
A weak (resp. strong) Morita equivalence between two Lie groupoids $\GGG \toto M$ and $\GGG' \toto M' $ is a triple
$(M'',p,p') $ where $M''$ is a manifold, $p: M'' \to M $ and $q:M'' \to M'$ are generalized surjective submersions
(resp. surjective submersions),
together with an isomorphism between the pull-back Lie groupoid $\GGG[p,M'']\toto M'' $ and the pull-back Lie groupoid  $ \GGG'[p',M''] \toto M'' $.
\end{defi}
In terms of commutative diagram, Morita equivalence of Lie groupoid can be visualized as follows
$$ \xymatrix{\GGG[p,M''] \ar@<1pt>[dr] \ar@<-1pt>[dr] &\approx                        & \GGG[p',M''] \ar@<1pt>[dl] \ar@<-1pt>[dl]\\
             \GGG \ar@<1pt>[d] \ar@<-1pt>[d]         &M'' \ar[dl]_{p} \ar[dr]^{p'}            & \GGG' \ar@<1pt>[d] \ar@<-1pt>[d]\\
             M                                       &                               & M'}$$
If there exists a Morita equivalence between two Lie groupoids, then they are said to be \emph{Mortita equivalent.}
\begin{examp} \label{ex:isom}
Any pair of isomorphic Lie groupoids are Morita equivalent.
\end{examp}

\begin{examp} \label{ex:pullback}
Every pull-back $\GGG[p,M]\toto M'$ of a given Lie groupoid $\GGG\toto M$ via a surjective submersion $p:M'\to M$ is Morita equivalent to the Lie groupoid $\GGG\toto M$, itself. which can be seen in the following diagram:
$$ \xymatrix{\GGG[p, M'] \ar@<1pt>[dr] \ar@<-1pt>[dr] \ar[rr]^{\phi} &                        & \GGG[p,M'][id_{M'}, M'] \ar@<1pt>[dl] \ar@<-1pt>[dl]\\
             \GGG \ar@<1pt>[d] \ar@<-1pt>[d]         &M' \ar[dl]_{p} \ar[dr]^{id_{M'}}            & \GGG[p,M'] \ar@<1pt>[d] \ar@<-1pt>[d]\\
             M                                       &                                & M'\ar[ll]^{p}}$$
where $\phi(m_1',\gamma,m_2'):=(m_1',(m_1',\gamma,m_2'),m_2').$
\end{examp}

Let us briefly recall the notion of \v{C}ech groupioid. Given a smooth manifold $M$ and an open cover $\{U_{i}\}$ we set the disjoint union $\coprod_{i}U_i$ to be objects. Let $U_{ij}$ stands for $U_i\cap U_j$, for every pair of indices $(i,j)$. We use the notation $x_{ij}$ for an element $x\in U_{ij}$ and see $U_{ij}$ as a subset of $M$. We set the disjoint union $\coprod_{i,j}U_{ij}$ to be the morphisms and we set the source and target as embeddings
\begin{equation*}
\begin{array}{c}
s,t:\coprod_{i,j}U_{ij}\to \coprod_{i}U_i\\
s(x_{ij})=x_i,
t(x_{ij})=x_j
\end{array}
\end{equation*}
and we define the multiplication to be $x_{ij}x_{jk}=x_{ik}$, where $x_{ij}\in U_{ij},x_{jk}\in U_{jk}$ are considering as the same element $x\in M$. The identity over an element $x_i\in U_i$ is $x_{ii}$. These dates define a Lie groupoid $\coprod_{i,j}U_{ij} \toto \coprod_{i}U_i$ which is called \v{C}ech groupioid.
\begin{examp}
Given a smooth manifold $M$ and  an open cover  $\{U_{i}\}$, the \v{C}ech groupoid $\coprod_{i,j}U_{ij} \toto \coprod_{i}U_i$ is the pull-back the Lie groupoid Lie groupoid $M \toto M$ via the trivial submersion $p(x_{ij})=x$ and hence Morita equivalent to the Lie groupoid $M \toto M$. The pair Lie groupoid $M \times M \toto M$ is strong Morita equivalent to a point. More generally, $P \times_{\varphi,M,\varphi} P \toto P$ is strong Morita equivalent to the trivial groupoid $M \toto M $ for every surjective submersion $\varphi: P \to M$.
\end{examp}


\begin{prop}\cite{BehrendXu,LaurentHabilitation}
Two Lie groupoids are weak Morita equivalent if and only if they are strong Morita equivalent.
\end{prop}
To state the next proposition, we need to make a small abuse of vocabulary.
Strictly speaking, equivalence relations are defined on sets, not on classes.
But of course, it can be made sense of equivalence relations on  a class:
it is just defined to be a category admitting the objects of the class as objects, and having one and only
one invertible arrow relating any two objects of the class that we wish ta make equivlent.
This allows to make sense of the next proposition:

\begin{prop}
Morita equivalence induces equivalence relation on the class of all Lie groupoids.
\end{prop}

By a differential stack, we mean a equivalence class of this equivalence relation on the class of all Lie groupoids:

\begin{defi}\cite{BehrendXu}
A differential stack is a Morita equivalence class of Lie groupoids.
\end{defi}

Morita equivalence preserve many properties of Lie groupoids.
For instance, it can be shown that for every pair $(\Gamma \toto M , \Gamma' \toto M')$ of Morita
 equivalent Lie groupoids:
  $$ 2{\rm dim}(M)-{\rm dim}(\Gamma) = 2{\rm dim}(M')-{\rm dim}(\Gamma')  ,$$
which allows to define the dimension of a stack as being the this number, computed
with the help of any representative of the differential stack.
Beyond this example is a general idea. In order to define something on differential stacks,
one defines it first on Lie groupoids, then we show that it is Morita invariant.
More precisely, if some set $R(\Gamma)$ is associated to Lie groupoids, if a Morita equivalence ${\mathcal M} $
between $\Gamma $ and $\Gamma' $ induces a map $f({\mathcal M}): R(\Gamma) \to R(\Gamma')$ in such a way that
$f({id})=id $ and $f({\mathcal M} \circ {\mathcal M}') = f({\mathcal M}) \circ f({\mathcal M})$, then we say that
the construction goes to the quotient at the differential stack level.


\subsection{\v{C}ech cohomology and non-Abelian cohomology}\label{sectionchech}

The purpose of this section is to introduce non-Abelian $1$-cocycle,
and to justify this notion by recalling first what \v{C}ech $2$-cocycles
are (since, unfortunately for the terminology, \v{C}ech $2$-cocycles are non-Abelian $1$-cocycles),
then we shall introduce twisted \v{C}ech $2$-cocycles, which is also a non-trivial example
of  non-Abelian $1$-cocycle. Then we will introduce Abelian $1$-cocycles.


For the \v{C}ech cohomology valued in a multiplicative Abelian group, $2$-cocycles are given by families
 $a : {\mathcal C}^\infty( \coprod_{i,j,k \in I}U_{ijk} , A) $
such that $$ a_{ijk}+ a_{ikl} - a_{ijl} - a_{jkl}=0$$

Let us now say a word on twisted \v{C}ech cocycles.
Assume $A$ is an additive Abelian group, and assume that we are given a $\HH$-principal bundle $P \to M$,
with $\HH$ a group acting on $A$ by Lie group morphisms of $A$, an action that we shall denote
with the help of the shorthand $ (h,a) \mapsto h(a) $ for all $h \in H, a \in A$.
Given an open cover $(U_i)_{i \in I}$ of $M$ and a cocycle $h_{ij}:U_{ij} \to H$
defining the principal bundle $P \to M$ and computed out of sections $\sigma_i: U_i \to P$ of $P \to M$, a complex can be obtained as follows:
chains are, as previously, maps (or continuous maps, or smooth maps, depending on the context)
from $ \coprod_{i_1,\dots,i_k} (U_{i_1 \dots i_k} $ to $ A $, and
the differential is given by
 $$ (\diff \omega)_{i_0,\dots,i_k} = h_{i_0i_1} (\omega_{i_1\dots,i_n}) +  \sum_{i=1 }^n (-1)^i \omega)_{\widehat{i_k}}$$
 with $ \widehat{i_k} = i_0, \dots,i_{k-1},i_{k+1},\dots,i_n$.
Equivalently, this can be understood as being the \v{C}ech complex
valued in $A$ on $P \to M $ associated to the cover $p^{-1}(U_i)$
are obtained restricting ourself to $H$-equivariant maps, maps that we identify with their pull-back through
$\sigma_i : U_i \to P$.
In particular, $2$-cocycles in this complex are given by families $a : C^\infty( \coprod_{i,j,k \in I}U_{ijk} , A) $
such that
$$ a_{ijk}+ a_{akl} - a_{ikl} = h_{ij}(a_{jkl})$$
where we recall that $ h_{ij} h_{jk} = h_{ik}$

For non-Abelian groups, \v{C}ech cohomology does not make sense, and defining it formally similar formulas is not relevant,
for it would not square to zero. However, it is possible, sometimes, to have \v{C}ech-like cocycles even when for non-Abelian groups.
For instance, $1$-cocycles still make sense, when defined as maps from $U_{ij} \to \G $ satisfying
 $$ \lambda_{ij} \lambda_{jk} = \lambda_{ik}$$
Identifying cocycles if and only if $ r_i \lambda_{ij} r_j^{-1} = \lambda_{ij}' $ with $r_i:U_i \to \G$, we get a notion of cohomology,
whose classes precisely describe $\G$ principal bundles.

More, generally, $1$-cocyles can be taken valued in a crossed-module.
A \emph{crossed module of Lie groups} (consult, for instance, \cite{BaezSchreiber}) is a quadruple $(\G, \HH, \rho,\jmath)$, where $\rho:\G \rightarrow \HH $
  and $\jmath : \HH \to \Aut (\G)$ are Lie group homomorphisms satisfying the next conditions, for all $ \fatg,\fatg' \in \G, h \in \HH$
\begin{enumerate}
\item $\rho \big(  h(\fatg) \big)=h \rho(\fatg) h^{-1}$
\item $ \rho (\fatg) \big(\fatg'\big)  =\fatg \fatg' \fatg^{-1} $
\end{enumerate}
with the understanding that $h(g) $, for every $ h\in \HH, g\in \G$ is a short hand for $\jmath (h) (g) $.
Notice that here we consider that the action of $\HH$ on $\G$ to be a left-action, which is not the usual convention, but is necessary to recover the formulas of the $\G \to \Aut (\G)$ case as they are stated in \cite{LaurentStienonXu,BL}.
 In order to avoid an easily done confusion between elements in $\G$ and in $ \HH$, we shall denote by bold letters,
 $\fatg,\fatg'$ elements of $\G$, and in ordinary letters $h,h'$ elements in $\HH $. Also, bold letters shall be used for
  $\G$-valued functions. Last, it is customary to denote a cross-module by $\G\stackrel{\rho}{\to}\HH$, forgetting to make explicit the morphism~$\jmath $.

\begin{examp}
\label{ex:Abelian}
$A \to 1$ with $A$ Abelian.
\end{examp}

\begin{examp}
\label{ex:twistedAbelian}
$ A \to H$ with $A$ Abelian and $H$ a Lie group that acts on  $A$ by group automorphisms.
\end{examp}

\begin{examp}
\label{ex:groups}
$ 1 \to H$ with $H$ an arbitrary Lie group.
\end{examp}

\begin{examp}
\label{ex:GtoAutG}
$ G \to Aut(G)$ with $G$ an arbitrary Lie group.
Theorem : For G a Lie group, Aut(G) is a Lie group.
\end{examp}


\begin{rem}
Cross-Modules of Lie groups induces crossed modules of Lie algebras.
\end{rem}

We now recall from \cite{Dedecker} the notion of non-Abelian $1$-cocycles valued in an arbitrary crossed module $\G \to \HH$.
\begin{defi}\label{def:nonab}
Let $\G\stackrel{\rho}{\to} \HH$ be a crossed module, and ${\mathcal U}=(U_i)_{i \in I} $ an open covering of a manifold $N$.
A non-Abelian $1$-cocycle w.r.t. ${\mathcal U} $ with values in $\G\to \HH$ is a pair $(\lambda, \gggg)  \in \mathcal C ^{\infty}( \coprod_{i,j \in I}U_{ij},\HH ) \times \mathcal C ^{\infty} (\coprod_{i,j,k \in J}U_{ijk},\G) $
 required to satisfy the following conditions:
   \begin{equation}\label{nonAbe}  \left\{
    \begin{array}{l}
   \rho(\gggg_{ijk})  \lambda_{ik} =  \lambda_{ij}
    \lambda_{jk}    \\
     \gggg_{ijk}\gggg_{ikl}=
    \lambda_{ij}(\gggg_{jkl})\gggg_{ijl}\\
    \gggg_{iik}=e
\end{array}   \right. \end{equation}
for all possible indices (here $\lambda_{ij}$ (resp.$\gggg_{ijk}$) stands for the restriction of $\lambda$ (resp. $\gggg$ ) to $ U_{ij} $ (resp.$ U_{ijk} $ )).
\end{defi}

\begin{rem} \label{rmk:lambdaii} Note that the first of the relations (\ref{nonAbe}), when $ i=j$, implies $\lambda_{ii}=\e $. Note that  we have used the notation $e$ for the nutral elements of both Lie groupoids $\G$ and $\HH$.
\end{rem}

\begin{examp}
\label{ex:cocycles_Abelian}
From the crossed-module of example \ref{ex:Abelian},
a Non-Abelian $1$-cocycle amounts to a \v{C}ech 2-cocycle.
\end{examp}

\begin{examp}
\label{ex:cocycles_twistedAbelian}
From the crossed-module of example \ref{ex:twistedAbelian},
a Non-Abelian $1$-cocycle amounts to a twisted Cech 2-cocycle.
$ A \to H$ with $A$ Abelian and $H$ a Lie group that acts on  $A$ by group automorphisms.
\end{examp}

\begin{examp}
\label{ex:cocycles_groups}
From the crossed-module of examp \ref{ex:groups}, a Non-Abelian $1$-cocycle amounts to a $H$-valued $1$-cocycle.
\end{examp}

It follows from the crossed module axioms that $\rho(\G)$ is a distinguished subgroup of $\HH$,
so that the quotient $ \HH/\rho(\G)$ is a group.  Every non-Abelian 1-cocycles induces an
$\HH/\rho (\G)$-principal bundle over the base manifold $M$ which we called the \emph{band}
of the non-Abelian $1$-cocycle.

\begin{rem}
When $\rho$ is injective. Then the second relation follows from the first one.
Moreover the first relation implies that the second one holds up to an element in the kernel of ~$\rho$.
\end{rem}











\section{Lie groupoids $\G \to \HH$-extensions}
\label{sec:GpdExt}

Let $\G \to \HH$ be a crossed-module of finite dimensional Lie groups. The purpose of this section is to give a complete description, purely in terms of Lie groupoids, of  $\G \to \HH $-gerbes over a given stack, and to check that, when the stack in question is simply a manifold $N$, our notion gives back an already known description \cite{BM,BaezSchreiber,Dedecker} in terms of non-Abelian cohomology.

\subsection{Definition of Lie groupoids $\G \to\HH$-extensions}
\label{subsec:defofextensions}

 In \cite{BL}-\cite{LaurentStienonXu} gerbes are described as Lie groupoids extensions (up to Morita equivalence of those). But this description mainly covers the case of the so-called $\G$-gerbes, i.e. gerbes valued in the crossed module $\G \to \Aut (\G)$. In order to describe, in purely Lie groupoid terms, $\G \to \HH$-gerbes, one needs to go further and to work with $\G$-extensions endowed with some $\HH$-principal bundle structure (up to Morita equivalence of those).

We first wish to introduce Lie groupoid extensions.

  \begin{defi}\cite{LaurentStienonXu}
  A Lie groupoid extension is a triple $(\mathcal{R},\GGG,\varphi)$, denoted by $\mathcal{R} \stackrel{\varphi}{\to}\GGG$ or simply by $\mathcal{R}\to \GGG$, when there is no risk of confusion,  where $\mathcal{R}\toto M$ and $\GGG \toto M$ are Lie groupoids over the (same) manifold $M$ and the map $\varphi: \mathcal{R}  \to \GGG $ is a groupoid morphism over the identity of $M$ such that $\varphi$ is a surjective submersion.
  \end{defi}

  The \emph{kernel} of a Lie groupoid extension $\mathcal{R} \stackrel{\varphi}{\to}\GGG$ is, by definition, the inverse image through $\varphi$ of the unit manifold of $\GGG$, \emph{i.e.} the set  $$K=\{r\in \mathcal{R} : \varphi(r) \in \epsilon(M) \}.$$ Since $\varphi$ is a surjective submersion, the kernel is a submanifold of $\mathcal{R}$. Also, since $\varphi$ is a groupoid homomorphism over the identity of $M$, $K$ is indeed a bundle of Lie groups (\emph{i.e.} it is a Lie groupoid whose source and target maps coincide). Notice that $K$ is normal in $\mathcal{R}$ in the sense that $r^{-1} \gpdproductR k \gpdproductR r \in K$ for all admissible $k \in K, r \in \mathcal{R}$.
  The previous assignment defines indeed a Lie groupoid action of $\mathcal{R} $ on $K \to M $, action that we shall denote by $\gpdactkernel $.

  \begin{defi}\cite{LaurentStienonXu}
   Let $\G$ be a Lie group. A Lie groupoid extension $\mathcal{R} \stackrel{\varphi}{\to}\GGG $ is called a Lie groupoid $\G$-extension if its kernel $K$ is locally trivial with typical fiber $\G$,\emph{i.e.} if every point $x\in M$ ($M$ being the base manifold of both $\mathcal{R}$ and $\GGG $) admits a neighborhood $U$ such that $K_{U}=\varphi^{-1}(\epsilon(U))$ is isomorphic to $\G\times U$.
  \end{defi}

  To a Lie groupoid $G$-extension $\mathcal{R} \stackrel{\varphi}{\to}\GGG$, we now associate an $\Aut(\G)$-principal bundle over the groupoid $ \mathcal{R}  \toto M$, called the band of the extension. We first recall the notion of $\HH$-principal bundle over a Lie groupoid. See \cite{LTX} for instance.

  \begin{defi}
  Let $\HH$ be a Lie group, and $ \mathcal{R}  \toto M$ a Lie groupoid.
  A principal $\HH$-bundle over $\mathcal{R}\toto M$ is an usual (right) principal $\HH$-bundle $P \stackrel{\pi}{\to} M$  together with a (left) action of the Lie groupoid $ \mathcal{R} \toto M$ on $P \stackrel{\pi}{\to} M$ such that the ${\mathcal R}$ and the $\HH$ actions commute, \emph{i.e.} denoting the action of the Lie groupoid $\mathcal{R}\toto M$ and the action of Lie group $\HH$ on $P \stackrel{\pi}{\to} M$, both by the same notation $\cdot$, then $(\gamma\cdot p)\cdot h=\gamma\cdot(p\cdot h)$, for all admissible $\gamma \in \mathcal R, p\in P, h\in \HH$.
  \end{defi}

    We also define morphisms between two principal bundles w.r.t.different groups over different Lie groupoids, as follow.

\begin{defi}    \label{defi:isomppalbundles}
A morphism from a principal $\HH$-bundle $P \stackrel{\pi}{\to} M$ over a Lie groupoid $\mathcal{R}\toto M$ to a  principal $\HH'$-bundle $P' \stackrel{\pi'}{\to} M'$ over a Lie groupoid $\mathcal{R'}\toto M'$ is triple a $(\Phi,\Psi,\jmath)$, where $\Phi : \mathcal R \rightarrow \mathcal R'$ is an morphism of Lie groupoids $\Psi: P \to P' $ is diffeomorphism and  and $\jmath: \HH \to \HH' $ be a Lie group morphism, such that:
    $$ \Psi ( \gamma \gpdactP p \cdot h ) = \Phi(\gamma) \gpdactP  \Psi (p) \cdot \jmath(h)   $$
    for all pair $(\gamma,p) \in \mathcal{R} \times_{t,M,\pi} P $ and all $h \in \HH$.
%
When $\mathcal{R}\toto M$ and $\mathcal{R'}\toto M'$ are identically the same Lie groupoid and the map $\Phi $ is the identity map, then the morphism $(\Phi,\Psi,\jmath)$ is called a \emph{morphism over the identity of ${\mathcal R} \toto M $}, and we simply denote it by the pair $(\Psi,\jmath)$.
  \end{defi}

The band \cite{Giraud,BM} is in general defined for the gerbe itself, but \cite{LaurentStienonXu} introduced a notion of band for Lie groupoid  $\G$-extensions that boils down to the band of the gerbe. By construction, it is the set of all Lie group morphisms from $\G$ to some fiber of its kernel. More precisely, let $\G$ be a Lie group and $\mathcal{R} \rightarrow \GGG$ be a $\G$-extension, we set
\begin{equation}
Band( \mathcal{R} \rightarrow \GGG ) :=\coprod_{m\in M} isom(\G,K_{m}).
\end{equation}
to be the set of all possible Lie group isomorphisms from $\G$ to some fiber $K$.
Recall from \cite{LaurentStienonXu} that
\begin{enumerate}
\item \label{actAutonBand}$Band( \mathcal{R} \rightarrow \GGG ) $ admits a natural manifold structure, for which the projection on $M$ is a smooth surjective submersion.
We let $Band_m( \mathcal{R} \rightarrow \GGG )$ stands for the fiber over a point  $m \in M $.
\item $\Aut (\G)$ acts (on the right) freely and transitively on the fibers of $Band( \mathcal{R} \rightarrow \GGG ) $ as follows: $b_{m}\cdot \rho:=b_{m}\circ \rho$, for  $\rho \in \Aut (\G), b_{m}\in Band_m( \mathcal{R} \rightarrow \GGG ) $.
\end{enumerate}
All these items together imply that $Band( \mathcal{R} \rightarrow \GGG ) \stackrel{\pi}{\to}M$, where $\pi$ is the obvious projection, is a (right) $\Aut(\G)$-principal bundle over the base manifold $M$. Moreover, $Band( \mathcal{R} \rightarrow \GGG ) $  is an $\Aut(\G)$-principal bundle over the Lie groupoid $\mathcal{R}\toto M$, when it is equipped with the left action of $\mathcal{R}\toto M$ on $Band( \mathcal{R} \rightarrow \GGG ) \stackrel{\pi}{\to}M$ defined by setting $r\gpdactband b_{m} $ to be the Lie group morphism from $\G$ to $K_{s(r)} $ given by
 \begin{equation}\label{def:def of gpdactband}
 g \mapsto r b_m(g) r^{-1},
 \end{equation}
 for all $r\in \mathcal{R}$ with $t(r)=m,b_m\in isom(\G,K_m)$


We now have all the tools required for defining the type of extension whose (to be defined in section \ref{sec:MoritaExt}) quotients shall define $\G \to \HH $-gerbes.

\begin{defi}\label{def: Azimi-extension}
Let $\G \stackrel{\rho}{\to} \HH$ be a crossed module, with action map $\jmath: \mathbf H \to \Aut (\G)$, and $\GGG \toto M $ a Lie groupoid. A $\G \to \mathbf H$-extension of $\GGG \toto M $ is a triple $(\mathcal{R} \rightarrow \GGG, P \to M, \chi)$, where:
 \begin{enumerate}
 \item $\mathcal{R} \rightarrow \GGG$ is a Lie groupoid $\G$-extension,
 \item  $P \to M $ is an $\HH$-principal bundle over the Lie groupoid $ \mathcal{R} \toto M$,
 \item $(\chi,\jmath)$ is a morphism over the identity of $ {\mathcal R} \toto M$ (see definition \ref{defi:isomppalbundles}) from the $\HH$-principal bundle $P \to M $ to the $\Aut(\G)$-principal bundle $Band( \mathcal{R} \rightarrow \GGG ) $,
 \end{enumerate}
      such that, for all $p\in P, g \in \G$:
  \begin{equation}
  \label{eq:chi} \begin{array}{rcl}
    p \cdot \rho (g) & =  & \chi(p) (g) \gpdactP p  \end{array}
    \end{equation}
    (recall that $\chi(p)$ belongs to $Band_{\pi(p)}(  \mathcal{R} \rightarrow \GGG)= Isom(\G,K_{\pi(p)}) $, so that  $\chi(p) (g)$ is an element in $ K_{\pi (p)} \subset {\mathcal{R}}$: it makes therefore sense to let it act on $p \in P$).
\end{defi}

It shall be convenient to draw the following diagram in order to represent $\G \to \HH$-extensions.
Below, it shall be understood that an arrow of the type $\mathcal{R} \xymatrix{  \ar@{{}{--}^{)}}[r] & P}$ means that
the groupoid $\mathcal{R} $ acts on $P$.
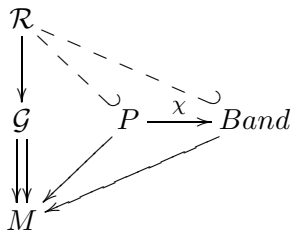
$$\xymatrix{  {\mathcal R} \ar@{{}{--}^{)}}[dr] \ar@{{}{--}^{)}}[drr]  \ar[d]  &   \\ {\GGG} \ar@<2pt>[d] \ar@<-2pt>[d] & \ar[dl] P \ar[r]^-{\chi }  & \ar[dll] Band \\  M\\} $$

Let $\G \to \HH$ be a crossed module, by an \emph{isomorphism between two $\G \to \HH$-extensions of a Lie groupoid $\GGG \toto M $}, namely
$(\mathcal{R} \rightarrow \GGG \toto M, P \to M, \chi)$
and
$(\mathcal{R}' \rightarrow \GGG \toto M, P' \to M, \chi')$, we mean an isomorphism $(\Phi,\Psi,id_\HH)$ of principal bundles over Lie groupoids (see definition \ref{defi:isomppalbundles}) such that the following diagram commutes:
\begin{equation}\label{diagram for isomorphism} \begin{CD}
P @>\Psi>> P'\\
@VV\chi V @VV\chi'V\\
Band(\mathcal{R}\rightarrow\GGG ) @>\bar{\Phi}>> Band(\mathcal{R'}\rightarrow\GGG)
\end{CD}
\end{equation}
where $\bar{\Phi}(\eta)(g)=\Phi(\eta(g))$, for $\eta \in Band(\mathcal{R}\rightarrow\GGG ), g \in \G $. For such an isomorphism we use the notation $(\Phi,\Psi)$ instead of $(\Phi,\Psi,id_\HH)$.




\begin{examp}
Given a \v{C}ech 2-cocycle, we have a Lie groupoid extension of the \v{C}ech groupoid defined as
a $1 \to \mathbf H$-extension of the \v{C}ech groupoid.
\end{examp}

\begin{examp}
When the crossed module is simply $ \{1\} \to \HH $, then $\G \to \HH$-extensions are nothing than $\HH$-principal bundles over Lie groupoids, and isomorphisms of $\G \to \HH$-extensions amount to isomorphisms of those.
\end{examp}

\begin{examp}\label{ex:G_extensions_as_Azimi_extensions}
For every $\G$-extension $\mathcal{R} \rightarrow \GGG \toto M$, the quadruple
 \begin{equation}\label{eq:quadruple} (\mathcal{R} \rightarrow \GGG , Band(\mathcal{R} \rightarrow \GGG  ) \to M, \id_{Band(\mathcal{R} \rightarrow \GGG )}) \end{equation}
 is a $\G \to \Aut (\G)$-extension.
Conversely, when the crossed module $\G \to \HH$ is $ \G \to \Aut (\G) $, then, for every $\G \to \HH$-extension $(\mathcal{R} \rightarrow \GGG , P \to M, \chi)$ , then $(\chi,Id_{\Aut (\G)})  $ is an isomorphism of principal bundles over the identity of  $\mathcal{R} \toto M$.
  In conclusion, the assignment of (\ref{eq:quadruple}) induces a one to one correspondence between $ \G \to \Aut (\G) $-extensions and $\G$-extensions. This correspondence is an equivalence of categories, for isomorphisms of $\G \to \Aut(\G)$-extension amount to isomorphisms of the corresponding $\G$-extension.
\end{examp}

 These examples lead to the guess that non-Abelian $1$-cocycle should induce crossed-modules extension.
 This is the content of the next section.

\subsection{The manifold case: $\G \to \HH $-valued non-abelian cocycles as $\G \to \HH$-extensions over  Lie groupoid}
\label{subsec:mfdcase}

Throughout the present section, we shall fix an open covering ${\mathcal U}= (U_i)_{i \in I}$ of a manifold~$N$.

 Our purpose is to show that $\G \to \HH$-extensions of the  \v{C}ech groupoid $\Cechs$ correspond to $\G \to \HH$-valued non-Abelian $1$-cohomology, computed with respect to ${\mathcal U}$. We first recall the notion of non-Abelian $1$-cocycles \cite{Dedecker,BM}, as introduced by Dedecker. Then, we show that these are in one to one correspondence with (a certain set of) $\G \to \HH$-extensions of the  \v{C}ech groupoid $\Cechs$. Proving that $\G \to \HH$-coboundaries correspond to isomorphisms of these extensions shall then yield to the desired conclusion.

\begin{defi} \label{definition of Adapted extension}
An adapted Lie groupoid $\G \to \HH$-extension of the \v{C}ech groupoid $\Cechs $ is a $\G \to \HH$-extension  $({\mathcal R} \stackrel{\varphi}{\to}  \Cechs, P \to \coprod_{i \in I} U_{i}, \chi) $
on which we impose the following constraints:
\begin{enumerate}
  \item ${\mathcal R}$ is the space $ \G \times  \coprod_{i,j \in I} U_{ij}$ and $ \varphi$ is the projection onto the second component;
    \item $P$ is the space $  \coprod_{i \in I} U_{i} \times \HH $, equipped with the trivial right $\HH$-action $ (x_i,h) \cdot h'= (x_i,hh')$ for all $h,h' \in \HH, x\in U_{i}$;
          \item The map $\chi : P \to Band ( {\mathcal R} \stackrel{\phi}{\to}  \Cechs)$ maps $ (x_i,h) \in P $ to the element of the band over $x_i$ given by $  g \mapsto (h (g),x_{ii}) $ for all $g \in \G, h \in \HH, x \in U_i $;
   \item the Lie groupoid product $\gpdproductR$ of ${\mathcal R}$ satisfies the relation $(g,x_{ii}) \gpdproductR (g',x_{ij}) =
     (gg', x_{ij}) $ for all~$x \in U_{ij}, g,g' \in \G, i,j\in I$.
    \end{enumerate}
\end{defi}
  Items 1 and 4 of the definition imply that the kernel of an adapted Lie groupoid $\G \to \HH$-extension of the \v{C}ech groupoid $\Cechs $ is the trivial bundle of group:
$K=\G \times \coprod_{i \in I} U_{ii}\simeq \G \times \coprod_{i, \in I} U_{i}$
so that the band $Band(\mathcal R \rightarrow \Cechs)$ is canonically isomorphic to
 $\Aut(\G) \times \coprod_{i \in I} U_{i} $.

 Given manifold $N$, an open covering ${\mathcal U}= (U_i)_{i \in I}$ of a manifold~$N$ and a crossed module $\G \to \HH $, adapted Lie groupoid $\G \to \HH$-extensions of the same \v{C}ech groupoid $\Cechs $ may only differ by two things: the Lie groupoid product on  $\mathcal{R}=\G \times \coprod_{i,j \in I} U_{ij}$ and the action of $ \mathcal{R}=\G \times \coprod_{i,j \in I} U_{ij}$ on $P:= \coprod_{i \in I} U_i \times \G $.

\begin{nota}\label{not:adaptextensions}
For $ {\mathcal U}$ an open covering of $N$, we shall denote adapted Lie groupoid $\G \to \HH$-extensions of $\Cechs $ as triples $({\mathcal U},\bullet, \star)$, where $\mathcal U$ refers to the open covering, $\bullet $ refers to the multiplication of the Lie groupoid $\mathcal{R}:= \G \times \coprod_{i,j \in I} U_{ij}$ and $\star$ refers to the action of $\mathcal{R}:= \G \times \coprod_{i,j \in I} U_{ij}$ on the principal bundle $P:=  \coprod_{i \in I} U_i \times \G $.
\end{nota}


\begin{rem}\label{rmk:remark_on_adapted}
For an adapted extension $(\mathcal U,\bullet, \star)$, the action of an element $(g, x_{ii}) $  in the kernel of ${\mathcal R} \stackrel{\phi}{\to}  \Cechs$ on an admissible element $(x_i,h) \in P$ is given by $(x_i,\rho(g)h)$. We prove it as follows.
First notice that
\begin{equation}
\begin{array}{rcll}
(x_{i},h) \cdot \rho(g) & = & (x_{i},h\rho(g)) & \hbox{by definition \ref{definition of Adapted extension}, item 2}\\
                 & = & (x_{i},h\rho(g)h^{-1}h) & \hbox{}\\
                 & = & (x_{i},\rho(h(g))h) & \hbox{by axioms of crossed module.}
\end{array}
\end{equation}
On the other hand:
\begin{equation}
\begin{array}{rcll}
(x_{i},h) \cdot \rho(g)
                 & = & \chi (x_{i}, h)(g) \star (x_{i},h) &\hbox{by (\ref{eq:chi}) in definition \ref{def: Azimi-extension}}\\
                 & = & (h(g),x_{ii}) \star (x_{i},h) & \hbox{by definition \ref{definition of Adapted extension}, item 3.}
\end{array}
\end{equation}
The result follows by substituting $h(g) $ by $g $ in the previous relations.
\end{rem}

We now prove the desired correspondence which generalizes \cite{BL}(recall that $\jmath: \HH \to \Aut (\G)$ is part of the crossed module structure, see section \ref{sectionchech}).

\begin{prop}\label{prop:cocycles}
Let $ {\mathcal U}=(U_i)_{i \in I}$ be an open covering of a manifold $N$,
and $\G \to \HH$ a crossed modules of Lie groups.
\begin{enumerate}
\item Let $(\mathcal U,\bullet, \star)$ be an adapted Lie groupoid $\G \to \HH$-extension of \v{C}ech groupoid $N[\mathcal U]$. We define $(\lambda,\gggg) \in  \mathcal C ^{\infty}( \coprod_{i,j \in I}V_{ij},\HH ) \times \mathcal C ^{\infty} (
\coprod_{i,j,k \in J}V_{ijk},\G) $ gluing together the family of maps $\lambda_{ij}:U_{ij} \to \HH$ and $\gggg_{ijk}:U_{ijk} \to \G $ defined by
    \begin{equation}
  \label{eq:cons}
      \left\{ \begin{array}{rcll}

                        (\e,x_{ij}) \star (x_{j},\e) &=& ( x_{i}, \lambda_{ij} ) & \forall i,j\in I, \forall x\in U_{ij} \\
                      (\e,x_{ij}) \bullet  (\e,x_{jk}) &=& (\gggg_{ijk}, x_{ik}) & \forall i,j,k\in I \forall x\in U_{ijk} .\\
       \end{array} \right.
  \end{equation}
  Then $(\lambda_, \gggg)$ is a non-Abelian $1$-cocycle.


 \item Given a non-Abelian $1$-cocycle  $(\lambda ,\gggg)$, we define:
   \begin{enumerate}
   \item a Lie groupoid structure $\bullet$ on $\mathcal R =\G \times \coprod_{i,j\in I} U_{ij} $ by:
     \begin{equation}\label{eq:recons}  \left\{ \begin{array}{rcl} (g,x_{ij}) \bullet (g',x_{jk}) &:=&
     (g\lambda_{ij}(g')\gggg_{ijk}, x_{ik}) \\ (g,x_{ij})^{-1} &:=& (
     \lambda_{ij}^{-1}(g^{-1}\gggg_{iji}^{-1}) , x_{ji})
 \end{array} \right. \end{equation}
 for all $g, g' \in \G, i,j \in I, x \in U_{ij}$, where $\lambda_{ij}(g')$ stands for $\jmath (\lambda _{ij})(g')$,
   \item a map $ \phi: \mathcal R \to N[\mathcal U] $ given by $ (g,x_{ij}) \mapsto x_{ij} $, for all $ g \in \G ,i,j \in I , x \in U_{ij} $,
   \item a structure of $H$-principal bundle $\star$ on $P:=\coprod_{i\in I} U_{i} \times \HH $ over the Lie groupoid $\mathcal R \toto \coprod_{i\in I} U_{i} $ given by
  \begin{equation}\label{defi: ppalstronP}
(g,x_{ij}) \star (x_{j},h)=(x_{i},\rho(g)\lambda_{ij} h),
\end{equation}
for all $g\in \G, h\in \HH, i,j\in I, x\in U_{ij}$,
   \item a map $\chi: P \to Band(\mathcal R \to N[\mathcal U])$ by $(x_{i},h)\mapsto (\jmath(h),x_{i}).$
   \end{enumerate}
   Then $(\mathcal U,\bullet , \star )$ is an adapted Lie groupoid $\G \to \HH$-extension of the \v{C}ech groupoid $N[\mathcal U]$

   \item the procedures in items 1 and 2 are inverse to each other.


 \end{enumerate}
 \end{prop}

The proof will go through a lemma.

\begin{lem} \label{key lemma}
Let $(\mathcal U,\bullet, \star) $ be an adapted Lie groupoid $\G \to \HH$-extension of the \v{C}ech groupoid $\Cechs $. Define the maps $\lambda_{ij}:U_{ij} \to \HH$ and $\gggg_{ijk}:U_{ijk} \to \G $ as in (\ref{eq:cons}). Then the following relation holds for all $i,j \in I, g\in \G$ and $x \in U_{ij}$:
 $$(e,x_{ij}) \bullet (g,x_{jj})=(\lambda_{ij}(g), x_{ij})$$
\end{lem}
\begin{proof}
First observe that,
\begin{equation*}
\begin{array}{rcll}
\chi((e,x_{ij}) \star (x_{j},e))(g) &=& \chi(x_{i},\lambda_{ij})(g) & \hbox{by (\ref{eq:cons}), i.e. definition of $\lambda_{ij}$}\\
                                  &=& (\lambda_{ij}(g),x_{ii}) & \hbox{by definition \ref{definition of Adapted extension} item 3},
\end{array}
\end{equation*}
for all $i,j \in I, g\in \G$ and $x \in U_{ij}$. On the other hand,
\begin{equation*}
\begin{array}{rcll}
                                 && \chi((e,x_{ij})\star (x_{j},e))(g)\\
                                  &=& ((e,x_{ij}) \gpdactband \chi(x_{j},e)) \, (g) & \hbox{$\chi$ is a morphism of ppal bundles over grpds}\\
                                  &=& (e,x_{ij}) \bullet ( g,x_{jj})\bullet (e,x_{ij})^{-1} &\hbox{by def of $\gpdactband$ i.e. (\ref{def:def of gpdactband})}
\end{array}
\end{equation*}
for all $i,j \in I, g\in \G$ and $x \in U_{ij}$.
%
Multiplying on the right of both sides of the last two relations by $ (e,x_{ij}) $ and using item $4$ in definition \ref{definition of Adapted extension} yield the desired relation.
%
\end{proof}

\begin{proof}(of proposition \ref{prop:cocycles}).\emph{1}) We first prove that the maps defined in item $1$ form a non-Abelian $1$-cocycle. The relation $\gggg_{iij}=e$ is obtained by putting $i=j$ in the second relation of (\ref{eq:cons}).
 By definition of groupoid action, we have
 \begin{equation}\label{associativity}
 ((e,x_{ij})\bullet (e,x_{jk}))\star (x_{k},e)=(e,x_{ij})\star((e,x_{jk})\star(x_ k,e))
\end{equation}
for all indices $i,j,k$ and all $x\in U_{ijk}$. The LHS of (\ref{associativity}) gives:
\begin{equation*}
\begin{array}{rcll}\label{LHS}
\hbox{LHS of} \,(\ref{associativity}) & = & (\gggg_{ijk},x_{ik})\star(x_{k},e) & \hbox{by (\ref{eq:cons}), i.e. def. of $\gggg_{ijk}$}              \\
                           & = & ((\gggg_{ijk},x_{ii}) \bullet (e,x_{ik}))\star(x_{k},e)  & \hbox{by def. \ref{definition of Adapted extension}, item 4}\\
                            & = & (\gggg_{ijk},x_{ii}) \star ((e,x_{ik}))\star(x_{k},e)  & \hbox{by axioms of groupoid action}\\
                           & = & (\gggg_{ijk},x_{ii}) \star (x_{i},\lambda_{ik}) &\hbox{by (\ref{eq:cons}), i.e. def. of $\lambda_{ik}$}  \\
                           & = & (x_{i},\rho(\gggg_{ijk})\lambda_{ik}) & \hbox{by remark \ref{rmk:remark_on_adapted}}
\end{array}
\end{equation*}
while the RHS of (\ref{associativity}) gives:
\begin{equation*}
\begin{array}{rcll}\label{RHS}
\hbox{RHS of}\, (\ref{associativity}) & = &(e,x_{ij})\star (x_{j},\lambda_{jk}) & \hbox{by (\ref{eq:cons}), i.e. def. of $\lambda_{jk}$} \\
 & = & (e,x_{ij})\star(x_{j},e)  \lambda_{jk} & \hbox{by def. \ref{definition of Adapted extension}, item 2,} \\
                                       & = & (x_{i},\lambda_{ij})\lambda_{jk} & \hbox{by def. of $\lambda_{ij}$} \\
                                       & = & (x_{i},\lambda_{ij}\lambda_{jk}) & \hbox{by (\ref{eq:cons}), i.e. def. \ref{definition of Adapted extension}, item 2}
\end{array}
\end{equation*}
Comparing these relations, we obtain the first condition of (\ref{nonAbe}).
 To show that the families $(\lambda_{ij})_{i,j \in I}$ and $(\gggg_{ijk})_{i,j,k \in I}$ satisfy the second condition of (\ref{nonAbe}), we write the associativity condition of the Lie groupoid multiplication of $\mathcal{R} $ as follows:
\begin{equation}\label{associativity 2}
((e,x_{ij})\bullet(e,x_{jk}))\bullet(e,x_{kl})=(e,x_{ij})\bullet((e,x_{jk})\bullet(e,x_{kl}))
\end{equation}
for all indices $i,j,k,l \in I$ and $x\in U_{ijkl}$. The LHS of (\ref{associativity 2}) amounts to:
\begin{equation*}
\begin{array}{rcll}
\hbox{LHS of } (\ref{associativity 2})& = & (\gggg_{ijk},x_{ik})\bullet(e,x_{kl}) & \hbox{by (\ref{eq:cons}), i.e. definition of $\gggg_{ijk}$}\\
                                     & = & ((\gggg_{ijk},x_{ii})\bullet(e,x_{ik}))\bullet(e,x_{kl}) & \hbox{by definition \ref{definition of Adapted extension}, item 4}\\
                                     & = & (\gggg_{ijk},x_{ii})\bullet((e,x_{ik})\bullet(e,x_{kl})) & \hbox{(by associativity of the gpd product)}\\
                                     & = & (\gggg_{ijk},x_{ii})\bullet(\gggg_{ikl},x_{il})   & \hbox{by (\ref{eq:cons}), i.e. definition of $\gggg_{ikl}$ }\\
                                     & = & (\gggg_{ijk}\gggg_{ikl},x_{il})            & \hbox{by definition \ref{definition of Adapted extension}, item 4},
\end{array}
\end{equation*}
while the RHS of (\ref{associativity 2}) gives
 \begin{equation*}
 \begin{array}{rcll}
  \hbox{RHS of } (\ref{associativity 2}) & = & (e,x_{ij})\bullet(\gggg_{jkl},x_{jl}) & \hbox{by  (\ref{eq:cons}), i.e. definition of $\gggg_{jkl}$ } \\
                                       & = & (e,x_{ij})\bullet(\gggg_{jkl},x_{jj}) \bullet (e,x_{jl})  &\hbox{by definition \ref{definition of Adapted extension}, item 4}\\
                                       & = & (\lambda_{ij}(\gggg_{jkl}),x_{ij}) \bullet (e,x_{jl})               & \hbox{by lemma \ref{key lemma}}\\
                                       & = & (\lambda_{ij}(\gggg_{jkl}),x_{ii}) \bullet (e,x_{ij})\bullet(e,x_{jl})     &\hbox{by definition \ref{definition of Adapted extension}, item 4}\\
                                       & = & (\lambda_{ij}(\gggg_{jkl}),x_{ii}) \bullet (\gggg_{ijl},x_{il})     & \hbox{by (\ref{eq:cons}), i.e. definition of $\gggg_{ijl}$ }\\
                                       & = & (\lambda_{ij}(\gggg_{jkl})\gggg_{ijl},x_{il})            & \hbox{by definition \ref{definition of Adapted extension}, item 4}
 \end{array}
 \end{equation*}
Comparing these relations, we obtain the second condition of (\ref{nonAbe}), which completes the proof of the first item.

\bigskip
\noindent
\emph{2}   We need to check that the multiplication\, $\bullet$ defined in (\ref{eq:recons}) is a Lie groupoid multiplication. We first prove the associativity.
Let $i,j,k \in I, x \in U_{ijkl}$ and $g,g',g'' \in \G$. Then
\begin{equation*}
\begin{array}{rcll} & & (g,x_{ij})\bullet ((g',x_{jk}) \bullet(g'',x_{kl})) & \\
 & = & (g,x_{ij}) \bullet (g'\lambda_{jk}( g'')\gggg_{jkl},x_{jl})   &  \hbox{by (\ref{eq:recons}), i.e. def. of $\bullet$}\\
                                    & = & (g\lambda_{ij}(g'\lambda_{jk}( g'')\gggg_{jkl})\gggg_{ijl}, x_{il}) & \hbox{by (\ref{eq:recons}), i.e. def. of$\bullet$}\\
                                    & = & (g\lambda_{ij}( g')\lambda_{ij}(\lambda_{jk}( g''))\lambda_{ij}(\gggg_{jkl})\gggg_{ijl}, x_{il}) & \hbox{by crossed-modules axioms}\\
                                    & = & (g\lambda_{ij}( g')(\rho(\gggg_{ijk})\lambda_{ik}( g''))\lambda_{ij} (\gggg_{jkl})\gggg_{ijl},x_{il}) & \hbox{by (\ref{nonAbe}) in definition \ref{def:nonab}}\\
                                    & = & (g\lambda_{ij}( g')\gggg_{ijk}\lambda_{ik}(g'')\gggg_{ijk}^{-1}\lambda_{ij} (\gggg_{jkl})\gggg_{ijl},x_{il}) &\hbox{by crossed module axioms}\\
                                    & = & (g\lambda_{ij} (g')\gggg_{ijk}\lambda_{ik}(g'') \gggg_{ikl},x_{il}) & \hbox{by (\ref{nonAbe}) in definition \ref{def:nonab}}\\
                                    & = & (g\lambda_{ij}( g')\gggg_{ijk},x_{ik})\bullet(g'',x_{kl}) & \hbox{by (\ref{eq:recons}), i.e. def. of $\bullet$}\\
                                    & = & ((g,x_{ij}) \bullet (g',x_{jk})) \bullet (g'',x_{kl}) & \hbox{by (\ref{eq:recons}), i.e. def. of $\bullet$}.
\end{array}
\end{equation*}
 It is routine to check that the henceforth defined multiplication admits the source (resp. target) map $s, (\hbox{resp.}t ): \mathcal{R} \to \coprod_{i \in I} U_i$ given by $(g,x_{ij})\mapsto x_i$ (resp $x_j$), admits the map $\epsilon : \coprod_{i \in I} U_i \to \mathcal R$ given by $x_i \mapsto (e,x_{ii})$ as an unit map, and admits an inverse given as in (\ref{eq:recons}). Altogether, these structural maps endow $\mathcal R$ with a structure of Lie groupoid, and eventually turn $\mathcal R\stackrel{\phi}{\to} \coprod_{i \in I} U_{ij}$ into a Lie groupoid $\G$-extension. It is also routine to check that (\ref{defi: ppalstronP}) gives a structure of principal $\HH$-bundle over the Lie groupoid $ \mathcal R \toto \coprod_{i\in I} U_{i} $. In order to check that
 $$(\mathcal R \to N[\mathcal U], P \to \coprod_{i\in I} U_{i}, \chi)$$
 is a $\G \to \HH$-extension, we are left with the task of showing that $(\chi, \jmath)$ is a morphism of principal-bundles over the identity of $\mathcal R$.
 One condition is obvious:
 $$ \chi((x_i,h)\cdot h')=\chi(x_i,hh')=(\jmath(hh'),x_i)=(\jmath(h)\jmath(h'),x_i)=(\jmath(h),x_i)\jmath(h')=\chi(x_i,h)\jmath(h')$$
  while the following proves that $p \cdot \rho(g)= \chi(p)(g)\star p$ for all $p\in P, g\in \G$, hence proves the claim:

 \begin{equation*}
 \begin{array}{rcll}
                                       & &\chi((x_i,h)\cdot h')(g)\star (x_i,h)&\hbox{}\\
                                       & = & (\jmath(h)(g),x_{ii})\star (x_i,h)& \hbox{by def. of $\chi$}\\
                                       & = & (h(g),x_{ii})\star (x_i,h) & \hbox{}\\
                                       & = & (x_i,\rho(h(g))\lambda_{ii}h) & \hbox{by (\ref{defi: ppalstronP}})\\
                                       & = & (x_i,h\rho(g)h^{-1}h)& \hbox{by crossed module axiom}\\
                                       & = & (x_i,h)\cdot \rho(g).&
 \end{array}
 \end{equation*}
 Now items 1-3 of definition \ref{definition of Adapted extension} hold by construction and item 4 holds because $\gggg_{iij}$ is assumed, in definition \ref{def:nonab}, to be equal to the neutral element $\e$ of $\G$. This completes the proof of the second item.


  \bigskip
\noindent
\emph{3})
   Next, we prove that items 1 and 2 in the proposition yield constructions which are inverse one to the other. For this purpose, we first notice that (\ref{eq:recons}) and (\ref{defi: ppalstronP}) hold for any adapted Lie groupoid $\G \to \HH$-extension, hence the construction of item 2 is injective. Assume that we are given a $\G \to \HH$-valued non-Abelian $1$-cocycle $(\lambda_{ij}, \gggg_{ijk})_{i,j,k \in I}$, then applying the procedure in item 2 we obtain an adapted Lie groupoid $\G \to \HH$-extension, to which we apply the construction in item $1$ to yield a  $\G \to \HH$-valued non-Abelian $1$-cocycle $(\lambda'_{ij}, \gggg'_{ijk})_{i,j,k \in I}$. We need to show that these two non-Abelian $1$-cocycles are equal. For this, observe that, by construction in item 2, we have $(x_i, \lambda'_{ij})= (e,x_{ij}) \star (x_j, e)$ while it follows from item 1 that $(e,x_{ij}) \star (x_j, e)=(x_i,\rho(e)\lambda_{ij}e)$. These two relations together prove that $\lambda_{ij}=\lambda'_{ij}$ for all $i,j \in I$. A  similar argument proves that $\gggg_{ijk}=\gggg'_{ijk}$, hence the claim. This implies that if two adapted Lie groupoid $\G \to \HH$-extensions of the \v{C}ech groupoid have the same $\G \to \HH$-valued non-Abelian $1$-cocycles associated with, they are equal. This proves the claim.
\end{proof}

Having made explicit a one to one correspondence between adapted $\G \to \HH$-extensions and non-Abelian $1$-cocycles, we now prove that, under this correspondence, isomorphisms of adapted $\G \to \HH$-extensions correspond to non-Abelian coboundaries, a notion that we now introduce, following \cite{BM} and\cite{Dedecker}.

\begin{defi}\label{defi:cobound}
Let $ {\mathcal U}=(U_i)_{i \in I}$ be an open covering of a manifold $N$
and $\G \to \HH$ a crossed module of Lie groups.
 A $\G \to \HH$-valued $1$-coboundary is a pair $(r,\vvvv) \in  \mathcal C ^{\infty}( \coprod_{i,j \in I}U_{ij},\HH ) \times \mathcal C ^{\infty} (
\coprod_{i,j,k \in J}U_{ijk},\G) $. We say that a $\G \to \HH$-valued $1$-coboundary  $ (r,\vvvv) $  relates two non-Abelian $1$-cocycles  $(\lambda,\gggg) $ and $(\lambda',\gggg') $ if
\begin{equation}\label{eq:relate}   \left\{ \begin{array}{rcll}
                   \lambda_{ij}'      &=&  \rho(\vvvv_{ij})   r_i
                    \lambda_{ij}  r_j^{-1} & (*) \\
\gggg_{ijk}'  \vvvv_{ik} &=& \lambda_{ij}' (\vvvv_{jk}) \vvvv_{ij}
                   r_i(\gggg_{ijk})  & (**) \\ \end{array}
\right. \end{equation}
for all possible indices. We  recall that $r_i, \vvvv_{ij}$ stand for the restriction of non-Abelian $1$-coboundary  $(r,\vvvv)$ to the intersection $U_{ij}$.
\end{defi}

The next proposition relates coboundaries and isomorphisms of adapted extensions which generalizes the results of \cite{BL} to arbitrary crossed-modules.

\begin{prop}\label{prop:coboundaries}  Let $(\mathcal U, \bullet,\star )$ and $(\mathcal U, \bullet',\star')$ be two adapted Lie groupoid $\G \to \HH$-extensions of $\Cechs$.
Let $(\lambda,\gggg) $ and $(\lambda',\gggg')$ be the $\G \to \HH$-valued  non-Abelian $1$-cocycles w.r.t. ${\mathcal U}$
  associated with the adapted Lie groupoid $(\mathcal U, \bullet,\star )$ and $(\mathcal U, \bullet',\star')$, respectively,(as in proposition \ref{prop:cocycles}).
Then the following construction defines a one to one correspondence between
 the set of isomorphisms of Lie groupoid $\G \to \HH$-extensions of $\Cechs$
from $(\mathcal U, \bullet,\star )$
to $(\mathcal U, \bullet',\star' )$,
and the set of $\G \to \HH$-valued $1$-coboundaries relating $(\lambda,\gggg) $ and $(\lambda',\gggg')$:
\begin{enumerate}
\item
Given an isomorphism $(\Phi_{\mathcal R}, \Phi_P) $ of adapted  Lie groupoid $\G \to \HH$-extensions of $\Cechs$ between $(\mathcal U, \bullet,\star )$ and $(\mathcal U, \bullet',\star')$, we define $r_i: U_i \to \HH$ and $\vvvv_{ij}: U_{ij} \to \G$ by:
  \begin{equation}\label{eq:isomclasses} \begin{array}{rcl} ( x_{i},r_i) &=&  \Phi_P ( x_{i},e)
  \\ (\vvvv_{ij}^{-1}, x_{ij})  &=& \Phi_{\mathcal R} (e, x_{ij})  \end{array} \end{equation}

\item Given a $\G \to \HH$-valued $1$-coboundary $(r , \vvvv) $ such that relates the non-Abelian $1$-cocycles
$(\lambda,\gggg) $ and $(\lambda',\gggg') $, define an isomorphism of $\g\to \HH$-extensions $(\Phi_{\mathcal{R}},\Phi_{P})$
between the corresponding adapted Lie groupoid $\G$-extensions $(\mathcal U, \bullet,\star )$ and $(\mathcal U, \bullet',\star')$ as follows
  \begin{equation} \label{eq:def_PhiX} \Phi_\mathcal R (g,x_{ij}) =   (r_i (g) \vvvv_{ij}^{-1}, x_{ij} ) \quad \quad \mbox{for all}\,\, i,j \in I,x \in U_{ij}, g \in \G.\end{equation}
  and an isomorphism between $P$ and $P'$ by:
    \begin{equation} \label{eq:def_PhiP}  \Phi_P (x_{i},h) = (x_i ,  r_i h),\quad \quad \mbox{for all}\,\, i\in I,x \in U_{i}, h \in \HH.\end{equation}
  \end{enumerate}
\end{prop}
\begin{proof}
1) First we prove that by following the construction in item 1, for a given  isomorphism of Lie groupoid $\G \to \HH$-extensions $(\Phi_{\mathcal R},\Phi_P, id_{\HH})$ between adapted Lie groupoid $\G \to \HH$-extensions $(\mathcal U, \bullet,\star )$
and $(\mathcal U, \bullet',\star' )$,  we obtain a $\G \to \HH$-valued $1$-coboundary. For this we need
 to prove that the pair $(r , \vvvv)$ obtained as in (\ref{eq:isomclasses})
satisfy relations (\ref{eq:relate}). We first prove the first of those relations, by
exploiting the fact that $(\Phi_{\mathcal R},\Phi_P, id_{\HH})$ is a morphism of principal bundles over Lie groupoids, which amounts to:
\begin{equation}\label{eq:cobound1}  \Phi_P (  (e,x_{ij}) \star (x_j,e) )  =
\Phi_\mathcal R (  (e,x_{ij}))  \star' \Phi_P ( (x_j,e) ), \quad \quad \mbox{for all}\,\, i,j \in I , x \in U_{ij}
\end{equation}
The LHS of (\ref{eq:cobound1}) is given by
\begin{equation*}
\begin{array}{rcll}
\Phi_P (  (e,x_{ij}) \star (x_j,e) ) & = & \Phi_P (x_{i},\lambda_{ij})  & \hbox{by (\ref{eq:cons}), i.e. definition of $\lambda_{ij}$}\\
                                              & = & \Phi_P (x_{i},e) \cdot \lambda_{ij} & \hbox{$\Phi_P$ being a $\HH$-ppal bundle morphism}\\
                                              & = & (x_i,r_i) \cdot \lambda_{ij} & \hbox{by (\ref{eq:isomclasses}), i.e. definition of $r_i$}\\
                                              & = & ( x_i , r_i\lambda_{ij}) & \hbox{by definition \ref{definition of Adapted extension}, item 2.} \\
 \end{array}
 \end{equation*}
The RHS of (\ref{eq:cobound1}) is given by
\begin{equation*}
\begin{array}{rcll}
                             & & \Phi_\mathcal R   (e,x_{ij})  \star' \Phi_P  (x_j,e) & \\
                            & = & (\vvvv_{ij}^{-1},x_{ij}) \star' (x_j, r_j) & \hbox{by  (\ref{eq:isomclasses}), i.e. def. of $\vvvv_{ij}^{-1}$  and $r_j$ }\\
                             & = & (( \vvvv_{ij}^{-1}, x_{ii}) \bullet'  (e,x_{ij})) \star' (x_j, e)\cdot r_j &  \hbox{by def. \ref{definition of Adapted extension}  item 2 and 4} \\
                              & = &( \vvvv_{ij}^{-1}, x_{ii}) \star' ( x_i, \lambda'_{ij}) \cdot r_j & \hbox{by  (\ref{eq:cons}), i.e. def. of $\lambda_{ij}'$}\\
                              & = & ( \vvvv_{ij}^{-1}, x_{ii}) \star' ( x_i, \lambda'_{ij} r_j) & \hbox{by def. \ref{definition of Adapted extension}  item 2}\\
                             & = & ( x_i, \rho(\vvvv_{ij}^{-1})\lambda'_{ij} r_j ) & \hbox{by (\ref{defi: ppalstronP})}
\end{array}
\end{equation*}

Comparing the LHS and RHS of (\ref{eq:cobound1}) yields
$$r_i\lambda_{ij} = \rho(\vvvv_{ij}^{-1})\lambda'_{ij} r_j \hbox{ $\Leftrightarrow $ } \lambda'_{ij}= \rho(\vvvv_{ij})r_i\lambda_{ij}  r_j^{-1}$$
 which is the first relation of (\ref{eq:relate}).
Before proving the second relation of (\ref{eq:relate}), we need to explore the consequences of the commutativity of the diagram displayed in (\ref{diagram for isomorphism}).
It follows from item 3 in definition \ref{definition of Adapted extension} that $ \chi((x_i,e))$ is the element in the band given by $\chi((x_i,e))(g)=  ( g, x_{ii} ) $, so that $ \overline{\Phi_\mathcal R}(\chi((x_i,e))) $ is by definition the element of the band given by $g \mapsto \Phi_\mathcal R((g,x_{ii})) $.
Now, $\Phi_P((x_i,e))= (r_i,e) $ by (\ref{eq:isomclasses}), i.e. definition of $r_i$, so that $\chi'( \Phi_P((x_i,e)) ) $ is the element of the band given by
$g \mapsto (r_i(g),x_{ii}) $, by item (3) of definition of adapted extensions again. The commutativity of diagram (\ref{diagram for isomorphism}) can therefore be expressed by meaning that the next relation holds for all $g \in \G$:
 \begin{equation}\label{eq:PhiX} \Phi_\mathcal R(g,x_{ii}) = (r_i(g),x_{ii}) .\end{equation}
 Exploiting the assumption that $\Phi_\mathcal R$ is a Lie groupoid morphism, we can derive a more general formula as follows
\begin{equation}  \label{eq:cobound2}
\begin{array}{rcll}
\Phi_\mathcal R(g,x_{ij})&=&\Phi_{\mathcal R}((g,x_{ii})\bullet (e,x_{ij}))& \hbox{by definition \ref{definition of Adapted extension} item 4}\\
   &=&\Phi_\mathcal R(g , x_{ii} )\bullet'\Phi_\mathcal R(e,x_{ij})& \hbox{$\Phi_{\mathcal R}$ being a Lie groupoid morphism}\\
                &=&(r_i(g) , x_{ii} )\bullet'\Phi_\mathcal R(e,x_{ij})& \hbox{by (\ref{eq:PhiX})}\\
                &=&(r_i(g) , x_{ii} )\bullet'(\vvvv_{ij}^{-1}, x_{ij})& \hbox{by (\ref{eq:isomclasses}) definition of $ \vvvv_{ij}$}\\
                &=&(r_i(g)\vvvv_{ij}^{-1}, x_{ij})& \hbox{by definition \ref{definition of Adapted extension} item 4.}
\end{array}
\end{equation}
Now, we derive the second of the relations (\ref{eq:relate}) by comparing the left and right hand sides of a relation following from the assumption that $\Phi_\mathcal R$ be a Lie groupoid morphism:
 \begin{equation}  \label{eq:cobound3}
 \Phi_\mathcal R( (e,x_{ij}) \bullet (e,x_{jk}) ) = \Phi_{\mathcal R}  ((e,x_{ij})) \bullet' \Phi_\mathcal R ((e,x_{jk})),
 \end{equation}
a computation that goes as follows:
\begin{equation*}
\begin{array}{rcll}
                            & & \Phi_\mathcal R ( (e,x_{ij}) \bullet (e,x_{jk}) )& \\
                            & = & \Phi_\mathcal R ( \gggg_{ijk} ,x_{ik} )  & \hbox{by (\ref{eq:cons}), i.e. definition of $\gggg_{ijk}$ }\\
                            & = &    (r_i(\gggg_{ijk})\vvvv_{ik}^{-1},x_{ik}) & \hbox{by (\ref{eq:cobound2})},  \\
\end{array}
\end{equation*}
while the RHS of (\ref{eq:cobound3}) is:
\begin{equation*}
\begin{array}{rcll}
 & & \Phi_\mathcal R  ((e,x_{ij})) \bullet' \Phi_\mathcal R ((e,x_{jk}))& \\
 & =&  ( \vvvv_{ij}^{-1}, x_{ij} ) \bullet' ( \vvvv_{jk}^{-1}, x_{jk} )  & \hbox{by (\ref{eq:isomclasses}), i.e. def. of $\vvvv_{ij}$ and $\vvvv_{jk} $ } \\
  & = & (  \vvvv_{ij}^{-1}  \lambda'_{ij}(\vvvv_{jk}^{-1}) \gggg_{ijk}', x_{ik})  & \hbox{by (\ref{eq:recons})}. \\
\end{array}
\end{equation*}
Comparing the LHS  and the RHS of (\ref{eq:cobound3}) leads to
$$ r_i(\gggg_{ijk})\vvvv_{ik}^{-1} = \vvvv_{ij}^{-1}  \lambda'_{ij}(\vvvv_{jk}^{-1}) \gggg_{ijk}'
\hbox{ $\Leftrightarrow $ } \gggg_{ijk}'\vvvv_{ik}=\lambda'_{ij}(\vvvv_{jk})\vvvv_{ij}r_i(\gggg_{ijk}) ,$$
which is precisely the second relation of (\ref{eq:relate}), and completes the proof of the first item.

\bigskip
\noindent
2) Second we prove that, given a $\G \to \HH$-valued $1$-coboundary, by following the construction in item 2, we get an isomorphism of adapted Lie groupoid $\G \to \HH$-extensions. In order to show that the triple $(\Phi_\mathcal R, \Phi_P, \id_{\HH}) $ with $\Phi_\mathcal R, \Phi_P$ as in (\ref{eq:def_PhiX}) and (\ref{eq:def_PhiP}), is an isomorphism of $\G \to \HH$-extensions, we need to check that
 \begin{enumerate}
 \item[(a)] $\Phi_\mathcal R : \mathcal R \to \mathcal R'$ is a morphism of Lie groupoids,
 \item[(b)] $\Phi_P : P \to P'$ is a morphism of principal bundles over Lie groupoids,
 \item[(c)] the following diagram commutes:
 $$ \xymatrix{P \ar [r]^{\Phi_P} \ar [d]^{\chi}& P'\ar [d]_{\chi'}\\
              Band((\mathcal R,\bullet)\to \Cechs)\ar [r]^{\bar{\Phi}_\mathcal R} & Band((\mathcal R,\bullet')\to \Cechs)} $$
              with $\bar{\Phi}_\mathcal R$ being defined as in (\ref{diagram for isomorphism})
 \end{enumerate}
We first check that condition (a) holds, i.e that $\Phi_\mathcal R( r \bullet r'  ) = \Phi_\mathcal R(r) \bullet' \Phi_\mathcal R(r') $
for arbitrary elements of the form $r=(g,x_{ij}) \in \mathcal R$ and $r'=(g',x_{jk}) \in \mathcal R$. On the one hand:
\begin{equation*}
\begin{array}{rcll}
\Phi_\mathcal R((g,x_{ij}) \bullet (g',x_{jk})) & = & \Phi_\mathcal R (g \lambda_{ij} ( g')  \gggg_{ijk}, x_{ik}) & \hbox{by (\ref{eq:recons}) in prop. \ref{prop:cocycles}}                      \\
 & = & ( r_i ( g \lambda_{ij} (g') \gggg_{ijk}) \vvvv_{ik}^{-1}, x_{ik}) & \hbox{by (\ref{eq:def_PhiX}), i.e. definition of $\Phi_\mathcal R$},
\end{array}
\end{equation*}
while on the other hand:
\begin{equation*}
\begin{array}{rcll}
  & & \Phi_\mathcal R( (g,x_{ij}) ) \bullet' \Phi_\mathcal R (g',x_{jk}) & \\
  & = &  (r_i(g)\vvvv_{ij}^{-1},x_{ij})   \bullet'  (r_j(g')\vvvv_{jk}^{-1},x_{jk}) & \hbox{by (\ref{eq:def_PhiX}), i.e. def. of $\Phi_\mathcal R$} \\
  & = & (   r_i(g)\vvvv_{ij}^{-1}  \lambda_{ij}'( r_j(g')\vvvv_{jk}^{-1}) \gggg_{ijk}' ,x_{ik}) & \hbox{by (\ref{eq:recons}) in prop.  \ref{prop:cocycles}.}             \\
\end{array}
\end{equation*}
Of course, $\Phi_\mathcal R$ is a Lie groupoid isomorphism if and only if both sides of the previous relations are equal for all $g,g' \in \G$, i.e. if and only if
\begin{equation}\label{eq:bothsides} r_i ( g \lambda_{ij} (g') \gggg_{ijk}) \vvvv_{ik}^{-1}  = r_i(g)\vvvv_{ij}^{-1}  \lambda_{ij}'( r_j(g')\vvvv_{jk}^{-1}) \gggg_{ijk}'\end{equation}
which reduces, multiplying both sides by $r_i(g^{-1})$, to require that, for all $g' \in \G $:
$$ r_i ( \lambda_{ij} (g') \gggg_{ijk}) \vvvv_{ik}^{-1}  =\vvvv_{ij}^{-1}  \lambda_{ij}'( r_j(g')\vvvv_{jk}^{-1}) \gggg_{ijk}'$$
or, equivalently:
$$r_i ( \lambda_{ij} (g') ) r_i( \gggg_{ijk} ) \vvvv_{ik}^{-1} =  \vvvv_{ij}^{-1}  \lambda_{ij}'( r_j(g')) \vvvv_{ij} \vvvv_{ij}^{-1} \lambda_{ij}'(\vvvv_{jk}^{-1})   \gggg_{ijk}' $$
 By the second relation in (\ref{eq:relate}), $ r_i( \gggg_{ijk} ) \vvvv_{ik}^{-1} =  \vvvv_{ij}^{-1} \lambda_{ij}'(\vvvv_{jk}^{-1})   \gggg_{ijk}' $,
 so that, eventually,  $\Phi_\mathcal R$ is a Lie groupoid isomorphism if and only if for all $g ' \in \G$
\begin{equation}\label{eq:*}
  r_i ( \lambda_{ij} (g') )= \vvvv_{ij}^{-1}  \lambda_{ij}'( r_j(g')) \vvvv_{ij}
\end{equation}
By axiom of crossed module the RHS of (\ref{eq:*}) is equal to
$\rho(\vvvv_{ij}^{-1})\lambda_{ij}'(r_j (g'))$
so that eventually $\Phi_\mathcal R$ is a Lie groupoid isomorphism if and only if
$$r_i \circ \lambda_{ij} \,(g')= \rho(\vvvv_{ij}^{-1}) \circ \lambda_{ij}' \circ r_j\, (g'); \forall g' \in \G ,$$
an equation which is obtained by applying $\jmath: \HH \to \Aut (\G)$ to the first relation in (\ref{eq:relate}), and is therefore true, here $\circ$ refers to the composition low of $\Aut(\G)$,
Hence, $\Phi_\mathcal R $ is a Lie groupoid isomorphism.

We wish now to check that condition (b) holds, i.e that $\Phi_P( r \star p  ) = \Phi_\mathcal R(r) \star' \Phi_P(p) $
for arbitrary elements $r=(g,x_{ij}) \in \mathcal R$ and $p=(x_i,h) \in P$. On the one hand, we compute:
\begin{equation} \label{eq:condition4}
\begin{array}{rcll}
                                          & & \Phi_P((g,x_{ij}) \star (x_{j},h)) & \\
                                          & = & \Phi_P(x_i,\rho(g)\lambda_{ij}h)  & \hbox{by (\ref{eq:recons}) in prop. \ref{prop:cocycles}}\\
                                          & = & (x_i, r_i \rho(g)\lambda_{ij}h )  & \hbox{by (\ref{eq:def_PhiP}), i.e. definition of $\Phi_P$},
\end{array}
\end{equation}
while on the other hand, we compute:
\begin{equation} \label{eq:condition5}
\begin{array}{rcll}
                       &   & \Phi_\mathcal R(g,x_{ij}) \star'\Phi_P(x_{j},h)) &\\
                       & = & (r_i(g)\vvvv^{-1}_{ij},x_{ij}) \star' (x_j, r_jh) & \hbox{by (\ref{eq:def_PhiX}-\ref{eq:def_PhiP}), i.e. def. of $\Phi_\mathcal R $ and  $\Phi_P$}\\
                       & = & (x_i, \rho(r_i(g)\vvvv^{-1}_{ij})\lambda'_{ij}r_j h) & \hbox{by (\ref{eq:recons}) in prop. \ref{prop:cocycles}}
\end{array}
\end{equation}
Equations (\ref{eq:condition4}) and (\ref{eq:condition4}),  together with $\rho(r_i(g)\vvvv^{-1}_{ij})\lambda'_{ij}r_j h= r_i \rho(g)\lambda_{ij}h$
(an immediate consequence of (\ref{eq:relate})), imply that:
$$ \Phi_P((g,x_{ij}) \star(h,x_{j}))= \Phi_\mathcal R(g,x_{ij}) \star'\Phi_P(h,x_{j}))$$
which completes the proof of (b). Condition (c) is a direct computation.

Last, we have to check that both constructions in item 1 and item 2 are inverse one to the other. It is easy to see that, applying the construction of item 2 and  then the construction of item 1 to a $\G \to \HH $-valued coboundary $(r_i,\vvvv_{ij}) $, one obtains $(r_i,\vvvv_{ij}) $ again.
Moreover,  two
  $(\Phi_\mathcal R,\Phi_P) $, $(\Phi_\mathcal R',\Phi_P') $
  isomorphisms of $\G \to \HH$-extensions  which correspond to the same coboundary $(r_i,\vvvv_{ij}) $ need to be equal. This follows from (\ref{eq:cobound2}), which clearly implies that $\Phi_{\mathcal R}=\Phi_\mathcal R'$, and from (\ref{eq:isomclasses}), which implies that $\Phi_P$ and $\Phi_{P'}$ coincide on every element in $P$ of the form $(x_{i},e) $, and are therefore equal since principal bundle morphisms that coincide on some global section coincide globally. This completes the proof.
\end{proof}
Let $ {\mathcal U}=(U_i)_{i \in I}$ be an open covering of a manifold $N$,
and $\G \to \HH$ be a crossed module of Lie groups.
It follows from proposition \ref{prop:coboundaries} that coboundaries define an equivalence relation on
the set of $\G \to \HH$-valued $1$-cocycles w.r.t. ${\mathcal U} $.  The quotient set obtained by this equivalence relation is called \emph{$\G \to \HH$-valued $1$-cohomology w.r.t. ${\mathcal U}$} and denoted by $H^1_{\mathcal U}( \G \to \HH)$.

The next corollary follows from propositions \ref{prop:cocycles} and
 \ref{prop:coboundaries}.

\begin{cor}\label{cor:1-1Cohomology-adapted}
Let $ {\mathcal U}=(U_i)_{i \in I}$ be an open covering of a manifold $N$,
and $\G \to \HH$ be a crossed module of Lie groups.
There is a one to one correspondence between the set $H^1_{\mathcal U}( \G \to \HH)$
and the set of all adapted Lie groupoid $\G\to \HH$-extensions of $\Cechs $ up to isomorphisms (of Lie groupoids $\G \toto \HH $-extensions of $\Cechs$).
\end{cor}

The notion of adapted extension may appear to be somewhat arbitrary. We wish to convince the reader that it is not,
 by showing the next proposition.

\begin{prop}\label{prop:everyextensionisadapted}
Let $ {\mathcal U}=(U_i)_{i \in I}$ be an open covering of a manifold $N$ such that $U_{ij}$
is a contractible open set for all $i,j \in I$,
and let $\G \to \HH$ be a crossed module of Lie groups.
Then every Lie groupoid $\G \to \HH$-extension of the \v{C}ech groupoid $\Cechs $ is isomorphic (as a Lie groupoid $\G \to \HH$-extension of the \v{C}ech groupoid $\Cechs$) to an adapted Lie groupoid $\G \to \HH$-extension $\Cechs$.
\end{prop}
\begin{proof}
Let $( \mathcal R \stackrel{\phi}{\to} \Cechs, P\to \coprod_{i \in I}U_i, \chi) $ be a Lie groupoid $\G \to \HH$-extension of $\Cechs$.
Since $\coprod_{i \in I}U_i $ is a disjoint union of contractible sets (since $U_{ij}$ is by assumption contractible for all $i,j \in I$, so is $U_i = U_{ii}$), there exists a global section $\sigma $ of the $\HH$-principal bundle $ P\to \coprod_{i \in I}U_i$.

Since $\chi: P \to Band (\mathcal R \stackrel{\phi}{\to} \Cechs)$ is by assumption a morphism of principal bundles over the identity of $ \coprod_{i \in I}U_i $, the map $\hat{\sigma} :=\chi \circ \sigma  $ is a global section of the $Aut(\G)$-principal bundle $Band( \mathcal R \stackrel{\phi}{\to} \Cechs )$.
In turn, a global section of the band amounts to a global trivialization of the kernel $K \to  \coprod_{i \in I}U_i$, by considering the group bundle isomorphism $ \tau_{K} :\G \times \coprod_{i \in I} U_{i}  \simeq K $ given by
 $ (g,x_i) \mapsto  \hat{\sigma}(x_i) (g) $.
 Since, by construction, $\hat{\sigma}(x_i) $ belongs to $Band_{x_i}= Isom(\G,K_{x_i}) $, it is clear that
 $\tau_K $ is, as expected, a group bundle isomorphism over the identity of $\coprod_{i \in I} U_{ii} $.

 Now, the surjective submersion $\phi: \mathcal R \to  \coprod_{i,j \in I} U_{ij} $ restricts to a surjective submersion
 from $\mathcal R  \backslash K $ to $\coprod_{i \neq j} U_{ij} $, and the fibers of this submersion are acted upon transitively and freely by $K $. Using $\tau_K $, we endow
 $$\mathcal R  \backslash K  \to \coprod_{i,j \in I \hbox{ s.t.} i \neq j} U_{ij}$$
with a structure of $\G$-principal bundle as follows: the outcome of the action of  $g \in \G $ on $r \in  \mathcal R  \backslash K $
is defined to be $ \tau_K (g, {s(r)} ) \bullet_{\mathcal R} r   $.
Every principal bundle over a disjoint union of contractible open sets is trivial, which means, in this case, that there is a global section
$\sigma_{1}:\coprod_{i \neq j} U_{ij} \to \mathcal R\backslash K $.
Then we define $\tau_{\mathcal R \backslash K } : \G \times \coprod_{i \neq j} U_{ij} \to \mathcal R \backslash K $  by
$$ (g,x_{ij}) \mapsto \tau_K(g,x_i) \bullet_{\mathcal R} \sigma_1(x_{ij})$$

 for all $i,j \in I, i \neq j$.
By construction $\tau_{\mathcal R \backslash K } $ is a group bundle morphism over the identity of $\coprod_{i \neq j }U_{ij} $.
Gluing $\tau_{ K}  $ and $ \tau_{\mathcal R \backslash K}$, we get a map (over the identity of $\coprod_{i,j \in I}U_{ij} $)
that we denote by $\tau : \G \times \coprod_{i,j \in I} U_{ij} \to \mathcal R $, namely:
 $$ \tau(g,x_{ii}):= \tau_K(g,x_i)  ,\forall i\in I $$ and $$ \tau(g,x_{ij}):= \tau_{\mathcal R  \backslash K } (  g,x_{ij}  ),  \forall i,j \in I  \,\, \hbox{with}\,\, i \neq j  $$
 The section $\sigma$ of $ P \to  \coprod_{i \in I} U_i$ also induces a map
   $\Psi_P :   \coprod_{i \in I} U_{i} \times \HH \simeq  P $ given by:
   \begin{equation} \label{eq:psiP}
    (x_{i}, h)  \mapsto \sigma (x_i) \cdot  h.
    \end{equation}
With the help of this pair of maps $\Psi_P $ and $\tau$, the structure of $\G \to \HH$-extensions on $ (\mathcal R \stackrel{\phi}{\to} \Cechs, P\to \coprod_{i \in I}U_i, \chi)$ is transported and induces a structure of $\G \to \HH$-extension on $( \G \times \coprod_{i,j \in I} U_{ij} \stackrel{\phi}{\to} \Cechs,  \coprod_{i} U_{i} \times \HH \to \coprod_{i \in I}U_i, \chi' )  $. Explicitly the induced Lie groupoid structure on $\G\times \coprod_{i,j\in I}U_{ij} \toto \coprod_{i\in I}U_{i}$ is given by:
$$(g,x_{ij})\bullet(g',x_{jk}):=\tau^{-1}(\tau(g,x_{ij})\bullet_{\mathcal R} \tau(g',x_{jk})),$$
for all $g,g' \in \G, i,j,k\in I, x\in U_{ijk},$
the induced action of Lie groupoid $G\times \coprod_{i,j\in I}U_{ij} \toto \coprod_{i\in I}U_{i}$ on $\coprod_{i\in I}U_{i}\times \HH$ is given by:
$$(g,x_{ij})\star(x_j,h):=\Psi_P^{-1}(\tau(g,x_{ij}) \bullet_{\mathcal R,P}\Psi_P(x_j,h)),$$
for all $g, \in \G,h\in \HH, i,j\in I, x\in U_{ij},$ and the induced principal bundle structure on $\coprod_{i\in I}U_{i}\times \HH \to \coprod_{i}U_{i}$ over the Lie groupoid $\G\times \coprod_{i,j\in I}U_{ij} \toto \coprod_{i,j\in I}U_{ij}$ is given by:
$$(x_i,h)\cdot h'=(x_i,hh'),$$
for all $h,h' \in \HH, i\in I, x\in U_{i},$
Last we define $ \chi':\coprod_{i\in I}U_{i} \times \HH \to Band(\G\times \coprod_{i,j\in I}U_{ij} \toto \coprod_{i,j\in I}U_{ij})$ by
$$ (x_i,h)\mapsto (x_i,j(h)),$$
for all $h \in \HH, i\in I, x\in U_{i},$
We claim that:
\begin{enumerate}
\item The extension $Ext_2:=(G\times \coprod_{i,j \in I}U_{ij} \to \coprod_{i,j \in I}U_{ij}, \coprod_{i\in I}U_{i} \times \HH \to \coprod_{i \in I}U_{i},\chi' )$ is a Lie groupoid $\G \to \HH$-extension.
\item The Lie groupoid $\G \to \HH$-extension $Ext_2$ is isomorphic to the  Lie groupoid $\G \to \HH$-extension $Ext_1:=(\mathcal R \to \coprod_{i,j \in I}U_{ij}, P \to \coprod_{i \in I}U_{i}, \chi)$.
\item The Lie groupoid $Ext_2$    is an adapted Lie groupoid $\G \to \HH$-extension.
\end{enumerate}
These claims complete the proof of the proposition.
For the proof of claim 1), it is enough to check that $(x_i,h)\cdot \rho(g)=\chi'(x_i,h)(g)\star(x_i,h)$ for all $x\in N,i\in I,h\in \HH, g\in \G,$ which goes as follows:
\begin{equation}
\begin{array}{rcll}
                        & &\chi'(x_i,h)(g)\star (x_i,h) & \\
                        &=& (h(g),x_{ii}) \star (x_i,h)& \hbox{by def of $\chi'$}\\
                        &=& \Psi_P^{-1}(\tau (h(g),x_{ii})\bullet_{\mathcal R,P}\Psi_P(x_i,h)) & \hbox{by def of $\star$}\\
                        &=& \Psi_P^{-1}(\chi \circ \sigma (x_i)(h(g)) \bullet_{\mathcal R,P}\sigma(x_i)\cdot h) &\hbox{by def of $\tau$ and def of $\Psi_P$}\\
                        &=& \Psi_P^{-1}( \chi(\sigma(x_i)\cdot h)(g)\bullet_{\mathcal R,P}\sigma(x_i)\cdot h)& \hbox{$\chi$ is morphism of ppal bundles}\\
                        &=& \Psi_P^{-1}(\sigma(x_i)\cdot h \cdot \rho(g))& \hbox{since $Ext_1$ is a $\G \to \HH$-extension}\\
                        &=& (x_i,h\cdot\rho(g))& \hbox{by def of $\Psi_P^{-1}$}\\
                        &=& (x_i,h)\cdot\rho(g)&.

\end{array}
\end{equation}
For the proof of claim 2), since  $\Psi_P$ and $\tau$ are clearly morphisms of principal bundles and Lie groupoids respectively, it is enough to prove that the following diagram commutes: \label{diag*}
$$\xymatrix{\coprod_{i\in I}U_i \times\HH \ar[d]^{\chi'} \ar[rrr]^{\Psi_P}&& & P \ar[d]\\
             Band(\G\times\coprod_{i,j \in I}U_{ij} \ar[r]& \coprod_{i,j \in I}U_{ij})\ar[r]^{\bar{\tau}} &Band(\mathcal R \ar[r] &\coprod_{i,j\in I}U_{ij} ) }$$
In turn, the commutativity of this diagram follows from:
\begin{equation*}
\begin{array}{rcll}
                     & & (\bar{\tau} \circ \chi' (x_i,h))(g) & \\
                     &=& \tau(\chi'(x_i,h)(g)) & \hbox{by def of $\bar{\tau}$}\\
                     &=& \tau(h(g),x_{ii}) & \hbox{by def of $\chi'$}\\
                     &=& (\chi\circ \sigma (x_i))(h(g)) & \hbox{by def of $\tau$}\\
                     &=& \chi \circ \sigma(x_i) \circ j(h)(g) &\hbox{ }\\
                     &=& \chi(\sigma(x_i)\cdot h)(g)& \hbox{since $\chi$ is a morphism of ppal bundles}\\
                     &=& \chi \circ \Psi_P(x_i, h)(g)& \hbox{by def of $\Psi_P $},

\end{array}
\end{equation*}
for all $i\in I , x_i \in U_i, h \in \HH, g\in \G$. \\
 Last we need to prove claim 3). For this it is enough to check that axiom 2 in Definition \ref{definition of Adapted extension} holds. Note that the other axioms in Definition \ref{definition of Adapted extension} hold by construction of $\chi$ and $\bullet$.
 \begin{equation}\label{eq:adap}
 (g,x_{ii})\bullet (g',x_{ij})=(gg',x_{ij})
 \end{equation}
  for all $x \in N, i,j \in I , g,g' \in \G$. This goes as follows:
\begin{equation*}
\begin{array}{rcll}
                     & & \hbox{LHS of (\ref{eq:adap})}& \\
                     &=& \tau^{-1}(\tau(g,x_{ii})\bullet_{\mathcal R}\tau(g',x_{ij}))&\hbox{by def of $\bullet$}\\
                     &=& \tau^{-1}(\chi \circ \sigma(x_i)(g)\bullet_{\mathcal R}\chi \circ\sigma(x_i)(g')\bullet_{\mathcal R}\sigma_1(x_{ij}))&\hbox{by def of $\tau$}\\
                     &=& \tau^{-1}(\chi \circ\sigma(x_i)(gg')\bullet_{\mathcal R} \sigma_1(x_{ij}))&\hbox{$\chi \circ\sigma(x_i)$ is a morphism of groups}\\
                     &=& \tau^{-1}(\tau(gg',x_{ij}))&\hbox{by def of $\tau$ }\\
                     &=& \hbox{RHS of (\ref{eq:adap})}&
\end{array}
\end{equation*}
so that condition (\ref{eq:chi}) in definition (\ref{def: Azimi-extension}) is satisfied. The other conditions are satisfied by construction.
\end{proof}

We can now state the conclusion of this section, which follows immediately from proposition \ref{prop:everyextensionisadapted} and corollary \ref{cor:1-1Cohomology-adapted}.

\begin{them}
Let $ {\mathcal U}=(U_i)_{i \in I}$ be an open covering of a manifold $N$ such that $U_{ij}$
is a contractible open set for all $i,j \in I$, and $\G \to \HH$ a crossed module of Lie groups.
There is a one to one correspondence between
\begin{enumerate}
\item[(i)] the set $H^1_{\mathcal U}( \G \to \HH)$,
\item[(ii)] adapted Lie groupoid $\G \to \HH$-extensions of $\Cechs $ up to isomorphisms of Lie groupoid $\G \to \HH$-extensions of $\Cechs $,
\item[(iii)] Lie groupoid $\G \to \HH$-extensions of $\Cechs $ up to isomorphisms of Lie groupoid $\G \to \HH$-extensions of $\Cechs $.
\end{enumerate}
\end{them}
\begin{proof}
The correspondence  between (i) and (ii) was already stated in corollary \ref{cor:1-1Cohomology-adapted}.
The correspondence between (iii) and (ii) comes from proposition \ref{prop:everyextensionisadapted} which states that every Lie groupoid $\G \to \HH$-extension
 of $\Cechs $ is isomorphic to an adapted one. Of course, a given extension can be isomorphic (as Lie groupoid $\G \to \HH$-extensions of $\Cechs $)
 to two different adapted $\G \to \HH$-extensions $\Cechs $, but both adapted extensions are then isomorphic (as Lie groupoid $\G \to \HH$-extensions of $\Cechs $), so that the assignment from (iii) to (ii) is well-defined and is one to one by construction. This concludes the proof of the theorem.
\end{proof}

\section{Morita equivalence of $\G \to \HH$-gerbes}\label{sec:MoritaExt}

Let $\G \to \HH$ be a crossed-module. We intend in this section to define, purely in terms of Lie groupoids, the notion of $\G \stackrel{\rho}{\to} \HH $-gerbes over a given Lie groupoid $\BBB \toto \BBB_{0}$, having in mind the case where $\BBB \toto \BBB_{0}$ is the trivial Lie groupoid $N \toto N$ associated to a manifold~$N$.

In view of the preceding section, it is reasonable to consider all the $\G \to \HH$-extensions of all the possible pull-back of $\BBB \toto \BBB_{0}$ with respect to surjective submersions. For instance, when $\BBB \toto \BBB_{0} $ is of the form $N \toto N$, with $N$ a manifold, this includes all the $\G \to \HH$-extensions of the \v{C}ech groupoids associated to an arbitrary open cover of $N$ (because the \v{C}ech groupoid $N[{\mathcal U}]  \toto \coprod_{i \in I}U_i$ is the pull-back groupoid of $N \toto N $ with respect to the natural inclusion maps $\imath: \coprod_{i \in I}U_i \to N $).
  But of course, we shall later have to take a quotient of that class. We do it by identifying
two $\G \to \HH$-extensions which are Morita equivalent in some sense described below.

\subsection{Definition of Morita equivalences of $ \G \to \HH$-extensions and $\G \to \HH $-gerbes}

 Let us first define
what the pull-back of a $ \G \to \HH$-extension is.

 Given a Lie groupoid extension $\mathcal{R} \stackrel{\phi}{\to} \GGG \toto M$ and a surjective submersion $p:M' \to M$, the functor of definition \ref{def:pback} applied to $\mathcal{R} \stackrel{\phi}{\to} \GGG $ yields a Lie groupoid extension
  $$\mathcal{R}[p] \stackrel{\phi[p]}{\longrightarrow} \GGG[p,M].$$
  It is routine to check that   $\mathcal{R}[p] \stackrel{\phi[p]}{\longrightarrow} \GGG[p,M] \toto M'$ is again a Lie groupoid extension.
    This construction still goes through under the weaker assumption that $p$ is a generalized surjective submersion for the Lie groupoid $\GGG\toto M $. Notice that $p$ is a generalized surjective submersion for the Lie groupoid $\GGG \toto M$
    if and only if it is a generalized surjective submersion for the Lie groupoid $\mathcal{R} \toto M$, so that we could say
    that this construction still goes through under the weaker assumption that $p$ be a generalized surjective submersion for the Lie groupoid ${\mathcal R}\toto M $.
     For all such maps $p:M'\to M$, we call the Lie groupoid extension $\mathcal{R}[p] \stackrel{\phi[p]}{\longrightarrow} \GGG[p,M] \toto M'$ the \emph{pull-back of the Lie groupoid extension $\mathcal{R} \stackrel{\phi}{\to} \GGG \toto M$ with respect to $p$}.

Having defined the pull-back of Lie groupoid extensions, we wish to define the pull-back of Lie groupoids $\G \to\HH $-extensions.
   This shall require to go through some technical considerations about the pull-back of the kernel and the band of a Lie groupoid $\G$-extension.

   There is a clear notion of pull-back for group bundles (resp. principal bundles): to say it in one word, given $P \stackrel{\pi}{\to} M $
   a group bundle (resp. principal bundle), and $p: M' \to M $ a smooth map, then the fibered product $P \times_{\pi,M,p} M' $ endows a natural structure of group bundle (resp. principal bundle).  To a Lie groupoid extension, we have associated in section \ref{subsec:defofextensions} a bundle of group, called the kernel, and, provided that the extension is a $\G$-extension, we have also constructed an $\Aut(\G) $-principal bundle, called the band. The next proposition claims that these two constructions behave well with respect to pull-back.

  \begin{prop}\label{prop:pullbackExtensions} Let $M,M'$ be smooth manifolds, $p : M' \to M $ be a surjective submersion,
  and $\mathcal{R}\stackrel{\phi}{\to} \GGG \toto M$ a Lie groupoid extension over the base manifold $M$.Then:
  \vspace{-8mm}
  \begin{enumerate}
   \item there is a canonical isomorphism between the kernel of Lie groupoid extension $\mathcal{R}[p] \stackrel{\phi[p]}{\longrightarrow} \GGG[p,M] \toto N$ and the pull-back of the kernel $K$ of the Lie groupoid extension $\mathcal{R}\stackrel{\phi}{\to} \GGG \toto M$ by the surjective submersion $p$,
   \item the pull-back of a Lie groupoid $\G$-extension by a surjective submersion is a Lie groupoid $\G$-extension,
   \item in the case of a $\G$-extension, there is a canonical isomorphism between the band of $\mathcal{R}[p] \stackrel{\phi[p]}{\lra} \GGG[p,M] \toto M'$ and the pull-back of the band of $\mathcal{R}\stackrel{\phi}{\to} \GGG \toto M$ by $p$.
   \end{enumerate}
    The same holds true when $p$ is a generalized surjective submersion.
  \end{prop}
  \begin{proof}
  The kernel of the pull-back of the Lie groupoid extension $\mathcal{R}[p] \stackrel{\phi[p]}{\lra} \GGG[p,M] \toto M'$, denoted by $K[p]$, is, as a set, given by $\{ (n,k,n) | n\in M' , k\in K_p(n) \} $, where $K \to M$ is the kernel of the Lie groupoid extension $\mathcal{R} \stackrel{\phi}{\to} \GGG \toto M$. As a bundle of group, $K[p]$ can be identified, therefore, with $M'\times_M K$. This proves the first item.

  In particular, the fiber $K[p]_{n}$ of the kernel $K[p]$ over a given point $n\in M'$ is isomorphic to
   $K_{p(n)} $, and, more generally, if $K$ is locally trivial with typical fiber $\G$, so is its pull-back $K[p]$. This  means precisely that the pull-back of a Lie groupoid $\G$-extension is again a Lie groupoid $\G$-extension.
   This proves the second item.

     The identification between $K[p] $ and the kernel $K'$ of $\mathcal{R}[p] \stackrel{\phi[p]}{\to} \GGG[p,M] \toto M'$ induces an identification between
     the set of all Lie group automorphisms from $\G$ to $K_m'$ and $Band_{p(m)}(\mathcal{R} \stackrel{\phi}{\to} \GGG) $ for all $ m \in M'$. All together, these identifications yield an identification $Band(\mathcal{R}[p] \stackrel{\phi[p]}{\to} \GGG[p,M])$ and  $M'\times_{M} Band(\mathcal{R} \stackrel{\phi}{\to} \GGG)$. This proves the last item.
   \end{proof}

   We are now able to define clearly the notion of pull-back of a $\G \to \HH$-extension $( \mathcal{R} \to \GGG, P \to M, \chi ) $.
Let $p: M' \to M $ be a (maybe generalized) surjective submersion. According to the second item in proposition \ref{prop:pullbackExtensions}, the pull-back extension $ \mathcal{R}[p] \stackrel{\phi[p]}{\to} \GGG[p,M] $ is again a $\G$-extension. Moreover, $p^* P = P \times_M M' \to M' $ is an $\HH$-principal bundle over $M'$, which is acted upon by  $\mathcal{R}[p] \toto M' $ as follows:
    $$ (n,r,n') \bullet (x,n') = (r \bullet x , n) ,$$
    for all $n,n' \in M', x\in P, r \in R $ subject to the constraints $p(n)=s(r),t(r)=p(n')=p(x) $.
   The map $\chi[p] : P \times_M M'  \to Band(  \mathcal{R} \to \GGG ) \times_M M' $ defined by $ (p,n) \to (\chi(p), n) $, composed with the canonical isomorphism between $Band(  \mathcal{R} \to \GGG ) \times_M M' $ and
    $ Band(  \mathcal{R}[p]  \stackrel{\phi[p]}{\to}  \GGG[p,M] ) $ of item 3 in proposition \ref{prop:pullbackExtensions}, satisfies all the requirements needed to guarantee that   $( \mathcal{R}[p] \stackrel{\phi[p]}{\to} \GGG[p,M], p^*P \to M', \chi[p] ) $ is a $\G \to \HH$-extension.

    \begin{defi} Let $( \mathcal{R}  \stackrel{\phi}{\to} \GGG, P \to M, \chi ) $ be a Lie groupoid $\G \to \HH$-extension.
Let $p: M' \to M $ be a (generalized) surjective submersion. We call the Lie groupoid $\G \to \HH$-extension defined in the lines above the \emph{pull-back of the Lie groupoid $\G \to \HH$-extension $(\mathcal{R} \stackrel{\phi}\to \GGG, P \to M, \chi)$ with respect to $p$} and we denote it by
$  ( \mathcal{R}[p] \stackrel{\phi[p]}\to \GGG[p,M], P[p] \to M[p], \chi[p] ) $.
    \end{defi}

Indeed, we need a notion which is slightly more subtle. Recall that our purpose is to define gerbes as being the quotient of a sub-class of all $\G \to \HH$-extensions by some relation. We can now be more precise, and define, given a Lie groupoid $\BBB \toto \BBB_0 $,
a \emph{$\G \to \HH$-extension over $\BBB \toto \BBB_0$} to be a quadruple
$(q, \mathcal{R}\stackrel{\phi}{\to} \BBB[q],P \to M,\chi)$ where:
\begin{enumerate}
\item $q: M \to \BBB_0$ is a surjective submersion,
\item $ ( \mathcal{R} \stackrel{\phi}{\to} \BBB[q],P \to M,\chi)$ in a $\G \to \HH$-extension of the pull-back groupoid
$ \BBB[q] \toto M  $ of $\BBB \toto \BBB_0  $ with respect to $q$.
\end{enumerate}

We define the pull-back of those.

\begin{defi}
The pull-back of a Lie groupoid $\G \to \HH$-extension $(q, \mathcal{R}\stackrel{\phi}{\to} \BBB[q],P \to M,\chi)$  over the Lie groupoid $\BBB \toto \BBB_0$ w.r.t the surjective submersion $p: M'\to M$ is the Lie groupoid $\G \to \HH$-extension $(q\circ p,Y \stackrel{\phi[p]}{\to} \BBB[q \circ p],p^* P \to M',\chi[p])$ over the Lie groupoid $\BBB \toto \BBB_0$.
\end{defi}

\begin{rem}
 The previous definition used implicitly the existence of a natural isomorphism $ \BBB[q][p] \simeq  \BBB[q \circ p]$:
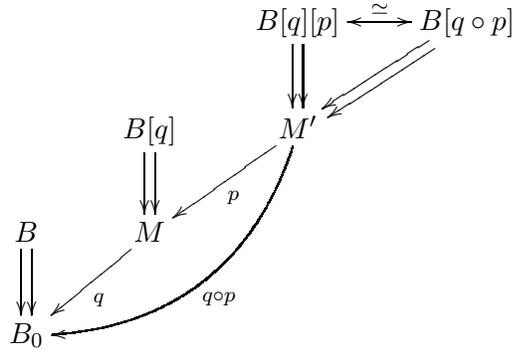
  $$ \xymatrix{  & &\BBB[q][p] \ar@<2pt>[d] \ar@<-2pt>[d]  \ar@{<->}^-{\simeq}[r] &  \BBB[q \circ p] \ar@<2pt>[dl] \ar@<-2pt>[dl]\\ & \BBB[q] \ar@<2pt>[d] \ar@<-2pt>[d] & M' \ar[dl]^p \ar@/^2pc/[ddll]^{q \circ p}\\  \BBB\ar@<2pt>[d] \ar@<-2pt>[d] & M\ar[dl]^q \\ \BBB_{0}  }  $$
  Indeed, the pull-back of the $\G \to \HH$-extension
  $( \mathcal{R}\stackrel{\phi}{\to} \BBB[q],P \to M,\chi)$ with respect to $p$ is a priori a $\G \to \HH$-extension
  of $\BBB[q][p]$. But in view of the
  isomorphism $ \BBB[q][p] \simeq  \BBB[q \circ p]$, it can be considered as a $\G \to \HH$-extension
   of $\BBB[q\circ p] $, and $(q \circ p, \mathcal{R}[p] \stackrel{\phi[p]}{\to} \BBB[q \circ p],p^* P \to M,\chi[p])$ is a $\G \to \HH$-extension over $\BBB \toto \BBB_0 $.
\end{rem}

We can now define the notion of Morita equivalence that we are interested in.

\begin{defi}
A Morita equivalence between two Lie groupoid $\G \to \HH$-extensions $(q,{\mathcal R} \stackrel{\phi}{\to} \BBB[q],P \to M,\chi)$ and
$(q',{\mathcal R} \stackrel{\phi}{\to} \BBB[q'],P \to M,\chi)$ over $\BBB \toto \BBB_0 $ is a triple
$(M'',p,p') $ where $M''$ is a manifold, $p: M'' \to M $ and $q:M'' \to M'$ are surjective submersions, such that:
\begin{enumerate}
\item the following diagram commutes:
 $$ 
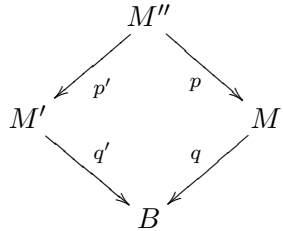
\xymatrix{ & M'' \ar[dr]_{p} \ar[dl]^{p'} & \\ M' \ar[dr]^{q'} & & M \ar[dl]_{q} \\ & B & \\ }  $$

 $$ q \circ p=q' \circ p', $$
\item the pull-back of the Lie groupoid $ \G \to \HH$-extension $(\mathcal R \stackrel{\phi}{\to} \BBB[q],P \to M,\chi)$ with respect to $p$ is isomorphic
to the pull-back of the Lie groupoid  $ \G \to \HH$-extension $(\mathcal R' \stackrel{\phi}{\to} \BBB[q'],P' \to M',\chi')$ with respect to $p'$
(notice that both pull-back Lie groupoid $\G \to \HH $-extensions are $\G \to \HH $-extensions of $\BBB[q'\circ p'] = \BBB[q \circ p]$).
\end{enumerate}
\end{defi}
In terms of commutative diagram, Morita equivalence of Lie groupoid $\G \to \HH$-extensions over $\BBB \toto \BBB_0 $
can be visualized as follows
$$ \xymatrix{ & p^* P \ar@/^3pc/[rrrr]_{\simeq}  \ar[ddrr] & \mathcal R[p] \ar@/_1pc/[rr]^{\simeq} \ar[dr] & & \mathcal R[p'] \ar[dl]& p'^* P' \ar[ddll] & \\
         P \ar[ddr] & \mathcal R \ar[d]& & \BBB [q \circ p]\approx \BBB [q'\circ p'] \ar@<1pt>[d] \ar@<-1pt>[d]  & & \mathcal R' \ar[d] & P' \ar[ddl]\\
          & \mathcal B[q] \ar@<1pt>[d] \ar@<-1pt>[d] & & M'' \ar@{.>}[dll]^p \ar@{.>}[drr]_{p'}  & & \mathcal B[q'] \ar@<1pt>[d] \ar@<-1pt>[d] &\\
           & M \ar[drr]^{q} & & \mathcal B \ar@<1pt>[d] \ar@<-1pt>[d] & & M' \ar[dll]_{q'} & \\
            & & & \mathcal B_{0} & & & }$$
          \begin{examp} \label{ex:isomOver}
A pair  $(q,\mathcal R \stackrel{\phi}{\to} \BBB[q],P \to M,\chi)$
and  $(q,\mathcal R' \stackrel{\phi}{\to} \BBB[q],P'\to M,\chi)$ of
$\G \to \HH$-extensions over $\BBB \toto \BBB_0 $ which are isomorphic over the identity of
$\BBB[q]$ are Morita equivalent.
\end{examp}
\begin{examp} \label{ex:pullbackOver}

Every Lie groupoid $\G \to \HH$-extension over a Lie groupoid is Morita equivalent to its pull-back with respect to a (generalized) surjective submersion.
\end{examp}

We can not say, strictly speaking, that Morita equivalence of $\G \to \HH$-extensions over a given Lie groupoid
$\BBB \toto \BBB_0 $ is an equivalence relation because $\G \to \HH$-extensions over $\BBB \toto \BBB_0 $ do not form a set. However, the axioms of equivalence relations remain satisfied, as shown in the next proposition

\begin{prop}
Let $\BBB \toto \BBB_0 $ be a Lie groupoid.
\begin{enumerate}
\item A $\G \to \HH$-extension over $\BBB \toto \BBB_0 $ is always Morita equivalent to itself.
\item Let $Ext_1,Ext_2$ be $\G \to \HH$-extensions over $\BBB \toto \BBB_0 $. $Ext_1$ is Morita equivalent to $Ext_2$
if and only if $ Ext_2$ is Morita equivalent to $Ext_1$.
\item Let $Ext_1,Ext_2,Ext_3$ be $\G \to \HH$-extensions over $\BBB \toto \BBB_0 $. If
$ Ext_1$ is Morita equivalent to $Ext_2$ and $ Ext_2$ is Morita equivalent to $Ext_3$, then $ Ext_1$ is Morita equivalent to $Ext_3$.
\end{enumerate}
\end{prop}
\begin{proof}
Only the third item merits some justification. If $M$, together with the surjective submersions $p,q$ give a
Morita equivalence between $Ext_1 $ and $Ext_2 $, while $M' $ together with the surjective submersions $p',q'$ give a
Morita equivalence between $Ext_2 $ and $Ext_3 $, then we introduce $ M'':= M \times_{q,M_2,p'} M' $
and equip it with the surjective submersions $(m,m') \to p(m) $ and $(m,m') \to q'(m') $ onto $M_1$
and $M_3 $ respectively, where, in the previous, $M_i, i=1,2,3 $ is the base manifold of the $\G \to \HH$-extension $Ext_i $.
A cumbersome but easy computation shows that the pull-back of $Ext_1 $ and $Ext_3 $ to $ M'' $ are isomorphic Lie groupoid $\G \to \HH$-extensions.
\end{proof}

This proposition allows one to give, at last, the following definition.

\begin{defi}
A $\G \to \HH$-gerbe over $ \BBB \toto \BBB_0$ is a Morita equivalence class of Lie groupoid $\G \to \HH$-extensions over $ \BBB \toto \BBB_0$.
\end{defi}

To justify this definition, we shall in subsection \ref{subsec:manifoldcase} show that, when the Lie groupoid $ \BBB \toto \BBB_0$ is simply a manifold
 Lie groupoid $ \BBB \toto \BBB_0$, $\G \to \HH$-gerbe are precisely the same thing as  $\G \to \HH$-valued non-Abelian $1$-cohomology.

\subsection{The manifold case: $\G \to \HH$-gerbes as non-Abelian $1$-cohomology}\label{subsec:manifoldcase}

The notion of  $\G \to \HH$ non-Abelian $1$-cohomology w.r.t. a given open covering was introduced in section \ref{subsec:mfdcase} .
As usual, $\G \to \HH$ non-Abelian $1$-cohomology is obtained by inductive limits of those.
More precisely, we proceed as follows. By a \emph{refinement} of an open cover ${\mathcal U}=(U_i)_{i \in I} $, we mean a pair $ ({\mathcal V},\sigma)$ made of an open cover $ {\mathcal V}=(V_j)_{j \in J}$ together with a map $\sigma: J \to I $ such that $ V_j \subset U_{\sigma (j)} $ for all $j \in J$.
Notice that $\sigma $ induces a map, again denoted by $\sigma $, from $\coprod_{k,l \in J}V_{kl} $ to
$\coprod_{i,j\in I}U_{ij} $ (resp. $\coprod_{k,l,m \in J}V_{klm} $ to
$\coprod_{i,j,k\in I}U_{ijk} $), obtained by mapping $x_{kl} \in V_{kl}  $
to $x_{\sigma(k)\sigma(l)} \in U_{\sigma(k)\sigma(l)} $ (using the notations of section \ref{notation:Cech}).
By the \emph{pull-back} of a non-Abelian $1$-cocycle $(\lambda,\gggg)   \in \mathcal C ^{\infty}( \coprod_{i,j \in I}V_{ij},\HH ) \times \mathcal C ^{\infty} (
\coprod_{i,j,k \in J}V_{ijk},\G)  $ w.r.t. ${\mathcal U}$, we mean
the pair of functions $(\sigma^* \lambda,\sigma^* \gggg )$ in $\mathcal C ^{\infty}( \coprod_{i,j \in I}V_{ij},\HH ) \times \mathcal C ^{\infty} (
\coprod_{i,j,k \in J}V_{ijk},\G)  $. Notice that, by construction, $ (\sigma^* \lambda)_{ij} = \left. \lambda_{\sigma(i)\sigma(j)} \right|_{V_{ij}} $
and
 $  (\sigma^* \gggg)_{ijk} =  \left.\gggg_{\sigma(i)\sigma(j)\sigma(k)} \right|_{V_{ijk}} $
 for all $i,j,k \in J $.

\begin{lem}
Let  $({\mathcal V},\sigma)$  be a refinement of ${\mathcal U} $.
The pull-back of a $\G \to \HH $-valued non-Abelian $1$-cocycle w.r.t. ${\mathcal U} $  is   a $\G \to \HH $-valued non-Abelian $1$-cocycle w.r.t ${\mathcal V} $. Moreover, two $\G \to \HH $-valued non-Abelian $1$-cocycles that differ by a coboundary have pull-back that differ by a coboundary again.
\end{lem}

We now identify two $\G \to \HH $-valued non-Abelian $1$-cocycles $(\lambda, \gggg  ) $ and $(\lambda', \gggg') $, defined on covering $ \mathcal U  $ and $\mathcal U' $ of $N$ respectively, if there exists a common refinement of both
$ \mathcal U  $ and $\mathcal U' $  such that the pull-back to that refinement of
 $(\lambda, \gggg  ) $ and $(\lambda', \gggg') $ differ by a coboundary. We denote by $H^1(\G \to \HH ) $ the set henceforth obtained and we call this set the \emph{$\G \to \HH $-valued non-Abelian $1$-cohomology on $N$}. In general, $H^1(\G \to \HH ) $ has no group structure.

We can now state the main result of this section.

\begin{them}\label{th:coho=gerbes}
Let $N$ be a manifold. There is a one to one correspondence between:
\begin{enumerate}
\item $\G \to \HH$-valued non-Abelian $1$-cohomology on $N$,
\item $\G \to \HH$ gerbes over $N \toto N $.
\end{enumerate}
\end{them}


The proof of the theorem requires two lemmas. For the first one, recall from proposition \ref{prop:cocycles}
that, given an open covering  $\mathcal U $ of $N$, there is a one to one correspondence between
non-Abelian $1$-cocycles and adapted extensions of the \v{C}ech groupoid $N[\mathcal U]$.

\begin{lem}\label{lem:lemma A}
Let  $({\mathcal V},\sigma) $ be a refinement of ${\mathcal U}$ and $(\lambda,\gggg)$ be a non-Abelian $1$-cocycle w.r.t.~${\mathcal U}$. Then the adapted Lie groupoid $\G \to \HH$-extension associated to the pull-back of the non-Abelian $1$-cocycle $(\lambda, \gggg) $ is isomorphic (as a $\G \to \HH$-extension) to the pull-back of the adapted Lie groupoid $\G \to \HH$-extension associated to $(\lambda, \gggg) $. This can be expressed as the commutativity (up to a a canonical isomorphism) of the diagram:

$$ \xymatrix{*\txt{non.Ab. $1$-cocycle \\ w.r.t ${\mathcal V}$ } \ar@{<->}[rrr]^{\txt {corresp. of \\prop. \ref{prop:cocycles}}} &&& *\txt{adapted ext. of \\ $N[ {\mathcal V}]$ }   \\
&&& \\
             *\txt{ non.Ab. $1$-cocycle \\ w.r.t ${\mathcal U}$ }\ar[uu]^{pull-back} \ar@{<->}[rrr]^{\txt {corresp. of \\prop. \ref{prop:cocycles}}} &&& *\txt{adapted ext. of \\ $N[ {\mathcal U}]$ } \ar[uu]_{pull-back}}$$

\end{lem}
\begin{proof}
Let  $(\mathcal U , \bullet , \star)$  (resp. $(\mathcal V , \bullet' , \star')$) be the adapted Lie groupoid $\G \to \HH$-extension associated to the $\G \to \HH$-valued non-Abelian $1$-cocycle $(\lambda,\gggg)$ (resp. $(\lambda',\gggg')$, the pull-back of $(\lambda,\gggg)$ w.r.t. $\sigma$). Let $\sigma $ stand for the map $ \coprod_{j \in J} V_j\to \coprod_{i\in I} U_i$ (resp. $ \coprod_{k,l \in J} V_{kl}\to \coprod_{i,j\in I} U_{ij}$). The pull-back Lie groupoid  $ (\G \times \coprod_{i,j \in I}U_{ij})[\sigma] $
is isomorphic to $\G \times \coprod_{i,j \in  J}V_{ij}$ through the isomorphism defined for all $i,j \in J, x \in V_{ij}, g \in \G$ by
 $$ \psi (x_{i} ,(g,x_{\sigma(i)\sigma(j)}),x_{j}):=(g,x_{ij}),$$
 and the map  $$\psi' (x_{i} ,(x_{\sigma(i)},h)):=(x_i, h)$$ is an isomorphism between the pull-back of $\coprod_{i \in I} U_i \times \HH $  through $\sigma $
  and $ \coprod_{j\in J} V_j \times \HH $. We leave it to reader to prove that $(\psi,\psi',id_{\HH})$ is an isomorphism of Lie groupoid $\G \to \HH$-extensions.
\end{proof}

The next lemma shall also have its importance. The reader can replace the Lie groupoid $\BBB \toto \BBB_0$ by $N \toto N $ for the sake of simplicity, since we shall only use the lemma in that case.


\begin{lem} \label{lem:lemma B}
Let $(q,{\mathcal R} \stackrel{\phi}{\to} \BBB[q],P \to M,\chi)$ be a Lie groupoid $\G \to \HH$-extension over $\BBB \toto \BBB_0$. Let $\tau: M' \to M $ be a map such that $ q \circ \tau $ is a surjective submersion.
Then:
\begin{enumerate}
\item $\tau $ is a generalized surjective submersion for both Lie groupoids $ \mathcal{R}\toto M $ and $\BBB[q] \toto M $,
\item $( q \circ \tau, \mathcal{R}[\tau] \stackrel{\phi[\tau]}{\to}\BBB[q \circ \tau ],\tau^* P \to M', \chi[\tau])$
is a Lie groupoid $\G \to \HH$-extension,
\item this Lie groupoid $\G \to \HH$-extension is Morita equivalent (over the identity of $\BBB \toto \BBB_0$) to $(q,\mathcal{R} \stackrel{\phi}{\to} \BBB[q],P \to M,\chi)$.
  \end{enumerate}
\end{lem}
\begin{proof}
We wish to show that the map $\xi: M' \times_{\tau,M,s} B[q]\to M $ given by $(b,m)\mapsto t(b)$ for all $b\in B[q], m\in M$ is a surjective submersion. Let $m\in M$, take $m'\in (q\circ \tau)^{-1}(q(m))$ (which is non-empty by assumption). Now since $t^{-1}(q(m))$ is not empty so $(m',(\tau(m'),b,m))$
projects on $m$ by $\xi $, where $b\in t^{-1}(q(m)) $. This proves the surjectivity. To check that $\xi $ is indeed a submersion, we have to think in terms of infinitesimal paths.
Let $(m' , (\tau(m'),b,m ))  $ be a point in $M' \times_{\tau,M,s} B[q]$, and $m \in M$ such that $ \xi(m' , (\tau(m'),b,m ))=m$. Let
 $m (\epsilon) $ be a path in $M$ starting from $m $. Since the target map (of Lie groupoid $\BBB \toto \BBB_0 $) is always a surjective submersion there exists a path $b(\epsilon)$ in $\BBB $ starting at $b$ such that $t(b(\epsilon))=q(m(\epsilon)) $ (for all $\epsilon $ small enough). Since $q \circ \tau $ is a surjective submersion by assumption, there exists also a path $m'(\epsilon) $ in $M'$ starting at $m'$
such that $ q \circ \tau (m'(\epsilon)) =s(b(\epsilon)) $ for $\epsilon  $ small enough. By construction, the path $ (m'(\epsilon),\tau(m'(\epsilon), b(\epsilon), m(\epsilon) ))  $ is a path in  $ M' \times_{\tau,M,s} B[q]\to M$ which starts at $(m' , (\tau(m'),b,m ))  $ and projects by $\xi$ onto $m(\epsilon) $, which completes the proof of the first item.

In view of the proof of lemma \ref{lem: generalized surjective submersion}, all the algebraic axioms of Lie groupoid $\G \to \HH$-extensions are satisfied
by $(\mathcal R[\tau] \stackrel{\phi}{\to}\BBB[q \circ \tau ],\tau^* P \to M',\tau^* \chi)$. Lemma \ref{lem:manifold} implies that the sets involved are manifolds. This completes the proof of the second item.

For the last item, the manifold that we shall consider to construct an explicit Morita equivalence is:
   $$  T= M' \times_{ \tau, M, t } \mathcal R $$
   equipped with the surjective submersions  $q_M' :T \to M'$ and $q_M :T \to M$ given by the projection on the first component and the target of the second component respectively.
    By construction, the following diagram commutes:
    $$ \xymatrix{ & T \ar[dr]_{q_{M'}} \ar[dl]^{q_{M}}& \\ M \ar[dr]^{q} & & M' \ar[ll]_{\sigma} \ar[dl]_{q \circ \tau} \\ & B& \\ }. $$
    This implies that $ \mathcal R[\tau][q_{M'}] \simeq \mathcal R[ \tau \circ q_{M'}] = \mathcal R[q_M] $ and also $$ p_{M'}^* \sigma^* P \simeq  (\sigma \circ q_{M'})^* P= p_M^* P.$$
     It is routine to check that this pair of isomorphisms form an isomorphism of Lie groupoid $\G \to \HH$-extensions between the pull-back of $( q \circ \tau, \mathcal R[\sigma] \stackrel{\phi}{\to}\BBB[q \circ \tau ],\sigma^* P \to M',\tau^* \chi , )$ w.r.t. $q_{M'}$ and the pull-back of $( q , \mathcal R \stackrel{\phi}{\to}\BBB[q ],P \to M',\chi , )$ w.r.t. $q_M$.

\end{proof}
We now prove theorem  \ref{th:coho=gerbes}.
\begin{proof}
According to the first item of proposition \ref{prop:cocycles},
 to an arbitrary $\G \to \HH$-valued non-Abelian $1$-cocycle $( \lambda,\gggg) $ with respect to an arbitrary open cover ${\mathcal U}$ corresponds an
 adapted Lie groupoid $\G \to \HH$-extension (which is by construction a Lie groupoid $\G \to \HH$-extension above the Lie groupoid $N \toto N $).

This assignment goes to the quotient to yield an assignment from  $\G \to \HH$-valued non-Abelian $1$-cohomology on $N$ to $\G \to \HH$-gerbes over $N \toto N $. This follows from the fact that the adapted Lie groupoid $\G \to \HH$-extensions associated to a
$\G \to \HH$-valued non-Abelian $1$-cocycle and a pull-back of it are Morita equivalent over the identity of $N \toto N $ by Lemma \ref{lem:lemma A}.
Also, by proposition \ref{prop:coboundaries},  the adapted extensions associated to two $\G \to \HH$-valued
non-Abelian $1$-cocycles that differ by a coboundary are isomorphic, hence Morita equivalent over the identity of $N \toto N$
by example \ref{ex:isomOver}.
Hence, the adapted Lie groupoid $\G \to \HH$-extensions associated to two
$\G \to \HH$-valued non-Abelian $1$-cocycles that define the same element in cohomology are  Morita equivalent over the identity of $N \toto N$, yielding a well-defined map from  $\G \to \HH$-valued non-Abelian $1$-cohomology on $N$ to $\G \to \HH$-gerbes over $N \toto N$, that we denote by $\Xi $.

 We first check that $\Xi $ is surjective. Let $(q,\mathcal R \stackrel{\phi}{\to} N[q],P \to M,\chi)$ be an arbitrary Lie groupoid $\G \to \HH$-extension over $N \toto N$. There exists an open cover $\mathcal U=(U_i)_{i \in I}$ of $N$ such that $q: M \to N $
 admits local sections $\sigma_i : U_i \to M $ for all $i \in I $, which, altogether, define a map $\sigma: \coprod_{i \in I} U_i \to M $.
 By lemma \ref{lem:lemma B}, $(q,\mathcal R \stackrel{\phi}{\to} N[q],P \to M,\chi)$  is Morita equivalent over the identity of $N \toto N$ to
 its pull-back with respect to $\sigma $. The pull-back being a Lie groupoid $\G \to \HH$-extension of the \v{C}ech groupoid is, by proposition \ref{prop:everyextensionisadapted}, isomorphic (hence Morita equivalent by example \ref{ex:isom}) to an adapted one. Hence $(q,\mathcal R \stackrel{\phi}{\to} N[q],P \to M,\chi)$ is Morita equivalent to an adapted Lie groupoid $\G \to \HH$-extension, which, by proposition  \ref{prop:cocycles}, comes from some non-Abelian $1$-cocycle. This proves that the assignment $\Xi $ is surjective.

 We then check that $ \Xi$ in injective. The proof is based on  the following general property of open covers.
 Assume that there is a commutative diagram as follows:
  $$  \xymatrix{  & M' \ar[dr]^p \ar[dl]_{p'} & \\ \coprod_{i \in I} U_i \ar[dr]_\imath & & \coprod_{j\in J} V_j \ar[dl]^\imath \\ & N &   }$$
  with $p$ and $p'$ surjective submersions (above, the symbol $ \imath$ stands for all the canonical inclusions, and  $(U_i)_{i \in I} $ and $(V_j)_{j \in J} $ are open covers of the manifold $N$). Then there is a common refinement $(W_k)_{k \in K} $ of $(U_i)_{i \in I} $ and $(V_j)_{j \in J} $ and a map $\tau: \coprod_{k \in K} W_k  \to M'$ such that the following diagram commutes:
   \begin{equation}\label{eq:commonrefinement}\xymatrix{ & \coprod_{k \in K} W_k \ar[d]^\tau \ar[ddr]^\imath \ar[ddl]_\imath & \\ & M' \ar[dr]_p \ar[dl]^{p'} & \\ \coprod_{i \in I} U_i \ar[dr]^\imath & & \coprod_{j\in J} V_j \ar[dl]_\imath \\ & N &   }
   \end{equation}
  where, again, we use the symbol $\imath $ to denote all the canonical inclusions.

Assume now two $\G \to \HH$-valued non-Abelian $1$-cocycles, defined w.r.t. open covers  $(U_i)_{i \in I} $ and $(V_j)_{j \in J} $ respectively, have adapted Lie groupoid $\G \to \HH$-extensions $Ext1 $  and $Ext2 $ (which are hence over the Lie groupoid $N\toto N$) associated with which are Morita equivalent.
By the very definition of Morita equivalence of $\G \to \HH$-extensions, this implies that there exists a manifold $M'$ together with surjective submersions $p :M' \to \coprod_{i\in I}U_i$ and $p' :M' \to \coprod_{j\in J}V_j$  such that $\iota \circ p = \iota \circ p' $ and such that the pull-backs $Ext1[p] $  and $Ext2[p'] $ of both extensions
to $M'$ are isomorphic. According to the discussion above, there exists a common refinement
 $\coprod_{k \in K} W_k$ of both open covers and a map  $\tau: \coprod_{k \in K} W_k \to M' $
   such that the diagram (\ref{eq:commonrefinement}) commutes. According to lemma \ref{lem:lemma B}, the pull-back of the adapted extension $Ext1 $
   on $\coprod_{k \in K} W_k$ is isomorphic to the pull-back of $Ext1[p] $ by $\sigma $. Similarly,  the pull-back of the adapted extension $Ext2 $ on $\coprod_{k \in K} W_k$ is isomorphic to the pull-back of $Ext2[p] $ by $\sigma $. Since $Ext1[p] $ and $Ext2[p] $ are isomoprhic, this implies that the pull-back of $Ext1 $
   and $Ext2 $ to $\coprod_{k \in K} W_k$ are isomorphic.
   According to Lemma \ref{lem:lemma A}, this means that the pull-back of both cocycles to $(W_k)_{k \in K} $ have corresponding adapted extensions that are isomorphic. By proposition \ref{prop:coboundaries}, it means that their pull-back to  $(W_k)_{k \in K} $ differ by a coboundary, i.e. that both cocycles define the same class in cohomology. This proves the injectivity of $\Xi $.
\end{proof}

\subsection{$\G \to \HH$-gerbes over differentiable stacks}

Recall that a Morita equivalence between two Lie groupoids $ \BBB \toto \BBB_0$
and $\BBB' \toto \BBB_0' $ is a quadruple ${\mathcal M}=(T,f,g,\Phi) $, with $T$ a manifold, $f,g$ surjective submersions from $T$ to $\BBB_0$
 and to $\BBB_0' $ respectively, and $\Phi$ a Lie groupoid isomorphism over the identity of $T$ between
 $\BBB[f] \toto T$ and $\BBB'[g] \toto T $.
 (Alternatively, Morita equivalence may be defined with the help of the notion of bi-modules, a description which happens to be
  equivalent to the previous one, see \cite{BehrendXu}.)
  Morita equivalent Lie groupoids often share similar properties, in particular which regards to cohomology.
  The next theorem shows that they also have the same gerbes over them.

\begin{them}\label{thm:morita_and_gerbes}
A Morita equivalence between two groupoids $ \BBB \toto \BBB_0$
and $\BBB' \toto \BBB_0' $ induces a one to one correspondence between:
\begin{enumerate}
\item $\G \to \HH$-gerbes over $ \BBB \toto \BBB_0 $,
\item $\G \to \HH$-gerbes over $ \BBB' \toto \BBB_0' $.
\end{enumerate}
\end{them}
\begin{proof}
Let ${\mathcal M}=(T,f,g,\Phi)$ be a Morita equivalence between the Lie groupoids $ \BBB \toto \BBB_0$ and $\BBB' \toto \BBB_0' $, i.e. $f:T \to \BBB_0$ and $g:T \to \BBB'_0$ are surjective submersions and $\Phi: \BBB[f] \to \BBB'[g]$ is an isomorphism of Lie groupoids between the the pull-back groupoids $ \BBB[f] \toto T$ and $\BBB'[g] \toto T $.
We intend to assign to an arbitrary Lie groupoid $\G \to \HH$-extension $Ext:=(q, X\stackrel{\phi}{\to} \BBB[q], P\to M, \chi)$ over $ \BBB \toto \BBB_0 $ a  Lie groupoid $\G \to \HH$-extension over $\BBB' \toto \BBB_0'$.
 To start with, we consider the set $M':= M\times_{q,\BBB_{0},f}T$. One checks easily that $M'$ is a manifold such that the projections $\alpha,\beta $ onto the first and second components are surjective submersions,
 and as well as the maps $q \circ \alpha : M' \to \BBB_0$ and $q':= g \circ \beta : M'\to \BBB_0' $.
Then we consider the pull-back of $Ext $ by $\alpha $, namely
 $$ ( \mathcal R[\alpha]\stackrel{\phi[\alpha]}{\to} \BBB[q \circ \alpha], \alpha^* P\to  M', \chi[\alpha]) .$$
We claim that there exists an isomorphism $\Phi': \BBB[q \circ \alpha] \to \BBB'[g \circ \beta]$, so that
$$ (g \circ \beta ,  \mathcal R[\alpha]\stackrel{\Phi ' \circ \phi[\alpha]}{\longrightarrow} \BBB'[g \circ \beta], \alpha^* P\to  M', \chi[\alpha]) $$
  is a Lie groupoid $\G \to \HH$-extension over $\BBB' \toto \BBB_0' $.
For all $((m_1,t_1),b,(m_2,t_2)) \in M' \times_{q \circ \alpha , \BBB_0 , s} B  \times_{t, \BBB_0,q \circ \alpha } M'$, we have the relations
  $$ f(t_1)=q(m_1), f(t_2)=q(m_2) , s(b)=q \circ \alpha (m_1,t_1)=q(m_1) , t(b)=q \circ \alpha (m_2,t_2)=q(m_2) $$
which impliy that $f(t_1)=s(b), f(t_2)=t(b)$, hence that $(t_1,b,t_2) \in \BBB[f]$. This allows us to set
  $$ \Phi'((m_1,t_1),b,(m_2,t_2)):=((m_1,t_1),b',(m_2,t_2))$$
where $b' \in \BBB$ is given by $\Phi(t_1,b,t_2)=(t_1,b',t_2)$.
By construction,  $((m_1,t_1),b',(m_2,t_2)) $ is in $B'[ g \circ \beta]$,
and it is routine to check that $\Phi'$ is an isomorphism of Lie groupoids.

 %
 We have therefore assigned  a Lie groupoid $\G \to \HH$-extension over
 $\BBB' \toto \BBB_0' $ to a Lie groupoid  $\G \to \HH$-extension over
 $\BBB \toto \BBB_0 $, as intended. In picture:
            $$
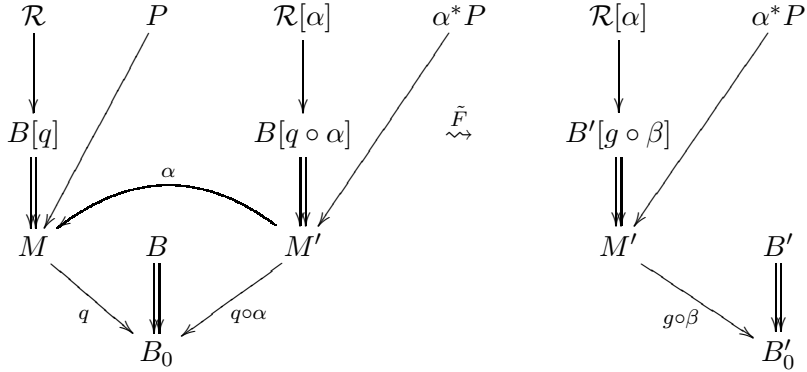
\xymatrix{\mathcal R \ar[d] & P \ar[ddl] & \mathcal R[\alpha] \ar[d] & \alpha ^* P \ar[ddl] & \mathcal R[\alpha] \ar[d] & \alpha ^* P \ar[ddl] \\
                      B[q]\ar@<1pt>[d] \ar@<-1pt>[d]& & B[q \circ \alpha] \ar@<1pt>[d] \ar@<-1pt>[d] & \stackrel{\tilde{F}}{\rightsquigarrow} & B'[g \circ \beta]\ar@<1pt>[d] \ar@<-1pt>[d]\\
                      M \ar[dr]_{q} & B \ar@<1pt>[d] \ar@<-1pt>[d] & M' \ar@/_2pc/[ll]_{\alpha}\ar[dl]^{q \circ \alpha} & & M' \ar[dr]_{g \circ \beta}& B' \ar@<1pt>[d] \ar@<-1pt>[d]\\
                         & B_{0} & & & & B'_{0}  } $$
                         $$\alpha(m,n)=m, q'(m,n)=g(n)$$
 The same construction could be done to assign a $\G \to \HH$-extension over
 $\BBB' \toto \BBB'_0 $ to a $\G \to \HH$-extension over
 $\BBB \toto \BBB_0 $.
Since the roles of $\BBB \toto \BBB_0 $ and
 $\BBB' \toto \BBB_0' $ can be exchanged, in order to check that both assignments induce a one to one correspondence between the corresponding gerbes, it is necessary and sufficient to check
 that:
\begin{itemize}
\item[(i)] Morita equivalent Lie groupoid $\G \to \HH$-extensions over the Lie groupoid $\BBB \toto \BBB_0 $ are mapped to Morita equivalent Lie groupoid $\G \to \HH$-extensions over the Lie groupoid $\BBB' \toto \BBB'_0 $ by the first assignment,
\item[(ii)] applying the first assignment, then the second one to a Lie groupoid $\G \to \HH$-extension $Ext $ over $\BBB \toto \BBB_0$ yields a $\G \to \HH$-extension which is Morita equivalent to $Ext $.
\end{itemize}
Let us check these two points. The second one is an immediate consequence of example \ref{ex:pullbackOver}, since a Lie groupoid $\G \to \HH$-extension over $\BBB \toto \BBB_0 $ is always Morita equivalent to its pull-back.
The first one is more involved. Let $Ext_1,Ext_2$ be two Lie groupoid $\G \to \HH$-extensions over $\BBB \toto \BBB_0 $, namely, to fix notations
 $$  Ext_i:=( q_i, \mathcal R_i\stackrel{\phi_i}{\to} \BBB[q_i], P_i\to M_i, \chi_i )  \hbox{ , } i=1,2, $$
 and let
 $$Ext_i':= (g \circ \beta ,  \mathcal R_i[\alpha]\stackrel{\Phi ' \circ \phi_i[\alpha]}{\longrightarrow} \BBB'[g \circ \beta], \alpha^* P\to  M'_i, \chi_i[\alpha]) \hbox{ , } i=1,2. $$
be the associated Lie groupoid $\G \to \HH$-extensions over $\BBB' \toto \BBB'_0 $ constructed as above.
Assume now that $Ext_i \hbox{ , } i=1,2.  $ are Morita equivalent. This means, first, that there is a commutative diagram of surjective submersions:
  $$\xymatrix{ & M \ar[dr]^{p_2} \ar[dl]_{p_1}& \\ M_1 \ar[dr]_{q_1} & & M_2 \ar[dl]^{q_2} \\ & \BBB_0& \\ }   $$
 This implies that the following is also a commutative diagram of surjective submersions:
  $$\xymatrix{ &  M\times_{\BBB_{0}}T\ar[dr]^{(p_2,id_T)} \ar[dl]_{(p_1,id_T)}& \\ M'_1 \ar[dr]_{g \circ \beta} & & M'_2 \ar[dl]^{g \circ \beta} \\ & \BBB_0& \\ }   $$
  where $M'_i:=M_i \times_{\BBB_{0}}T$ and the fibred product $M\times_{\BBB_{0}}T$ are considered w.r.t. the maps $q_1 \circ p_1 (=q_2 \circ p_2): M \to \BBB_0$ and $f: T \to \BBB_0$.

  To show that both pull-back of $Ext'_i$ w.r.t. $(p_i,id_T)$ with $i=1,2$ are isomorphic as $\G \to \HH $-extensions, it suffices to show that $\mathcal R_i[\alpha \circ (p_i,id_T)] \toto M'_i$ with $i=1,2$ are isomorphic Lie groupoids. But this is a consequence of the existence of isomorphism between $\mathcal R_i[p_i] \toto M_i \hbox{ , } i=1,2$ which is in turn a consequence of the assumption of $(M,p_1,p_2,\phi)$ being a Morita equivalence between $Ext_1, Ext_2$.

\end{proof}

Let us denote by $F({\mathcal M}) $ the correspondence between $\G \to \HH$-gerbes associated to a given Morita equivalence ${\mathcal M} $.
It is easy to check that, given ${\mathcal M}_1$ and ${\mathcal M}_2$ composable Morita equivalences, the relation  $F({\mathcal M}_1 ) \circ F( {\mathcal M}_2)
 = F({\mathcal M}_1 \circ {\mathcal M}_2) $ holds when the composition ${\mathcal M}_1 \circ {\mathcal M}_2$ of Morita equivalences,
 is defined as in \cite{BehrendXu}.
According to \cite{BehrendXu}, Lie groupoids up to Morita equivalences are one possible description  of differential stacks,
so that theorem \ref{thm:morita_and_gerbes} makes sense of the notion of $\G \to \HH$-gerbes valued in a differential stack.

\chapter{Nijenhuis forms on $L_\infty$-algebras}

This chapter is devoted to the equivalent of Nijenhuis tensors on $L_\infty$-algebras.
We refer to the Introduction, section \ref{Int:1stglanceatLinfty}, for a discussion motivating this concept.
The way we shall handle this problem requires to re-interpret $L_\infty$-structures in terms of Poisson elements
for the Richardson-Nijenhuis bracket, a notion that we now study.
\section{Richardson-Nijenhuis bracket and $L_\infty$-structures.}\label{section:3.1}

The purpose of this section is to introduce the Richardson-Nijenhuis bracket on arbitrary graded vector spaces,
and to give some of its properties. The Richardson-Nijenhuis bracket of vector-valued forms
on a vector space was introduced in \cite{NR}. Then, it was extended for vector-valued forms on a manifold and was used, for instance,
 in \cite{KMS} which will be the most important reference in this section.

\subsection{Graded symmetric spaces}\label{subsection:3.1.1}

We first recall and fix notations for symmetric algebras on graded vector spaces. Let $E$ be a graded vector space over a field $\mathbb{K}=\mathbb{R}$ or $\mathbb{C}$, that is a vector space of the form
$$ \oplus_{i \in {\mathbb Z}} E_i. $$
For a given $i\in \mathbb{Z}$, the vector space $E_i$ is called the component of degree $i$,
elements of $E_i$ are called \emph{homogeneous elements of degree $i$}, and elements in the union $\cup_{i\in \mathbb{Z}} E_i$
are called the \emph{homogeneous elements}. Notice that only zero can be homogeneous of two different degrees,
which allows to denote by $|X|$ the degree of a non-zero homogeneous element $X$, a convention that we will often use in an implicit manner. Given a graded vector space $E=\oplus_{i\in \mathbb{Z}}E_i$ and an integer $p$, one may shift all the degrees
by $p$ to get a new grading on the vector space $E$. We use the notation $E[p]$ for the graded vector space $E$ after shifting the degrees
 by $p$, that is the graded vector space whose  component of degree $i$, is $E_{i+p}$.

  We denote by $\otimes E$ the tensor algebra of $E$, with product given by concatenation. The symmetric space of $E$, denoted by $S(E)$, is by definition, the quotient space of the tensor algebra $\otimes E$ by the two-sided idea
l $I \subset \otimes E$  generated by elements of the type $X\otimes Y-(-1)^{|X||Y|}Y\otimes X$, with $X$ and $Y$ arbitrary homogeneous elements in $E$.
For a given $k\geq 0$, $S^k(E)$ is the image of $\otimes^{k} E$ through the quotient map $\otimes E \mapsto \frac{\otimes E}{I} = S(E)$.
Of course, the decomposition
   $$ S(E) = \oplus_{k \geq 0} S^k (E) $$
    holds and $S^0(E)$ is simply the field $\mathbb{K}$.
  Moreover, when all the components in the graded space $E$ are of finite dimension, the dual of $S^k(E)$ is isomorphic to $S^k(E^*)$, for all $k\geq 0$. In this case,  there is a one to one correspondence between
\begin{enumerate}
\item graded symmetric $k$-linear maps on the graded vector space $E$,
\item linear maps from the space $S^k(E)$ to $E$,
 \item $S^k(E^*)\otimes E$.
\end{enumerate}
Elements of the space $S^k(E^*)\otimes E$ are called \emph{symmetric vector valued  $k$-forms}. Notice that $S^0(E^*)\otimes E$, the space of vector valued zero-forms, is isomorphic to the space $E$.

 Having the decomposition $S(E)=\oplus_{k\geq 0}S^k(E)$, every element in $S(E)$ is the sum of finitely many elements in $S^k(E)$, $k \geq 0$.
 We absolutely need to consider also infinite sums, which is often referred to, in the literature, as
 taking the completion of $S(E)$. By a \emph{formal sum}, we mean a sequence  $\phi: \mathbb{N}\bigcup \{0\} \to S(E) $ mapping an integer $k$
 to an element $a_k \in S^k(E)$: we shall, by a slight abuse of notation, denote by $\sum_{k=0}^\infty a_k$
 such an element. We denote the set of all formal sums by $\tilde{S}(E)$.  The algebra structure on $S(E)$ extends in an unique manner to $\tilde{S}(E)$.
 For two formal sums $a=\sum_{k=0}^\infty a_k$ and $b=\sum_{k=0}^\infty b_k$ we define $a+b$ to be $\sum_{k=0}^\infty (a_k+b_k)$,
 while the product of $a$ and $b$ is the infinite sum $\sum_{k=0}^\infty c_k$ with
   $ c_k = \sum_{i=0}^k a_i \cdot b_{k-i}  $ (with $\cdot$ being the product of $S(E)$).

When all the components in the graded space $E$ are of finite dimension, there is a one to one correspondence between
\begin{enumerate}
\item collections indexed by ${k\geq 0}$ of graded symmetric $k$-linear maps on the graded vector space $E$,
\item collections  indexed by ${k\geq 0}$ of linear maps from $S^k(E)$ to $E$,
 \item $\tilde{S}(E^*)\otimes E$.
\end{enumerate}
Elements of the space $\tilde{S}(E^*)\otimes E$ are called \emph{ symmetric vector valued forms} and shall be written as infinite sums $\sum_{i\geq 0} K_i$ with $K_i \in S^i(E^*)\otimes E$. For all  $m\geq 0$,  $\tilde{S}^{k\geq m}(E^*)\otimes E$ denotes the space of all symmetric vector valued $k$-forms, with $k\geq m$.

\subsection{Richardson-Nijenhuis bracket}\label{subsection:3.1.2}

We first define the insertion operator, following \cite{KMS} as the guideline, but with different sign conventions.

\begin{defi}\label{def:insertion}
Let $E$ be a graded vector space, $E= \oplus_{i\in \mathbb{Z}}E_i$. The insertion operator of a symmetric vector valued $k$-form $K$ is denoted by $\iota_{K}$ and is an operator
$$ \iota_K : S(E^*)\otimes E \to S(E^*)\otimes E$$
defined by:
\begin{equation}\label{equ:insertion}
\iota_{K}L(X_1,...,X_{k+l-1})=\sum_{\sigma \in Sh(k,l-1)}\epsilon(\sigma)L(K(X_{\sigma(1)},...,X_{\sigma(k)}),...,X_{\sigma(k+l-1)}),
\end{equation}
for all $L \in S^{l}(E^*)\otimes E $, $l\geq 0$ and $X_1,\cdots X_{k+l-1} \in E$, where $\epsilon(\sigma)$ is the Koszul sign.
\end{defi}
If $L$ is an element in $S^0(E^*)\otimes E \simeq E$, then (\ref{equ:insertion}) should be understood as meaning that
$\iota_{K}L=0$, for all vector valued forms $K$ and
\begin{equation}\label{eq:insertionof0}
\iota_{L}K (X_1,...,X_{k-1})=K(L,X_1,...,X_{k-1}),
\end{equation}
for all vector valued $k$-form $K$.

We shall extend the previous definition of the insertion operator, allowing $L$ and $K$ to be  symmetric vector valued forms, as follows
\begin{equation}\label{inifinitsums}
\begin{array}{rcl}
\iota_{K_1+K_2+K_3+\dots}(L_1+L_2+L_3+\dots )
                                             &=& \iota_{K_1}L_1+\iota_{K_1}L_2+\dots\\
                                             && +\iota_{K_2}L_1+\iota_{K_2}L_2+\dots \\
                                             &&+ \dots\\
                                             && +\iota_{K_k}L_1+\iota_{K_k}L_2+\dots\\
                                             &&+\dots\\

\end{array}
\end{equation}
with $K_i,L_i \in S^{i}(E^*)\otimes E , i\geq 0.$
The Equation (\ref{inifinitsums}) above makes sense since, for all  $m \geq 0$, the component in $S^{m} (E^*) \otimes E $ of the right hand side
is only a finite sum.

The insertion operator $\imath_{K}L $ can also be defined using an expression similar to (\ref{equ:insertion})
but where $L$ is now an element in $S^l(E^*)$, i.e. a linear form on $S^l(E) $.
This procedure extends to infinite sums, and one obtains for every $K \in \tilde{S}(E^*)\otimes E$
a linear operator on $ \tilde{S}(E^*)$.


\begin{lem}\label{lem:insertion1-1}
The insertion operator $ \iota_K: \tilde{S}(E^*) \to \tilde{S}(E^*) $, with $K \in \tilde{S}(E^*)\otimes E$,
is equal to zero if and only if $K=0.$
\end{lem}

Insertion operation allows us to define the Richardson-Nijenhuis bracket.
\begin{defi}
 Given a symmetric vector valued $k$-form $K\in S^{k}(E^*)\otimes E$ and a symmetric vector valued $l$-form $L\in S^{l}(E^*)\otimes E$, the Richardson-Nijenhuis bracket of $K$ and $L$ is the symmetric vector valued $(k+l-1)$-form $[K,L]_{_{RN}}$, defined as
\begin{equation*}
[K,L]_{_{RN}}=\iota_{K}L-(-1)^{\bar{K}\bar{L}}\iota_{L}K,
\end{equation*}
where $\bar{K}$ is the degree of $K$ as a graded map, that is $K(X_1,\cdots, X_k)\in E_{1+\cdots+k+\bar{K}},$ for all $X_i \in E_i$.
\end{defi}
For an element $X \in E$, $\bar{X}=|X|$, that is, the degree of a vector valued $0$-form, as a graded map, is just its degree as an element of $E$.
In \cite{LMS}, the authors defined a multi-graded Richardson-Nijenhuis bracket, in a graded vector space, but with a different approach of ours.
\begin{them}\label{theo:RNisgradedLiebracket}
The space $\tilde{S}(E^*) \otimes E$, equipped with the Richardson-Nijenhuis bracket, is a graded Lie algebra.
\end{them}
\begin{proof}The space $S(E^*)$ of symmetric forms on $E$, is a graded algebra, where the product is concatenation and the grading is defined as follows. For every $k\geq 0,$ and every symmetric $k$-form $\alpha \in S^k(E^*)$, $\bar{\alpha}$ is the degree of $\alpha$ as a graded map, that is $\alpha(X_1,\cdots, X_k)\in E_{1+\cdots+k+\bar{\alpha}},$ for all $X_i \in E_i$. Now let $S^{\bar{k}}$ denotes the space of all symmetric forms of degree $k$. Then
\begin{equation*}
S(E^*)=\oplus_{k\geq 0}S^{\bar{k}}(E^*).
\end{equation*}
Let $Der(E)$ be the space of graded derivations on $S(E^*)$. We consider the graded commutator of two derivations $D_1$ and $D_2$ given by
\begin{equation*}
[D_1,D_2]=D_1\circ D_2-(-1)^{\bar{D_1}\bar{D_2}}D_2\circ D_1,
\end{equation*}
where for a graded derivation $D$, $\bar{D}$ stands for degree of $D$. By definition,
\begin{equation*}
D(\alpha)\in S^{\overline{\bar{D}+\bar{\alpha}}}(E^*),
\end{equation*}
for all $\alpha \in S(E^*)$. It is known that the space of all graded derivations on a graded vector space $E$, together with the graded commutator, is a graded Lie algebra.
Let $K$ be a vector valued $k$-form, $L$ be a vector valued $l$-form and $\omega$ be a $1$-form on the graded vector space $E$. It is clear from Definition \ref{def:insertion} that $ \iota_K$ and $ \iota_L$ are derivations of degree $\bar{K}$ and $\bar{L}$ respectively, of the space $\tilde{S}(E^*)$. Therefore we have
\begin{equation*}
\begin{array}{rcl}
                                             & & ([\iota_K,\iota_L]\omega)(X_1,\cdots,X_{k+l-1})\\
                                             &=&((\iota_K \circ \iota_L -(-1)^{\bar{K}\bar{L}}\iota_L \circ \iota_K)\omega)(X_1,\cdots,X_{k+l-1})\\
                                             &=&\sum_{\sigma \in Sh(k,l-1)}\epsilon(\sigma)\iota_L \omega (K(X_{\sigma(1)},\cdots,X_{\sigma(k)}),\cdots,X_{\sigma(k+l-1)})\\
                                             & &-(-1)^{\bar{K}\bar{L}}\sum_{\sigma \in Sh(l,k-1)}\epsilon(\sigma)\iota_K \omega (L(X_{\sigma(1)},\cdots,X_{\sigma(l)}),\cdots,X_{\sigma(k+l-1)})\\
                                             &=&\sum_{\sigma \in Sh(k,l-1)} \epsilon(\sigma)\omega (L(K(X_{\sigma(1)},\cdots,X_{\sigma(k)}),\cdots,X_{\sigma(k+l-1)}))\\
                                             & &-(-1)^{\bar{K}\bar{L}}\sum_{\sigma \in Sh(l,k-1)}\epsilon(\sigma) \omega (K(L(X_{\sigma(1)},\cdots,X_{\sigma(l)}),\cdots,X_{\sigma(k+l-1)})).\\
\end{array}
\end{equation*}
On the other hand
\begin{equation*}
\begin{array}{rcl}
                   & &(\iota_{[K,L]_{_{RN}}}\omega)(X_1,\cdots,X_{k+l-1})\\
                   &=& \omega([K,L]_{_{RN}}(X_1,\cdots,X_{k+l-1}))\\
                   &=&\omega(\iota_KL-(-1)^{\bar{K}\bar{L}}\iota_LK)(X_1,\cdots,X_{k+l-1})\\
                   &=&\omega\big((\sum_{\sigma \in Sh(k,l-1)}\epsilon(\sigma)L(K(X_{\sigma(1)},\cdots,X_{\sigma(k)}),\cdots,X_{\sigma(k+l-1)}))\\
                   & &-(-1)^{\bar{K}\bar{L}}(\sum_{\sigma \in Sh(l,k-1)}\epsilon(\sigma)K(L(X_{\sigma(1)},\cdots,X_{\sigma(l)}),\cdots,X_{\sigma(k+l-1)}))\big).
 \end{array}
\end{equation*}
Computations above show that
\begin{equation}\label{LiemorphismDR-N}
\iota_{[K,L]_{_{RN}}}=[\iota_{K},\iota_{L}]
\end{equation}
on $1$-forms. Since derivations are completely determined on $1$-forms, (\ref{LiemorphismDR-N}) holds on the space of forms.
  Equation (\ref{LiemorphismDR-N}) together with Lemma \ref{lem:insertion1-1} show that the Richardson-Nijenhuis bracket is a graded Lie bracket such that the insertion operator is a Lie morphism from the graded Lie algebra of vector valued forms $(\tilde{S}(E^*) \otimes E,\, \left[.,.\right]_{_{RN}})$ to the graded Lie algebra of derivations on forms $(Der(E),\, \left[.,.\right])$.
  \end{proof}

The following, which we use in a slightly implicit manner, in the sequel, is a direct consequence of definitions together with Theorem \ref{theo:RNisgradedLiebracket}.
\begin{cor}
The space $\tilde{S}^{\geq 1} (E^*) \otimes E$ is a sub-graded Lie algebra of the graded Lie algebra $(\tilde{S}(E^*) \otimes E, \left[.,.\right]_{_{RN}})$. Also, vector valued $0$-forms, that is, elements of $ S^0(E^*)\otimes E$, form an abelian sub-graded Lie algebra of the graded Lie algebra $(\tilde{S}(E^*) \otimes E, \left[.,.\right]_{_{RN}})$.
\end{cor}

Note that for a symmetric vector valued $k$-form $K$, $K(X_1,\dots ,X_k)$ can be recovered by taking successive Richardson-Nijenhuis brackets with $X_1,\dots,X_k$, seen as vector valued $0$-forms,
a point of view that may be useful in computations.

\begin{prop}\label{DrivedBt}
For a vector valued $k$-form $K\in S^k(E^*)\otimes E$ we have:
\begin{equation*}
K(X_1,\cdots,X_k):=[X_k,\cdots,[X_2,[X_1,K]_{_{RN}}]_{_{RN}}\cdots]_{_{RN}}
\end{equation*}
for all $X_1,\cdots,X_k\in E$.
\end{prop}
\begin{proof}
The proof can be done using induction  and follows from the fact that for a vector valued $k$-form $K$
\begin{equation*}
[X_1,K]_{_{RN}}(X_2,\cdots ,X_k)=K(X_1,X_2,\cdots ,X_k)
\end{equation*}
as a consequence of (\ref{eq:insertionof0}), for all $X_1, X_2,\cdots, X_k \in E$.
\end{proof}
We will use the following easy lemma several times:
\begin{lem}\label{lem:[mu[k,mu]]=0}
Let $L$ be a vector valued form of odd degree and $K$ be a vector valued form of any degree such that $[L,L]_{_{RN}}=0$, then
\begin{equation}
[L,[K,L]_{_{RN}}]_{_{RN}}=[[L,K]_{_{RN}},L]_{_{RN}}=0.
\end{equation}
\end{lem}
\begin{proof}
The proof only uses the Jacobi identity for the Richardson-Nijenhuis bracket and it goes as follows
\begin{equation*}
\begin{array}{rcl}
[L,[K,L]_{_{RN}}]_{_{RN}}&=&[[L,K]_{_{RN}},L]_{_{RN}}+(-1)^{\bar{K}\bar{L}}[K,[L,L]_{_{RN}}]_{_{RN}}\\
                      &=&(-1)^{\bar{K}\bar{L}+\bar{L}(\bar{K}+\bar{L})}[L,[K,L]_{_{RN}}]_{_{RN}}\\
                      &=&(-1)^{\bar{L}^2}[L,[K,L]_{_{RN}}]_{_{RN}}.

\end{array}
\end{equation*}
Therefore, if $\bar{L}$ is an odd number, then $[L,[K,L]_{_{RN}}]_{_{RN}}=[[L,K]_{_{RN}},L]_{_{RN}}=0$.
\end{proof}

\subsection{ Characterization of $L_\infty$-structures in terms of Richardson-Nijenhuis bracket}
In the Introduction we gave the definition of an  $L_{\infty}$-algebra as a graded vector space together with a collection of graded symmetric multi-linear maps of degree $1$, satisfying a certain graded Jacobi identity (see Definition \ref{def:SymLinfty}).
In this section we will give a criterion that determines whether a vector valued form is an $L_{\infty}$-structure on a given graded vector space.
  First let us recall the definition of what is called a \emph{curved $L_{\infty}$-algebra}.
\begin{defi}\label{def:curvedSymLinfty}
A curved $L_{\infty}$-algebra is a graded vector space $E$ together with a family of symmetric vector valued forms $(l_i)_{i\geq 0}$ such that
\begin{enumerate}
\item $\bar{l_i}=1$ for all $i\geq 0$ ,
\item $l_1(l_0)=0$,
\item The graded Jacobi identity holds for the family $(l_i)_{i\geq 0}$, i.e.,
\begin{equation}
l_{n+1}(l_0,X_1,\cdots, X_{n})+\sum_{i+j=n+1}\sum_{\sigma \in Sh(i,j-1)} \epsilon(\sigma) l_j(l_i(X_{\sigma(1)},\cdots,X_{\sigma(i)}),\cdots,X_{\sigma(n)})
=0;
\end{equation}
for all $n\geq 1$ and all $X_1,\cdots,X_n$, where $\epsilon(\sigma)$ is the Koszul sign.
\end{enumerate}
The vector valued $0$-form $l_0$ of degree $1$, which is in fact an element in $E_1$, is called the curvature of the curved $L_{\infty}$-algebra.
\end{defi}

\begin{rem}
When in a curved $L_{\infty}$-algebra $l_0=0$, we get the notion of $L_{\infty}$-algebra.
\end{rem}
\begin{rem}
Equation (\ref{L_infty-def}) in Definition \ref{def:SymLinfty} (the definition of $L_\infty$-algebra) is equivalent to the following
\begin{equation*}
\sum_{i+j=n+1}\iota_{l_i}l_j=0,
\end{equation*}
for all $n\geq 1$.
\end{rem}
The following statement appears, in a more or less implicit form, in
\cite{Roytenberg}.

\begin{them}\label{th:linftyRN}
Let $E=\oplus_{i\in \mathbb{Z}}E_i$ be a graded vector space,
 $(l_i)_{i \geq 1}: \otimes^{i}E \to E$ be a family of symmetric vector valued forms on $E$ and $\mu= \sum_{i\geq 1} l_i$.
\begin{enumerate}
\item If the symmetric vector valued forms $(l_i)_{i \geq 1}$  define an $L_{\infty}$-structure on $E$, then $[\mu,\mu]_{_{RN}}=0$.
\item If  for each $i\geq 1$, the degree of  $l_i$ is $+1$ and $[\mu,\mu]_{_{RN}}=0$, then $(E,(l_i)_{i \geq 1})$ is an $L_{\infty}$-algebra.
\end{enumerate}
\end{them}
\begin{proof} First observe that if $\bar{l_i}=1$ for all $i\geq 1$, then
\begin{equation}\label{mumu}
[\mu,\mu]_{_{RN}}=\sum_{n\geq 1}(\sum_{i+j=n+1}[l_i,l_j]_{_{RN}})=2\sum_{n\geq 1}(\sum_{i+j=n+1}\iota_{l_i}l_j).
\end{equation}
1. Assume that the symmetric vector valued forms $(l_i)_{i \geq 1}$  define an $L_{\infty}$-structure on $E$. Then, $\bar{l_i}=1$ for all $i\geq 1$ and $\sum_{i+j=n+1}\iota_{l_i}l_j=0$, for all $n\geq 1$. Thus, Equation (\ref{mumu}) implies that $[\mu,\mu]_{_{RN}}=0$.

2. Assume that $\bar{l_i}=1$ for all $i\geq 1$ and $[\mu,\mu]_{_{RN}}=0$. Then Equation (\ref{mumu}) implies that $\sum_{i+j=n+1}\iota_{l_i}l_j=0$, for all $n\geq 1$, which means that $(E,(l_i)_{i \geq})$ is an $L_{\infty}$-algebra.

%
%
\end{proof}
The same statement remains true for curved $L_\infty$-algebras.

\begin{them}\label{th:curvedlinftyRN}
Let $E=\oplus_{i\in \mathbb{Z}}E_i$ be a graded vector space, $(l_i)_{i \geq 0}: \otimes^{i}E \to E$ be a family of symmetric vector valued forms on $E$ and $\mu= \sum_{i \geq 0} l_i$.
\begin{enumerate}
\item If the symmetric vector valued forms $(l_i)_{i \geq 0}$  define a curved $L_{\infty}$-structure on $E$, then $[\mu,\mu]_{_{RN}}=0$.
\item If  for each $i\geq 0$, the degree of  $l_i$ is $+1$ and $[\mu,\mu]_{_{RN}}=0$, then $(E,(l_i)_{i \geq 0})$ is a curved $L_{\infty}$-algebra.
\end{enumerate}
\end{them}
\section{Multiplicative $L_{\infty}$-structures}
In this section we introduce the concept of multiplicative $L_{\infty}$-structures and  classify all multiplicative $L_{\infty}$-structures on $\Gamma(\wedge A)[2]$,
for $A \to M $ an arbitrary vector bundle over a manifold $M$.

\subsection{Richardson-Nijenhuis bracket on multi-derivations}

There is an important sub-graded Lie algebra of $(\tilde{S}(E^*) \otimes E,\,\left[.,.\right]_{{RN}})$, when $E$ itself is equipped with a graded
commutative associative algebra structure on $E[2]$, denoted by $\wedge$, that is, a bilinear operation such that for all $X \in E_i, Y \in E_j, Z\in E_k$
\begin{itemize}
 \item $  X\wedge Y \in E_{i+j+2}$,
 \item $(X\wedge Y) \wedge Z = X\wedge (Y \wedge Z)$,
 \item $ X \wedge Y = (-1)^{|X||Y|} Y\wedge X$,
\end{itemize}
 where $|X|=i+2$ and $|Y|=j+2$.
%
\begin{defi}
Let $E$ be a graded vector space equipped with an associative graded commutative algebra structure, that is a graded symmetric bilinear map $\wedge$ of degree zero which is associative. An element $D\in S^{d}(E^*)\otimes E$ is called a multi-derivation vector valued $d$-form, if
\begin{equation}\label{eqdef:multi-derivation}
\begin{array}{rcl}
&&D(X_1, \cdots,X_{i-1}, Y\wedge Z, X_{i+1},\cdots,X_d)\\
&=&(-1)^{|Z|(|X_{i+1}|+\cdots+|X_{d}|)}D(X_1, \cdots,X_{i-1}, Y, X_{i+1},\cdots,X_d)\wedge Z\\
&&+(-1)^{|Y|(|X_{1}|+\cdots+|X_{i-1}|+\bar{D})}Y\wedge D(X_1, \cdots,X_{i-1}, Z, X_{i+1},\cdots,X_d),\\
\end{array}
\end{equation}
for all $X_1, \cdots,X_d,Y,Z \in E$, where $\bar{D}$ is degree of $D$ as a graded map.
\end{defi}
\begin{rem}
Graded commutativity of the product $\wedge$ implies that the Equation (\ref{eqdef:multi-derivation}) is equivalent to
\begin{equation}
\begin{array}{rcl}
&&D(X_1,\cdots,X_{d-1},Y\wedge Z)\\
&=&D(X_1, \cdots,X_{d-1},Y)\wedge Z +(-1)^{|Y||Z|} D(X_1,\cdots,X_{d-1},Z)\wedge Y.\\
\end{array}
\end{equation}
\end{rem}

 The space of all multi-derivations is denoted by $MultiDer(E)$. Elements of $S^{1}(E^*)\otimes E$ are simply called \emph{derivations}. By definition, $E \subset MultiDer(E)$ and we have the following:
\begin{prop}\label{prop:Multider}
    $MultiDer(E)$ is a sub-graded Lie algebra of  $(\tilde{S}(E^*) \otimes E,\, \left[.,.\right]_{_{RN}})$.
 \end{prop}
 We will use the following lemmas in the proof of Proposition \ref{prop:Multider}.
 \begin{lem}\label{lem:multiDer1}
 Let $D_1$ and $D_2$ be two derivations. Then $[D_1,D_2]_{_{RN}}$ is also a derivation.
 \end{lem}
 \begin{proof}
 We use the fact that the space of derivations of an associative graded commutative algebra, equipped with the graded commutator bracket, is a graded Lie algebra. We have
 \begin{equation*}
 \begin{array}{rcl}
 [D_1,D_2]_{_{RN}}&=& D_2\circ D_1-(-1)^{\bar{D_1}\bar{D_2}}D_1\circ D_2\\
                  &=&-(-1)^{\bar{D_1}\bar{D_2}}[D_1,D_2],
 \end{array}
 \end{equation*}
 where $\left[.,.\right]$ is the graded commutator on the space of derivations of the graded associative commutative algebra $(E,\, \wedge)$. This proves that $[D_1,D_2]_{_{RN}}$ is a derivation.

 \end{proof}

 \begin{lem}\label{lem:multiDer2}
 If $D\in S^{d}(E^*)\otimes E $ is a multi-derivation vector valued form, then for all $X\in E$, $[X,D]_{_{RN}}$ is a multi-derivation vector valued $(d-1)$-form.
 \end{lem}
 \begin{proof}
 It is a direct consequence of
 \begin{equation*}
 [X,D]_{_{RN}}(X_1,\cdots,X_{d-2}, Y\wedge Z)=D(X,X_1,\cdots,X_{d-2}, Y\wedge Z),
 \end{equation*}
 for all elements $Y,Z,X_1,\cdots,X_{d-2} \in E$ which comes from the definition of Richardson-Nijenhuis bracket.
 \end{proof}
 \begin{proof}(of Proposition \ref{prop:Multider})
 Let $D, D^{\prime}$ be two multi-derivation vector valued $d$- and $d'$-forms respectively. We show that $[D,D^{\prime}]_{_{RN}}$ is a multi-derivation vector valued $(d+d'-1)$-form. We use the induction on the number $n=d+d'-1$. Lemmas \ref{lem:multiDer1} and \ref{lem:multiDer2} prove the case $n=1$. Assume by induction that $[D,D^{\prime}]_{_{RN}}$ is a multi-derivation vector valued $(d+d'-1)$-form and let $D_1$ and $D_2$ be two multi-derivation vector valued $d_1$- and $d_2$-forms respectively, such that $d_1+d_2-1=n+1$. Using Proposition \ref{DrivedBt} we have
 \begin{equation}\label{eq:firstinprop}
 [D_1,D_2]_{_{RN}}(X_1,\cdots,X_{d_1+d_2-2},Y\wedge Z)=[Y\wedge Z,[X_{d_1+d_2-2},\cdots,[X_1,[D_1,D_2]_{_{RN}}]_{_{RN}}\cdots]_{_{RN}}]_{_{RN}},
 \end{equation}
 for all $X_1,\cdots,X_{d_1+d_2-2},Y, Z\in E$. Using Jacobi identity for the vector valued forms $D_1,D_2$ and $X_1$, Equation (\ref{eq:firstinprop}) can be rewritten as
 \begin{equation}
 \begin{array}{rcl}
 &&[D_1,D_2]_{_{RN}}(X_1,\cdots,X_{d_1+d_2-2},Y\wedge Z)\\
 &=&[Y\wedge Z,[X_{d_1+d_2-2},\cdots,[[X_1,D_1]_{_{RN}},D_2]_{_{RN}}\cdots]_{_{RN}}]_{_{RN}}\\
 &&+(-1)^{\bar{D_1}\bar{X_1}}[Y\wedge Z,[X_{d_1+d_2-2},\cdots,[D_1,[X_1,D_2]_{_{RN}}]_{_{RN}}\cdots]_{_{RN}}]_{_{RN}}.\\
 \end{array}
 \end{equation}
 By Lemma \ref{lem:multiDer2}, $[X_1,D_1]_{_{RN}}$ and $[X_1,D_2]_{_{RN}}$ are multi-derivation vector valued $(d_1-1)$- and $(d_2-1)$-forms respectively, and hence using the assumption of induction, $[[X_1,D_1]_{_{RN}},D_2]_{_{RN}}$ and $[D_1,[X_1,D_2]_{_{RN}}]_{_{RN}}$ are multi-derivation vector valued $n$-forms. Therefore
  \begin{equation}
 \begin{array}{rcl}
 &&[D_1,D_2]_{_{RN}}(X_1,\cdots,X_{d_1+d_2-2},Y\wedge Z)\\
 &=&[[X_1,D_1]_{_{RN}},D_2]_{_{RN}}(X_2,\cdots,X_{d_1+d_2-2},Y\wedge Z)\\
 &&+(-1)^{\bar{D_1}\bar{X_1}}[D_1,[X_1,D_2]_{_{RN}}]_{_{RN}}(X_2,\cdots,X_{d_1+d_2-2},Y\wedge Z)\\
 &=&[[X_1,D_1]_{_{RN}},D_2]_{_{RN}}(X_2,\cdots,X_{d_1+d_2-2},Y)\wedge Z\\
 &&+(-1)^{|Y||Z|}[[X_1,D_1]_{_{RN}},D_2]_{_{RN}}(X_2,\cdots,X_{d_1+d_2-2},Z)\wedge Y\\
 &&+(-1)^{\bar{D_1}\bar{X_1}}[D_1,[X_1,D_2]_{_{RN}}]_{_{RN}}(X_2,\cdots,X_{d_1+d_2-2},Y)\wedge Z\\
 &&+(-1)^{\bar{D_1}\bar{X_1}}(-1)^{|Y||Z|}[D_1,[X_1,D_2]_{_{RN}}]_{_{RN}}(X_2,\cdots,X_{d_1+d_2-2},Z)\wedge Y\\
 &=&[D_1,D_2]_{_{RN}}(X_1,\cdots,X_{d_1+d_2-2},Y)\wedge Z\\
 &&+(-1)^{|Y||Z|}[D_1,D_2]_{_{RN}}(X_1,\cdots,X_{d_1+d_2-2},Z)\wedge Y.\\
 \end{array}
 \end{equation}
 This completes the induction and hence the proof of proposition.
 \end{proof}

\subsection{Examples around Lie algebroids}\label{subsection:Example around Lie algebroids}
We start by introducing the notion of multiplicative $L_{\infty}$-algebra.
\begin{defi}
An $L_\infty$-structure $\mu=\sum_{i=1}^\infty l_i$ on a graded vector space $E$ equipped with
a graded commutative product $\wedge : E_i \times E_j \to E_{i+j}$
is called multiplicative if all the multi-linear brackets $l_i$
are multi-derivations.
\end{defi}
The notion of multiplicative $L_\infty$-algebra will offer us an opportunity to review the notions of
 Lie algebroid, pre-Lie algebroid, Lie bialgebroid, and quasi-Lie bialgebroid,
then we interpret those in terms of $L_\infty$-structures.

A \emph{pre-Lie algebroid structure on a vector bundle $A\to M$ over a manifold $M$}, is a pair $(\rho, \left[.,.\right])$ with $\rho: A \to TM $
a vector bundle morphism over the identity of $M$, called \emph{anchor map}, and $[ .,.] $
a skew-symmetric bilinear endomorphism of $\Gamma(A) $
subject to the so-called Leibniz identity:
   $$ [X,fY] = f[X,Y] + (\rho(X) f) \, Y,$$
   for all $X,Y \in \Gamma(A)$ and all $f \in C^\infty (M) $.
   When, moreover, $ \left[.,.\right]$ is a Lie algebra bracket, the pair $(\left[.,.\right],\, \rho)$ is called a
   \emph{Lie algebroid structure on $A\to M$}.

Let $(\left[.,.\right],\, \rho)$ be a pre-Lie algebroid structure on a vector bundle $A \to M$. Set $E_i:=\Gamma(\wedge^{i+1} A)$ and $ E=\oplus_{i\geq -1}E_i$ , where $E_{-1}=\Gamma(\wedge^0A)=\mathcal{C}^{\infty}(M)$. The \emph{Schouten-Nijenhuis} bracket on $E$ is defined as follows:
 \begin{equation}\label{Schouten-Nijenhuis1}
 \left[ X,Y\right]_{_{SN}}= \sum_{i=1 }^{p} \sum_{j=1 }^{q} (-1)^{i+j} [X_i,Y_j]  \widehat{X_i} \wedge \widehat{Y_j}
 \end{equation}
and
\begin{equation}\label{Schouten-Nijenhuis2}
\left[X , f\right]_{_{SN}} = \sum_{i=1}^{p} (-1)^{i+1}\rho(X_i)f  \widehat{X_i},
\end{equation}
for all $X=X_1\wedge\dots\wedge X_p \in E_{p-1}, Y=Y_1\wedge\dots\wedge Y_q \in E_{q-1},p,q\geq 1$, and  $f \in {\mathcal C}^\infty (M)$, where $\widehat{X_i}$ stands for $X_1\wedge \cdots\wedge X_{i-1}\wedge X_{i+1}\wedge \cdots X_p$. From Equations \ref{Schouten-Nijenhuis1} and \ref{Schouten-Nijenhuis2}, we obtain the following well known properties of the Schouten-Nijenhuis bracket on $E$:
\begin{equation} \label{eq:defGerstenhaber}
\begin{array}{rcl}
\left[ X,f \right]_{_{SN}} &=& \rho(X) f  \\
\left[ X,Y \right]_{_{SN}} & = & \left[ X,Y \right]  \\
\left[ P,Q \right]_{_{SN}} &=& -(-1)^{pq} \left[ Q,P\right]_{_{SN}}\\
\left[ P,Q\wedge R \right]_{_{SN}}&=&\left[ P,Q \right]_{_{SN}}\wedge R +(-1)^{qr}\left[ P,R \right]_{_{SN}}\wedge Q,
\end{array}
 \end{equation}
 for all $P\in \Gamma(\wedge^{p+1}A), Q\in \Gamma(\wedge^{q+1}A),R \in \Gamma(\wedge^{r+1}A)$ and $f\in {\mathcal C}^\infty (M)$.
 It follows from these properties that the Schouten-Nijenhuis bracket is a graded skew symmetric bracket of degree zero on $E=\oplus_{i\geq -1}E_i$. From the above relations, one sees that the Schouten-Nijenhuis bracket entirely encodes the initial pre-Lie algebroid structure $(\rho,\left[.,.\right])$ on $A \to M$. Also, it is known that a pre-Lie algebroid structure $(\rho,\left[.,.\right]) $ is a Lie algebroid structure on the vector bundle $A \to M,$ if and only if $\left[.,.\right]_{_{SN}}$ is a graded Lie algebra bracket on $E=\Gamma(\wedge A)[1]$. It is also classic that a pre-Lie algebroid $(A,\,\left[.,.\right],\,\rho)$ is entirely determined by its "de Rham"-like or "Chevalley-Eilenberg"-like differential, that is, the derivation $\diff^A$ of $\Gamma(\wedge A^*) $ given by:
  \begin{equation}\label{def:d^A} \diff^{A} \omega (X_0 , \dots, X_k) := \sum_{i=0}^k (-1)^i \rho(X_i) \omega(\widehat{X_i})+\sum_{0\leq i<j\leq k}(-1)^{i+j}
 \omega([X_i,X_j], \widehat{X_{i,j}}),
 \end{equation}
   for all $X_0,\dots,X_k \in \Gamma(A) , \omega \in \Gamma(\wedge^k A^*)$, where $\widehat{X_{i}}$ and $\widehat{X_{i,j}}$ stand for
   \begin{equation*}
   X_1,\cdots, X_{i-1}, X_{i+1}, \cdots, X_{k}\quad \mbox{and} \quad X_1,\cdots, X_{i-1}, X_{i+1}, \cdots, X_{j-1}, X_{j+1}, \cdots, X_{k}
    \end{equation*}
    respectively. Notice that in the above expression, we have implicitly identified elements of $\Gamma(\wedge^k A^*) $ with skew-symmetric $k$-linear maps from $\Gamma(A)\times \cdots\times\Gamma(A)$ to ${\mathcal C}^{\infty}(M)$. The "Chevalley-Eilenberg"-like differential $\diff^{A}$ squares to zero, if and only if, $(A,\,\left[.,.\right],\,\rho)$ is Lie algebroid.

    The discussion above leads to the conclusion that there are two ways to see Lie algebroids as $L_\infty$-structures: the first one will make it an $L_{\infty}$-structure on $\Gamma(\wedge A) $, and the second one will make it an $L_\infty$-structure on $\Gamma(\wedge A^*) $. More precisely:

\begin{prop}\label{prop:algebroid}
   Let  $A \to M $ be a vector bundle and $A^* \to M $ its dual. There is a one to one correspondence between:
   \begin{enumerate}
   \item[(i)] pre-Lie algebroid structures $(\rho, \left[.,.\right]) $ on $A \to M$,
   \item[(ii)] binary multi-derivations of $\Gamma(\wedge A)[2] $ of degree $1$,
   \item[(iii)] unary multi-derivations of $\Gamma(\wedge A^*)[2] $ of degree $1$.
   \end{enumerate}
   The one to one correspondence above restricts to a one to one correspondence between:
   \begin{enumerate}
   \item[($i'$)] Lie algebroid structures $(\rho, \left[.,.\right]) $ on $A \to M$,
   \item[($ii'$)] multiplicative $L_{\infty}$-structures on $\Gamma(\wedge A)[2] $ given by a binary bracket,
   \item[($iii'$)] multiplicative $L_{\infty}$-structures on $\Gamma(\wedge A^*)[2] $ given by a unary bracket.
   \end{enumerate}
   \end{prop}
\begin{proof}
Let $(\left[.,.\right],\, \rho)$ be a pre-Lie algebroid structure on the vector bundle $A\to M$. Then, the Schouten-Nijenhuis bracket $ \left[.,.\right]_{_{SN}}$ is a graded skew-symmetric bracket of degree zero on  $\Gamma(\wedge A)[1]$. Therefore the bilinear map
\begin{equation*}
l_2^{^{\left[.,.\right]}}: \Gamma(\wedge A)[2]\times \Gamma(\wedge A)[2]\to \Gamma(\wedge A)[2]
\end{equation*}
defined as
\begin{equation}\label{eq:braToL2}
l_2^{^{\left[.,.\right]}}(X,Y)=(-1)^{|X|}[X,Y]_{_{SN}},
\end{equation}
with $X\in \Gamma(\wedge^{|X|+1}A)$,
is a graded symmetric vector valued $2$-form of degree $1$. Hence, the forth equation in  (\ref{eq:defGerstenhaber}) shows that $l_2^{^{\left[.,.\right]}}$ is a binary multi-derivation of $\Gamma(\wedge A)[2] $. This shows the correspondence from $(i)$ to $(ii)$ in the proposition.  Moreover, if $(\rho, \left[.,.\right]) $ is a Lie algebroid structure on the vector bundle $A\to M$, then $l_2^{^{\left[.,.\right]}}$ is, in addition, a graded symmetric Lie algebra bracket and hence a multiplicative $L_{\infty}$-structure on $\Gamma(\wedge A)[2] $. This shows the restricted correspondence from $(i')$ to $(ii')$ in the proposition.\\
The correspondence from $(i)$ to $(iii)$ is obtained by associating to a pre-Lie algebroid $(A,\,\left[.,.\right],\,\rho)$, its "Chevalley-Eilenberg"-like differential $\diff^{A}$. While the correspondence from $(iii)$  to $(i)$ holds because, as we already mentioned, from $\diff^A$ we may define the pre-Lie algebroid structure $(\rho, \, \left[.,.\right])$. The anchor is obtained by $\rho(X)f=\diff^Af(X)$ with $X\in E$ and $f\in \mathcal{C}^{\infty}(M)$ and the bracket is directly obtained from (\ref{def:d^A}). The restricted correspondence between $(i')$ and $(iii^\prime)$ in the proposition follows from the fact that the "Chevalley-Eilenberg"-like differential squares to zero, if and only if, $(A,\,\left[.,.\right],\,\rho)$ is Lie algebroid.\\
Let us now look at the converse correspondence. Given a binary multi-derivation $l_2$ on $\Gamma(\wedge A)[2]$, of degree $1$, Equation (\ref{eq:braToL2}) defines a skew-symmetric bracket $\left[.,.\right]_{_{SN}}$ of degree zero, on $\Gamma(\wedge A)[1]$. For sections $X,Y\in \Gamma(A)$, set $[X,Y]=[X,Y]_{_{SN}}$; then $\left[.,.\right]$ is a skew symmetric bilinear bracket on $\Gamma(A)$. Now, if $f,g \in \mathcal{C}^{\infty}(M)$ and $X\in \Gamma(A)$ we have, since $l_2$ is a multi-derivation,
\begin{equation*}
l_2(X,fg)=l_2(X,f)g+fl_2(X,g)
\end{equation*}
which means that $l_2(X,.)|_{\mathcal{C}^{\infty}(M)}$ is a derivation on $\mathcal{C}^{\infty}(M)$. So, there exists a map $\rho:\Gamma(A)\to \Gamma(TM)$ such that
\begin{equation*}
[X,f]=\rho(X)f,\,\,\,\mbox{for all}\,\,\, X\in \Gamma(A) \,\,\,\mbox{and}\,\,\,f\in \mathcal{C}^{\infty}(M).
\end{equation*}
Moreover, $\rho$ is $\mathcal{C}^{\infty}$-linear:
\begin{equation*}
\rho(hX)f=[hX,f]_{_{SN}}=h[X,f]_{_{SN}}=h\rho(X)f,
\end{equation*}
for all $f,h \in \mathcal{C}^{\infty}(M)$ and $X\in \Gamma(A)$. Hence there exists a vector bundle morphism, denoted by the same letter, $\rho:A\to TM$. Finally, because $l_2$ is a multi-derivation, the condition
\begin{equation*}
[X,fY]=f[X,Y]+\rho(X)f.Y
\end{equation*}
holds for all $X,Y\in \Gamma(A)$ and $f\in \mathcal{C}^{\infty}(M)$, so that pair $(\rho, \left[.,.\right])$ is a pre-Lie algebroid structure on $A\to M$. This shows the correspondence from $(ii)$ to $(i)$. For the correspondence from $(ii')$ to $(i')$, it is enough to notice that if $[l_2,l_2]_{_{RN}}=0$ holds, then the bracket $\left[.,.\right]$ on $\Gamma(A)$ satisfies the Jacobi identity.
\end{proof}
We recall the following:

 \begin{defi}\label{def:bialgebroid}
 A pair $(A,A^*)$ of vector bundles over manifold $M$, in duality, is said to be a Lie bialgebroid if $A$ and $A^*$ are both equipped with Lie algebroid structures $(\rho_A,\left[.,.\right]^A)$ and $(\rho_{A^*},\left[.,.\right]^{A^*})$ respectively, satisfying

 \begin{equation}\label{eq:bialgebra}
 \diff^{A^*} [P,Q]^A_{_{SN}}= [\diff^{A^*} P,Q]^A_{_{SN}} + (-1)^{p-1}[P,\diff^{A^*} Q]^A_{_{SN}}
 \end{equation}
  for all $P \in  \Gamma(\wedge^p A),Q \in  \Gamma(\wedge^q A)$, where $\left[.,.\right]^A_{_{SN}}$ stands for the Schouten-Nijenhuis bracket on $\Gamma(\wedge A)$ corresponding to the bracket $\left[.,.\right]^A$, and $\diff^{A^*}$ stands for the de Rham differential of $A^*$.
  \end{defi}


  In view of (\ref{eq:braToL2}), (\ref{eq:bialgebra}) amounts to:
    $$ \diff^{A^*} l_2^{[,.,]} (P,Q)= - l_2^{[,.,]}(\diff^{A^*} P,Q) + (-1)^{p-1}l_2(P,\diff^{A^*} Q),$$
  which is equivalent to
    $$ [l_2^{[,.,]}, \diff^{A^*}]_{_{RN}} =0.$$

  Since, in view of Proposition \ref{prop:algebroid}, there is a one to one correspondence between
  Lie algebroid structures on $A $ and binary multi-derivations $l_2 $
  with $[l_2,l_2]_{_{RN}}=0 $ and a one to one correspondence between
  Lie algebroid structures on $A^* $ and unary multi-derivations $l_1 $
  with $[l_1,l_1]_{_{RN}}=0 $, we arrive at the following conclusion (which is already well-known, but stated in terms in Differential Graded Lie Algebras):

  \begin{prop}
   Let  $A \to M $ be a vector bundle, there is a one to one correspondence between:
   \begin{enumerate}
   \item[(i)] Lie bialgebroid structures on $A$,
   \item[(ii)] multiplicative $L_{\infty}$-structures on $\Gamma(\wedge A)[2] $
   which are the sum of an unary and a binary multi-derivations.
   \end{enumerate}
  \end{prop}

 In fact, this result has been extended by Roytenberg \cite{Roytenberg} to the case of quasi-Lie bialgebroids. First, let us recall what a quasi-Lie bialgebroid is.

\begin{defi}\label{def:quasi-Liebialgebroid}
Let $(A, A^*)$ be a pair of vector bundles over $M$ in duality, such that $A$ is equipped with a pre-Lie algebroid structure $(\rho_A,\left[.,.\right]^A)$ and $(A^*,\rho_{A^*},\left[.,.\right]^{A^*})$ is a Lie algebroid. Let $\diff^A$ and $\diff^{A^*}$ be the associated de Rham pre-differential and de Rham differential respectively, and  $\omega \in \Gamma(\wedge^3 A^*)$ be a $ \diff^{A}$-closed $3$-form.
 We say that the triple $(A,A^*,\omega)$ forms a quasi-Lie bialgebroid when
  \begin{equation}\label{eq:quasi1}
  \diff^{A^*} [P,Q]_{_{SN}}^A= [\diff^{A^*} P,Q]_{_{SN}}^A + (-1)^{p-1}[P,\diff^{A^*} Q]_{_{SN}}^A,
  \end{equation}
\begin{equation}\label{eq:quasi2}
\begin{array}{c}
(-1)^{p+q}[P,[Q,R]_{_{SN}}^A]_{_{SN}}^A +c.p.+\diff^{A^*} \underline{\omega} (P,Q,R) + \underline{\omega} (\diff^{A^*} P,Q,R)\\
 +(-1)^{p} \underline{\omega} (P,\diff^{A^*} Q,R) +(-1)^{p+q}\underline{\omega} (P,Q,\diff^{A^*} R)=0
\end{array}
\end{equation}
 for all $P \in  \Gamma(\wedge^p A),Q \in  \Gamma(\wedge^q A)$ and $R\in \Gamma(\wedge^r A)$, where
   $$\underline{\omega} (P,Q,R)=\sum_{i,j,k=1}^{p,q,r}(-1)^{i+j+k+|Q|+1}\omega(P_i,Q_j,R_k)\widehat{P_i}\wedge \widehat{Q_j}\wedge \widehat{R_k},$$
   with $\widehat{P_i}=P_1\wedge \cdots \wedge P_{i-1}\wedge P_{i+1}\wedge \cdots \wedge P_p$, etc, $p\geq 1,q\geq 1$ and $r\geq 0$.
\end{defi}

Notice that a quasi-Lie bialgebroid $(A, A^*,\omega)$ with $\omega=0$ is a Lie bialgebroid.
\begin{rem} The previous definition is not exactly the one given by Roytenberg in \cite{RoytenbergWeinstein} and \cite{Roytenberg}, but they are equivalent.
Roytenberg's definition includes
\begin{equation}\label{eq:quasiRoy1}
[[X,Y]^A,Z]^A+c.p.=\diff^{A^*}\omega (X,Y,Z) +\omega (\diff^{A^*} X,Y,Z) - \omega(X,\diff^{A^*} Y,Z) +\omega(X,Y,\diff^{A^*} Z)
\end{equation}
and
\begin{equation}\label{eq:quasiRoy2}
\rho_A[X,Y]^A=[\rho_{A}(X),\rho_{A}(Y)]+\rho_{A^*}\omega(X,Y),
\end{equation}

with $X,Y,Z\in \Gamma(A)$, while in Definition \ref{def:quasi-Liebialgebroid} these two are encoded in Equation (\ref{eq:quasi2}). Starting with Equation (\ref{eq:quasi2}) and taking $P=X,Q=Y,R=Z$, with $X,Y,Z\in \Gamma(A)$, we get Equation (\ref{eq:quasiRoy1}). Then taking  $P=X,Q=Y,R=f$, with $X,Y\in \Gamma(A)$ and $f\in \mathcal{C}^{\infty}(M)$, gives
\begin{equation*}
[X,[Y,f]_{_{SN}}^A]_{_{SN}}^A+[Y,[f,X]_{_{SN}}^A]_{_{SN}}^A+[f,[X,Y]_{_{SN}}^A]_{_{SN}}^A=-\underline{\omega}(X,Y,\diff^{A^*}f),
\end{equation*}
which implies
\begin{equation*}
\rho_{A}(X)\rho(Y)f-\rho_{A}(Y)\rho_{A}(X)f-\rho_{A}([X,Y])f=-\iota_{\diff^{A^*}f}\omega(X,Y).
\end{equation*}
Hence
\begin{equation*}
\rho_{A}([X,Y]^A)=[\rho_{A}(X),\rho_{A}(Y)]+\rho_{A^*}(\omega(X,Y)),
\end{equation*}
where $\omega(X,Y)=\iota_{X\wedge Y}\omega$. Using the fact that $\left[.,.\right]^A$ is a multi-derivation with respect to $\wedge$, it is easy to see that from Equation (\ref{eq:quasiRoy1}) we get Equation (\ref{eq:quasi2}), with $p,q,r\geq 1$. It is also clear that from Equation (\ref{eq:quasiRoy2}) we get Equation (\ref{eq:quasi2}), with $p=q=1$ and $r=0$.
\end{rem}
\begin{prop}
   Let  $A \to M $ be a vector bundle, there is a one to one correspondence between:
   \begin{enumerate}
   \item[(i)] quasi-Lie bialgebroid structures on $A$,
   \item[(ii)] multiplicative $L_\infty$-structures on $\Gamma(\wedge A)[2]$ which are the sum of unary, binary and 3-ary multi-derivations.
   \end{enumerate}
  \end{prop}
  \begin{proof} Let $(A,A^*, \omega)$ be a quasi-Lie bialgebroid and set $E_i:=\Gamma(\wedge^{i+2}A)$ and $E=\oplus_{i\geq -2} E_i$.
  Let $l_1:= \diff^{A^*}, l_2(P,Q):=(-1)^{p}[P,Q]_{_{SN}}^A$, with $P\in E_{p}$ and $l_3:=\frac{1}{2}\underline{\omega}$. First notice that $(A^*,\rho_{A^*},\left[.,.\right]^{A^*})$ is a Lie algebroid, if and only if, $\diff^{A^*}$ squares to zero which is equivalent to $[l_1,l_1]_{_{RN}}=0$. Equation (\ref{eq:quasi1}) is equivalent to $[l_1,l_2]_{_{RN}}=0$, Equation (\ref{eq:quasi2}) is equivalent to $[l_2,l_2]_{_{RN}}+2[l_3,l_1]_{_{RN}}=0$. We prove that $\omega$ is $\diff^A$-closed, if and only if,  $[l_2,l_3]_{_{RN}}=0$. It follows from definitions that
  \begin{equation*}
  \diff^A\omega =[l_2,l_3]_{_{RN}}|_{E_{-1}\times E_{-1}\times E_{-1}\times E_{-1}}.
  \end{equation*}
   This proves that if $[l_2,l_3]_{_{RN}}=0$, $\omega$ is $\diff^A$-closed. Conversely assume that $\omega$ is $\diff^A$-closed. Assume by induction (on $n=p+q+r+s$) that $[l_2,l_3]_{_{RN}}(P,Q,R,S)=0$, where $P\in E_p, Q\in E_q,R\in E_r$ and $S\in E_s$. Since both $l_2$ and $l_3$ are multi-derivations,  $[l_2,l_3]_{_{RN}}$ is a multi-derivation. Therefore
   $[l_2,l_3]_{_{RN}}(P\wedge X,Q,R,S)=[l_2,l_3]_{_{RN}}(P,Q\wedge X,R,S)=[l_2,l_3]_{_{RN}}(P,Q,R\wedge X,S)=[l_2,l_3]_{_{RN}}(P,Q,R,S\wedge X)=0$ for all $X\in E_{-1}.$ This proves that $[l_2,l_3]_{_{RN}}(P,Q,R,S)=0$ when, $P\in E_p, Q\in E_q,R\in E_r$ and $S\in E_s$ and $p+q+r+s=n+1$. This completes the induction. Notice that by definition of $\underline{\omega}$ and $l_3 $ we have $[l_3,l_3]_{_{RN}}=0$, since $\underline{\omega}$ takes value in $\mathcal{C}^{\infty}(M)$ and $\underline{\omega}$ vanishes on $E_{-2}\times E_{-1}\times E_{-1}$.
   We conclude that the triple $(A,A^*,\omega)$ forms a quasi-Lie bialgebroid, if and only if, $[l_1+l_2+l_3,l_1+l_2+l_3]_{_{RN}}=0$.
  \end{proof}


\section{Nijenhuis vector valued forms}
Having an $L_\infty$-structure on a graded vector space, that is, by Theorem \ref{th:linftyRN}, a symmetric vector valued form $\mu=\sum_{i=0}^{\infty} l_i$ with $\bar{l_i}=1$ and $[\mu,\mu]_{_{RN}}=0$, one may deform the structure $\mu$ by a given symmetric vector valued form $\mathcal{N}$ yielding the vector valued form $\mu^{\mathcal N}=[{\mathcal N},\mu]_{_{RN}}$. In this section we define a Nijenhuis vector valued form with respect to a given vector valued form $\mu$ and deformation of $\mu$ by a Nijenhuis vector valued form. Then we show that deforming an $L_{\infty}$-structure by a Nijenhuis vector valued form, one gets an $L_{\infty}$-structure.
\subsection{Definition and properties}
 We start by fixing some notations.
\begin{defi}\label{def:Nijenhuis}
Let $E$ be a graded vector space and $\mu$ be a symmetric vector valued form on $E$ of degree $1$. A vector valued form ${\mathcal N}$ of degree zero is called
\begin{enumerate}
  \item weak Nijenhuis vector valued form with respect to $\mu$ (or simply weak Nijenhuis with respect to $\mu$, when there is no risk of confusion), if
       \begin{equation}\label{weakN}
       \left[\mu,\left[{\mathcal N},\left[{\mathcal N},\mu\right]_{_{RN}}\right]_{_{RN}}\right]_{_{RN}}=0,
       \end{equation}
  \item Nijenhuis vector valued form with respect to $\mu$ (or simply Nijenhuis with respect to $\mu$, when there is no risk of confusion), if there exists a vector valued form ${\mathcal K}$ of degree $0$, such that
       \begin{equation}\label{strongN}
       \left[{\mathcal N},\left[{\mathcal N},\mu\right]_{_{RN}}\right]_{_{RN}}=\left[{\mathcal K},\mu\right]_{_{RN}}\,\,\,\hbox{and }\,\,\,\left[{\mathcal N},{\mathcal K}\right]_{_{RN}}=0.
       \end{equation}
  Such a ${\mathcal K} $ is called a square of ${\mathcal N} $. If ${\mathcal N}$ contains an element of the underlying graded vector space, that is, ${\mathcal N}$ has a component which is a vector valued zero form, then ${\mathcal N}$ is called Nijenhuis vector valued form with curvature.
\end{enumerate}
\end{defi}


The following is an immediate consequence of Definition \ref{def:Nijenhuis}.
\begin{prop}
Let $\mu$ be a  (curved) $L_\infty$-structure on a graded vector space $E$ and ${\mathcal N}$ be a vector valued form on $E$. If ${\mathcal N}$ is a Nijenhuis vector valued form with respect to $\mu$, then ${\mathcal N}$ is a weak Nijenhuis vector valued form with respect to $\mu$.
\end{prop}
\begin{proof}
Let ${\mathcal N}$ be a Nijenhuis vector valued form with respect to a (curved) $L_\infty$-structure $\mu$ with square $\mathcal{K}$. By the Lemma \ref{lem:[mu[k,mu]]=0} we get
\begin{equation}
[\mu,[{\mathcal K},\mu]_{_{RN}}]_{_{RN}}=0.
\end{equation}
This implies that
 \begin{equation*}
       [\mu,[{\mathcal N},[{\mathcal N},\mu]_{_{RN}}]_{_{RN}}]_{_{RN}}=0,
       \end{equation*}
       which completes the proof.

\end{proof}

\begin{rem}\label{rem:whynotsquares}
It would be of course tempting to choose ${\mathcal K} = \iota_{\mathcal N}{\mathcal N} $, having in mind what happens for manifolds, and the fact
that  $\iota_{ {\mathcal N}}{{\mathcal N}} = {\mathcal N}^2 $ for vector valued $1$-forms.
However, it is not what examples show to be a reasonable definition.

Also, it is easy to see that if ${\mathcal N}$ is a vector valued $2$-form we do not have
  $$ [\iota_{\mathcal N} {\mathcal N} , {\mathcal N}]_{_{RN}} =0,$$
 which says $\iota_{\mathcal N} {\mathcal N}  $ is not a good candidate for the square, except maybe for vector valued $1$-forms.
\end{rem}

The properties of the Richardson-Nijenhuis bracket allow one to
write the following definition.

\begin{defi}
Given a symmetric vector valued form $\mu$ of degree $1$ and a symmetric vector valued form ${\mathcal N} $
 of degree $0$ on a graded vector space, we call $[{\mathcal N},\mu]_{_{RN}}$ the deformed bracket
 of $\mu$ by ${\mathcal N} $ and  denote it by $\mu^{\mathcal N}$.
 More generally,  $\mu^{{\mathcal N},\dots , {\mathcal N} }$
 (with $k$ copies of ${\mathcal N}$) stands for
  $$ [{\mathcal N},\dots, [{\mathcal N}, \mu]_{_{RN}}]_{_{RN}}$$
 (with $k$ copies of ${\mathcal N}$).
\end{defi}

\begin{prop}\label{prop:deformedisLinfty}
Let $(E,\mu)$ be a (curved)  $L_\infty$-algebra and ${\mathcal N}$ be a symmetric vector valued form on $E$. Then ${\mathcal N}$ is weak Nijenhuis with respect to $\mu$ if and only if $\mu^{\mathcal N} $ is a (curved) $L_{\infty}$-algebra.
\end{prop}
\begin{proof}

Using the Jacobi identity for $\mu, {\mathcal N}$ and $[{\mathcal N}, \mu]_{_{RN}}$ we get
\begin{equation}\label{k1}
[\mu,[{\mathcal N},[{\mathcal N},\mu]_{_{RN}}]_{_{RN}}]_{_{RN}}=[[\mu,{\mathcal N}]_{_{RN}},[\mu,{\mathcal N}]_{_{RN}}]_{_{RN}}+[{\mathcal N},[\mu,[{\mathcal N},\mu]_{_{RN}}]_{_{RN}}]_{_{RN}}
\end{equation}
and hence Lemma \ref{lem:[mu[k,mu]]=0} implies that
\begin{equation}\label{k3}
[\mu,[{\mathcal N},[{\mathcal N},\mu]_{_{RN}}]_{_{RN}}]_{_{RN}}=[[\mu,{\mathcal N}]_{_{RN}},[\mu,{\mathcal N}]_{_{RN}}]_{_{RN}}.
\end{equation}
Since $\mu$ is of degree $1$, the degree of $[{\mathcal N},\mu]_{_{RN}}$ is $1$, if and only if, the degree of $\mathcal{N}$ is zero. Therefore, (\ref{k3}) means that ${\mathcal N}$ is weak Nijenhuis with respect to $\mu$ if and only if $[{\mathcal N},\mu]_{_{RN}}$ is a (curved) $L_{\infty}$-algebra.
\end{proof}

\begin{cor}
Let $(E,\mu)$ be a (curved) $L_\infty$-algebra and ${\mathcal N}$ be a vector valued form on $E$. If ${\mathcal N}$ is Nijenhuis with respect to $\mu$, then $\mu^{\mathcal N} $ is a (curved) $L_{\infty}$-algebra.
\end{cor}

Nijenhuis tensors, in the classical case, allow to construct hierarchies of compatible structures. With some restrictions, the same
phenomena shall appear here.

\begin{defi}
Two $L_\infty$-structures $\mu_1$ and $\mu_2$ on a graded vector space
are said to be compatible, if they commute with respect to the Richardson Nijenhuis bracket.
\end{defi}

\begin{examp}
Given a weak Nijenhuis vector valued form   ${\mathcal N}$ with respect to an $L_{\infty}$-structure $\mu$,
$\mu$ and $\mu^{{\mathcal N},{\mathcal N}}$ are compatible.
\end{examp}

\begin{prop}
Given two $L_\infty$-structures $\mu_1$ and $\mu_2$ on a graded vector space $E$,
the following conditions are equivalent:
\begin{enumerate}
\item $\mu_1$ and $\mu_2$ are compatible,
\item $a\mu_1 + b\mu_2 $ is an $L_\infty$-structure on $E$, for all $a,b \in {\mathbb R}$.
\end{enumerate}
 \end{prop}
\begin{proof}This is a direct consequence of the following
\begin{equation*}
[a\mu_1+b\mu_2,a\mu_1+b\mu_2]_{_{RN}}=a^2[\mu_1,\mu_1]_{_{RN}}+2ab[\mu_1,\mu_2]_{_{RN}}+b^2[\mu_2,\mu_2]_{_{RN}}.
\end{equation*}
\end{proof}

Weak Nijenhuis vector valued forms do not, in general, give hierarchies in any sense.
However, Nijenhuis vector valued forms do.

\begin{them}\label{theo:Hierarchy}
Let ${\mathcal N}$ be a Nijenhuis vector valued form with respect to a  (curved) $L_\infty$-structure ${\mu}$ with square ${\mathcal K}$, then :
\begin{enumerate}
\item For all integers $k \geq 1$, $\mu^{{\mathcal N}, \dots, {\mathcal N}}$ (with $k$ copies of ${\mathcal N}$) is a (curved) $L_\infty$-structure denoted by $\mu_k$ and ${\mathcal N}$ is Nijenhuis of square ${\mathcal K}$ with respect to $\mu_k$.
\item For all integers $k,l \geq 1$, $\mu_k$ and $\mu_l$ are compatible.
\end{enumerate}
\end{them}
\begin{proof}
1) Assume, by induction, that ${\mathcal N}$ is Nijenhuis with respect to $\mu_k$ with square ${\mathcal K}$. Then we have
\begin{equation*}
[{\mathcal N},[{\mathcal N},\mu_k]_{_{RN}}]_{_{RN}}=[{\mathcal K},\mu_k]_{_{RN}},
\end{equation*}
that implies
\begin{equation}\label{J1}
[{\mathcal N},[{\mathcal N},[{\mathcal N},\mu_k]_{_{RN}}]_{_{RN}}]_{_{RN}}=[{\mathcal N},[{\mathcal K},\mu_k]_{_{RN}}]_{_{RN}}.
\end{equation}
Applying the Jacobi identity on the right hand side of (\ref{J1}) and using the fact that  ${\mathcal N}$ and  ${\mathcal K}$ commute with respect to the Richardson-Nijenhuis bracket, we get
\begin{equation*}
[{\mathcal N},[{\mathcal N},\mu_{k+1}]_{_{RN}}]_{_{RN}}=[{\mathcal K},\mu_{k+1}]_{_{RN}}.
\end{equation*}
This means that ${\mathcal N}$ is Nijenhuis with respect to $\mu_{k+1}$, with square ${\mathcal K}$.\\
%
2) Let $k$ and $l$ be two positive integers with $k\geq l$. If $k=l$, then 1) together with Theorem \ref{th:linftyRN} implies that $\mu_k$ and $\mu_l$ commute with respect to the Richardson-Nijenhuis bracket. For the case $k>l$, assume by induction that $\mu_k$ commutes with $\mu_n$ for all integers $k\geq n\geq l$. Then by the Jacobi identity
\begin{equation}\label{ind}
[\mu_{k+1},\mu_{l}]_{_{RN}}=[[{\mathcal N},\mu_{k}]_{_{RN}},\mu_{l}]_{_{RN}}=
[{\mathcal N},[\mu_{k},\mu_{l}]_{_{RN}}]_{_{RN}}-[\mu_{k},[{\mathcal N},\mu_{l}]_{_{RN}}]_{_{RN}},
\end{equation}
while both terms in the right hand side of (\ref{ind}) vanish by assumption of induction. This completes the induction and shows that $\mu_k$ and $\mu_l$ commute for all $k$ and $l$.
\end{proof}


Assume that ${\mathcal N}$ is Nijenhuis with respect to $\mu$ with square $ {\mathcal K}$. Then, for every non-negative integer $k$, we have  $[{\mathcal K},\mu_k]=\mu_{k+2} $. This implies that
\begin{equation}
[\mu_k,[\mathcal K,[\mathcal K,\mu_k]_{_{RN}}]_{_{RN}}]_{_{RN}}=[\mu_k, \mu_{k+4}]_{_{RN}}=0,
\end{equation}
which means that $ {\mathcal K}$ is weak Nijenhuis with respect to $\mu_k$. However $ {\mathcal K}$ may not be a Nijenhuis with respect to $\mu_k$.
This contrasts with  Nijenhuis tensors on manifolds and on Lie algebras,
where it is true that $N$ being Nijenhuis implies that $N^2$ is Nijenhuis.

\subsection{Generalities on Nijenhuis forms on $L_{\infty}$-algebras}

Let $\mu=\sum_{n\geq 1} l_n$ be an $L_\infty$-structure on a graded vector space $E=\oplus_{i\in \mathbb{Z}}E_i$.
This implies that $l_1$ is a map of degree $+1$ squaring to zero
defining therefore a cohomology denoted by $H^*(\mu)$.

Since $[l_2,l_1]_{_{RN}}=0$, it is classical that $l_2$ goes to the quotient
to yield a symmetric bilinear map $ \widetilde{l_2}:H^*(\mu) \otimes H^*(\mu) \to  H^*(\mu)$.
The relation $[l_2,l_2]_{_{RN}}+2[l_3,l_1]_{_{RN}}=0$ implies that
 $l_2$, which is \emph{not} a Lie bracket on $E$,  induces a Lie bracket on $H^*(\mu)$.

Now, there is a well-known notion of Nijenhuis tensor on a Lie algebra, see Subsection \ref{sec:deformationofLinfty}. It is natural to see whether a Nijenhuis vector valued form on an $L_\infty$-algebra will give a Nijenhuis
tensor on the cohomology of $ H^*(\mu)$. The answer is negative in general due to the fact that, in full generality,
a Nijenhuis form does not induce anything at the cohomological level.
Indeed, the only thing that can be proved is the following:

 \begin{lem}
 Let ${\mathcal N}= \sum_{k \geq 1} N_k$ be a Nijenhuis vector valued form with respect to an $L_{\infty}$-structure
  $\mu=\sum_{n\geq 1} l_n$, on a graded vector space $E=\oplus_{i \in \mathbb{Z}} E_i$, with square ${\mathcal K}= \sum_{k \geq 1} K_k$. Then
   $$ 2 N_1 l_1  N_1 = l_1 (N_1^2 -K_1) + (N_1^2 +K_1)l_1.$$
 \end{lem}
\begin{proof}
The relation above is simply obtained by considering the vector valued $1$-form component in $\mu^{{\mathcal N},\mathcal N} = \mu^{\mathcal K}.$
\end{proof}
This relation only proves that $N_1 l_1 N_1 $ maps $l_1$-closed elements in $ E_k $ into
exact elements in $E_{k+1}$. However, in many examples we shall consider in the sequel, we have
 $l_1 N_1 + N_1 l_1 = \lambda l_1 $
 and
 $l_1 K_1 + K_1 l_1 = \lambda' l_1$,
for some constants $\lambda,\lambda'$ (i.e. $N_1$ and $K_1$ are  quasi chain maps) or the same kind of  relations with different signs.
This implies that $N_1$ defines an endomorphism $H^*(N_1)$ of degree $0$ of $ H^*(\mu)$,
and the following holds true:

\begin{prop}\label{prop:chomologyNijenhuis}
Let ${\mathcal N}= \sum_{k \geq 1} N_k$ be a Nijenhuis vector valued form with respect to an $L_{\infty}$-structure $\mu=\sum_{n\geq 1} l_n$ on a graded vector space $E$, with $N_1$ and $K_1$ being quasi-chain maps. Then the induced maps
  $H^*(N_1), H^*(K_1) : H^*(\mu) \to H^*(\mu) $ satisfy:
   $$[H^*(N_1),[H^*(N_1),\widetilde{l_2}]_{_{RN}}]_{_{RN}}=[H^*(K_1),\tilde{l_2}]_{_{RN}},$$
   where $\widetilde{l_2}$ is the induced Lie bracket on the cohomology level $H^*(\mu)$.
\end{prop}
%
\begin{lem}\label{lem:bikhod}
If $N_1$ is a quasi chain map, then for all $X,Y\in Ker(l_1)$
\begin{enumerate}
\item $[N_1,l_1]_{_{RN}}(X)=0$,
\item $[N_1,l_1]_{_{RN}}l_2(X,Y)=0$.
\end{enumerate}
\end{lem}
\begin{proof}
The first item follows only from the definitions of Richardson-Nijenhuis bracket and quasi chain map, while for the second item we also need to use the fact that if $X,Y\in Ker(l_1),$ then $[l_1,l_2]_{_{RN}}=0$ implies that $l_1l_2(X,Y)=0.$
 \end{proof}
 \begin{proof}(of Proposition \ref{prop:chomologyNijenhuis})\\
 Let $X,Y\in E$, with $l_1(X)=l_1(Y)=0$. Then, using the graded Jacobi identity of the Richardson-Nijenhuis bracket for $l_1,l_2,N_1$ and using the fact that $[l_1,l_2]_{_{RN}}=0,$ we get
  \begin{equation*}
  [l_1,[N_1,l_2]_{_{RN}}]_{_{RN}}(X,Y)=[[l_1,N_1]_{_{RN}},l_2]_{_{RN}}(X,Y)
  \end{equation*}
  or
  \begin{equation*}
  \begin{array}{rcl}
  &&l_2^{N_1}(l_1X,Y)+l_2^{N_1}(X,l_1Y)-l_1(l_2^{N_1}(X,Y))\\
  &=&l_2([N_1,l_1]_{_{RN}}X,Y)
  +l_2(X,[N_1,l_1]_{_{RN}}Y)-[N_1,l_1]_{_{RN}}(l_2(X,Y)).
  \end{array}
  \end{equation*}
  This and  Lemma \ref{lem:bikhod} imply that
  \begin{equation*}
  l_1(l_2^{N_1}(X,Y))=0,
  \end{equation*}
   which means that $l_2^{N_1}$ induces a map on the cohomology level $H^*(\mu)$. As a consequence, the following diagram commutes :

%
%

\begin{equation}\label{diagram:cohomology}
  \xymatrix@R+2em@C+10em{
  l_2 \ar[r]^{deforming}_{N_1,N_1} \ar[d]_{quotient} &l_2^{N_1,N_1}  \ar[d]^{quotient} \\
   \widetilde{l_2}\ar[r]^ {deforming}_{H^*(N_1),H^*(N_1)}& \widetilde{l_2}^{H^*(N_1),H^*(N_1)}=\widetilde{l_2^{N_1,N_1}}
  }
 \end{equation}
and

\begin{equation}\label{1}
\widetilde{l_2}^{H^*(K_1)}=\widetilde{l_2^{K_1}}.
\end{equation}
Let $\mu_2^{{\mathcal N},\mathcal N}$ denote the sum of components which are vector valued $2$-forms, appearing in $\mu^{{\mathcal N},\mathcal N}$. Then
\begin{equation}\label{2}
\mu_2^{{\mathcal N},\mathcal N}=[N_1,[N_1,l_2]_{_{RN}}]_{_{RN}}+[N_1,[N_2,l_1]_{_{RN}}]_{_{RN}}+[N_2,[N_1,l_1]_{_{RN}}]_{_{RN}}.
\end{equation}
Applying the Jacobi identity for $[N_1,[N_2,l_1]_{_{RN}}]_{_{RN}}$ in (\ref{2}) gives
\begin{equation}\label{3}
\mu_2^{{\mathcal N},\mathcal N}=l_2^{N_1,N_1}+[[N_1,N_2]_{_{RN}},l_1]_{_{RN}}+2[N_2,[N_1,l_1]_{_{RN}}]_{_{RN}}.
\end{equation}
If we apply (\ref{3}) to $(l_1)$-closed elements $X,Y$, using Lemma \ref{lem:bikhod},  we get:
\begin{equation}\label{4}
\begin{array}{rcl}
\mu_2^{{\mathcal N},\mathcal N}(X,Y)&=&l_2^{N_1,N_1}(X,Y)\\
                                    &-&l_1[N_1,N_2]_{_{RN}}(X,Y)+[N_1,N_2]_{_{RN}}(l_1(X),Y)+[N_1,N_2]_{_{RN}}(X,l_1(Y))\\
                                    &+&2[N_1,l_1]_{_{RN}}N_2(X,Y)+N_2([N_1,l_1]_{_{RN}}X,Y)+N_2(X,[N_1,l_1]_{_{RN}}Y)\\
                                    &=&l_2^{N_1,N_1}(X,Y)-l_1[N_1,N_2]_{_{RN}}(X,Y)+2l_1N_1N_2(X,Y)-2N_1l_1N_2(X,Y)\\
                                    &=&l_2^{N_1,N_1}(X,Y)\\
                                    &-&l_1[N_1,N_2]_{_{RN}}(X,Y)+4l_1N_1N_2(X,Y)-2\lambda l_1N_2(X,Y).\\
\end{array}
\end{equation}
On the other hand, the sum of components which are vector valued $2$-forms appearing in $\mu^{\mathcal K}$, will be denoted by $\mu_2^{\mathcal K}$ and is of the type
\begin{equation}\label{5}
l_2^{K_1}+l_1^{K_2}.
\end{equation}
If we apply (\ref{5}) to $(l_1)$-closed elements $X,Y$  we get:
\begin{equation}\label{6}
\mu_2^{\mathcal K}(X,Y)=l_2^{K_1}(X,Y)+l_1K_2(X,Y).
\end{equation}
Now since $\mathcal{N}$ is Nijenhuis with respect to $\mu$ with square $\mathcal{K}$, we have $\mu^{\mathcal{N},\mathcal{N}}=\mu^{\mathcal{K}}$. This together with (\ref{4}) and (\ref{6}) imply that
\begin{equation*}
\widetilde{l_2^{N_1,N_1}}=\widetilde{l_2}^{H^*(K_1)}.
\end{equation*}
Since diagram (\ref{diagram:cohomology}) commutes, we get
\begin{equation}\label{}
\widetilde{l_2}^{H^*(N_1),H^*(N_1)}=\widetilde{l_2}^{H^*(K_1)}.
\end{equation}
\end{proof}

In general, the sum of two Nijenhuis vector valued forms with respect to a vector valued form $\mu$ on a graded vector space $E$, is not a Nijenhus vector valued form with respect to $\mu$.
But if two Nijenhuis vector valued forms, with respect to a vector valued form $\mu$ on a graded space, are compatible in the sense of  next definition, then their sum is going to be Nijenhuis.
\begin{defi}
Let $\mu$ be a vector valued form on a graded vector space $E$ and ${\mathcal N}_1,{\mathcal N}_2$ be Nijenhuis vector valued forms, with squares ${\mathcal K}_1,{\mathcal K}_2$, respectively, with respect to $\mu$. We say that the Nijenhuis forms ${\mathcal N}_1,{\mathcal N}_2$ are compatible if:
$$\left[\mathcal{N}_1,\left[\mathcal{N}_2,\mu\right]_{_{RN}}\right]_{_{RN}}+\left[\mathcal{N}_2,\left[\mathcal{N}_1,\mu\right]_{_{RN}}\right]_{_{RN}}=0$$
and
$$\left[\mathcal{N}_1, \mathcal{K}_2\right]_{_{RN}}+\left[\mathcal{N}_2, \mathcal{K}_1\right]_{_{RN}}=0.$$
More generally, given an $n$-tuple of Nijenhuis vector valued forms ${\mathcal N}_i$, with respect to a vector valued form $\mu$, with squares ${\mathcal K}_i$, we say that the ${\mathcal N}_i$´s are compatible if they are pairwise compatible.
\end{defi}
\begin{prop}
Let ${\mu}$ be a vector valued form on a graded vector space $E$. For every $n$-tuple of compatible Nijenhuis vector valued forms ${\mathcal N}_i$ with respect to $\mu$, of squares ${\mathcal K}_i$ and every constants $a_1,\cdots,a_n$,
$\sum_{i=1}^n a_i {\mathcal N}_i $ is a Nijenhuis vector valued form with respect to $\mu$, of square $\sum_{i=1}^n a_i^2 {\mathcal K}_i$.
 \end{prop}
\begin{proof}
The statement is a direct consequence of the following computations:
\begin{equation*}
\begin{array}{rcl}
&&\left[\sum_{i=1}^{n}a_i{\mathcal N}_i,\left[\sum_{i=1}^{n}a_i{\mathcal N}_i,\mu\right]_{_{RN}}\right]_{_{RN}}\\
&=& \sum_{i=1}^{n}a_{i}^2\left[{\mathcal N}_i,\left[{\mathcal N}_i,\mu\right]_{_{RN}}\right]_{_{RN}}+\frac{1}{2}\sum_{ \substack{i,j=1 \\ i\not=j} }^{n}a_{i}a_{j}\left(\left[{\mathcal N}_i,\left[{\mathcal N}_j,\mu\right]_{_{RN}}\right]_{_{RN}}+\left[{\mathcal N}_j,\left[{\mathcal N}_i,\mu\right]_{_{RN}}\right]_{_{RN}}\right)\\
&=&\sum_{i=1}^{n}a_{i}^2\left[{\mathcal N}_i,\left[{\mathcal N}_i,\mu\right]_{_{RN}}\right]_{_{RN}}
\end{array}
\end{equation*}
and
\begin{equation*}
\left[ \sum_{i=1}^na_i\mathcal{N}_i,\sum_{i=1}^na_i\mathcal{K}_i \right]_{_{RN}}=\sum_{i=1}^na^2_i\left[ \mathcal{N}_i,\mathcal{K}_i \right]_{_{RN}}+\frac{1}{2}\sum_{\substack{i,j=1\\i\not=j}}^na_ia_j\left(\left[\mathcal{N}_i,\mathcal{K}_j\right]_{_{RN}}+\left[\mathcal{N}_j,\mathcal{K}_i\right]_{_{RN}} \right).
\end{equation*}
\end{proof}

\section{Nijenhuis forms on GLA and DGLA}

\subsection{The Euler map}

The Euler map, the map that simply counts the degree, is always a Nijenhuis vector valued $1$-form, with respect to any vector valued form of degree $1$ (or with respect to any vector valued form which is the sum of vector valued forms of degrees $1$), with square itself.
\begin{lem}\label{Euler}
Let $E=\oplus_{i\in \mathbb{Z}}E_i$ be a graded vector space. Let $S:E\rightarrow E$ be the Euler map, that is $S(X)=-|X|X$, for all homogeneous elements $X\in E$ of degree $|X|$. Then,
for every symmetric vector valued form $\alpha$, we have
\begin{equation*}
[S,\alpha]_{_{RN}}=\bar{\alpha}\,\alpha.
\end{equation*}
\end{lem}
\begin{proof}First notice that the degree of $S$ as a graded map is zero. Assume that $\alpha$ is a symmetric vector valued $k$-form, then
\begin{equation*}
\begin{array}{rcl}
                               & &[S,\alpha]_{_{RN}}(X_1,\cdots,X_k)\\
                               &=&(\iota_S\alpha -\iota_{\alpha}S)(X_1,\cdots,X_k)\\
                               &=&\alpha(SX_1,X_2,\cdots,X_k)+\cdots+\alpha(X_1,,\cdots,X_{k-1},SX_{k})-S\alpha(X_1,\cdots,X_k)\\
                               &=&-(|X_1|+\cdots+|X_k|)\,\alpha(X_1,\cdots,X_k)+(|X_1|+\cdots+|X_k|+ \bar{\alpha})\,\alpha(X_1,\cdots,X_k)\\
                               &=&\bar{\alpha}\,\alpha(X_1,\cdots,X_k).
\end{array}
\end{equation*}
\end{proof}
The next proposition is the first example of Nijenhuis vector valued form and will be used to construct  more examples.
\begin{prop}
Let $E=\oplus_{i\in \mathbb{Z}}E_i$ be a graded vector space and $S:E\rightarrow E$ be the Euler map introduced in the Lemma \ref{Euler}. Let $\mu$ be a vector valued form of degree $1$ on the graded vector space $E=\oplus_{i\in \mathbb{Z}}E_i$, that is all the components of $\mu$  are of degree $1$. Then $S$ is Nijenhuis with respect to $\mu$ with square itself.
\end{prop}
\begin{proof}
Let $\mu=\sum_{i=1}^{\infty}l_i$. Applying Lemma \ref{Euler} to each $l_i, 1\leq i\leq \infty$, and taking the sum we get:
\begin{equation*}
[S,\mu]_{_{RN}}=\sum_{i=1}^{\infty}[S,l_i]_{_{RN}}=\sum_{i=1}^{\infty}l_i=\mu.
\end{equation*}
Therefore
\begin{equation*}
[S,[S,\mu]_{_{RN}}]_{_{RN}}=[S,\mu]_{_{RN}}.
\end{equation*}
Since $\bar{S}=0$, Lemma \ref{Euler} implies that $[S,S]_{_{RN}}=0$ and this completes the proof.
\end{proof}

Of course, the result can be enlarged as follows for every $\mu$-cocycle, that is, a vector valued form $\eta$ such that $[\mu,\eta]_{_{RN}}=0$.
\begin{examp} Let $\mu=\sum_il_i$ be a vector valued form of degree $+1$ on a graded vector space. Then for every element $\alpha$ of degree $0$ in $ S(E^*) \otimes E $ with $[\mu,\alpha]_{_{RN}} =0$, $S+\alpha $ is Nijenhuis with respect to $\mu$, with square $S$.\end{examp}
\subsection{Lie algebra, GLA and DGLA}\label{subsection:Lie algebra,GLAandDGLA}

Recall that a \emph{symmetric graded  Lie algebra (symmetric GLA)} is a $\mathbb{Z}$-graded vector space $E=\oplus E_i$ endowed with a binary graded symmetric bracket $\mu=\left[.,.\right]$ of degree $1$, satisfying the graded Jacobi identity i.e.
\begin{equation}\label{JacobiLie}
[X,[Y,Z]]=(-1)^{|X|+1}[[X,Y],Z]+(-1)^{(|X|+1)(|Y|+1)}[Y,[X,Z]],
\end{equation}
for all $X,Y,Z \in E$. Such a symmetric GLA will be denoted by $(E=\oplus_{i\in \mathbb{Z}}E_i,\, \left[.,.\right] )$. We will sometimes use the notation $l_2$ for the bracket $\left[.,.\right]$ and we will  say that $\left[.,.\right]$ (or $l_2$) is a GLA structure on the graded vector space $E$. Note that when the graded vector space is concentrated on degree $-1$, that is, all the vector spaces $E_i$ are zero, except $E_{-1}$, then (\ref{JacobiLie}) is the usual Jacobi identity and we get a Lie algebra with symmetric bracket. We would like to remark that (\ref{JacobiLie}) can be written as
\begin{equation}
\mu(\mu(X,Y),Z)+(-1)^{|Y||Z|}\mu(\mu(X,Z),Y)+(-1)^{|X|(|Y|+|Z|)}\mu(\mu(Y,Z),X)=0,
\end{equation}
for all $X,Y,Z \in E$, or
\begin{equation}
\frac{1}{2}(\iota_{\mu}\mu+\iota_{\mu}\mu)(X,Y,Z)=\frac{1}{2}[\mu,\mu]_{_{RN}}(X,Y,Z)=0.
\end{equation}
This means that a \emph{symmetric graded Lie algebra is simply an $L_{\infty}$-algebra such that all the multi-brackets are zero except the binary one.} And a \emph{Lie algebra is an $L_{\infty}$-algebra concentrated in degree $-1$, that is an $L_{\infty}$-algebra on a graded vector space concentrated on degree $-1$ for which all the brackets are zero except the binary bracket.} As we mentioned in the Introduction, applying the d\'{e}calage isomorphism to a symmetric GLA, we get the notion of a (GLA) graded Lie algebra.

 Recall from \cite{YKS} that a \emph{Nijenhuis tensor on a graded Lie algebra} $(E,\mu=\left[.,.\right])$ is a $(1,1)$-tensor $N:E\rightarrow E$ such that the Nijenhuis torsion of $\mu$ is identically zero, that is
\begin{equation}\label{Torsion}
T_{\mu}N(X,Y)=\mu(NX,NY)-N(\mu(NX,Y)+\mu(X,NY)-N(\mu(X,Y)))=0,
\end{equation}
for all $X,Y\in E$.
For a binary bracket $\mu=\left[.,.\right]$ the deformed bracket by $N$ is denoted by $\left[.,.\right]_N$ and is given by $[X,Y]_N=[NX,Y]+[X,NY]-N[X,Y]$. It has been shown in \cite{YKS} that if $N$ is Nijenhuis on a Lie algebra $(E,\left[.,.\right])$, then $(E,\left[.,.\right]_N)$ is also a Lie algebra and $N$ is a morphism of Lie algebras.
Also it has been shown that $N$ is Nijenhuis if and only if deforming the original bracket of the Lie algebra twice by $N$ is equivalent with deforming it once by $N^2$, that is $([X,Y]_N)_N=[X,Y]_{N^2}$. This can be stated using the notion of Richardson-Nijenhuis bracket on the space of vector valued forms on a graded vector space $E$, as follows:
\begin{equation}
[N,[N,\mu]_{_{RN}}]_{_{RN}}=[N^2,\mu]_{_{RN}}.
\end{equation}

Nijenhuis structures in this usual and traditional sense are of course Nijenhuis structures in our sense also.

\begin{prop}
Let $\left[.,.\right]$ be a graded skew-symmetric Lie algebra structure on a graded vector space $\mathfrak g={\oplus \mathfrak g}_i $.
Assume that $\left[.,.\right]^{\prime}$ is the correspondent symmetric bracket by the use of d\'{e}calage isomorphism on the graded vector space $E=\oplus E_i$, with $E_i:={\mathfrak g}_{i+1} $. Then $N$ is a Nijenhuis $(1-1)$-tensor on the graded skew-symmetric Lie algebra $(\mathfrak g,\left[.,.\right])$ if and only if it is a Nijenhuis vector valued $1$-form with respect to the symmetric bracket $\mu=\left[.,.\right]^{\prime}$ seen as symmetric vector valued $2$-form, with square $N^2$.
\end{prop}

A \emph{symmetric differential graded Lie algebra (symmetric DGLA)} is an $L_\infty$-algebra $(E=\oplus_{i\in \mathbb{Z}} E_i,\mu=\sum_{i=1}^{\infty}l_i)$, with all the brackets, except $l_1$ and $l_2$, being zero. In other words, a symmetric DGLA is a symmetric GLA $(\oplus_{i\in \mathbb{Z}} E_i,\left[.,.\right])$ endowed with a differential $\diff$, that is, a linear map $\diff :\oplus E_i\to\oplus E_i $ of degree $1$ squaring to zero, satisfying the compatibility condition
\begin{equation}\label{def:DGLA}
\diff[X,Y]+[\diff(X),Y]+(-1)^{|X|}[X,\diff(Y)]=0,
\end{equation}
for all $X,Y \in E.$ Such a symmetric DGLA will be denoted by $(E=\oplus_{i\in \mathbb{Z}} E_i,\, \diff, \left[.,.\right])$. Sometimes we will use the notation $l_1$ for the differential $\diff$ and the notation $l_2$ for the bracket $\left[.,.\right]$ and we will say that $\mu=l_1+l_2$ is a DGLA structure on the graded vector space $E=\oplus_{i\in \mathbb{Z}} E_i$ or simply on the graded vector space $E.$
\begin{defi}\label{def:curvedDGLA}
A curved symmetric DGLA is a $\mathbb{Z}$-graded vector space $E=\oplus_{i\in \mathbb{Z}} E_i$ together with an element $C \in E_{1}$, a linear map $\diff: E \to E$ of degree $1$ and a bilinear map $\left[.,.\right]:E \times E \to E$ of degree $1$ such that:
\begin{enumerate}
\item $\left[.,.\right]$ is a graded symmetric Lie bracket,
\item $\diff(C)=0$,
\item $[C,X]+\diff^2X=0$, for all $X \in E$,
\item $\diff [X,Y]+[\diff X,Y]+(-1)^{|X|}[X,\diff Y]=0$, for all $X,Y \in E$.
\end{enumerate}
 Such a curved symmetric DGLA will be denoted by $(E=\oplus_{i\in \mathbb{Z}} E_i,\, C, \diff, \left[.,.\right])$. Sometimes, we will use the notation $l_1$ for the linear map $\diff$ and the notation $l_2$ for the bracket $\left[.,.\right]$ and we will say that $\mu=C+l_1+l_2$ is a curved DGLA structure on the graded vector space $E=\oplus_{i\in \mathbb{Z}} E_i$ or simply on the graded vector space $E.$
\end{defi}
A \emph{Maurer Cartan element} in a symmetric DGLA $(E,\diff,\left[.,.\right])$ is an element $e \in E_0$ such that
\begin{equation} \label{MaurerCartanElementsymm}
\diff(e)-\frac{1}{2}[e,e]=0.
\end{equation}
A \emph{Maurer Cartan element} in a symmetric  curved DGLA $(E,\, C,\diff, \left[.,.\right])$ is an element $e \in E_0$ such that
\begin{equation} \label{MaurerCartanElementsymmCurved}
(\diff(e)-l_0)-\frac{1}{2}[e,e]=0.
\end{equation}
\begin{rem}
In literature the Maurer Cartan element for a skew-symmetric $L_{\infty}$-algebra $(E=\oplus_{i\in \mathbb{Z}}E_i,\,\mu=\sum_{i\geq 1}l_i)$ is defined to be an element $e\in E_{1}$ such that
\begin{equation}
\sum_{i\geq 1}\frac{1}{i!}l_i(e,\cdots,e)=0.
\end{equation}
Of course after using the d\'{e}calage isomorphism we get the following formula for the symmetric $L_{\infty}$-algebra $((E=\oplus_{i\in \mathbb{Z}}E_i)[1],\,\mu=\sum_{i\geq 1}l'_i)$:
\begin{equation}
\sum_{i\geq 1}(-1)^{\frac{i-1}{2}}\frac{1}{i!}l'_i(e,\cdots,e)=0,
\end{equation}
where $e$ is an element in $(E[1])_{0}=E_{1}$.
\end{rem}
A \emph{Poisson element} in a (curved) $L_{\infty}$-algebra $(E,\mu=\sum_{i\geq 0}l_i)$ is an element $\pi \in E_0$, such that $l_2(\pi,\pi)=0$.

The next propositions provide examples of Nijenhuis vector valued forms on symmetric graded and symmetric differential graded Lie algebras.
\begin{prop}\label{deformationofGLA}
Let $\mu=l_0+l_2$ be a curved symmetric graded Lie algebra structure on a graded vector space $E=\oplus_{i\in \mathbb{Z}}E_i$ and $\pi\in E_0$. Then $\pi+S$ is a Nijenhuis vector valued form with curvature, with respect to $\mu$ and with square $2\pi+S$ if, and only if, $\pi$ is a Poisson element.
\end{prop}
\begin{proof}
The proof is a direct consequence of the following computations:
\begin{equation*}
[\pi+S,l_0+l_2]_{_{RN}}=l_2(\pi,.)+l_0+l_2,
\end{equation*}

\begin{equation*}
\begin{array}{rcl}
[\pi+S,[\pi+S,l_0+l_2]_{_{RN}}]_{_{RN}}&=&l_2(\pi,\pi)+l_0+2l_2(\pi,.)+l_2\\
                                       &=&l_2(\pi,\pi)+[2\pi+S,l_0+l_2]_{_{RN}}
\end{array}
\end{equation*}
and
\begin{equation*}
[\pi+S,2\pi+S]_{_{RN}}=2[\pi,\pi]_{_{RN}}+[\pi,S]_{_{RN}}+2[S,\pi]_{_{RN}}+[S,S]_{_{RN}}=0.
\end{equation*}
Notice that
\begin{equation*}
[\pi,S]_{_{RN}}=[S,\pi]_{_{RN}}=0,
\end{equation*}
since $\pi \in E_0$.
\end{proof}
\begin{cor}\label{cor1poissononGLA}
Let $\mu=l_2$ be a symmetric graded Lie algebra structure on a graded vector space $E=\oplus_{i\in \mathbb{Z}}E_i$ and $\pi\in E_0$. Then
\begin{enumerate}
\item $\pi+S$ is a Nijenhuis vector valued form with curvature, with respect to $\mu$ and with square $2\pi+S$ if, and only if, $\pi$ is a Poisson element.
\item Assuming 1. is satisfied, the deformed structure $[\pi+S,l_2]_{_{RN}}$ is a DGLA structure on $E$, with the same binary bracket $l_2=\left[.,.\right]$ and $l_1=\diff=[\pi,.]$.
\end{enumerate}
\end{cor}
\begin{proof}
1. It is a direct consequence of Proposition \ref{deformationofGLA} by setting $l_0=0$.\\
%
2. Let $X,Y \in E$. Applying the graded Jacobi identity in the symmetric graded Lie algebra $(E,l_2=\left[.,.\right])$, for the vectors $\pi,X,Y$ we get
\begin{equation*}
[[\pi,X],Y]+(-1)^{|X||Y|}[[\pi,Y],X]+[[X,Y],\pi]=0,
\end{equation*}
or
\begin{equation*}
\diff [X,Y]+[\diff X,Y]+(-1)^{|X|}[X,\diff Y]=0.
\end{equation*}
This proves that $l_2=\left[.,.\right]$ and $l_1=\diff=[\pi,.]$ is a DGLA structure on $E$.
\end{proof}

\begin{prop}\label{deformationofDGLA}
Let $\mu=C+l_1+l_2$ be a curved symmetric DGLA structure on a graded vector space $E=\oplus_{i\in \mathbb{Z}}E_i$ and $\pi\in E_0$. Then $\pi+S$ is a Nijenhuis vector valued form  (with curvature $\pi \in E_0$) with respect to $\mu$ and with square $2\pi+S$ if, and only if, $\pi$ is a Poisson element.
\end{prop}
\begin{proof}
The proof is a direct consequence of the following computations:
\begin{equation*}
[\pi+S,C+l_1+l_2]_{_{RN}}=C+l_1(\pi)+(l_2(\pi,.)+l_1)+l_2,
\end{equation*}

\begin{equation*}
\begin{array}{rcl}
[\pi+S,[\pi+S,C+l_1+l_2]_{_{RN}}]_{_{RN}}&=&[\pi+S,C+l_1+l_2+l_1(\pi)+l_2(\pi,.)]_{_{RN}}\\
                                 &=&C+l_1+l_2+2l_1(\pi)+2l_2(\pi,.)+l_2(\pi,\pi)\\
                                 &=&[2\pi+S,C+l_1+l_2]_{_{RN}}+l_2(\pi,\pi)
\end{array}
\end{equation*}
and
\begin{equation*}
[\pi+S,2\pi+S]_{_{RN}}=2[\pi,\pi]_{_{RN}}+[\pi,S]_{_{RN}}+2[S,\pi]_{_{RN}}+[S,S]_{_{RN}}=0.
\end{equation*}
\end{proof}
\begin{cor}\label{cor2poissononDGLA}
Let $l_1+l_2$ be a symmetric DGLA structure on a graded vector space $E=\oplus_{i\in \mathbb{Z}}E_i$ and $\pi\in E_0$. Then
\begin{enumerate}
\item $\pi+S$ is a Nijenhuis vector valued form  (with curvature $\pi \in E_0$) with respect to $\mu$ and with square $2\pi+S$ if, and only if, $\pi$ is a Poisson element.
\item Assuming $1.$ is satisfied, the deformed structure  $l_1(\pi)+(l_2(\pi,.)+l_1)+l_2$ is a curved DGLA.
\end{enumerate}
\end{cor}
\begin{proof}
1. It follows from Proposition \ref{deformationofDGLA} by setting $C=0$.\\
2. First notice that $l_1(\pi)\in E_1$. We need to prove: (see Definition \ref{def:curvedDGLA})
\begin{enumerate}
\item[(a)] $l_2$ is a graded symmetric Lie bracket,
\item[(b)] $(l_1+l_2(\pi,.))l_1(\pi)=0$,
\item[(c)] $l_2(l_1(\pi),X)+(l_1+l_2(\pi,.))^2(X)=0$,
\item[(d)] $l_1l_2(X,Y)+l_2(\pi,l_2(X,Y))+l_2(l_1(X)+l_2(\pi,X),Y)\\+(-1)^{|X|}l_2(X,l_1(Y)+l_2(\pi,Y))=0$,
\end{enumerate}
for all $X,Y\in E$.\\
(a) follows from definition of symmetric DGLA.\\
Let us prove (b). Since $\mu=l_1+l_2$ is a symmetric DGLA structure on $E$, we have
\begin{equation}\label{inja1}
l_1^2(\pi)=0
\end{equation}
and
\begin{equation}\label{inja2}
l_1(l_2(\pi,\pi))+2l_2(l_1(\pi),\pi)=0,
\end{equation}
where the last Equation follows from Equation (\ref{def:DGLA}) by setting $X=Y=\pi$ and $\left[.,.\right]=l_2, \diff=l_1$. Using the assumption that $\pi$ is a Poisson element, from Equation (\ref{inja2}) we get
\begin{equation} \label{inja3}
l_2(l_1(\pi),\pi)=0.
\end{equation}
Equations (\ref{inja1}) and (\ref{inja3}) imply that
\begin{equation}\label{inja4}
(l_1+l_2(\pi,.))(l_1(\pi))=0,
\end{equation}
which proves (b).\\
To prove (c) observe that since $(E,\mu=l_1+l_2)$ is a symmetric DGLA, we have
\begin{equation}\label{fench}
\begin{array}{rcl}
l_1^2(X)&=&0,\\
l_2(l_2(X,Y),Z)+l_2(l_2(Y,Z),X)+l_2(l_2(Z,X),Y)&=&0,\\
l_2(l_1(X),Y)+(-1)^{|X|}l_2(X,l_1(Y))+l_1l_2(X,Y)&=&0,
\end{array}
\end{equation}
for all $X,Y\in E$. On the other hand, using the fact that $l_2(\pi,\pi)=0$, we get by the Jacobi identity
\begin{equation}\label{fencheh}
l_2(\pi,l_2(\pi,X))=0.
\end{equation}
Setting $Y=\pi$, in the third equation in (\ref{fench}) and using Equation (\ref{fencheh}) we get
\begin{equation*}
l_2(l_1(\pi),X)+l_1^2(X)+l_1(l_2(\pi,X))+l_2(\pi,l_1(X))+l_2(\pi,l_2(\pi,X))=0,
\end{equation*}
which is equivalent to (c).\\
For (d) set $Z=\pi$ in the second equation in (\ref{fench}) and add with the last equation in (\ref{fench}). Then,
\begin{equation*}
l_1l_2(X,Y)+l_2(\pi,l_2(X,Y))+l_2(l_1(X)+l_2(\pi,X),Y)+(-1)^{|X|}l_2(X,l_1(Y)+l_2(\pi,Y))=0,
\end{equation*}
which is equivalent to (d).
%
\end{proof}


%
%
\begin{prop}\label{prop:MaurerCartanDGLA}
Let $\mu=C+l_1+l_2$ be a curved symmetric DGLA structure on a graded vector space $E=\oplus_{i\in \mathbb{Z}}E_i$ and $\pi\in E_0$. Then $Id_E+\pi$ is a Nijenhuis vector valued form with curvature $\pi$, with respect to $\mu$ and with square itself if, and only if, $\pi$ is a Maurer-Cartan element of the curved DGLA $(E,\mu)$.
\end{prop}
\begin{proof} First notice that
\begin{equation*}
\begin{array}{rcl}

                                   & &[\pi+Id_{E},[\pi+Id_E,C+l_1+l_2]_{_{RN}}]_{_{RN}}\\
                                   &=&[\pi+Id_E,(l_1(\pi)-C)+l_2(\pi,.)+l_2]_{_{RN}}\\
                                   &=&l_2(\pi,\pi)+l_2(\pi,.)-l_1(\pi)+C+l_2\\
                                   &=&-C-2((l_1(\pi)-C)-\frac{1}{2}l_2(\pi,\pi))+l_1(\pi)+l_2(\pi,.)+l_2\\
                                   &=&-2((l_1(\pi)-C)-\frac{1}{2}l_2(\pi,\pi))+[\pi+Id_E,C+l_1+l_2]_{_{RN}}.
\end{array}
\end{equation*}
This together with the fact that $[\pi+Id_E,\pi+Id_E]_{_{RN}}=0$ imply that $Id_E+\pi$ is Nijenhuis vector valued form with respect to $\mu$ if, and only if,  $\pi$ is a Maurer-Cartan element of the curved DGLA $(E,\mu)$.
\end{proof}
\begin{cor}\label{cor3poissononDGLA}
Let $\mu=l_1+l_2$ a DGLA structure on a graded vector space $E=\oplus_{i\in \mathbb{Z}}E_i$ and $\pi\in E_0$. Then
\begin{enumerate}
\item $Id_E+\pi$ is Nijenhuis vector valued form with curvature $\pi$, with respect to $\mu$ and with square itself if, and only if, $\pi$ is a Maurer-Cartan element of the  DGLA $(E,\mu)$,
\item Assuming $1.$ is satisfied, the deformed structure is the curved symmetric DGLA $l_1(\pi)+ l_2(\pi,.) +l_2$.
\end{enumerate}
\end{cor}
\begin{proof}
$1$ is a direct consequence of Proposition \ref{prop:MaurerCartanDGLA}. For item $2$ we need to prove
\begin{enumerate}
\item[(a)] $l_2(\pi,l_1(\pi))=0,$
\item[(b)] $l_2(l_1(\pi),X)+l_2(\pi,l_2(\pi,X))=0$,
\item[(c)] $l_2(\pi,l_2(X,Y))+l_2(l_2(\pi,X),Y)+(-1)^{|X|}l_2(X,l_2(\pi,Y))=0$.
\end{enumerate}
Since $\mu=l_1+l_2$ is a DGLA structure on the graded vector space $E$, we get
\begin{equation*}
l_1(l_2(\pi,\pi))+l_2(l_1(\pi),\pi)+l_2(\pi,l_1(\pi))=0.
\end{equation*}
This together with the fact that $\pi$ is a Maurer-Cartan element, hence of degree zero, gives (a). Since $\mu=l_1+l_2$ is a DGLA structure on the graded vector space $E$, we get
\begin{equation*}
l_2(l_2(\pi,\pi),Z)+l_2(l_2(\pi,Z),\pi)+l_2(l_2(Z,\pi),\pi)=0,
\end{equation*}
for all $Z\in E$ or
\begin{equation*}
l_2(l_2(\pi,\pi),Z)+2l_2(l_2(\pi,Z),\pi)=2(\frac{1}{2}l_2(l_2(\pi,\pi),Z)+l_2(\pi,l_2(\pi,Z)))=0,
\end{equation*}
for all $Z\in E$. This together with the fact that $\frac{1}{2}l_2(\pi,\pi)=l_1(\pi)$, gives (b). While (c) follows from the fact that
\begin{equation*}
[l_2,l_2]_{_{RN}}(X,Y,\pi)=0,
\end{equation*}
for all $X,Y\in E.$
\end{proof}
\section{Lie $n$-algebras and $n$-plectic manifolds}
Symmetric Lie $n$-algebras are particular cases of symmetric $L_{\infty}$-algebras. We give the definition and shortly discuss about Nijenhuis form on Lie $n$-algebras. Then we will give a certain class of Nijenhuis forms on Lie $n$-algebras associated to the so called $n$-plectic manifolds, a notion that we  recall.
\subsection{Nijenhuis forms on Lie $n$-algebras}\label{subsection:Lienalgebras}
A $\mathbb{Z}$-graded vector space $E=\oplus_{i\in \mathbb{Z}}E_i$ is said to be \emph{concentrated in degrees $p_1,\cdots p_k$}, with $p_1,\cdots, p_k \in \mathbb{Z}$, if $E_{p_1},\cdots,E_ {p_k}$ are the only non-zero components of $E$. Let us start with the definition of symmetric Lie $n$-algebra.
\begin{defi}
A symmetric Lie $n$-algebra is a symmetric $L_{\infty}$-algebra whose underlying graded vector space is concentrated on degrees $-n,\cdots,-1$.
\end{defi}
\begin{rem}\label{degreereason}
Note that by degree reasons, the only non-zero symmetric vector valued forms (multi-brackets) are $l_1,\cdots,l_{n+1}$.
\end{rem}
\begin{prop}\label{khodemoun}
 Let $(E=E_{-n}\oplus \cdots\oplus E_{-1},\, \mu=l_1+\cdots +l_{n+1})$ be a Lie $n$-algebra. Let $\frac{n+3}{2}\leq k\leq n+1 $ and $N$ be any symmetric vector valued $k$-form of degree zero on $E$. Then $S+N$ is a Nijenhuis vector valued form with respect to $\mu=l_1+\cdots +l_{n+1}$, with square $S+2N$ and the deformed Lie $n$-algebra structure on $E$ is of the form $$l_1+\cdots+l_{k-1}+(l_{k}+[N,l_1]_{_{RN}})+\cdots+(l_{n+1}+[N,l_{n-k+2}]_{_{RN}}).$$
\end{prop}
\begin{proof} By Remark \ref{degreereason}, any vector valued $(m+k-1)$-form, with  $m\geq n-k+3$, is identically zero and any vector valued $(2k+m-2)$-form, with $m\geq 1$ is identically zero, because from the conditions $\frac{n+3}{2}\leq k\leq n+1 $ and $m\geq 1$ we get $2k+m-2\geq n+2.$ Hence,
\begin{enumerate}\label{mesal}
\item $[N,l_m]_{_{RN}}=0$ for all $m\geq n-k+3$ and
\item $[N,[N,l_m]_{_{RN}}]_{_{RN}}=0$ for all  $m\geq 1$.
\end{enumerate}
These imply that
\begin{equation}\label{mesall}
[S+N,\mu]_{_{RN}}=\mu+[N,l_1]_{_{RN}}+\cdots+[N,l_{n-k+2}]_{_{RN}}
\end{equation}
and
\begin{equation}\label{mesalll}
[S+N,[S+N,\mu]_{_{RN}}]_{_{RN}}=\mu+2[N,l_1]_{_{RN}}+\cdots+2[N,l_{n-k+2}]_{_{RN}}.
\end{equation}
On the other hand using Lemma \ref{Euler} we have
\begin{equation}\label{mesallll}
[S+N,S+2N]_{_{RN}}=[S,S]_{_{RN}}+[S,N]_{_{RN}}+[N,N]_{_{RN}}=0.
\end{equation}
Equations (\ref{mesalll}) and (\ref{mesallll}) show that $S+N$ is a Nijenhuis vector valued form with respect to $\mu=l_1+\cdots+l_{n+1}$, with square $S+2N$. and Equation (\ref{mesall}) shows that the deformed Lie $n$-algebra by $S+N$ is $$l_1+\cdots+l_{k-1}+(l_{k}+[N,l_1]_{_{RN}})+\cdots+(l_{n+1}+[N,l_{n-k+2}]_{_{RN}}).$$
\end{proof}

\begin{prop}\label{corlast}
Let $(E=E_{-n}\oplus \cdots\oplus E_{-1}, \mu=l_1+\cdots +l_{n+1})$ be a Lie $n$-algebra. Let $N_1,\cdots,N_l$ be a family of symmetric vector valued $k_1,\cdots,k_l$-forms, respectively, of degree zero on $E$, with $\frac{n+3}{2}\leq k_1\leq \cdots\leq k_l \leq n+1$. Then $S+\sum_{i=1}^l N_i$ is a Nijenhuis vector valued form with respect to $\mu$, with square $S+2\sum_{i=1}^l N_i$. The deformed Lie $n$-algebra structure is
\begin{equation*}
\begin{array}{rcl}
\left[S+\sum_{i=1}^l N_i,\mu\right]_{_{RN}}&=&\mu+\left[\sum_{i=1}^l N_i,l_1\right]_{_{RN}}+\cdots+\left[\sum_{i=1}^l N_i,l_{n-k_l+2}\right]_{_{RN}}+\\
       &+&\left[\sum_{i\not=l} N_i,l_{n-k_l+3}\right]_{_{RN}}+\cdots+\left[\sum_{i\not=l} N_i,l_{n-k_{l-1}+2}\right]_{_{RN}}+\\
       &+&\left[\sum_{i\not=l,l-1} N_i,l_{n-k_l+3}\right]_{_{RN}}+\cdots+\left[\sum_{i\not=l,l-1} N_i,l_{n-k_{l-1}+2}\right]_{_{RN}}+\\
       &+&\cdots+\\
       &+&\left[ N_1,l_{n-k_2+3}\right]_{_{RN}}+\cdots+\left[ N_1,l_{n-k_1+2}\right]_{_{RN}}.
\end{array}
\end{equation*}
\end{prop}
\begin{proof}
Let $1\leq i,j\leq l$. By Remark \ref{degreereason}, any vector valued $(m+k_i-1)$-form, with  $m\geq n-k_i+3$,  is identically zero and any vector valued $(k_i+k_j+m-2)$-form, with $m\geq 1$ is identically zero, because from the conditions $\frac{n+3}{2}\leq k_1\leq \cdots\leq k_l \leq n+1$ and $m\geq 1$ we get $k_i+k_j+m-2\geq n+2.$ Hence,
\begin{enumerate}\label{mesal1}
\item $\left[N_i,l_m\right]_{_{RN}}=0$ for all $m\geq n-k_i+3$ and
\item $\left[N_i,\left[N_j,l_m\right]_{_{RN}}\right]_{_{RN}}=0$ for all  $m\geq 1$.
\end{enumerate}
This implies that
\begin{equation*}
\begin{array}{rcl}
\left[S+\sum_{i=1}^l N_i,\mu\right]_{_{RN}}&=&\mu+\left[\sum_{i=1}^l N_i,l_1\right]_{_{RN}}+\cdots+\left[\sum_{i=1}^l N_i,l_{n-k_l+2}\right]_{_{RN}}+\\
       &+&\left[\sum_{i\not=l} N_i,l_{n-k_l+3}\right]_{_{RN}}+\cdots+\left[\sum_{i\not=l} N_i,l_{n-k_{l-1}+2}\right]_{_{RN}}+\\
       &+&\left[\sum_{i\not=l,l-1} N_i,l_{n-k_l+3}\right]_{_{RN}}+\cdots+\left[\sum_{i\not=l,l-1} N_i,l_{n-k_{l-1}+2}\right]_{_{RN}}+\\
       &+&\cdots+\\
       &+&\left[ N_1,l_{n-k_2+3}\right]_{_{RN}}+\cdots+\left[ N_1,l_{n-k_1+2}\right]_{_{RN}}
\end{array}
\end{equation*}
 and
 \begin{equation*}
 \left[S+\sum_{i=1}^l N_i,\left[S+ \sum_{i=1}^l N_i,\mu\right]_{_{RN}}\right]_{_{RN}}=\mu+2\left[\sum_{i=1}^l N_i,\mu\right]_{_{RN}}=\left[S+2\sum_{i=1}^l N_i,\mu\right]_{_{RN}}.
 \end{equation*}
 Note that it follows from the condition $\frac{n+3}{2}\leq k_1\leq \cdots\leq k_l \leq n+1$ and $m\geq 1$ that for $1\leq i,j\leq l$, we have $k_i+k_j-1\geq n+2$. Hence, $\left[N_i,N_j\right]_{_{RN}}=0$ for all $1\leq i,j\leq l$, which implies that
 \begin{equation*}
 \left[S+\sum_{i=1}^l N_i,S+2\sum_{i=1}^l N_i\right]_{_{RN}}=0.
 \end{equation*}
\end{proof}
\begin{rem}
In Proposition \ref{corlast} one may replace each vector valued $k_i$-form $N_i$ by family of (infinite) symmetric vector valued $k_i$-forms.
\end{rem}
\subsection{Nijenhuis forms on $n$-plectic manifolds}
In this subsection we construct Nijenhuis vector valued forms on a certain type of Lie n-algebras, those which are determined by n-plectic manifolds. Let us recall some definitions from \cite{CRoger}.
\begin{defi}\label{def:nplecticmanifold}
An $n$-plectic manifold, is a manifold $M$ equipped with a non-degenerate, closed $(n+1)$-form $\omega$ and is denoted by $(M,\omega)$.
\end{defi}
An $(n-1)$-form $\alpha$ on an $n$-plectic manifold $(M,\omega)$ is called \emph{Hamiltonian form} if there exists a smooth vector field $\chi_{\alpha}$ on $M$, called \emph{Hamiltonian vector field} associated to $\alpha$, such that $\diff \alpha=-\iota_{\chi_{\alpha}}\omega$. The space of all Hamiltonian forms on an $n$-plectic manifold $(M,\omega)$ is denoted by $\Omega_{Ham}^{n-1}(M)$.\\\\
For two Hamiltonian forms $\alpha, \beta$ on an $n$-plectic manifold $(M,\omega)$, we define a bracket $\{.,.\}$ by \begin{equation}\label{bracketonhamiltonians}
\{\alpha,\beta\}:=\iota_{\chi_{\alpha}}\iota_{\chi_{\beta}}\omega.
\end{equation}
In fact the space of Hamiltonian forms is closed under the bracket $\{.,.\}$.
\begin{prop}\label{prop:weldefined}
Let $\alpha, \beta$ be Hamiltonian forms on an $n$-plectic manifold $(M,\omega)$, with associated Hamiltonian vector fields $\chi_{\alpha},\chi_{\beta}$ respectively. Then $\{\alpha,\beta\}$ is a Hamiltonian form with associated Hamiltonian vector field $[\chi_{\alpha},\chi_{\beta}]$.
\end{prop}
\begin{proof}
see \cite{CRoger}.
\end{proof}
Following \cite{CRoger} we may associate to an $n$-plectic manifold $(M,\omega)$ a symmetric Lie $n$-algebra.
\begin{them}\label{Crogernalgebrafromnplectic}
Let $(M,\omega)$ be an $n$-plectic manifold. Set
$$E_i = \begin{cases}
                      \Omega^{n-1}_{Ham}(M), & \mbox{if } \,\,\,\,\,i=-1,\\
                       \Omega^{n+i}(M), & \mbox{if }  -n\leq i\leq -2
        \end{cases}$$
and $E=\oplus_{i=-n}^{-1}E_i$. Let the collection $\{l_k :E\times\cdots\times E\to E\,\,;\,\,\mbox{with}\,\,k \,\,\mbox{copies of}\,\,E, 1\leq k\leq \infty\}$ of symmetric multi-linear maps be defined as
$$l_1(\alpha) = \begin{cases}
                      (-1)^{|\alpha|}\diff \alpha, & \mbox{if } \,\,\,\,\,\alpha \not\in E_{-1},\\
                       0, & \mbox{if } \,\,\,\,\, \alpha \in E_{-1},
        \end{cases}$$
$$l_{k}(\alpha_1,\cdots,\alpha_k) = \begin{cases}
                      0, & \mbox{if } \alpha_i \not\in E_{-1} \,\,\,\mbox{for some}\,\,\,\, 0\leq i \leq k, \\
                       (-1)^{\frac{k}{2}+1}\iota_{\chi_{\alpha_{1}}}\cdots \iota_{\chi_{\alpha_{k}}}\omega, & \mbox{if }  \alpha_i \in E_{-1} \,\,\,\mbox{for all}\,\,\,\,\,\ \,\,0\leq i \leq k \,\,\,\mbox{and}\,\,\,\, k \,\,\mbox{is even},\\
                       (-1)^{\frac{k-1}{2}}\iota_{\chi_{\alpha_{1}}}\cdots \iota_{\chi_{\alpha_{k}}}\omega, & \mbox{if }  \alpha_i \in E_{-1} \,\,\,\mbox{for all}\,\,\,\,\,\ \,\,0\leq i \leq k \,\,\,\mbox{and}\,\,\,\, k \,\,\mbox{is odd},\\
        \end{cases}$$
        for $k \geq 2$, where $\chi_{\alpha_{i}}$ is the unique Hamiltonian vector field associated to $\alpha_i$.
        Then $(E,(l_k)_{k \geq 1})$ is a symmetric Lie $n$-algebra.
\end{them}
\begin{proof}
In \cite{CRoger}, an $L_{\infty}$-algebra is defined to be a graded vector space $L$ equipped with a collection
${l_k:L^{{\otimes}^k}\to L}$
of skew-symmetric maps, with $\bar{l_k}=k-2$, satisfying a relation so called graded Jacobi identity. However, by translations of degrees in the graded vector space as $L_i \to L_{-i}$, it is equivalent to say an $L_\infty$-algebra is a graded vector space $L$ equipped with a collection
$\{l_k:L^{{\otimes}^k}\to L\}$
of skew-symmetric maps, with $\bar{l_k}=2-k$, satisfying a relation so called graded Jacobi identity. Now, after translating in degrees in the graded vector space as $L_i \to L_{-i}$, it is enough to shift the degrees of the graded vector space in Theorem $3.14.$ in \cite{CRoger}, by $1$, and use the d\'{e}calage isomorphism to get the desired result.
\end{proof}
\begin{rem}
In Theorem \ref{Crogernalgebrafromnplectic}, for all $k\not= 2$, $l_k$ is a map of degree $1$ which guaranties that it is well-defined because the outcome
can never be a $(n-1)$-form. While for the case $k=2$, $l_2(\alpha_1,\alpha_2)=\{\alpha_1,\alpha_2\}$ for all $\alpha_1,\alpha_2 \in E_{-1}$, where $\{.,.\}$ is given by Equation (\ref{bracketonhamiltonians}) and by Proposition \ref{prop:weldefined} this bracket is again in $E_{-1}$, so that $l_2$ is a well-defined bilinear map of degree $1$.
\end{rem}
In the next proposition we give an example of a Nijenhuis vector valued form, with respect to the $L_\infty$-algebra (Lie $n$-algebra) structure associated to a given $n$-plectic manifold, which is the sum of a symmetric vector valued $1$-form with a symmetric vector valued $i$-form, with $i=2,\cdots,n$.

\begin{prop}\label{thm:NijenN}
Let $(M,\omega)$ be an $n$-plectic manifold with the associated symmetric Lie $n$-algebra structure $\mu=l_1+\cdots+l_{n+1}$. For any $n$-form $\eta$ on the manifold $M$,  and any $i=2,\dots,n$, define $\widetilde{\eta}_i$ to be the symmetric vector valued $i$-form of degree zero given by
\begin{equation}\label{tilde}
\widetilde{\eta_i}(\beta_1,\cdots,\beta_i)=\begin{cases}
                                            \iota_{\chi_{\beta_1}}\cdots\iota_{\chi_{\beta_i}} \eta, & \mbox{if}\,\,\,\,\beta_i \in E_{-1} \,\,\,\mbox{for all}\,\,\,\,\,\ \,\,2\leq i \leq n,\\
                                           0, & \mbox{otherwise},
                                           \end{cases}
\end{equation}
where $\chi_{\beta_1},\cdots,\chi_{\beta_n}$ are the Hamiltonian vector fields of $\beta_1,\cdots,\beta_n,$ respectively. Then
 \begin{enumerate}
 \item $S+\widetilde{\eta}_i$ is a Nijenhuis vector valued form with respect to $\mu$ of square $S+2\widetilde{\eta}_i$. The deformed structure is
     \begin{equation*}
     [S+\widetilde{\eta_i},\mu]_{_{RN}}=\mu+[\widetilde{\eta_i},l_1]_{_{RN}}+[\widetilde{\eta_i},l_2]_{_{RN}}.
     \end{equation*}
 \item Moreover, for every pair $\eta,\xi$ of $n$-forms on the manifold $M$ and every $i,j \in \{2,\dots,n\}$,
 the Nijenhuis forms $S+\widetilde{\eta}_i$  and $S+\widetilde{\xi}_j$ are compatible.
 \end{enumerate}
\end{prop}
Proof of Proposition \ref{thm:NijenN} is based on the following lemma.
\begin{lem}\label{lem5ghesmati} For all $2\leq i\leq n$, and all homogeneous elements $\alpha_1,\cdots,\alpha_{i}\in E$ we have:\\
\begin{enumerate}
\item[(1)] $\widetilde{\eta_i}$ vanishes on $\oplus_{i=-n}^{-2} E_i$ and takes values in $ E_{-i},$
\item[(2)] $\widetilde{\eta_i}(l_1(\alpha_1),\alpha_2,\cdots,\alpha_{i})=0,$
\item[(3)] $[\widetilde{\eta_i},l_m]_{_{RN}}=0 ;\mbox{for all} \,\,\, m\geq 3,$
\item[(4)] $[\widetilde{\eta_i},l_1]_{_{RN}}= \diff \circ \widetilde{\eta_i},$
\item[(5)] $[\widetilde{\eta_i},l_2]_{_{RN}}=-\iota_{l_2}\widetilde{\eta_i},$
\item[(6)] $[\widetilde{\eta_i}[\widetilde{\eta_i},l_1]_{_{RN}}]_{_{RN}}=0,$
\item[(7)] $[\widetilde{\eta_i}[\widetilde{\eta_i},l_2]_{_{RN}}]_{_{RN}}=0.$

\end{enumerate}
\end{lem}
\begin{proof} (of the lemma)
Item (1) holds by definition of $\widetilde{\eta_i}$. Let $\alpha_1\in E_{-2}$ and $l_1(\alpha_1)$ be a Hamiltonian form with the associated Hamiltonian vector field $\chi_{l_1(\alpha_1)}$, then we have
\begin{equation*}
\iota_{\chi_{l_1(\alpha_1)}}\omega=-\diff (l_1(\alpha_1))=-\diff^2 \alpha_1=0.
 \end{equation*}
 Hence $\chi_{l_1(\alpha_1)}=0$, since $\omega$ is non-degenerate. This proves item (2). Since $\widetilde{\eta_i}$ takes value in $E_{-i}$ and since $i\geq 2$, definition of $l_m$ implies that
\begin{equation}\label{item3-1}
l_m(\widetilde{\eta_i}(\alpha_1,\cdots,\alpha_{i}),\cdots,\alpha_{m+i-1})=0,
 \end{equation}
 for all $m\geq 3$ and $\alpha_1,\cdots,\alpha_{i+m-1}\in E$. Since $l_m$ takes value in $E_{-m+1}$, we have
 \begin{equation}\label{item3-2}
\widetilde{\eta_i}( l_m(\alpha_1,\cdots,\alpha_{m}),\cdots,\alpha_{m+i-1})=0,
 \end{equation}
 for all $m\geq 3$ and $\alpha_1,\cdots,\alpha_{m+i-1}\in E.$
 Equations (\ref{item3-1}) and (\ref{item3-2}) imply item (3). From item (2) and definition of $\widetilde{\eta_i}$ we get item (4). Since $\widetilde{\eta_i}$ takes value in $E_{-i}$ we have $\iota_{\widetilde{\eta_i}}l_2=0$, hence $[\widetilde{\eta_i},l_2]_{_{RN}}=-\iota_{l_2}\widetilde{\eta_i}$. The same argument as in the proof of item (1) shows that the Hamiltonian vector field associated to a closed Hamiltonian form is zero. Hence
 \begin{equation}\label{item6-1}
 \begin{array}{rcl}
  &&\iota_{[\widetilde{\eta_i},l_1]_{_{RN}}}\widetilde{\eta_i}(\alpha_1,\cdots,\alpha_{2i-1})\\
  &=&\widetilde{\eta_i}\big(\diff (\widetilde{\eta_i}(\alpha_1,\cdots,\alpha_{i})),\cdots,\alpha_{2i-1}\big)\\
  &=&0,
  \end{array}
 \end{equation}
 for all $\alpha_1,\cdots,\alpha_{2i-1}\in E.$ Since $i\geq 2$, item (4) and the fact that $\widetilde{\eta_i}$ takes value in $E_{-i}$ imply that
 \begin{equation}\label{item6-2}
 \iota_{\widetilde{\eta_i}}[\widetilde{\eta_i},l_1]_{_{RN}}(\alpha_1,\cdots,\alpha_{2i-1})=0.
 \end{equation}
 Equations (\ref{item6-1}) and (\ref{item6-2}) imply item (6).
Since $\widetilde{\eta_i}$ does not take value in $E_{-1}$,\\ $l_2(\widetilde{\eta_i}(\alpha_1,\cdots,\alpha_{i}),\alpha_{i+1})=0$, for all $\alpha_1,\cdots,\alpha_{i+1}\in E$. Hence, using item (5) we get
 \begin{equation}\label{item7-1}
 \iota_{\widetilde{\eta_i}}[\widetilde{\eta_i},l_2]_{_{RN}}=0.
 \end{equation}
 A similar argument shows that
  \begin{equation}\label{item7-2}
 \iota_{[\widetilde{\eta_i},l_2]_{_{RN}}}\widetilde{\eta_i}=0.
 \end{equation}
 Equations (\ref{item7-1}) and (\ref{item7-2}) give item (7).

\end{proof}
\begin{proof}(of Proposition \ref{thm:NijenN})
Let $\eta$ and $\xi$ be two $n$-forms on the manifold $M$ and $i,j\geq 2$. It follows from definitions of $\widetilde{\eta_i}$, $\widetilde{\xi_j}$ and item $1$ in the Lemma \ref{lem5ghesmati} that both $\iota_{\widetilde{\eta_i}}\widetilde{\xi_j}$ and $\iota_{\widetilde{\xi_j}}\widetilde{\eta_i}$ vanish. This implies that
\begin{equation}\label{haminzeperty}
[S+\widetilde{\eta_i},S+\widetilde{\xi_j}]_{_{RN}}=0,
\end{equation}
for all $i,j \geq 2$.
Let $i\geq 2$. Item $3$ in Lemma \ref{lem5ghesmati}  implies that
\begin{equation}\label{themnplectic}
[S+\widetilde{\eta_i},\mu]_{_{RN}}=\mu+[\widetilde{\eta_i},l_1]_{_{RN}}+[\widetilde{\eta_i},l_2]_{_{RN}}.
\end{equation}
Applying $[S+\widetilde{\eta_i},.]_{_{RN}}$ to both sides of Equation (\ref{themnplectic}) and using items $3,6,7$ in the Lemma \ref{lem5ghesmati} we get
\begin{equation*}
[S+\widetilde{\eta_i},[S+\widetilde{\eta_i},\mu]_{_{RN}}]_{_{RN}}=\mu+2[\widetilde{\eta_i},l_1]_{_{RN}}+2[\widetilde{\eta_i},l_2]_{_{RN}}.
\end{equation*}
Therefore, again by using item $3$ in the Lemma \ref{lem5ghesmati} we have
\begin{equation}\label{eta-i}
[S+\widetilde{\eta_i},[S+\widetilde{\eta_i},\mu]_{_{RN}}]_{_{RN}}=[S+2\widetilde{\eta_i},\mu]_{_{RN}}.
\end{equation}
 Equation (\ref{haminzeperty}), with $\widetilde{\xi_j}=\widetilde{\eta_i}$ together with Equation (\ref{eta-i}) prove item $1$ in the proposition. Now, item (2) follows directly from Equation (\ref{haminzeperty}).
\end{proof}
From Proposition \ref{thm:NijenN} we immediately get the following result.
\begin{them}\label{them:Nijenhuisonnplectic}
Let $\eta$ be an arbitrary $n$-form on an $n$-plectic manifold $(M,\omega)$. Let $(E=E_{-n}\oplus \cdots\oplus E_{-1}, \mu=l_1+\cdots+l_{n+1})$ be the Lie $n$-algebra associated to $(M,\omega)$. For each $2\leq i\leq n$ define the maps $\widetilde{\eta_i}$ as in (\ref{tilde}).
 Then $\mathcal{N}:=S+\sum_{i=2}^{n}\widetilde{\eta_i}$ is a Nijenhuis vector valued form with respect to the the Lie $n$-algebra structure $\mu=l_1+\cdots+l_{n+1}$, associated to the $n$-plectic manifold $(M,\omega)$, with square $S+2\sum_{i=2}^{n}\widetilde{\eta_i}$. Moreover, the deformed structure $ [\mathcal{N}, \mu]_{_{RN}} =\sum_{i=1}^{n+1} l^{\mathcal N}_{i} $, with $l^{\mathcal N}_{i}$ being the component in the vector valued form  $[\mathcal{N}, \mu]_{_{RN}}$ which is a vector valued $i$-form, is given by:
\begin{equation*}
l^{\mathcal N}_{i}=\begin{cases}
                                l_1, &  \mbox{for} \,\,\, i=1,\\
                                l_i+\diff \circ \widetilde{\eta_i}-\iota_{l_2} \widetilde{\eta_{i-1}}, & \mbox{for} \,\,\, i\geq 2.
                   \end{cases}
\end{equation*}
\end{them}
A special case of the previous theorem is considered in the next proposition.
\begin{prop}
Let $(M,\omega)$ be an $n$-plectic manifold and $\alpha$ be a Hamiltonian form on $(M,\omega)$. For each $2\leq i\leq n$ define the maps $\widetilde{\alpha_i}$ as
\begin{equation*}
\widetilde{\alpha_i}(\beta_1,\cdots,\beta_i)=\begin{cases}
                                             \iota_{\chi_{\alpha}}\iota_{\chi_{\beta_1}}\cdots\iota_{\chi_{\beta_i}}\omega, & \mbox{if} \,\,\,\ \beta_k \in E_{-1}\,\,\, \mbox{for all} \,\,\,1\leq k\leq i\\
                                             0, & \mbox{otherwise}
                                             \end{cases}
\end{equation*}where $ \chi_{\alpha},\chi_{\beta_1},\cdots,\chi_{\beta_i}$ are the unique Hamiltonian vector fields associated to the Hamiltonian forms $\alpha,\beta_1,\cdots,\beta_i$ respectively. Then $S+\sum_{i=2}^{n}\widetilde{\alpha_i}$ is a Nijenhuis vector valued form with respect to the the Lie $n$-algebra structure $\mu=l_1+\cdots+l_{n+1}$, associated to the $n$-plectic manifold $(M,\omega)$.
\end{prop}
Theorem \ref{them:Nijenhuisonnplectic} can be easily generalized if, instead of taking one $n$-form on the manifold $M$, we take a family of $n$-forms on $M$.
\begin{them}\label{nplecticlasttheorem}
Let $(\eta^j)_{j\geq 1}$ be a family of $n$-forms on an $n$-plectic manifold $(M,\omega)$. Let $(E=E_{-n}\oplus \cdots\oplus E_{-1}, \mu=l_1+\cdots+l_{n+1})$ be the Lie $n$-algebra associated to $(M,\omega)$. For each $2\leq i\leq n$ define the vector valued $i$-forms $\widetilde{\eta^j_i}$ as
\begin{equation*}
\widetilde{\eta^j_i}(\beta_1,\cdots,\beta_i)=\begin{cases}
                                             \iota_{\chi_{\beta_1}}\cdots\iota_{\chi_{\beta_i}}\eta^j, & \mbox{if} \,\,\,\ \beta_k \in E_{-1}\,\,\, \mbox{for all} \,\,\,1\leq k\leq i\\
                                             0, & \mbox{otherwise}
                                             \end{cases}
\end{equation*}where $ \chi_{\beta_1},\cdots,\chi_{\beta_i}$ are the unique Hamiltonian vector fields associated to the Hamiltonian forms $\beta_1,\cdots,\beta_i$ respectively. Then $\mathcal{N}:=S+\sum_{j\geq 1}\sum_{i=2}^{n}\widetilde{\eta^j_i}$ is a Nijenhuis vector valued form with respect to the the Lie $n$-algebra structure $\mu=l_1+\cdots+l_{n+1}$, associated to the $n$-plectic manifold $(M,\omega)$.
\end{them}

\subsection{Generalities on Lie $2$-algebras and crossed modules}\label{subsection:Liealgebrasandcros}


In this subsection we shall consider Lie $2$-algebras, that is,  graded vector spaces $E$ which are concentrated in degrees $-1$ and $-2$, equipped with $L_\infty$-structures.

As we have already remarked, by degree reasons, a Lie $2$-algebra structure on a graded vector space $E=E_{-2}\oplus E_{-1}$ has to be of the form $l_1+l_2+l_3$, with $l_1, l_2,l_3$ being symmetric vector valued $1$-form, $2$-form and $3$-form, respectively, of degree $+1$.

We will discuss about Courant algebroids in Section \ref{subsection:3.3.5}. So, in view of that, it is convenient to   introduce the following notations for Lie $2$-algebras. Assume that $\mu=l_1+l_2+l_3$ is a Lie $2$-algebra structure on a graded vector space $E=E_{-2}\oplus E_{-1}$. Then, the  binary bracket $l_2$ has two parts:
\begin{equation*}
l_{2}'':E_{-2}\times E_{-1}\rightarrow E_{-2}\,\,\,\mbox{and}\,\,\, l_{2}':E_{-1}\times E_{-1}\rightarrow E_{-1}.
\end{equation*}
 We denote
 \begin{equation}
 \begin{array}{rcl}
 l_1(f)        &\mbox{ by} & \partial f,\\
 l_2'(X,Y)     &\mbox{ by} &[X,Y]_{2},\\
 l_{2}''(f,X)  & \mbox{ by}& \chi(X)(f),\\
 l_3(X,Y,Z)    &\mbox{ by} &\omega(X,Y,Z),\\
 \end{array}
 \end{equation}
 for all $X,Y,Z \in E_{-1}, f\in E_{-2}$, with $\chi:E_{-1}\rightarrow End(E_{-2})$.
In general, to every vector valued form $\mu=l_1+l_2+l_3$, with $l_i$ being a vector valued $i$-form, $i=1,2,3$, on a graded vector space $E$ concentrated in degrees $-2$ and $-1$, one may correspond a quadruple $(\partial,\chi, \left[.,.\right]_2,\omega)$ as in above. The following proposition shows how these kind of quadruples are related to Lie $2$-algebras.
\begin{prop}\label{prop:quadrubleLie2algebra}
 Let $E=E_{-2}\oplus E_{-1}$ be a graded vector space. For $k=1,2,3$ let $l_k\in S^k(E^*)\otimes E$ and denote
 \begin{equation*}
 \begin{array}{rcl}
 \partial f&:=&l_1(f),\\
 \left[X,Y\right]_2&:=&l_2(X,Y),\\
 \chi(X)f&:=&l_2(X,f),\\
 \omega(X,Y,Z)&:=&l_3(X,Y,Z),
 \end{array}
 \end{equation*}
 for all $X,Y,Z \in E_{-1}$ and $f\in E_{-2}$. Then $\mu=l_1+l_2+l_3$ is a Lie $2$-algebra structure on $E$ if, and only if, for the quadruple $(\partial,\chi, \left[.,.\right]_2,\omega)$ the following relations hold for all $X,Y,Z,W \in E_{-1}$ and $f,g\in E_{-2}$:
 \begin{itemize}
 \item [(i)] $\chi(\partial f)g+\chi(\partial g)f=0$
 \item[(ii)] $[X,\partial f]_2=\partial\chi(X)f,$
 \item[(iii)] $\chi([X,Y]_2)f+\chi(Y)\chi(X)f-\chi(Y)\chi(X)f+\omega(X,Y,\partial f)=0,$
 \item[(iv)] $[[X,Y]_2,Z]_2+c.p.= \partial \omega (X,Y,Z),$
 \item[(v)]
 $
\chi(W)\omega(X,Y,Z)-\chi(Z)\omega(X,Y,W)+\chi(Y)\omega(X,Z,W)-\chi(X)\omega(Y,Z,W)=\\
-\omega([X,Y]_2,Z,W)+\omega([X,Z]_2,Y,W)-\omega([X,W]_2,Y,Z)\\
-\omega([Y,Z]_2,X,W)+\omega([Y,W]_2,X,Z)-\omega([Z,W]_2,X,Y).
$
 \end{itemize}
 \end{prop}
\begin{proof} Let $X,Y,Z,W\in E_{-1}$ and $f,g\in E_{-2}$. Then, the result follows directly from the followings:
$[l_1,l_2]_{_{RN}}(f,g)=0$ is equivalent to (i).
$[l_1,l_2]_{_{RN}}(X,f)=0$ is equivalent to (ii).
$(2[l_1,l_3]_{_{RN}}+[l_2,l_2]_{_{RN}})(X,Y,f)=0$ is equivalent to (iii).
$(2[l_1,l_3]_{_{RN}}+[l_2,l_2]_{_{RN}})(X,Y,Z)=0$ is equivalent to (iv).
$[l_2,l_3]_{_{RN}}(X,Y,Z,W)=0$ is equivalent to (v). Note that, for degree reasons, all the missing cases, for example the case $(2[l_1,l_3]_{_{RN}}+[l_2,l_2]_{_{RN}})(X,f,g)$ are identically zero.
\end{proof}

In view of Proposition \ref{prop:quadrubleLie2algebra}, we may denote a Lie $2$-algebra structure $\mu=l_1+l_2+l_3$ on a graded vector space $E=E_{-2}\oplus E_{-1}$, also, by a quadruple $(\partial,\chi, \left[.,.\right]_2,\omega)$ with $\chi:E_{-1}\to End(E_{-2})$ satisfying all the relations (i)-(iv) in Proposition \ref{prop:quadrubleLie2algebra}.

%

Now, we consider some special cases of Lie $2$-algebras. \\
A Lie $2$-algebra $(E=E_{-2}\oplus E_{-1},\,\mu=l_1+l_2+l_3)$, with $l_2=l_3=0$ and $l_1$ invertible is called \emph{trivial Lie $2$-algebra}. We may also consider the notions of string Lie algebras and crossed module of Lie algebras as special cases of Lie $2$-algebras. Before talking about string Lie algebras let us introduce the notion of Chevalley-Eilenberg differential of a vector valued form on a graded vector space $E$ concentrated on degrees $-2$ and $-1$. Next lemma gives a motivation for the definition of Chevalley-Eilenberg differential.
  \begin{lem}\label{Koskhlane}
  Let $(E=E_{-2}\oplus E_{-1},\,\mu=l_1+l_2+l_3)$ be a Lie $2$-algebra with corresponding quadruple $(\partial,\chi, \left[.,.\right]_2,\omega)$, in view of Proposition \ref{prop:quadrubleLie2algebra}. Let $\eta \in S^k(E^*)\otimes E$ be a vector valued $k$-form of degree $k-2$. Then,
\begin{equation}\label{eq:CEilengerg}
[\eta,l_2]_{_{RN}} (X_0 , \dots, X_k) = \sum_{i=0}^k (-1)^i \chi(X_i) \eta(\widehat{X_i})+\sum_{0\leq i<j\leq k}(-1)^{i+j}
 \eta([X_i,X_j], \widehat{X_{i,j}}),
 \end{equation}
   for all $X_0,\dots,X_k \in E_{-1}$ , where $\widehat{X_{i}}$ and $\widehat{X_{i,j}}$ stand for
   \begin{equation*}
   X_1,\cdots, X_{i-1}, X_{i+1}, \cdots, X_{k}\quad and \quad X_1,\cdots, X_{i-1}, X_{i+1}, \cdots, X_{j-1}, X_{j+1}, \cdots, X_{k}
    \end{equation*}
    respectively.
    \end{lem}
    \begin{proof}
    For degree reasons $\eta$ has to be of the form
$\eta:E_{-1}\times\cdots\times E_{-1}\to E_{-2}$, with $k$ copies of $E_{-1}$. Hence, Equation (\ref{eq:CEilengerg}) is a direct consequence of the definition of Richardson-Nijenhuis bracket.
    \end{proof}
    Considering Lemma \ref{Koskhlane} and regarding to the Equation (\ref{def:d^A}) we may define the Chevalley-Eilenberg differential as follows:
    \begin{defi}
    Let $E=E_{-2}\oplus E_{-1}$ be a graded vector space concentrated on degrees $-2$ and $-1$, $S_{k}(E)\subset S^{k}(E^*)\otimes E$ be the subspace of all symmetric vector valued $k$-forms of degree $k-2$ and $S^{\bullet}(E):=\oplus_{k\geq 1}S_{k}(E)$. Let $\chi:E_{-1} \to End(E_{-2})$ be a representation of vector spaces and $\left[.,.\right]:E_{-1} \times E_{-1}\to E_{-1}$ be a graded symmetric bilinear map. Then the Chevalley-Eilenberg differential $\diff^{CE}:S^{\bullet}(E) \to S^{\bullet+1}(E)$ is defined by
\begin{equation}\label{eq:CEilenberggeneral}
\diff^{CE}\eta(X_0 , \dots, X_k) := \sum_{i=0}^k (-1)^i \chi(X_i) \eta(\widehat{X_i})+\sum_{0\leq i<j\leq k}(-1)^{i+j}
 \eta([X_i,X_j], \widehat{X_{i,j}}),
 \end{equation}
   for all $X_0,\dots,X_k \in E_{-1}$ , where $\widehat{X_{i}}$ and $\widehat{X_{i,j}}$ stand for
   \begin{equation*}
   X_1,\cdots, X_{i-1}, X_{i+1}, \cdots, X_{k}\quad and \quad X_1,\cdots, X_{i-1}, X_{i+1}, \cdots, X_{j-1}, X_{j+1}, \cdots, X_{k}
    \end{equation*}
    respectively.
    \end{defi}
In general, the operator $\diff^{CE}$ does not square to zero. However, according to Lemma \ref{Koskhlane} it can be written as
  $$  {\diff}^{CE}  = [.,l_2 ]_{_{RN}},$$
  and using the graded Jacobi identity, we have that $\diff^{CE}$ squares to zero, if and only if, $[l_2,l_2]_{_{RN}}=0.$ Here  $l_2(X,f)=\chi(X)f$ and $l_2(X,Y)=[X,Y],$ for all $X,Y\in E_{-1}, f\in E_{-2}.$

  Now, let $(E=E_{-2}\oplus E_{-1},\,\mu=l_1+l_2+l_3)$ be a Lie $2$-algebra with corresponding quadruple $(\partial,\chi, \left[.,.\right]_2,\omega)$, in view of Proposition \ref{prop:quadrubleLie2algebra}. If $l_1=0$, then for all $X,Y,Z \in E_{-1}$
     \begin{equation*}
        0=[l_2,l_2]_{_{RN}}(X,Y,Z)=2([[X,Y]_2,Z]_2+c.p.)
     \end{equation*}
       which means that $\left[.,.\right]_2$ is a Lie bracket on $E_{-1} $. And for all  $X,Y \in E_{-1},f \in E_{-2}$
    \begin{equation*}
       [l_2,l_2]_{_{RN}}(X,Y,f)=\chi[X,Y]_2(f)-\chi(Y)\chi(X)f+\chi(X)\chi(Y)f=0
    \end{equation*}
       which means that $\chi $ is a representation of $E_{-1}$ on $E_{-2}$.
      Also, the condition $[l_2,l_3]_{_{RN}}=0$ means that $\omega$  is a Chevalley-Eilenberg-closed 3-form
      of this Lie algebra $E_{-1}$ valued in $E_{-2}$. We call \emph{string Lie algebras}
      this kind of Lie $2$-algebras.

  Next we want to explain how  a crossed module of Lie algebras can be seen as a Lie $2$-algebra. Let us first, recall the definition of a crossed module of Lie algebras from \cite{FWagemann}:
\begin{defi}
A crossed module of Lie algebras $(\mathfrak g,\, \left[.,.\right]^{\mathfrak g})$ and $(\mathfrak h,\, \left[.,.\right]^{\mathfrak h})$ is a homomorphism  $\partial:\mathfrak g \to \mathfrak h$ together with an action $\chi$ of $\mathfrak h$ on $\mathfrak g$ by derivation, that is a linear map $\chi: \mathfrak h \to Hom(\mathfrak g,\mathfrak g)$ such that
\begin{equation}\label{crossedmodule1}
\partial(\chi(h)g)=[h,\partial(g)]^{\mathfrak h};\,\,\,\,
\text{for all} \,\,\,g \in\mathfrak g   \,\,\text{and all} \,\,\,  h \in \mathfrak h
\end{equation}
and
\begin{equation}\label{crossedmodule2}
\chi(\partial(g_1))g_2= [g_1,g_2]^{\mathfrak g};\,\,\,\, \text{for all} \,\,\,g_1 ,g_2\in \mathfrak g.
\end{equation}
\end{defi}
Such a crossed module will be denoted by $(\mathfrak g,\,\mathfrak h,\,\partial,\,\chi)$. The next proposition shows how one can see a crossed module as a Lie $2$-algebra.
\begin{prop}
Let $(E=E_{-2}\oplus E_{-1},\,\mu=l_1+l_2+l_3)$ be a Lie $2$-algebra with corresponding quadruple $(\partial,\chi, \left[.,.\right]_2,\omega)$, in view of Proposition \ref{prop:quadrubleLie2algebra}.
If $l_3=0$, then $(E_{-2}, E_{-1}, \partial, \chi)$ is a crossed module of Lie algebras.
\end{prop}
\begin{proof}
Since $(E=E_{-2}\oplus E_{-1},\,\mu=l_1+l_2+l_3)$ is a Lie $2$-algebra we have $[l_2,l_2]_{_{RN}}=0$
which means that $\left[.,.\right]_2$ is a Lie bracket on $E_{-1} $.
   The remaining condition, i.e., $[l_1,l_2]_{_{RN}}=0$ allow us to define
   a crossed-module structure on $E_{-2} \mapsto E_{-1}$ as follows:
      Define a bracket $\left[.,.\right]$ on $E_{-2}$ as $[f,g]:= \chi(\partial (f))g $ for all $f,g \in E_{-2}$. From item (i) in Proposition \ref{prop:quadrubleLie2algebra} we get $ [f,g]=-[g,f]$. While applying $[l_1,l_2]_{_{RN}}$ to $ (\partial (f),g)$ we obtain
     \begin{equation}\label{partialpartial}
         \begin{array}{rcl}
             l_2(\partial f,\partial g)&=&\partial l_2(\partial f,g)\\
                                       &=&\partial(\chi(\partial f)g)\\
                                       &=&\partial[f,g],
        \end{array}
     \end{equation}
     for all $f,g\in E_{-2}$. On the other hand, for $f,g,h\in E_{-2}$ we have:
     \begin{equation*}
     \begin{array}{rcl}
     0&=&\frac{1}{2}\left[l_2,l_2\right]_{_{RN}}(\partial f,\partial g,h)\\
      &=&l_2(l_2(\partial f, \partial g),f)+l_2(l_2(\partial f,h),\partial g)+l_2(l_2(\partial g,h)\partial f)\\
      &=&\chi(l_2(\partial f, \partial g)h)+\chi(\partial g)\chi(\partial f).h+\chi(\partial f)\chi(\partial g)h.
      \end{array}
     \end{equation*}
 Hence,
 \begin{equation}\label{chichi}
 \chi[\partial f, \partial g]_2h=\chi(\partial f)\chi(\partial g)h-\chi(\partial g)\chi(\partial f)h,
 \end{equation}
 for all $f,g,h\in E_{-2}$.
 Using (\ref{chichi}) and (\ref{partialpartial}), the Jacobi identity for the bracket $\left[.,.\right]$ on $E_{-2}$ can be proved as follows:
        \begin{equation*}
           \begin{array}{rcl}
                          \left[[f,g],h\right]&=&\chi \partial([f,g])h \\
                                              &=&\chi [\partial(f),\partial(g)]_2h \\
                                              &=&\chi(\partial f)\chi(\partial g)h-\chi(\partial g)\chi(\partial f)h\\
                                              &=&\chi(\partial f)[g,h]-\chi(\partial g)[f,h]\\
                                              &=&[f,[g,h]]-[g,[f,h]]
           \end{array}
        \end{equation*}
 or
 \begin{equation*}
 [f,[g,h]]=[[f,g],h]+[g,[f,h]].
 \end{equation*}
  Therefore, $(E_{-2},\, \left[.,.\right])$ is a Lie algebra, $\partial$ is a Lie algebra morphism and $\chi $ is a representation of $E_{-1}$ on $E_{-2}$.
   Applying $[l_1,l_2]_{_{RN}}$ to $(X,f)$, for $X\in E_{-1}$ and $f\in E_{-2}$, we get
        \begin{equation*}
          l_1(l_2(X,f))=l_2(X,l_1f)
        \end{equation*}
        or
        \begin{equation*}
           \partial(\chi(X)f)=[X,\partial f]_2.
        \end{equation*}
   Also by definition of $[f,g]$ we have $$\chi(\partial(f))g= [f,g].$$
  Hence, $(E_{-2},\,,E_{-1},\,\partial,\,\chi)$ is a crossed module.
\end{proof}
\subsection{Nijenhuis forms on Lie $2$-algebras}
Proposition \ref{khodemoun} of Subsection \ref{subsection:Lienalgebras} provides the construction of Nijenhuis forms on Lie $n$-algebras. However, for the case n=2, that proposition does not give the possibility of having a Nijenhuis vector valued $2$-form. We intend to give an example of Nijenhuis vector valued form with respect to a Lie $2$-algebra structure $\mu$ on a graded vector space $E_{-2}\oplus E_{-1}$ which is not purely a $1$-form, i.e. not just a collection of maps from $E_i$ to $E_i$. As we have mentioned before, elements of degree $0$ in $\tilde{S}(E^*) \otimes E $ are necessarily of the form $N + \alpha$ with $N: E \to E$ a linear endomorphism preserving the degree and $\alpha : E \times E \to E$
a symmetric vector valued $2$-form (or $\alpha : E_{-1} \times E_{-1} \to E_{-2}$ a skew-symmetric $E_{-2}$-valued $2$-form on $E_{-2}$). 
\begin{them}\label{thm:alphalpha}
Let $\mu=l_1+l_2+l_3$ be a Lie $2$-algebra structure on a graded vector space $E=E_{-2}\oplus E_{-1}$ and $\alpha$ be a symmetric vector valued $2$-form of degree $0$. Then  $S+\alpha$ is a Nijenhuis vector valued form with respect to $\mu$ with square of $S+2\alpha $ if and only if
 $$ \alpha(l_1 \alpha (X,Y),Z ) +c.p.=0,$$
 for all $ X,Y,Z \in E_{-1}$.
\end{them}
\begin{proof}
By degree reason $[\alpha,[\alpha, l_1]_{_{RN}}]_{_{RN}}$ is of the form  $$[\alpha,[\alpha, l_1]_{_{RN}}]_{_{RN}}: E_{-1}\otimes E_{-1}\otimes E_{-1} \to E_{-2}$$ and
\begin{equation}\label{eq:Salpha0}
   \begin{array}{rcl}
       [\alpha,[\alpha, l_1]_{_{RN}}]_{_{RN}}(X,Y,Z)&=& [\alpha, l_1]_{_{RN}}(\alpha(X,Y),Z)+c.p.-\alpha([\alpha, l_1]_{_{RN}}(X,Y),Z)+c.p.\\
                                       &=& -2\alpha(l_1(\alpha(X,Y)),Z)+c.p\\

   \end{array}
\end{equation} for all $X,Y,Z \in E_{-1}$.
Again by degree reasons $[\alpha,[\alpha, l_2]_{_{RN}}]_{_{RN}}$ and $[\alpha, l_3]_{_{RN}}$ ($4$-forms of degree $+1$) are identically zero. So
\begin{equation}\label{eq:Salpha}
    \begin{array}{rcl}
                                            &&[S+\alpha,[S+\alpha,l_1+l_2+l_3]_{_{RN}}]_{_{RN}} \\
                                            &=& [S+\alpha,l_1+l_2+l_3+[\alpha, l_1]_{_{RN}}+[\alpha, l_2]_{_{RN}}]_{_{RN}}\\
                                            &=&l_1+l_2+l_3+2[\alpha, l_1]_{_{RN}}+2[\alpha, l_2]_{_{RN}}+[\alpha,[\alpha, l_1]_{_{RN}}]_{_{RN}}\\
                                            &=&l_1+l_2+l_3+[2\alpha, l_1]_{_{RN}}+[2\alpha, l_2]_{_{RN}}+[\alpha,[\alpha, l_1]_{_{RN}}]_{_{RN}}\\
                                            &=&[S+2\alpha,l_1+l_2+l_3]_{_{RN}}+[\alpha,[\alpha, l_1]_{_{RN}}]_{_{RN}}.

    \end{array}
\end{equation}
On the other hand Lemma \ref{Euler} implies that
\begin{equation}\label{Salpha2}
[S+\alpha,S+2\alpha]_{_{RN}}=[S,S]_{_{RN}}+[S,\alpha]_{_{RN}}+[\alpha,\alpha]_{_{RN}}=0.
\end{equation}
Equations (\ref{eq:Salpha0}), (\ref{eq:Salpha}) and (\ref{Salpha2}) show that $S+\alpha$ is a Nijenhuis vector valued form with respect to $\mu$ with square $S+2\alpha$, if and only if, $$\alpha(l_1(\alpha(X,Y)),Z)+c.p=0.$$
\end{proof}
\begin{cor}
Let $(E,\mu)$ be a string Lie algebra (see Subsection \ref{subsection:Liealgebrasandcros}.) Then for every vector valued $2$-form $\alpha$ of degree zero, $S+\alpha$ is a Nijenhuis vector valued form with respect to $\mu$.
\end{cor}
In the following example we see an application of Theorem \ref{thm:alphalpha} to a trivial Lie $2$-algebra.

\begin{examp}
Let $\mathfrak{g}$ be a vector space and $\left[.,.\right]_{\mathfrak{g}}$ be a skew-symmetric bilinear map on $\mathfrak{g}$. Let $E_{-1}:=\{-1\}\times\mathfrak{g}$, $E_{-2}:=\{-2\}\times\mathfrak{g}$ and let $\partial:E_{-2}\to E_{-1}$ be given by $(-2,x)\mapsto (-1,x)$. Define $\alpha:E_{-1}\times E_{-1}\to E_{-2}$ to be vector valued $2$-form on the graded vector space $E=E_{-2}\oplus E_{-1}$ as $((-1,x),(-1,y))\mapsto (-2,[x,y]_{\mathfrak{g}})$. Then $S+\alpha$ is Nijenhuis with respect to $\partial$, if and only if, $\left[.,.\right]_{\mathfrak{g}}$ is a Lie bracket, which is a direct consequence of Theorem \ref{thm:alphalpha}.
\end{examp}
Let us see now what the deformed Lie $2$-algebra structure looks like concretely.
The definitions of $\partial, \chi,\left[.,.\right] $ and $\omega$ are intended to make the
interpretation of the following result clear:

\begin{prop}\label{prop:deformed}
Let $\mu$ be a Lie $2$ algebra structure on a graded vector space $E_{-2}\oplus E_{-1}$, corresponding to the quadruple $(\partial, \left[.,.\right]_2,\chi,\omega)$, in view of Proposition \ref{prop:quadrubleLie2algebra}. Let $\alpha$ be a symmetric vector valued $2$-form of degree zero on $E$ and let ${\mathcal N}= S+\alpha$. If the deformed structure $\mu^{\mathcal N}$ is associated to a quadruple $(\partial', [,]'_2,\chi',\omega')$, then
\begin{equation}\label{deformedLie2algebra}
\begin{array}{rcl}
\partial'(f)&=& \partial(f),\\
\left[X,Y\right]'_{2} &=&[X,Y]_2+\partial \alpha(X,Y), \\
\chi'(X)f   &=& \chi(X)f -\alpha(\partial(f),X),\\
\omega'(X,Y,Z)&=& \omega(X,Y,Z)+{\diff}^{CE}\alpha(X,Y,Z),
\end{array}
\end{equation}
for all $X,Y,Z\in E_{-1}$ and $f\in E_{-2}$.
\end{prop}
\begin{proof}
The statement follows from the following easy relations:

\begin{enumerate}
\item $[S+\alpha,\mu]_{_{RN}}=l_1+(l_2+[\alpha,l_1]_{_{RN}})+(l_3+[\alpha,l_2]_{_{RN}})$,
\item $[\alpha,l_1]_{_{RN}}(X,Y)=l_1\alpha(X,Y);\,\,\, \text{for all}\,\, X,Y \in E_{-1}$,
\item $[\alpha,l_1]_{_{RN}}(f,X)=-\alpha(l_1(f),X)\,\,\, \text{for all}\,\, X \in E_{-1},f \in E_{-2}$,
\item $[\alpha,l_2]_{_{RN}}={\diff}^{CE}\alpha.$
\end{enumerate}

\end{proof}

Indeed, in this case, the Nijenhuis transformations that we consider are invertible, it suffices for this to consider
the deformation by $S-\alpha$. Combining Theorems \ref{thm:alphalpha} and \ref{theo:Hierarchy} we get the following proposition.
\begin{prop}\label{prop:inverseNijenhuis}
Let $\mu=l_1+l_2+l_3$ be a Lie $2$-algebra structure on a graded vector space $E=E_{-2}\oplus E_{-1}$.
Let $\alpha$ be a vector valued $2$-form of degree zero such that $ \alpha(l_1 \alpha (X,Y),Z ) +c.p.=0,$
 for all $ X,Y,Z \in E_{-1}$.
Let $\mu_k$ stands for the structure $\mu$ after $k$ times deformations, that is $\mu_k=[S+\alpha,[S+\alpha,\cdots ,[S+\alpha,\mu]_{_{RN}} \cdots]_{_{RN}}]_{_{RN}}$ with $k$ copies of $S+\alpha$.
Then $S+\alpha$ is a Nijenhuis vector valued form with respect to all the terms of the hierarchy of successive deformations $\mu_k$, with square of $S+2\alpha $.
Moreover, $[S-\alpha,[S+\alpha,\mu]] =\mu $.
\end{prop}
As we have seen previously, Lie $2$-algebras $(E_{-2}\oplus E_{-1},\, l_1+l_2+l_3)$ with $l_1=0$ are called string Lie algebras. They are in one to one correspondence with Lie algebra structures on ${\mathfrak g}:=E_{-1}$
together with a representation of the Lie algebra ${\mathfrak g}$ on the vector space $V:=E_{-2}$ and a Chevalley-Eilenberg $3$-cocycle $\omega$ for this representation,
so that string Lie algebras can be denoted as triples $({\mathfrak g},V,\omega)$ with $\mathfrak g $  a Lie algebra,
$V$  a representation and $\omega$  a Chevalley-Eilenberg 3-cocycle.
 In this case, the deformation by $S+\alpha $
just amounts to change the $3$-cocycle $\omega$ by $\omega + \diff^{CE} \alpha$.
So that, for string Lie algebras, adding up a coboundary, i.e. changing $({\mathfrak g},V,\omega)$
into  $({\mathfrak g},V,\omega+ \diff^{CE}\alpha)$  can be seen as a Nijenhuis transformation by
$S+\alpha$.

Let us now investigate Lie $2$-algebras structures for which $\chi=0$.
There may be quite a few such Lie $2$-algebras but, we are going to show that, after a Nijenhuis transformation of the form
$S+\alpha$, such Lie $2$-algebras will be decomposed as a direct sum of
a string Lie algebra with a trivial Lie $2$-algebra. \\
A \emph{Lie sub-$2$-algebra of a Lie $2$-algebra $(E=E_{-2}\oplus E_{-1}, \mu=l_1+l_2+l_3)$} is a Lie $2$-algebra $(E^\prime=E^\prime_{-2}\oplus E^\prime_{-1}, \mu'=l'_1+l'_2+l'_3)$ with $E^\prime_{-2}\subset E_{-2}$ and $E^\prime_{-1}\subset E_{-1}$ vector sub-spaces,
\begin{equation*}
l_1^{\prime}=l_1|_{E^{\prime}},\,\,l_2^{\prime}=l_2|_{E^{\prime}\times E^{\prime}}\,\,\,\mbox{and}\,\,\,l_3^{\prime}=l_3|_{E^{\prime}\times E^{\prime}\times E^{\prime}}.
\end{equation*}
\begin{prop}\label{prop:strict+trivial}
 Given a Lie $2$-algebra structure $l_1+l_2+l_3$ on a graded vector space $E=E_{-2}\oplus E_{-1}$ and corresponding quadruple $(\partial, \left[.,.\right]_2, \chi,\omega )$, with $\chi=0$, there exists a Nijenhuis transformation
of the form $S+\alpha$ with $\alpha$ a vector valued $2$-form of degree zero, such that the deformed bracket $[S+\alpha,l_1+l_2+l_3] $
is the direct sum of a strict $2$-algebra with a trivial $L_\infty$-algebra.
\end{prop}
\begin{proof}
We set $E_{-1}^{t}:={\rm Im} (\partial)$, $E_{-2}^s := {\rm Ker}(\partial) $ and we choose two subspaces $E_{-2}^t \subset E_{-2}$
and $ E_{-1}^{s} \subset E_{-1}$ such that the sums $E_{-2}^t\oplus E_{-2}^s=E_{-2}$ and $E_{-1}^t\oplus E_{-1}^s=E_{-1}$ are direct sums. Since $\chi=0$, by item (i) in Proposition \ref{prop:quadrubleLie2algebra}, the bracket $[ .,. ]_2$ vanishes on $ E_{-1}^{t}$ so that there exists a unique
skew-symmetric bilinear map $\alpha: E_{-1} \times E_{-1} \to E^{t}_{-2}$ such that
\begin{equation}
 \partial \alpha (X,Y) = -pr_{E_{-1}^{t}} ([X,Y]_2),  \hbox{ for all}\,\,\, X,Y \in E_{-1},
 \end{equation}
where $pr_{E_{-1}^{t}} $ stands for the projection on ${E_{-1}^{t}}$ with respect to $E_{-1}^{s}$.
Note that $\alpha(X,Y)=0$ if $X$ or $Y$ belong to $ E_{-1}^t$, therefore we have $\alpha(\partial \alpha (X,Y),Z)=0$, for all $X,Y,Z \in E_{-1}$. Hence by Theorem \ref{thm:alphalpha}
$S+\alpha$
is Nijenhuis of square $S+2\alpha $. We claim that, for the deformed bracket $l_1'+l_2'+l_3':=[S+\alpha, l_1+l_2+l_3]_{_{RN}} $,
$(E_{-1}^s \oplus E_{-2}^s,l_1'+l_2'+l_3')$ and $(E_{-1}^t \oplus E_{-2}^t ,l_1'+l_2'+l_3')$ are  Lie sub-$2$-algebras of $(E=E_{-2}\oplus E_{-1}, \mu=l_1+l_2+l_3)$, that the first one is a string Lie algebra
while the second one is a trivial Lie algebra, and that their direct sum is isomorphic to $(E_{-2}\oplus E_{-1},l_1'+l_2'+l_3')$.

Let $(\partial^\prime,\left[.,.\right]^\prime,\chi^\prime,\omega^\prime)$ stands for the corresponding quadruple associated to the deformed structure $l_1'+l_2'+l_3'$. From $[l_1,l_2]_{_{RN}}=0$ we get $l_2(l_1f,X)=0$, for all $f\in E_{-2}$. This means that $l_2$ vanishes on $E_{-1}^t$. Also, since  $\alpha(X,Y)=0$ if $X$ or $Y$ belongs to $ E_{-1}^t$, by Equation (\ref{deformedLie2algebra}),  we have that $\chi^\prime=0$ and $\left[.,.\right]^\prime=0$ and hence $l_2^\prime$ vanishes on ${E_{-1}^{t}}$.
From $[l_1,l_3]_{_{RN}}=0$ we get $\omega(X,Y,Z)$ vanishes  for all $X \in {E_{-1}^{t}}$, so by Equation (\ref{deformedLie2algebra}) the restriction of $l_3' $ to $E_{-1}^t $ vanishes. Since the restriction of $l_1$ to $ E_{-2}^t$ is a bijection onto its image, the
restriction of $l_1'+l_2'+l_3'$ to $E_{-1}^t \oplus E_{-2}^t $ is a Lie $2$-sub-algebra and it is a trivial Lie $2$-algebra.

Next we prove that $(E_{-2}^s\oplus E_{-1}^s, l'_1+l'_2+l'_3)$ is a Lie sub-$2$-algebra with $l_1'(E_{-2}^s)=0$ and hence is a string Lie algebra. Let $X,Y \in E_{-1}^s$. Then by Equation (\ref{deformedLie2algebra}) we have
\begin{equation*}
l_2'(X,Y)=[X,Y]_2+\partial\alpha(X,Y)=[X,Y]_2-pr_{E_{-1}^{t}} ([X,Y]_2).
\end{equation*}
This implies that
\begin{equation}\label{strictpartofl-2}
l_2'(X,Y)\in E_{-1}^s.
\end{equation}
Let $X,Y,Z\in E_{-1}^s$. Then, we have $l_1'(X)=l_1'(Y)=l_1'(Z)=0$. Hence, from $(2[l_1',l_3']_{_{RN}}+[l_2',l_2']_{_{RN}})(X,Y,Z)=0$ we get
\begin{equation}\label{strictpartofl-3}
l_1'l_3'(X,Y,Z)=l_2'(l_2'(X,Y),Z).
\end{equation}
Using Relation (\ref{strictpartofl-2}), the right hand side of Equation (\ref{strictpartofl-3}) belongs to $E_{-1}^s$, while according to the definition of $E_{-1}^t$, the left hand side of Equation (\ref{strictpartofl-3}) belongs to $E_{-1}^t$ and since $E_{-1}=E_{-1}^t\oplus E_{-1}^s$ is a direct sum, both sides of Equation (\ref{strictpartofl-3}) should be zero. This implies that
\begin{equation}\label{strictpartlast}
l_3'(X,Y,Z)\in E_{-2}^s.
\end{equation}
 Relation (\ref{strictpartofl-2}) and Equation (\ref{strictpartlast}) show that $(E_{-2}^s\oplus E_{-1}^s, l'_1+l'_2+l'_3)$ is a Lie $2$-sub-algebra. Also, By definition of $E_{-2}^s$, we have $l_1'(E_{-2}^s)=0$. This completes the proof.
%
\end{proof}

Next, it is interesting to see that Lie algebras themselves can be seen as Nijenhuis forms.

 We start by noticing that any vector valued $2$-form of degree zero on a graded vector space $E_{-2}\oplus E_{-1}$ is of the form
\begin{equation}\label{khouk}
\alpha(X,Y)=\begin{cases}-\alpha(Y,X),&\mbox{if}\,\, X,Y \in E_{-1},\\
                          0,         & otherwise.
            \end{cases}
\end{equation}
This together with the fact that $\alpha$ always takes value in $E_{-2}$, imply that
\begin{equation}\label{sag}
\alpha(\alpha(X,Y),Z)+c.p.=0.
\end{equation}
Equations (\ref{khouk}) and (\ref{sag}) mean that \emph{any symmetric vector valued $2$-form $\alpha$ on an arbitrary graded vector space $E_{-2}\oplus E_{-1}$ is a Lie algebra} (not a graded Lie algebra).
 However, there is also a way to get a Lie bracket on a graded vector space $E=E_{-2}\oplus E_{-1}$ from a Nijenhuis form with respect to a Lie $2$-algebra structure $\mu=l_1+l_2+l_3$ on the vector space $E$ as follows:
 \begin{prop}
 Let $(E=E_{-2}\oplus E_{-1},\, \mu=l_1+l_2+l_3)$ be a Lie $2$-algebra, with corresponding quadruple $(\partial,\left[.,.\right],\chi,\omega)$. Let $\alpha$ be a vector valued $2$-form of degree zero. Define a bilinear map $\tilde{\alpha}$ by
 \begin{equation}
 \tilde{\alpha}(X,Y)=\begin{cases}\alpha(X,Y) & X,Y \in E_{-1},\\
                                  \alpha(\partial X,Y) & X\in E_{-2},Y \in E_{-1},\\
                                  \alpha(X,\partial Y) & X \in E_{-1},Y\in E_{-2},\\
                                  \alpha(\partial X,\partial Y) & X,Y \in E_{-2}.
                     \end{cases}
  \end{equation}
  Then $S+\alpha$ is Nijenhuis vector valued $2$-form with respect to $\mu$, with square $S+2\alpha$ if, and only if, $(E,\tilde{\alpha})$ is a Lie algebra.
  \end{prop}
  \begin{proof}
  By definition, $\tilde{\alpha}$ is a skew-symmetric bilinear map on the vector space $E$ and we have
  \begin{equation}\label{jacobialpha}\begin{array}{rcl}
  \tilde{\alpha}(\tilde{\alpha}(X,Y),Z)+c.p.&=&\alpha(\partial\alpha(X,Y),Z)+c.p.,\\
  \tilde{\alpha}(\tilde{\alpha}(f,Y),Z)+c.p.&=&\alpha(\partial\alpha(\partial f,Y),Z)+c.p.,\\
  \end{array}
  \end{equation}
  for all $X,Y,Z \in E_{-1}$ and $f\in E_{-2}$. Hence, Theorem \ref{thm:alphalpha} together with (\ref{jacobialpha}) imply that $\tilde{\alpha}$ is a Lie bracket on the vector space $E$ if, and only if, $S+\alpha$ is a Nijenhuis form with respect to $\mu$, with square $S+2\alpha$.
 \end{proof}
 Last, we give a result involving weak Nijenhuis forms on a Lie $2$-algebra.
 \begin{prop}
  Let $\partial: E_{-2}\to E_{-1}$ be a Lie $2$-algebra structure on the graded vector space $E=E_{-2}\oplus E_{-1}$, that is, a Lie $2$-algebra structure $\mu=l_1+l_2+l_3$ with $l_1=\partial$ and $l_2=l_3=0$, on $E$. Let $\alpha$ be a symmetric vector valued $2$-form of degree zero on the graded vector space $E$. If $S+\alpha$ is a weak Nijenhuis vector valued form with respect to $\partial$, then $E_{-1}$ is a Lie algebra with a representation on $E_{-2}$.
  \end{prop}
  \begin{proof}
  $S+\alpha$ is a weak Nijenhuis vector valued form with respect to $\partial$ if, and only if, $[S+\alpha,\partial]_{_{RN}}$ is an $L_{\infty}$-structure on the graded vector space $E$ if, and only if,
 \begin{equation*}
 [[S+\alpha,\partial]_{_{RN}},[S+\alpha,\partial]_{_{RN}}]_{_{RN}}=0
 \end{equation*}
 if, and only if,
 \begin{equation}\label{partial}
 [[\alpha,\partial]_{_{RN}},[\alpha,\partial]_{_{RN}}]_{_{RN}}=0.
 \end{equation}
 Therefore, $S+\alpha$ is a weak Nijenhuis vector valued form with respect to $\partial$ if, and only if,
 \begin{equation}\label{eq:weakexample1}
 \partial\alpha(\partial\alpha(X,Y),Z)+c.p.(X,Y,Z)=0
 \end{equation}
 and
 \begin{equation}\label{eq:weakexample2}
 \alpha(\partial\alpha(X,Y),\partial f)+c.p.(X,Y,\partial f)=0,
 \end{equation}
 for all $X,Y,Z \in E_{-1}, f\in E_{-2}$. Equation (\ref{eq:weakexample1}) means that $[X,Y]:=\partial\alpha(X,Y)$ defines a Lie bracket on  $E_{-1}$ since clearly it is skew-symmetric. While if we denote $ X\cdot f:=\alpha(X,\partial f)$, then (\ref{eq:weakexample2}) can be rewritten as
  \begin{equation*}
  [X,,Y]\cdot f= X\cdot(Y\cdot f)-Y\cdot(X\cdot f),
  \end{equation*}
  which means that $\cdot:E_{-1}\times E_{-2}\to E_{-2}$ is a representation of $E_{-1}$.
 \end{proof}

\begin{rem}
A notion of Nijenhuis operator on a Lie $2$-algebra independently appeared in \cite{LSZ} while the present manuscript was about to be completed.
This notion does not match our notion, although there are some similarities, and it is certainly also interesting. First, in \cite{LSZ}, a Nijenhuis operator $N$ is necessarily a vector valued $1$-form. Second, the deformed structure that they consider matches what we also called the deformed structure. In fact, their idea is to require
that $N$ itself has to be a Lie $2$-algebra morphism from the deformed structure to the original one, generalizing a property that holds true for Lie algebra and it is also based on a quite natural generalization of the usual Nijenhuis operator.

 If $\mathcal N=(N_0,N_1)$ is a Nijenhuis operator, in the sense of Definition 3.2. in \cite{LSZ},  then the following conditions hold:
\begin{equation*}
\begin{array}{rcl}
\left[\mathcal N, \left[\mathcal N, l_1\right]_{_{RN}} \right]_{_{RN}}&=& \left[\mathcal{N}^2, l_1\right]_{_{RN}}\\
\left[\mathcal N, \left[\mathcal N, l_2\right]_{_{RN}} \right]_{_{RN}}&=& \left[\mathcal{N}^2, l_2\right]_{_{RN}}\\
\left[\mathcal N, \left[\mathcal N, l_3\right]_{_{RN}} \right]_{_{RN}}&=& \left[\mathcal{N}^2, l_3\right]_{_{RN}},
\end{array}
\end{equation*}
which means that $\mathcal N$ is a Nijenhuis vector valued form, in our sense, with square $\mathcal N  ^ 2$.
\end{rem}

\section{Courant structures}\label{subsection:3.3.5}

According to Roytenberg and Weinstein \cite{LWX}, a Lie $2$-algebra can be associated to an arbitrary Courant algebroid,
Lie $2$-algebra that encodes entirely the Courant structure in question. Now, various authors, \cite{CGM, GrB, YKSBrazil, CJP, AJC} have investigated Nijenhuis
tensors on Courant structures. We shall check that
such Nijenhuis tensors give examples of Nijenhuis vector valued forms on the corresponding Lie $2$-algebra. But these Nijenhuis tensors are always $1$-forms,
while, for degree reasons, the corresponding Lie $2$-algebra may admit Nijenhuis vector valued forms which are sum of vector valued $1$-forms
with  vector valued $2$-forms. Surprisingly, the Lie $2$-algebra structures deformed by those will never be Lie $2$-algebras constructed out of  Courant structures,
except in some degenerated cases, and except, of course, if the Nijenhuis tensor is of the type studied in \cite{CGM, GrB, YKSBrazil, CJP, AJC}.
\subsection{Review of definitions}
Let us recall some basic definitions. First, we recall the original definition of Courant algebroids.
\begin{defi}\label{def:Courantskew}\cite{LWX}
A Courant algebroid is a vector bundle $E\to M$ equipped  with a non-degenerate inner product $\langle . , . \rangle$, a skew-symmetric bracket $\left[.,.\right]:\Gamma(E) \times \Gamma(E)\to \Gamma(E)$, and a bundle map $\rho : E\to TM$ such that the following properties are satisfied:
\begin{enumerate}
  \item For any $X,Y,Z \in \Gamma(E), J(X,Y,Z)=\mathcal{D}T(X,Y,Z)$;
  \item for any $X,Y \in \Gamma(E), \rho[X,Y]=[\rho X,\rho Y]$;
  \item for any $X,Y \in \Gamma(E)$ and $f\in \mathcal{C}^{\infty}(M)$,\\  \begin{equation*}[X,fY]=f[X,Y]+(\rho(X)f)Y-\frac{1}{2}\langle X,Y\rangle \mathcal{D}f;\end{equation*}
  \item $\rho \circ \mathcal{D}=0$, i.e., for any $f,g \in \mathcal{C}^{\infty}(M),\langle \mathcal{D}f,\mathcal{D}g\rangle =0;$
  \item for any $X,Y_1,Y_2 \in \Gamma(E)$,\\ $$\rho(X)\langle Y_1,Y_2\rangle =\langle [X,Y_1]+\frac{1}{2}\mathcal{D}\langle X,Y_1\rangle ,Y_2\rangle +\langle Y_1,[X,Y_2]+\frac{1}{2}\mathcal{D}\langle X,Y_2\rangle \rangle ,$$
\end{enumerate}
    where $J$ is the Jacobiator of the bracket $\left[.,.\right]$ i.e.
\begin{equation*}
J(X,Y,Z)=[[X,Y],Z]+c.p.,
\end{equation*}
 $T$ is the function on the base manifold  $M$ defined by:
\begin{equation}\label{T(X,Y,Z)}
T(X,Y,Z)=\frac{1}{6}\langle [X,Y],Z\rangle +c.p.,
\end{equation}
and $\mathcal{D}:\mathcal{C}^{\infty}(M) \to \Gamma(E)$ is the map defined by
\begin{equation}\label{matcallD}
\langle \mathcal{D}f,X\rangle =\rho (X)f,
\end{equation}
for all $X,Y,Z \in \Gamma(E),f \in \mathcal{C}^{\infty}(M)$. We will denote a Courant algebroid by $(E,\left[.,.\right],\rho, \langle .,.\rangle )$.
\end{defi}

 It is also common to define a Courant algebroid in terms of a non-skew-symmetric bracket which is called \emph{Dorfman bracket}. We also recall this definition.
\begin{defi}\label{def:CourantDorfman}\cite{RoytenbergPh.D.}
  A Courant algebroid consists of a vector bundle $E\to M$, a bilinear map $\circ:\Gamma(E)\times \Gamma(E)\to \Gamma(E)$, a bundle map $\rho:E\to TM$ and a non-degenerate inner product satisfying the following axioms:
  \begin{enumerate}
    \item $X\circ (Y\circ Z)=(X\circ Y)\circ Z+Y\circ (X\circ Z)$ ,
    \item $\rho(X\circ Y)=[\rho (X),\rho (Y)]$ ,
    \item $X\circ fY= f(X\circ Y)+(\rho(X)f) Y$ ,
    \item $X \circ X= \frac{1}{2}\mathcal{D}\langle X,X\rangle $,
    \item $\rho(X)\langle Y,Z\rangle =\langle X \circ Y, Z\rangle +\langle Y,X \circ Z\rangle $;
  \end{enumerate}
  for all $X,Y,Z \in \Gamma(E),f \in \mathcal{C}^{\infty}(M)$, where $\mathcal{D}:\mathcal{C}^{\infty}(M)\to \Gamma(E) $ is given by (\ref{matcallD}).
 \end{defi}
\begin{rem}\label{remarkondefinitionsofCourant}
  It is shown in \cite{RoytenbergPh.D.} that Definitions \ref{def:Courantskew} and \ref{def:CourantDorfman} are equivalent. The relation between brackets is given by $[X,Y]=\frac{1}{2}(X \circ Y - Y\circ X)$ and $X\circ Y=[X,Y]+\frac{1}{2}\mathcal{D}\langle X,Y\rangle $. Uchino has shown in \cite{Uchino} that some of the axioms in  Definitions \ref{def:Courantskew} and \ref{def:CourantDorfman} follow from the other ones. Later, in the same vein Grabowski and Marmo \cite{Gra} have simplified Definition \ref{def:CourantDorfman}, while Kosmann-Schwarzbach \cite{YKSquasi} obtained the following definition of Courant algebroid.
  \end{rem}
  \begin{defi}\label{def:Courantdorfmanshort}
  A Courant algebroid is a vector bundle $E\to M$ together with a non-degenerate inner product  $\langle .,. \rangle$, a bundle map $\rho:E\to TM$ and a bilinear operator $\circ: \Gamma(E) \times \Gamma(E) \to \Gamma(E)$, such that the following axioms hold:
  \begin{enumerate}
    \item  $(\Gamma(E),\circ)$ is a Leibniz algebra i.e $X\circ (Y\circ Z)=(X\circ Y)\circ Z+Y\circ (X\circ Z) $ ,
    \item $\rho(X)\langle Y,Z\rangle = \langle X \circ Y, Z\rangle +\langle Y,X\circ Z\rangle $,
    \item $\rho(X)\langle Y,Z\rangle = \langle X,Y \circ Z\rangle +\langle X,Z\circ Y\rangle $;
  \end{enumerate}
  for all $ X,Y,Z \in \Gamma(E)$.
  \end{defi}
  When item $(1)$ does not hold, then the quadruple $(E, \circ, \rho, \langle . , . \rangle )$ is called a \emph{pre-Courant algebroid}(\cite{CJP}).

  The next proposition is stated in \cite{YKSquasi}, for Courant algebroids. Since the proof does not use the fact of $\circ$ being a Leibniz bracket, the result also holds for pre-Courant algebroids.
  \begin{prop}\label{prop:fmiadbiroonupto}
  For every pre-Courant algebroid $(E, \circ, \rho, \langle .,. \rangle)$ we have
  \begin{equation}\label{LeibnizforpreCourantholds}
  X\circ f Y=f (X\circ Y) + (\rho(X)f)Y,
  \end{equation}
  for all $X,Y \in \Gamma(E), f\in \mathcal{C}^{\infty}(M).$
  \end{prop}
  \begin{proof}
  Let $X,Y \in \Gamma(E)$ and $f\in \mathcal{C}^{\infty}(M).$ By the second item in Definition \ref{def:Courantdorfmanshort} we have
  \begin{equation*}
  \rho(X)\langle fY,Z\rangle = \langle X \circ fY, Z\rangle +\langle fY,X\circ Z\rangle,
  \end{equation*}
   hence,
  \begin{equation*}
  f\rho(X)\langle Y,Z\rangle+\langle Y,Z\rangle\rho(X)f=\langle X \circ fY, Z\rangle +f\langle Y,X\circ Z\rangle.
  \end{equation*}
  or
   \begin{equation*}
  f\left(\langle X \circ Y, Z\rangle +\langle Y,X\circ Z\rangle \right)+\langle Y,Z\rangle\rho(X)f=\langle X \circ fY, Z\rangle +f\langle Y,X\circ Z\rangle.
  \end{equation*}
  And using the non-degeneracy of the inner product $\langle .,.\rangle$, the desired result follows after simplification.
  \end{proof}
  Note that for Courant algebroids, Equation (\ref{LeibnizforpreCourantholds}) is exactly the third item of Definition \ref{def:CourantDorfman}.
  \begin{cor}\label{cor:Anchoruniqe}
  Let $(E, \circ, \rho, \langle .,. \rangle)$ and $(E, \circ', \rho', \langle .,. \rangle)$ be two pre-Courant algebroids. If $\circ =\circ'$, then $\rho=\rho'$.
  \end{cor}
  \begin{proof}
  Assume that $(E, \circ, \rho, \langle .,. \rangle)$ and $(E, \circ, \rho', \langle .,. \rangle)$ are both pre-Courant algebroids. By Proposition \ref{prop:fmiadbiroonupto} we have
  \begin{equation*}
   (\rho(X)f)Y=(\rho'(X)f)Y,
  \end{equation*}
  for all $X,Y\in \Gamma(M), f\in \mathcal{C}^{\infty}(M)$, which implies that $\rho=\rho'$.
  \end{proof}

\subsection{Nijenhuis tensors on Courant algebroids}
   We intend to define Nijenhuis deformation of Courant structures. Let $(E,\circ, \rho,\langle .,.\rangle)$ be a Courant algebroid. For a given endomorphism $N:\Gamma(E) \to \Gamma(E)$, the deformed Dorfman bracket by $N$ is a bilinear operation defined as:
  $$X\circ^{N}Y:=NX\circ Y+X\circ NY -N(X\circ Y),$$
  for all $X,Y \in \Gamma(E)$. The deformation of $\rho$ by $N$ is the map given by $\rho^N(X)=\rho(NX)$, $X\in \Gamma(E)$. The Nijenhuis torsion of $N$, with respect to the Dorfman bracket $\circ$, is defined as:
  \begin{equation*}
  T_{\circ}N(X,Y):=NX\circ NY-N(X\circ^{N}Y),
  \end{equation*}
  for all $X,Y \in \Gamma(E)$.
  A direct computation shows that
  $$T_{\circ}N=\frac{1}{2}(\circ^{N,N}-\circ^{N^{2}}).$$
   All endomorphisms $N$ of $\Gamma(E)$ that will be considered here are $\mathcal{C}^{\infty}(M)$-linear, that is to say they are $(1,1)$-tensors, that is smooth sections of endomorphisms of $E$.

 According to \cite{CGM}, for every vector bundle $E\to M$, if $(\Gamma(E),\circ)$ is a Leibniz algebra and $N: \Gamma(E) \to \Gamma(E)$ is any endomorphism  whose Nijenhuis torsion vanishes, then the pair $(\Gamma(E),\circ^N)$ is a Leibniz algebra.
 The difficulty, is that for a given Courant algebroid $(E, \circ, \rho, \langle .,.\rangle )$ and a given
 $(1,1)$-tensor $N$, $(E, \circ^N, \rho^N, \langle .,.\rangle )$ may fail to be a pre-Courant algebroid, even if the Nijenhuis torsion of $N$ vanishes.
 Indeed,  from \cite{CGM} we have the following:

 \begin{them}\label{them:Grobowski}
 If $N$ is a $(1,1)$-tensor on a pre-Courant algebroid $(E, \circ, \rho, \langle .,.\rangle )$, then the quadruple $(E , \circ^N, \rho^N, \langle .,.\rangle )$ is a pre-Courant algebroid if, and only if,
 \begin{equation}
 X \circ (N+N^*)Y=(N+N^*)(X \circ Y) \hbox{ and } (N+N^*)(Y \circ Y)=((N+N^*) Y) \circ Y
 \end{equation}
 for all $X,Y \in \Gamma (E)$, where $N^*$ stands for the transpose of $N$, with respect to $\langle.,.\rangle.$
 \end{them}
 \begin{rem}
 In fact, Theorem \ref{them:Grobowski} is slightly different from Theorem $4$ in \cite{CGM}, because there, the authors start from a Courant
 algebroid. But the same proof is still valid for the case of pre-Courant.
 \end{rem}
 A \emph{Casimir function} or simply a Casimir on a Courant algebroid $(E, \circ, \rho, \langle . , .\rangle )$ is a function $f\in \mathcal{C}^{\infty}(M)$ such that $\rho(X)f=0 $ for all $X\in \Gamma(E)$. It is easy to check that $f$ is a Casimir if, and only if, $\mathcal{D}f=0$. Also, if $f$ is a Casimir, then condition $3$ in Definition  \ref{def:CourantDorfman} implies that $f (X \circ Y)=X \circ (fY)$, for all $X,Y \in \Gamma(E)$. Moreover, $f$ being Casimir, a simple computation gives $(f X) \circ Y = f (X \circ Y)$, for all $X,Y \in \Gamma(E)$, so we have that for a Casimir functions $f$,
  \begin{equation}\label{eq:CasimirGetsOut}(f X) \circ Y = f (X \circ Y) = X \circ (fY) \end{equation}
 for every pair of sections $X,Y \in \Gamma(E)$.

 The next lemma is a slight generalization of a result in \cite{CGM}.
  \begin{lem}\label{pre}
  Given a pre-Courant algebroid $(E, \circ, \rho, \langle .,.\rangle )$ and a map $N:\Gamma(E) \to \Gamma(E)$,
  if $N+N^*=\lambda Id_{\Gamma(E)}$, for some Casimir function $\lambda \in \mathcal{C}^{\infty}(M)$, then $(E, \circ^N, \rho^N, \langle .,.\rangle )$ is a pre-Courant algebroid.
  \end{lem}
  \begin{proof}
  This lemma is a direct consequence of Theorem \ref{them:Grobowski} together with (\ref{eq:CasimirGetsOut}).
  \end{proof}
  \begin{them}\label{Theml^N}
Let $ (E , \circ, \rho, \langle .,.\rangle ) $ be a Courant algebroid and $N:\Gamma(E)\to \Gamma(E)$  be a $(1-1)$-tensor whose Nijenhuis torsion vanishes, such that
 \begin{equation}
    N+N^*=\lambda Id_{\Gamma(E)},
 \end{equation}
  with $\lambda$ being a Casimir function. Then $ (E, \circ^N, \rho^N, \langle . ,. \rangle ) $ is a Courant algebroid.
  \end{them}
\begin{proof}
 Note that $(E,\circ)$ is a Leibniz algebra which implies that $(E,\circ^N)$ is also a Leibniz algebra since the Nijenhuis torsion of $N$ vanishes. This together with Lemma \ref{pre} prove the theorem.
 \end{proof}
   \begin{rem}
  For a (pre-)Courant algebroid $(E,\circ,\rho,\langle .,.\rangle)$, and a $(1-1)$-tensor $N$ with $N+N^*=\lambda Id_{\Gamma(E)}$ and $\lambda$ a Casimir function, we have
  \begin{equation*}
  \rho^{N}(X)f=\rho(NX)f=\langle NX,\mathcal D f\rangle=\langle X,N^*\mathcal D f\rangle=\langle X,(-N+\lambda Id_{\Gamma(E)})\mathcal D f\rangle,
  \end{equation*}
for all $X\in \Gamma(E), f\in \mathcal{C}^{\infty}(M)$.  This means that the
  operator $ \mathcal{D}^N : \mathcal{C}^\infty(M) \to \Gamma(E)$ associated with the (pre-)Courant algebroid $(E,\circ^N,\rho^N,\langle .,.\rangle)$, is given by
  \begin{equation} \label{deformedD}
  \mathcal {D} ^N=(-N+\lambda Id_{\Gamma(E)})\circ\mathcal D.
  \end{equation}
  \end{rem}
\subsection{Courant algebroids as Lie $2$-algebras and Nijenhuis forms}
  For a given vector bundle $E\to M$, equipped with a skew-symmetric bilinear map $\left[.,.\right]: \Gamma(E)\times \Gamma(E)\to \Gamma(E)$, a bundle map $\rho: \Gamma(E) \to TM$ over the identity of $M$ and a symmetric bilinear form $\langle .,.\rangle $ define
  \begin{enumerate}
  \item A graded vector space as $V=\mathcal{C}^{\infty}(M)\oplus \Gamma(E)$, where the elements of $\mathcal{C}^{\infty}(M)$ have degree $-2$ and the elements of $\Gamma(E)$ have degree $-1$.
  \item  Symmetric vector valued forms  on $V$ by:
  \end{enumerate}
  \begin{equation}\label{Lie2fromCourant}
  \begin{array}{rcl}
  l_1f &=& \mathcal{D} f \,\,\,\,\mbox{for all} \,\,\,\, f\in \mathcal{C}^{\infty}(M),\\
  l_2(X,Y)&=& [X,Y]\,\,\,\,\mbox{for all} \,\,\,\, X,Y\in \Gamma(E),\\
  l_2(X,f)&=& \frac{1}{2}\langle X,\mathcal{D}f\rangle \,\,\,\,\mbox{for all} \,\,\,\, X\in \Gamma(E),f\in \mathcal{C}^{\infty}(M),\\
  l_3(X,Y,Z)&=&T(X,Y,Z)\,\,\,\,\mbox{for all} \,\,\,\, X,Y,Z\in \Gamma(E),\\
  \end{array}
\end{equation}
where $ \mathcal{D} : \mathcal{C}^\infty(M) \to \Gamma(E)$ is the operator associated to the quadruple $(E,\left[.,.\right],\rho,\langle .,.\rangle)$ defined as in (\ref{matcallD}), T is given by (\ref{T(X,Y,Z)})
and
$$l_1=0\,\,\mbox{on degree} -1, l_2=0 \,\,\mbox{on degree} -4 \,\,\mbox{ and } l_3=0 \,\,\mbox{on degree} < -3.$$
By construction, $l_1,l_2$ and $l_3$ are of degree $+1$ on the graded vector space $V$. Given a pre-Courant algebroid $(E,\left[.,.\right],\rho,\langle .,.\rangle )$ and the symmetric vector valued forms $l_1, l_2$ and $l_3$ on $V=\mathcal{C}^{\infty}(M)\oplus \Gamma(E)$ constructed as above, the pair $(V=\mathcal{C}^{\infty}(M)\oplus \Gamma(E),\mu=l_1+l_2+l_3 )$ is called \emph{pre-Lie $2$-algebra} associated to the pre-Courant algebroid $(E,\left[.,.\right],\rho,\langle .,.\rangle )$.
It makes sense therefore to ask when $\mu=l_1+l_2+l_3$ is a  Lie $2$-algebra structure on the graded vector space $V$.
An answer can be found in \cite{RoytenbergWeinstein}:

\begin{prop}\cite{RoytenbergWeinstein}\label{Roytenbergassociation}
If $(E,\left[.,.\right],\rho,\langle .,.\rangle )$ is Courant algebroid, then the pair $(V,l_1+l_2+l_3)$, constructed in above, is a symmetric Lie $2$-algebra.
\end{prop}

We call this symmetric  Lie $2$-algebra the symmetric Lie $2$-algebra \emph{associated to the Courant algebroid $(E,\left[.,.\right],\rho,\langle .,.\rangle )$} or simply Lie $2$-algebra associated to the Courant algebroid $(E,\left[.,.\right],\rho,\langle .,.\rangle )$, when there is no risk of confusion.

\begin{rem}
In fact, the previous proposition is not exactly the result in \cite{RoytenbergWeinstein},
   who define Lie 2-algebras with graded skew-symmetric brackets,
   but it suffices to shift the degrees by $1$ and to apply the d\'ecalage isomorphism to the original Lie $2$-algebra that appears in \cite{RoytenbergWeinstein} to get the one introduced here.
\end{rem}
Using (\ref{Lie2fromCourant}) and Remark \ref{remarkondefinitionsofCourant} to pass from the skew-symmetric bracket to the Dorfman bracket, the binary bracket $l_2$ and $3$-ary bracket $l_3$ in (\ref{Lie2fromCourant}) for a quadruple $(E,\circ,\rho,\langle .,.\rangle)$ are given by
   \begin{equation}\label{Lie2fromDorfman}
   \begin{array}{rcl}
    l_2(X,Y)&=& \frac{1}{2}(X\circ Y-Y\circ X)\,\,\,\,\mbox{for all} \,\,\,\, X,Y\in \Gamma(E),\\
    l_3(X,Y,Z)&=&\frac{1}{12}\langle (X\circ Y-Y\circ X),Z\rangle +c.p.\,\,\,\,\mbox{for all} \,\,\,\, X,Y,Z\in \Gamma(E).
   \end{array}
   \end{equation}
Starting with a Nijenhuis tensor for a Courant algebroid we construct, in the next proposition, a Nijenhuis form for the Lie $2$-algebra associated to the Courant algebroid, according to Proposition \ref{Roytenbergassociation}. First we need the following lemma.
\begin{lem}\label{leml^N}
Let $ (E , \circ, \rho, \langle .,.\rangle ) $ be a pre-Courant algebroid with the associated symmetric pre-Lie-$2$  algebra structure $\mu=l_1+l_2+l_3$, on the graded vector space $V=\mathcal{C}^{\infty}(M)\oplus \Gamma(E)$. Let $N:\Gamma(E)\to \Gamma(E)$  be a $(1-1)$-tensor such that
 \begin{equation}
    N+N^*=\lambda Id_{\Gamma(E)},
\end{equation}
  with $\lambda$ being a Casimir function. Then, the pre-Lie $2$-algebra structure associated to the pre-Courant algebroid $ (E, \circ^N, \rho^N, \langle . ,. \rangle ) $ is $\left[{\mathcal N},l_1+l_2+l_3\right]_{_{RN}}$, with ${\mathcal N}$  defined as
  \begin{equation*}
{\mathcal N}|_{\Gamma(E)}= N \,\,\,\mbox{and}\,\,\,{\mathcal N}|_{\mathcal{C}^{\infty}(M)}=\lambda Id_{\mathcal{C}^{\infty}(M)}.
\end{equation*}

\end{lem}
\begin{proof}
 Let us denote the pre-Lie $2$-algebra associated to the pre-Courant algebroid $ (E , \circ^N, \rho^N, \langle . ,. \rangle ) $ by $l_1^{N}+l_2^{N}+l_3^{N}$. By (\ref{Lie2fromCourant}), (\ref{Lie2fromDorfman}) and (\ref{deformedD}), we have, for all $f\in\mathcal{C}^{\infty}(M) $ and for all $X,Y,Z \in \Gamma(E)$:
\begin{equation}\label{l1N}
l_1^{N}f=\mathcal{D}^Nf=\lambda \mathcal{D}f-N\mathcal{D}f=l_1(\mathcal {N}f)-\mathcal{N}l_1(f)= [\mathcal{N},l_1]_{_{RN}}(f),
\end{equation}
\begin{equation}\label{l2N}
\begin{array}{rcl}
l_2^{N}(X,Y)&=&\frac{1}{2}(X\circ^N Y-Y\circ^N X)\\
            &=&\frac{1}{2}(NX\circ Y-Y\circ NX+X\circ NY-NY\circ X-N(X\circ Y-Y\circ X))\\
            &=&l_2(NX,Y)+l_2(X,NY)-Nl_2(X,Y)\\
            &=&[\mathcal{N},l_2]_{_{RN}}(X,Y),
\end{array}
\end{equation}
\begin{equation}\label{l22N}
\begin{array}{rcl}
l_2^{N}(X,f)&=&\frac{1}{2}\langle X, \mathcal{D}^N f\rangle\\
            &=&\frac{1}{2}\langle X,(- N+\lambda I)\mathcal{D} f\rangle\\
            &=&\frac{1}{2}\langle X,N^* \mathcal{D} f\rangle\\
            &=&\frac{1}{2}\langle NX, \mathcal{D} f\rangle\\
            &=&l_2(NX,f)+\lambda l_2(X,f)-\lambda l_2(X,f)\\
            &=&l_2(\mathcal{N}X,f)+l_2(X,\mathcal{N}f)-\mathcal{N}l_2(X,f)\\
            &=&[\mathcal{N},l_2]_{_{RN}}(X,f)
\end{array}
\end{equation}
and
\begin{equation}\label{l3N}
\begin{array}{rcl}
l_3^{N}(X,Y,Z)&=&\frac{1}{6}\langle l_2^N(X,Y), Z\rangle+c.p.(X,Y,Z)\\
            &=&\frac{1}{6}\langle l_2(NX,Y)+l_2(X,NY)-Nl_2(X,Y), z\rangle+c.p.(X,Y,Z)\\
            &=&\frac{1}{6}\langle l_2(NX,Y)+l_2(X,NY)+(N^*-\lambda I)l_2(X,Y), Z\rangle+c.p.(X,Y,Z)\\
            &=&\frac{1}{6}(\langle l_2(NX,Y),Z\rangle+\langle l_2(X,NY), Z\rangle+\langle l_2(X,Y), NZ\rangle-\lambda\langle l_2(X,Y), Z\rangle)\\
            &&+c.p.(X,Y,Z)\\
            &=&\frac{1}{6}(\langle l_2(NX,Y),Z\rangle +c.p.(NX,Y,Z)\\
            &+&\langle l_2(X,NY), Z\rangle+c.p.(X,NY,Z)\\
            &+&\langle l_2(X,Y),N Z\rangle+c.p.(X,Y,NZ)\\
            &-&\lambda \langle l_2(X,Y), Z\rangle+c.p.(X,Y,Z))\\
            &=&l_3(\mathcal N X,Y,Z)+l_3(X,\mathcal N Y,Z)+l_3(X,Y,\mathcal N Z)-\mathcal N l_3(X,Y,Z)\\
            &=&[\mathcal{N},l_3]_{_{RN}}(X,Y,Z),
\end{array}
\end{equation}
where we used $N+N^*=\lambda Id_{\Gamma(E)}$ in (\ref{l22N}) and (\ref{l3N}).
\end{proof}
For the case of a Courant algebroid we have the following result.
\begin{cor}\label{cornew}
Let $ (E , \circ, \rho, \langle .,.\rangle ) $ be a Courant algebroid with the associated symmetric Lie-$2$  algebra structure $\mu=l_1+l_2+l_3$, on the graded vector space $V=\mathcal{C}^{\infty}(M)\oplus \Gamma(E)$. Let $N:\Gamma(E)\to \Gamma(E)$  be a $(1-1)$-tensor such that
 \begin{equation}
 \begin{cases}
    N+N^*=\lambda Id_{\Gamma(E)},\\
    (\Gamma(E), \circ^N)\,\,\,\mbox{is a Leibniz algebra,}
    \end{cases}
\end{equation}
  with $\lambda$ being a Casimir function. Then, the Lie $2$-algebra structure associated to the Courant algebroid $ (E, \circ^N, \rho^N, \langle . ,. \rangle ) $ is $[{\mathcal N},l_1+l_2+l_3]$, with ${\mathcal N}$  defined as in (\ref{eq:defmathcalN}).
\end{cor}
\begin{prop}\label{prop:NijenhuisonCourant}
Let $ (E , \circ, \rho, \langle .,.\rangle ) $ be a Courant algebroid with the associated symmetric Lie-$2$  algebra structure $\mu=l_1+l_2+l_3$, on the graded vector space $V=\mathcal{C}^{\infty}(M)\oplus \Gamma(E)$. Let $N:\Gamma(E)\to \Gamma(E)$  be a  $(1-1)$-tensor whose Nijenhuis torsion with respect to the bracket $\circ$ vanishes and
satisfies the following conditions
 \begin{equation}\label{con:CourantNijenhuisconditions}
  \begin{cases}
  N+N^*=\lambda Id_{\Gamma(E)},\\
  N^2+(N^2)^*=\gamma Id_{\Gamma(E)}
  \end{cases}
  \end{equation}
  with $\lambda,\gamma$ being Casimir functions. Define   ${\mathcal N}$ and ${\mathcal K}$ as
\begin{equation}\label{eq:defmathcalN}
{\mathcal N}|_{\Gamma(E)}= N \,\,\,\mbox{and}\,\,\,{\mathcal N}|_{\mathcal{C}^{\infty}(M)}=\lambda Id_{\mathcal{C}^{\infty}(M)},
\end{equation}
\begin{equation}\label{eq:defmathcalK}
{\mathcal K}|_{\Gamma(E)}= N^2 = \lambda N + \frac{\gamma - \lambda^2}{2} Id_{\Gamma(E)} \,\,\,\mbox{and}\,\,\,{\mathcal K}|_{\mathcal{C}^{\infty}(M)}=\gamma Id_{\mathcal{C}^{\infty}(M)}.
\end{equation}
Then $\mathcal N$ is a Nijenhuis vector valued $1$-form with respect to $\mu$, with square $\mathcal K$.
\end{prop}
\begin{proof}
Applying Corollary \ref{cornew} for the Courant algebroid $ (E , \circ, \rho, \langle .,.\rangle ) $, the Nijenhuis $(1-1)$-tensor $N$ and the vector valued $1$-form $\mathcal{N}$, twice, we get
\begin{equation}\label{mah}
l_1^{N,N}+l_2^{N,N}+l_3^{N,N}=[\mathcal{N},[\mathcal{N},l_1+l_2+l_3]_{_{RN}}]_{_{RN}},
\end{equation}
where $l_1^{N,N}+l_2^{N,N}+l_3^{N,N}$ stands for the Lie $2$-algebra structure associated to the Courant algebroid $ (E , \circ^{N,N}, \rho^{N,N}, \langle .,.\rangle )$.
Applying again Corollary \ref{cornew} for the Courant algebroid $ (E , \circ, \rho, \langle .,.\rangle ) $, the $(1-1)$-tensor $N^2$ and the vector valued $1$-form $\mathcal{K}$, we get

\begin{equation}\label{khorshid}
l_1^{N^2}+l_2^{N^2}+l_3^{N^2}=[\mathcal{K},l_1+l_2+l_3]_{_{RN}},
\end{equation}
where $l_1^{N^2}+l_2^{N^2}+l_3^{N^2}$ stands for the Lie $2$-algebra structure associated to the Courant algebroid $ (E , \circ^{N^2}, \rho^{N^2}, \langle .,.\rangle )$.
On the other hand, since the Nijenhuis torsion of $N$ vanishes, the Courant algebroids $(E,\circ^{N,N},\rho^{N,N},\langle.,.\rangle)$ and $(E,\circ^{N^2},\rho^{N^2},\langle.,.\rangle)$ are the same. Therefore  (\ref{mah}) and (\ref{khorshid}) imply that
\begin{equation*}
[\mathcal{N},[\mathcal{N},l_1+l_2+l_3]_{_{RN}}]_{_{RN}}=[\mathcal{K},l_1+l_2+l_3]_{_{RN}}.
\end{equation*}
An easy computation shows that $[\mathcal{N},\mathcal{K}]_{_{RN}}$ vanishes both on functions and on sections of $E$.
\end{proof}

%
%
%
Since the Lie $2$-algebra structure entirely encodes the Courant algebroid structure,
there was a hope that we could, given a Courant structure, find a Nijenhuis deformation
by a Nijenhuis tensor which is the sum of a vector valued $1$-form and a vector valued $2$-form of the corresponding Lie $2$-algebra structure,
and prove, eventually, that the Lie $2$-algebra structure obtained by this procedure comes from a Courant structure.
But this fails, at least when the anchor is not identically zero,
for the following reason. First, notice that every  $C^\infty(M)$-linear vector valued form of degree $0$ on $E_{-2}\oplus E_{-1}$, where $E_{-2}:=\mathcal{C}^{\infty}(M)$ and $E_{-1}:=\Gamma(E)$,
is the sum of a $2$-form $\alpha$, a $(1-1)$-tensor $N$ and an endomorphism of $C^\infty(M)$
of the form $F \mapsto \lambda F $ for some smooth function $\lambda$,
hence we should denote them as a sum $\lambda+ N + \alpha $.
We will use Definition \ref{def:Courantskew} for Courant algebroid in the next theorem.
\begin{them} \label{them:NijenhuisCourant}
Let $ (\rho, \left[.,.\right], \langle . ,. \rangle ) $ be a Courant structure on a vector bundle $E\to M$ with the associated Lie $2$-algebra structure $l_1+l_2+l_3$ on the graded vector space $V=E_{-2}\oplus E_{-1}$, where $E_{-2}:=\mathcal{C}^{\infty}(M)$ and $E_{-1}:=\Gamma(E)$.
Let ${\mathcal N} = \lambda + N + \alpha$ be a $C^\infty(M)$-linear vector valued form of degree $0$ on $V$.
Assume also that $\rho$ is not equal to zero on a dense subset of the base manifold.
If $ [{\mathcal N}, l_1+l_2+l_3]_{_{RN}}$ is the Lie $2$-algebra associated to a Courant structure, with the same scalar product, then
\begin{enumerate}
\item $\lambda$ is a Casimir,
\item $\alpha=0,$
\item $ N+N^* = \lambda Id_{\Gamma(E)}.$
%
\end{enumerate}
In this case, the Courant structure that $ [{\mathcal N}, l_1+l_2+l_3]_{_{RN}}$ is associated with, is $(\left[.,.\right]_N, \rho^N,\langle . ,. \rangle).$
\end{them}
\begin{proof}
Let $\mu=l_1+l_2+l_3$ and let us denote the component of $i$-form in $[\mathcal{N},\mu]_{_{RN}}$ by $[\mathcal{N},\mu]^i_{RN}$, $i=1,2$. Then for all $X,Y \in \Gamma(E),f\in \mathcal{C}^{\infty}(M)$ we have
\begin{equation*}\label{equation1intheorem}
\begin{array}{rcl}
[\mathcal{N},\mu]^1_{RN}(f)&=&([\lambda,l_1]_{_{RN}}+[N,l_1]_{_{RN}})(f)\\
                           &=&l_1(\lambda f)-Nl_1(f)\\
                           &=&\lambda l_1( f)+fl_1(\lambda)-Nl_1(f).
\end{array}
\end{equation*}
   The first equation in (\ref{Lie2fromCourant}) implies that if $[\mathcal{N},\mu]_{_{RN}}$ is a Lie 2-algebra associated to a Courant algebroid, then $[\mathcal{N},\mu]^1_{RN}$ has to be a derivation, and this happens if and only if $l_1(\lambda)=0$, and so we get
\begin{equation}\label{equation1intheorem-1}
[\mathcal{N},\mu]^1_{RN}(f)=(\lambda Id_{\Gamma(E)}-N) l_1( f).
\end{equation}
On the other hand
\begin{equation}\label{equation2intheorem}
\begin{array}{rcl}
[\mathcal{N},\mu]_{_{RN}}^2(X,f)&=&([\lambda,l_2]_{_{RN}}+[N,l_2]_{_{RN}}+[\alpha,l_1]_{_{RN}})(X,f)\\
                           &=&l_2(X,\lambda f)-\lambda l_2(X,f)+l_2(NX,f)-\alpha(X,l_1(f))\\
                           &=&\frac{1}{2}\lambda\langle X,l_1(f)\rangle-\frac{1}{2}\lambda\langle X,l_1(f)\rangle+\frac{1}{2}\langle NX,l_1(f)\rangle-\alpha(X,l_1(f))\\
                           &=&\frac{1}{2}\langle NX,l_1(f)\rangle-\alpha(X,l_1(f)),
\end{array}
\end{equation}
and the same computations for $(f,X)$ instead of $(X,f)$ gives
\begin{equation}\label{eq:174A}
[\mathcal{N},\mu]_{_{RN}}^2(f,X)=\frac{1}{2}\langle NX,l_1(f)\rangle-\alpha(l_1(f),X).
\end{equation}
Since $[\mathcal{N},\mu]_{_{RN}}^2(X,f)=[\mathcal{N},\mu]_{_{RN}}^2(f,X)$, from (\ref{equation2intheorem}) and (\ref{eq:174A}) we get $\alpha(X,l_1(f))=0$, for all $X\in \Gamma(E), f\in \mathcal{C}^{\infty}(M).$ So
\begin{equation}\label{eq:deformedon(X,f)}
[\mathcal{N},\mu]_{_{RN}}^2(X,f)=\frac{1}{2}\langle NX,l_1(f)\rangle.
\end{equation}
For $X,Y\in \Gamma(E)$, we have
\begin{equation}\label{equation2intheorem-2}
\begin{array}{rcl}
[\mathcal{N},\mu]_{_{RN}}^2(X,Y)&=&([\lambda,l_2]_{_{RN}}+[N,l_2]_{_{RN}}+[\alpha,l_1]_{_{RN}})(X,Y)\\
                           &=&l_2(NX,Y)+l_2(X,NY)-Nl_2(X,Y)+l_1\alpha(X,Y).\\

\end{array}
\end{equation}
Third item in the Definition \ref{def:Courantskew} implies that if $[\mathcal{N},\mu]_{_{RN}}$ is a Lie $2$-algebra associated to a Courant algebroid, then we must have:
\begin{equation}\label{mainbody}
[\mathcal{N},\mu]_{_{RN}}^2(X,fY)=f[\mathcal{N},\mu]_{_{RN}}^2(X,Y)+2[\mathcal{N},\mu]_{_{RN}}^2(X,f).Y-\frac{1}{2}\langle X,Y\rangle[\mathcal{N},\mu]^1_{RN}(f).
\end{equation}
But using (\ref{equation1intheorem-1}), (\ref{eq:deformedon(X,f)}) and (\ref{equation2intheorem-2}),
\begin{equation}\label{eq:LHS}
\begin{array}{rcl}
                              &&[\mathcal{N},\mu]_{_{RN}}^2(X,fY)\\
                              &=&l_2(NX,fY)+l_2(X,NfY)-Nl_2(X,fY)+l_1\alpha(X,fY)\\
                              &=&fl_2(NX,Y)+2l_2(NX,f)Y-\frac{1}{2}\langle NX,Y\rangle l_1(f)\\
                              &+&fl_2(X,NY)+2l_2(X,f)NY-\frac{1}{2}\langle X,NY\rangle l_1(f)\\
                              &-&fNl_2(X,NY)-2l_2(X,f)NY+\frac{1}{2}\langle X,Y\rangle Nl_1(f)\\
                              &+&fl_1\alpha(X,Y)+\alpha(X,Y)l_1(f)\\
                              &=&f(l_2(NX,Y)+l_2(X,NY)-Nl_2(X,Y)+l_1\alpha(X,Y))+2l_2(NX,f)Y\\
                              &-&\frac{1}{2}\langle X,(N+N^*)Y\rangle l_1(f)+\frac{1}{2}\langle X,Y\rangle Nl_1(f)+\alpha(X,Y)l_1(f)\\
\end{array}
\end{equation}
and
\begin{equation}\label{eq:RHS}
\begin{array}{rcl}
&&f[\mathcal{N},\mu]_{_{RN}}^2(X,Y)+2[\mathcal{N},\mu]_{_{RN}}^2(X,f).Y-\frac{1}{2}\langle X,Y\rangle[\mathcal{N},\mu]^1_{RN}(f)\\
&=&f(l_2(NX,Y)+l_2(X,NY)-Nl_2(X,NY)+l_1\alpha(X,Y))+2l_2(NX,f).Y\\
&-&\frac{1}{2}\langle X,Y\rangle(\lambda Id_{\Gamma(E)}-N)l_1(f).
\end{array}
\end{equation}
Now Equations (\ref{mainbody}), (\ref{eq:LHS}) and (\ref{eq:RHS}) show that
\begin{equation*}\label{akhari}
\frac{1}{2}\langle X,(N+N^*-\lambda Id_{\Gamma(E)})Y\rangle l_1(f)=\alpha(X,Y)l_1(f),
\end{equation*}
for all $X,Y\in \Gamma(E)$ and all $f\in \mathcal{C}^{\infty}(M)$. Since $\alpha$ is skew-symmetric and $\langle . ,(N+N^*-\lambda Id).\rangle$ is symmetric on $\Gamma(E) \times \Gamma(E)$ and since the anchor is not zero everywhere, which implies that $l_1(f)$ is not always zero, we have $\alpha=0$ and $N+N^*-\lambda Id=0$.
\end{proof}
\begin{cor}\label{thewrongcorollary}
Let $(E,\circ, \rho,\langle.,.\rangle)$ be a Courant algebroid  with anchor $\rho$ being different from $0$ on a dense subset of $E$, with the associated Lie $2$-algebra structure $\mu$ on the graded vector space $\mathcal{C}^{\infty}(M)\oplus \Gamma(E)$. Then, there is a one to one correspondence between:
 \begin{enumerate}
 \item[(i)] quadruples $(N,K,\lambda,\gamma)$ with $N,K$ being $(1-1)$-tensors and $\lambda,\gamma$ being Casimir functions satisfying the following conditions:
    $$ \begin{cases}
     \circ^{N,N}=\circ^K,\\
     N K - K N =0,\\
     N+N^*=\lambda Id_{\Gamma(E)} ,\\
     K+K^* = \gamma Id_{\Gamma(E)},\\
     (\Gamma(E), \circ^N)\,\,\, \mbox{and}\,\,\, (\Gamma(E), \circ^K) \,\,\,\mbox{are Leibniz algebras}.
     \end{cases}$$
 \item[(ii)] Nijenhuis vector valued forms ${\mathcal N} $ with respect to $\mu$, with square ${\mathcal K} $ such that the deformed brackets $[ {\mathcal N}, \mu]_{_{RN}} $ and $[ {\mathcal K}, \mu]_{_{RN}} $ are Lie $2$-algebras associated to  Courant structures with the same scalar product.
 \end{enumerate}
\end{cor}
\begin{proof}
Given a quadruple $(N,K,\lambda,\gamma)$ satisfying  conditions in item $(i)$,  we define vector valued $1$-forms $\mathcal{N}$ and $\mathcal{K}$ on the graded vector space $\mathcal{C}^{\infty}(M)\oplus \Gamma(E)$ as $\mathcal{N}(f)=\lambda f,\mathcal{K}(f)=\gamma f,  \mathcal{N}(X)=NX,\mathcal{K}(X)=KX$, for all $X\in \Gamma(E), f\in \mathcal{C}^{\infty}(M)$. We prove that $\mathcal{N}$ is a Nijenhuis vector valued form with respect to $\mu$, with square $\mathcal{K}$. First notice that using Corollary \ref{cor:Anchoruniqe}, the assumption $\circ^{N,N}=\circ^{K}$ implies that $(E,\circ^{N,N}, \rho^{N,N}, \langle .,.\rangle)$ and $(E,\circ^{K}, \rho^{K}, \langle .,.\rangle)$ determines the same pre-Courant algebroids, hence, they have the same associated pre-Lie $2$-algebras. On the other hand, using Lemma \ref{leml^N} twice, the pre-Lie $2$-algebra associated to the pre-Courant algebroid $(E,\circ^{N,N}, \rho^{N,N}, \langle .,.\rangle)$ is
$[\mathcal{N},[\mathcal{N},\mu]_{_{RN}}]_{_{RN}}$ and by using again Lemma \ref{leml^N}, the pre-Lie $2$-algebra associated to the pre-Courant algebroid $(E,\circ^{K}, \rho^{K}, \langle .,.\rangle)$ is $[\mathcal{K},\mu]_{_{RN}}$. Hence,
\begin{equation}\label{injast}
[\mathcal{N},[\mathcal{N},\mu]_{_{RN}}]_{_{RN}}=[\mathcal{K},\mu]_{_{RN}}.
\end{equation}
Also, using the assumption $NK-KN=0$ we get
\begin{equation}\label{keiman}
[\mathcal{N},\mathcal{K}]_{_{RN}}=0.
\end{equation}
Equations \ref{injast} and \ref{keiman} show that $\mathcal{N}$ is a Nijenhuis vector valued $1$-form with respect to $\mu$, with square $\mathcal{K}$. By Corollary \ref{cornew}, $[\mathcal{N},\mu]_{_{RN}}$ is a Lie $2$-algebra associated to the Courant algebroid $(E,\circ^N, \rho, \langle .,.\rangle)$ and $[\mathcal{K},\mu]_{_{RN}}$ is a Lie $2$-algebra associated to the Courant algebroid $(E,\circ^K, \rho, \langle .,.\rangle)$.

Conversely, assume that $\mathcal{N}$ is a Nijenhuis vector valued form with respect to $\mu$, with square $\mathcal{K}$ such that $[\mathcal{N},\mu]_{_{RN}}$ and $[\mathcal{K},\mu]_{_{RN}}$ are Lie $2$-algebras associated to  Courant algebroids. Then by Theorem \ref{them:NijenhuisCourant} $\mathcal{N}$ is of the form $\lambda+N$ with $N+N^*=\lambda Id_{\Gamma(E)}$ and $\mathcal{K}$ is of the form $\gamma+K$ with $K+K^*=\gamma Id_{\Gamma(E)}$. (Note that as we discussed in the proof of Theorem \ref{them:NijenhuisCourant}, $\lambda, \gamma$ are multiples of identity on $\mathcal{C}^{\infty}(M)$ which we denote both, the map and the coefficient, by the same $\lambda, \gamma$.) Also, it implies that the Courant algebroid which is associated to the Lie $2$-algebra $[\mathcal{N}, \mu]_{_{R,N}}$ (respectively, $[\mathcal{K}, \mu]_{_{R,N}}$) is $(E,\circ^N,\rho^N,\langle .,.\rangle)$ (respectively, $(E,\circ^K,\rho^K,\langle .,.\rangle)$ ), which means that $(\Gamma(E), \circ^N)$ and $(\Gamma(E), \circ^K)$ are Leibniz algebras. And  Since $\mathcal{N}$ is a Nijenhuis vector valued form with respect to $\mu$, with square $\mathcal{K}$, we have
\begin{equation}\label{famus}
[\mathcal{N},[\mathcal{N},\mu]_{_{RN}}]_{_{RN}}=[\mathcal{K},\mu]_{_{RN}}
\end{equation}
and
\begin{equation}\label{famusmus}
[\mathcal{N},\mathcal{K}]_{_{RN}}=0.
\end{equation}
Applying both sides of Equation (\ref{famus}) on a pair of sections $X,Y\in \Gamma(E)$ we get $X\circ^{N,N}Y=X\circ^K Y$ which implies that $\circ^{N,N}=\circ^K$. While Equation (\ref{famusmus}) implies that $KN-NK=0$.
\end{proof}

According to the proof of Proposition $2.6.5$ in \cite{RoytenbergPh.D.} the first properties in both Definitions \ref{def:Courantskew} and \ref{def:CourantDorfman}  are equivalent. This implies that if $(E,\circ,\rho,\langle.,.\rangle)$ is a pre-Courant algebroid, then the associated skew-symmetrized bracket $\left[.,.\right]$ satisfies the third property in Definition \ref{def:Courantskew}. Also, it comes from the definition of the map $\mathcal{D}$, given in (\ref{matcallD}), associated to a pre-Courant algebroid $(E,\circ,\rho,\langle.,.\rangle)$ that it is a derivation. These two arguments are enough to restate Theorem \ref{them:NijenhuisCourant} as follow
\begin{them} \label{them:NijenhuisCourant2}
Let $ (\circ,\rho, \langle . ,. \rangle ) $ be a Courant structure on a vector bundle $E\to M$, with the associated symmetric Lie $2$-algebra structure $l_1+l_2+l_3$ on the graded vector space $V=E_{-2}\oplus E_{-1}$, where $E_{-2}:=\mathcal{C}^{\infty}(M)$ and $E_{-1}:=\Gamma(E)$.
Let ${\mathcal N} = \lambda + N + \alpha$ be a $C^\infty(M)$-linear vector valued form of degree $0$ on $V$.
Assume also that $\rho$ is not equal to zero on a dense subset of the base manifold.
If $ [{\mathcal N}, l_1+l_2+l_3]_{_{RN}}=l'_1+l'_2+l'_3 $, where the vector valued forms $l'_1,l'_2,l'_3$ are obtained from a pre-Courant algebroid, with the same scalar product, by the construction given in (\ref{Lie2fromCourant}), then
\begin{enumerate}
\item $\lambda$ is a Casimir,
\item $\alpha=0,$
\item $ N+N^* = \lambda Id_{\Gamma(E)}.$
%
\end{enumerate}
In this case, the Courant structure that $ [{\mathcal N}, l_1+l_2+l_3]_{_{RN}}$ is associated with, is $(\left[.,.\right]_N, \rho^N,\langle . ,. \rangle).$
\end{them}
And this leads to the next result:
\begin{cor}\label{thewrongcorollary2}
Let $(E,\circ, \rho,\langle.,.\rangle)$ be a Courant algebroid  with anchor $\rho$ being different from $0$ on a dense subset of $E$, with the associated Lie $2$-algebra structure $\mu=l_1+l_2+l_3$ on the graded vector space $\mathcal{C}^{\infty}(M)\oplus \Gamma(E)$. Then, there is a one to one correspondence between:
 \begin{enumerate}
 \item[(i)] quadruples $(N,K,\lambda,\gamma)$ with $N,K$ being $(1-1)$-tensors and $\lambda,\gamma$ being Casimir functions satisfying the following conditions:
    \begin{equation}\label{yekikamtar}
     \begin{cases}
     \circ^{N,N}=\circ^K,\\
     N K - K N =0,\\
     N+N^*=\lambda Id_{\Gamma(E)} ,\\
     K+K^* = \gamma Id_{\Gamma(E)}.\\
     \end{cases}
     \end{equation}
 \item[(ii)] Nijenhuis vector valued forms ${\mathcal N} $ with respect to $\mu$, with square ${\mathcal K} $ such that the deformed bracket is of the form $[ {\mathcal N}, \mu]_{_{RN}}=l'_1+l'_2+l'_3 $ and $l'_1,l'_2,l'_3$ are constructed by the procedure in (\ref{Lie2fromCourant}) obtained from a pre-Courant algebroid, with the same scalar product.
 \end{enumerate}
\end{cor}
\begin{proof}
Let $\mathcal{N}$ be a Nijenhuis vector valued form with respect to the Lie $2$-algebra structure $\mu=l_1+l_2+l_3$, with square $\mathcal{K}$ and assume that $\left[\mathcal{N},\mu\right]_{_{RN}}$ is obtained from a pre-Courant algebroid. Let
\begin{equation*}
\mathcal{N}|_{\Gamma(E)}=N,\quad \mathcal{N}|_{\mathcal{C}^{\infty}(M)}=\lambda Id_{\mathcal{C}^{\infty}(M)},\quad \mathcal{K}|_{\Gamma(E)}=K \quad\mbox{and}\quad\mathcal{K}|_{\mathcal{C}^{\infty}(M)}=\gamma Id_{\mathcal{C}^{\infty}(M)}.
\end{equation*}
By Theorem \ref{them:NijenhuisCourant2}, $N+N^*=\lambda Id_{\Gamma(E)}$ and $(E,\circ^N, \rho^N,\langle.,.\rangle)$ is a pre-Courant algebroid( it is, in fact, the pre-Courant algebroid which $\left[\mathcal{N},\mu\right]_{_{RN}}$ is obtained from!). Hence, by Lemma \ref{pre}, $(E,\circ^{N,N}, \rho^{N,N},\langle.,.\rangle)$ is a pre-Courant algebroid. Now, Lemma \ref{leml^N} implies that $\left[\mathcal{K},\mu\right]_{_{RN}}=\left[\mathcal{N},\left[\mathcal{N},\mu\right]_{_{RN}}\right]_{_{RN}}$ is obtained from the pre-Courant algebroid $(E,\circ^{N,N}, \rho^{N,N},\langle.,.\rangle)$, by the construction given in (\ref{Lie2fromCourant}). Therefore, by Theorem \ref{them:NijenhuisCourant2}, $K+K^*=\gamma Id_{\Gamma(E)}$. While, $\left[\mathcal{N},\mathcal{K}\right]_{_{RN}}=0$ implies that $NK-KN=0$ and $\left[\mathcal{N},\left[\mathcal{N},\mu\right]_{_{RN}}\right]_{_{RN}}=\left[\mathcal{K},\mu\right]_{_{RN}}$ implies that $\circ^{N,N}=\circ^{K}$.

Conversely, assume that we are given a quadruple $(N,K,\lambda,\gamma)$ satisfying the properties in (\ref{yekikamtar}). By Lemma \ref{pre}, $(E,\circ^{N},\rho^{N},\langle .,.\rangle)$ is a pre-Courant and by Lemma \ref{leml^N}, the pre-Lie $2$-algebra structure associated to the pre-Courant algebroid $(E,\circ^{N},\rho^{N},\langle .,.\rangle)$ is $\left[\mathcal{N},\mu\right]_{_{RN}}$. Similar arguments prove that the pre-Lie $2$-algebra structure associated to the pre-Courant algebroid $(E,\circ^{N,N},\rho^{N,N},\langle .,.\rangle)$ is $\left[\mathcal{N},\left[\mathcal{N},\mu\right]_{_{RN}}\right]_{_{RN}}$ and the pre-Lie $2$-algebra structure associated to the pre-Courant algebroid $(E,\circ^{K},\rho^{K},\langle .,.\rangle)$ is $\left[\mathcal{K},\mu\right]_{_{RN}}$. Now, the assumption $\circ^{N,N}=\circ^{K}$ and Lemma \ref{cor:Anchoruniqe} imply that $(E,\circ^{N,N},\rho^{N,N},\langle .,.\rangle)$ and  $(E,\circ^{K},\rho^{K},\langle .,.\rangle)$ are the same pre-Courant algebroids, therefore, we have $\left[\mathcal{N},\left[\mathcal{N},\mu\right]_{_{RN}}\right]_{_{RN}}=\left[\mathcal{K},\mu\right]_{_{RN}}$. It follows from the assumption $NK-KN=0$ that $\left[\mathcal{N},\mathcal{K}\right]=0$. Hence, $\mathcal{N}$ is a Nijenhuis vector valued form with respect to the Lie $2$-algebra structure $\mu$, with square $\mathcal{K}$.
\end{proof}

\section{Lie algebroids}
We saw in Subsection \ref{subsection:Example around Lie algebroids} that given a vector bundle $A\to M$ over a manifold $M$, there is a one to one correspondence between multiplicative $L_{\infty}$-structures on the graded vector space $\Gamma(\wedge A)[2]$, with only binary brackets (GLA-structures whose Lie brackets are derivations) and Lie algebroid structures on the vector bundle $A\to M$. On the other hand, in Subsection \ref{subsection:Lie algebra,GLAandDGLA} we studied examples of Nijenhuis vector valued $1$-forms (or at most Nijenhuis vector valued forms which were the sum of $1$-forms with $0$-forms) with respect to GLA-structures. In this section we study examples of Nijenhuis vector valued forms which are the sum of vector valued $1$-forms with vector valued $2$-forms with respect to multiplicative GLA-structures on the graded vector space $\Gamma(\wedge A)[2]$. Throughout this section,  $(A,\left[.,.\right],\rho)$ stands for an arbitrary Lie algebroid over a manifold $M$.
\subsection{Extensions by derivation}
Let us first fix some notations.\\
 Let $N$ be a $(1-1)$-tensor on a Lie algebroid $(A,\left[.,.\right],\rho)$. By \emph{extension of $N$ by derivation} on the graded vector space $\Gamma(\wedge A)[2]$ we mean a linear map denoted by $\underline{N}$ which is defined as $\underline{N}(f):=0$, for all $f\in \mathcal{C}^{\infty}(M)$, and
\begin{equation*}
\underline{N}(P):=\sum_{i=1}^{p}P_1\wedge \cdots \wedge P_{i-1}\wedge N(P_i)\wedge P_{i+1}\wedge \cdots \wedge P_p,
\end{equation*}
for all homogeneous multi-sections $P=P_1\wedge \cdots \wedge P_{p}\in \Gamma(\wedge A)[2]$. It follows from definition that for any  $(1-1)$-tensor field $N$, on a Lie algebroid $(A,\left[.,.\right],\rho)$, the extension of  $N$ by derivation is a multi-derivation on the graded vector space $\Gamma(\wedge A)[2]$ hence a symmetric vector valued $1$-form on $\Gamma(\wedge A)[2]$ and it is of degree zero. In general, for a $k$-form on $A$,$\kappa \in\Gamma(\wedge^kA^*)$, the extension of $\kappa$ by derivation is a $k$-linear map denoted by $\underline{\kappa}$  given by
\begin{equation*}
\underline{\kappa}(P_1,\cdots,P_k):=\sum_{i_1,\cdots, i_k=1}^{p_1,\cdots, p_k}(-1)^{\spadesuit}\kappa(P_{1,i_1},\cdots,P_{k,i_k})\widehat{P_{1,i_1}}\wedge\cdots\wedge\widehat{P_{k,i_k}},
\end{equation*}
for all homogeneous multi-sections $P_i=P_{i,1}\wedge\cdots\wedge P_{i,p_i} \in \Gamma(\wedge^{p_i} A)$, with $i=1,\cdots,k,$ where $1\leq i_j\leq p_j$ for all $1\leq j\leq k$,
\begin{equation*}
\widehat{P_{j,i_j}}=P_{j,1}\wedge\cdots\wedge P_{j,i_j-1}\wedge P_{j,i_j+1}\wedge\cdots\wedge P_{j,p_j}\in \Gamma(\wedge^{p_j-1} A)
\end{equation*}
and
\begin{equation*}
\spadesuit=2p_1+3p_2+\cdots+(k+1)p_k+i_1+\cdots +i_k+1.
\end{equation*}
It follows from definition that the extension of a $k$-form by derivation is a multi-derivation on the graded vector space $\Gamma(\wedge A)[2]$ equipped with the graded commutative associative product $\wedge$ and that it is a symmetric vector valued $k$-form  of degree $k-2$ on the graded vector space $\Gamma(\wedge A)[2]$. The fact that every multi-derivation on the graded vector space $\Gamma(\wedge A)[2]$ is uniquely determined on the space of sections $\Gamma(A)$ implies the following lemma which we will use it in the next subsection.
\begin{lem}\label{underlinescommut}
Let $(A, \left[.,.\right], \rho)$ be a Lie algebroid, $\alpha \in \Gamma(\wedge^k A^*)$ be a $k$-form and $\beta \in \Gamma(\wedge^l A^*)$ be an $l$-form. Then,
\begin{equation*}
\left[\underline{\alpha},\underline{\beta}\right]_{_{RN}}=0.
\end{equation*}
\end{lem}
\begin{proof}
The facts that $\underline{\alpha}$ (respectively $\underline{\beta}$) is a vector valued $k$-form (respectively $l$-form) of degree $k-2$ (respectively $l-2$), imply that $\left[\underline{\alpha},\underline{\beta}\right]_{_{RN}}$ is a vector valued $(k+l-1)$-form of degree $k+l-4$ on the graded vector space $\Gamma(\wedge A)=\oplus_{i\geq 0}\Gamma(\wedge^i A)$. Therefore, for all $l,k\geq 0$ the restriction of $\left[\underline{\alpha},\underline{\beta}\right]_{_{RN}}$ to the space of sections is zero and hence we have $\left[\underline{\alpha},\underline{\beta}\right]_{_{RN}}=0$, because $\left[\underline{\alpha},\underline{\beta}\right]_{_{RN}}$ is a multi-derivation and it is determined uniquely on the space of sections.
\end{proof}

According to Proposition \ref{prop:algebroid}, for a given Lie algebroid $(A,\,\left[.,.\right],\rho)$, the bracket $l_2^{\left[.,.\right]}$ given by
\begin{equation}\label{multiplicativestructure}
l_2^{\left[.,.\right]}(P,Q)=(-1)^{p-1}[P,Q]_{_{SN}}, \,\,\, P\in \Gamma(\wedge^{p}A) , Q\in \Gamma(\wedge^{q}A),
\end{equation}
 defines a multiplicative GLA-structure on the graded vector space $\Gamma(\wedge A)[2]$. Note that this is a one to one correspondence, meaning that if the bracket $\left[.,.\right]$ is not Lie, then $l_2^{\left[.,.\right]}$ is not a multiplicative GLA-structure on the graded vector space $\Gamma(\wedge A)[2]$. Now let $N$ be a $(1-1)$-tensor field on a Lie algebroid $(A,\,\left[.,.\right],\rho)$. Then one can deform the bracket $\left[.,.\right]$ by $N$ as
\begin{equation*}
\left[X,Y\right]_N=\left[NX,Y\right]+\left[X,NY\right]-N\left[X,Y\right],
\end{equation*}
for all $X,Y\in \Gamma(A)$. Using this bracket, of course we may consider $l_2^{\left[.,.\right]_N}$ in the same way as in Equation (\ref{multiplicativestructure}), where the Schouten-Nijenhuis bracket corresponding to the deformed bracket $\left[.,.\right]_N$ is the unique extension by derivation, of the bracket $\left[.,.\right]_N$, on the space of sections, to the space of multi-sections. Note that the bracket $l_2^{\left[.,.\right]_N}$ is not necessarily a multiplicative GLA-structure.  On the other hand, since $l_2^{\left[.,.\right]}$ is a symmetric vector valued $2$-form of degree $1$ and $\underline{N}$ is a (symmetric) vector valued $1$-form of degree zero, we may also talk about the deformation of $l_2^{\left[.,.\right]}$ by $\underline{N}$. The following lemma shows the relation between the deformed multiplicative GLA-structure $\left[\underline{N},l_2^{\left[.,.\right]}\right]_{_{RN}}$ with $l_2^{\left[.,.\right]_N}$.
\begin{lem}\label{lem:deformingOidsN}
 Let $N$ be a $(1,1)$-tensor field on a Lie algebroid $(A,\left[.,.\right],\rho)$.
 Then we have:
  \begin{equation}\label{eq:underlineN}
   \left[\underline N , l_2^{\left[.,.\right]} \right]_{_{RN}} = l_2^{\left[.,.\right]_N}.
   \end{equation}
  \end{lem}
\begin{proof}
The proof follows directly from the fact that the Schouten-Nijenhuis bracket on $\Gamma(\wedge A) $
    associated to the bracket $\left[.,.\right]_N$  is given by
 $$ \left[P,Q\right]'_{_{SN}}=\left[\underline N P,Q\right]_{_{SN}} + \left[P, \underline N Q\right]_{_{SN}} - \underline N \left[P,Q\right]_{_{SN}},$$
 for all $P,Q \in \Gamma(\wedge A)$. (see \cite{YKS}.)
\end{proof}
We will need the following lemma for our next purpose.
\begin{lem}\label{lem:commutatorAlgebroidContraction}
Let $(A,\left[.,.\right], \rho)$ be a Lie algebroid, with the associated de Rham differential $\diff^A$, and associated
multiplicative GLA-structure $l_2^{\left[.,.\right]}$ on the graded vector space $\Gamma(\wedge A)[2]$. Then,
  \begin{equation*}
  \left[\underline{\alpha}, l_2^{\left[.,.\right]} \right]_{_{RN}}= \underline{\diff^A \alpha},
  \end{equation*}
  for all $\alpha \in \Gamma(\wedge^n A^*)$.
\end{lem}
\begin{proof}
We shall prove the statement for $n=2$. A similar proof can be done for any $n\geq 1$. First note that $\left[\underline{\alpha}, l_2^{\left[.,.\right]} \right]_{_{RN}}$ is a vector valued $3$-form of degree $1$ on the graded vector space $\Gamma(\wedge A)[2]$. This implies that the restriction of $\left[\underline{\alpha}, l_2^{\left[.,.\right]} \right]_{_{RN}}$ on the sections $\Gamma(A)$ is of the form:
\begin{equation*}
\left[\underline{\alpha}, l_2^{\left[.,.\right]} \right]_{_{RN}}|_{\Gamma(A) \times \Gamma(A) \times \Gamma(A)}:\Gamma(A)\times \Gamma(A)\times \Gamma(A)\to {\mathcal C}^{\infty}(M)
\end{equation*}
and it is zero on all other possible terms. On the other hand by Proposition \ref{prop:Multider}, $\left[\underline{\alpha}, l_2^{\left[.,.\right]} \right]_{_{RN}}$ is a multi-derivation and hence its restriction to the sections $\Gamma(A)$ is a ${\mathcal C}^{\infty}(M)$-linear map. Therefore  $\left[\underline{\alpha}, l_2^{\left[.,.\right]} \right]_{_{RN}}\in \Gamma(\wedge^3A^*)$. Next, we show that
\begin{equation}\label{restriction to sections}
\left[\underline{\alpha}, l_2^{\left[.,.\right]} \right]_{_{RN}}|_{\Gamma(A) \times \Gamma(A) \times \Gamma(A)}=\diff^A\alpha
\end{equation}
and this together with the fact that $\left[\underline{\alpha}, l_2^{\left[.,.\right]} \right]_{_{RN}}$ is a multi-derivation will imply that \begin{equation*}
\left[\underline{\alpha}, l_2^{\left[.,.\right]} \right]_{_{RN}}=\underline{\diff^A\alpha},
\end{equation*}
 by the uniqueness of extension by derivation of $\diff^A\alpha$ to the graded vector space $\Gamma(\wedge A)[2]$. But a direct computation shows that
\begin{equation*}
\left[\underline{\alpha}, l_2^{\left[.,.\right]} \right]_{_{RN}}(X,Y,Z)=\left[\alpha(X,Y),Z\right]_{_{SN}}-\alpha(\left[X,Y\right]_{_{SN}},Z)+c.p.,
\end{equation*}
for all $X,Y,Z\in \Gamma(A)$. Hence, Equation (\ref{eq:defGerstenhaber}) together with the definition of $\diff^A$ imply that
\begin{equation*}
\left[\underline{\alpha}, l_2^{\left[.,.\right]} \right]_{_{RN}}(X,Y,Z)=\rho(Z)\alpha(X,Y)-\alpha(\left[X,Y\right],Z)+c.p.=\diff\alpha^A(X,Y,Z).
\end{equation*}
\end{proof}
\subsection{Nijenhuis forms on multiplicative GLA associated to Lie algebroids}
Let $(A,\left[.,.\right],\rho)$ be a Lie algebroid and $N:\Gamma(A)\to \Gamma(A)$ be a $(1-1)$-tensor field. Then, similar to the case of Lie algebras, the Nijenhuis torsion of $N$ with respect to the Lie bracket $\left[.,.\right]$ , denoted by $T_{\left[.,.\right]}N$, is defined similar to the Equation (\ref{Torsion}) and again a direct computation shows that
\begin{equation*}
T_{\left[.,.\right]}N(X,Y)=\frac{1}{2}(\left[X,Y\right]_{N,N}-\left[X,Y\right]_{N^2}),
\end{equation*}
for all $X,Y \in \Gamma(A)$.
A $(1-1)$-tensor field $N$ on a Lie algebroid $(A,\left[.,.\right],\rho)$ is said to be \emph{Nijenhuis} if the Nijenhuis torsion of $N$, with respect to the Lie algebroid bracket $\left[.,.\right]$, vanishes. As a consequence of Lemma \ref{lem:deformingOidsN}, we have the following proposition:
\begin{prop}\label{prop:NijenhuisAsNijenhuis}
For every Nijenhuis tensor field $N$ on a Lie algebroid $(A,\left[.,.\right],\rho)$,
the extension of $N$ by derivation, $\underline N $, is a Nijenhuis vector valued $1$-form with respect to the multiplicative GLA-structure $l_2^{\left[.,.\right]}$ on the graded vector space $\Gamma(\wedge A)[2]$, with square ${\underline{(N^2)}}$.
\end{prop}
\begin{proof}
 Applying Lemma \ref{lem:deformingOidsN} twice, for the tensor field $N$ and the bracket $l_2^{\left[.,.\right]}$, we get $\left[\underline N ,\left[\underline N , l_2^{\left[.,.\right]}\right]_{_{RN}}\right]_{_{RN}} = l_2^{\left[.,.\right]_{N,N}} $. While applying Lemma \ref{lem:deformingOidsN} for the tensor field $N^2$ and the bracket $l_2^{\left[.,.\right]}$ we get $\left[\underline {N^2} , l_2^{\left[.,.\right]}\right]_{_{RN}} = l_2^{\left[.,.\right]_{N^2}}$. Since $N$ is a Nijenhuis $(1-1)$-tensor field, we have $l_2^{\left[.,.\right]_{N,N}}=l_2^{\left[.,.\right]_{N^2}}$ which implies that $\left[\underline N ,\left[\underline N , l_2^{\left[.,.\right]}\right]_{_{RN}}\right]_{_{RN}} = \left[\underline {N^2} , l_2^{\left[.,.\right]}\right]_{_{RN}}$.
Also, ${\underline{(N^2)}} $ and ${\underline{N}} $ commute with respect to the Richardson-Nijenhuis bracket.
\end{proof}
We have now all the tools to prove the next proposition where we obtain a Nijenhuis vector valued form which is the sum of a vector valued $1$-form with a vector valued $2$-form.
\begin{prop}
Let $(A, \left[.,.\right], \rho)$ be a Lie algebroid, with de Rham differential $\diff^A$ and associated
multiplicative GLA-structure $l_2^{\left[.,.\right]}$. Then, for every section $\alpha \in \Gamma(\wedge^2 A^*)$,
$S+\underline{\alpha} $ is a Nijenhuis vector valued form with respect to $l_2^{\left[.,.\right]}$, with square $ S+2\underline{\alpha}$.
The deformed structure is $l_2^{\left[.,.\right]}+ \underline{\diff^A\alpha} $.
\end{prop}
\begin{proof}
As a direct consequence of Lemma \ref{lem:commutatorAlgebroidContraction} we have
\begin{equation}\label{eq:deformedstructurealpha}
\left[S+\underline{\alpha},l_2^{\left[.,.\right]}\right]_{_{RN}} = l_2^{\left[.,.\right]}+ \underline{\diff^A\alpha}.
\end{equation}
Hence, Lemma \ref{underlinescommut} implies that
\begin{equation*}
\left[S+\underline{\alpha},\left[S+\underline{\alpha},l_2^{\left[.,.\right]}\right]_{_{RN}}\right]_{_{RN}}=l_2^{\left[.,.\right]}+ 2\underline{\diff^A\alpha}.
\end{equation*}

The fact that $ \left[S+\underline{\alpha},S+2\underline{\alpha}\right]_{_{RN}}=0$ follows immediately from Lemma \ref{lem:commutatorAlgebroidContraction}.
\end{proof}

\subsection{$\Omega N$, Poisson-Nijenhuis and $\Pi \Omega $-structures on Lie algebroids}
Next, for a given Lie algebroid $(A,\left[.,.\right], \rho)$, we give a Nijenhuis vector valued form which is the sum of a vector valued $1$-form with a vector valued $2$-form, with respect to the associated GLA-structure $l_2^{\left[.,.\right]}$, by the help of a Nijenhuis $(1-1)$-tensor field $N$ on $A$.
\begin{prop}\label{omega_N}
Let $(A,\left[.,.\right], \rho)$ be a Lie algebroid, with the associated de Rham differential $\diff^A$ and with the associated multiplicative GLA structure $l_2^{\left[.,.\right]}$ on the graded vector space $\Gamma(\wedge A)[2]$, let $N$ be a Nijenhuis $(1-1)$-tensor field on the Lie algebroid and $\alpha\in \Gamma(\wedge^2 A^*)$ be a $2$-form such that  $\alpha_N:\Gamma(A)\times\Gamma(A)\to\Gamma(A)$ given by
$$\alpha_{_N}(X,Y)=\alpha(NX,Y)$$
is skew-symmetric and therefore a $2$-form on $A$. Then,
\begin{enumerate}
\item $\left[\underline{N},\underline{\alpha}\right]_{_{RN}}=2 \underline{\alpha_{_N}},$ \emph{(unaffected by the condition of $N$ being Nijenhuis)}
\item $\left[\underline{N}+\underline{\alpha},l_2^{\left[.,.\right]}\right]_{_{RN}}=l_2^{\left[.,.\right]_N}+\underline{\diff^A\alpha}$
\item $\left[\underline{N}+\underline{\alpha},\left[\underline{N}+\underline{\alpha},
      l_2^{\left[.,.\right]}\right]_{_{RN}}\right]_{_{RN}}=
      \left[\underline{N^2},l_2^{\left[.,.\right]}\right]_{_{RN}}-2\underline{\diff^A\alpha_{_{N}}}+2\left[\underline{N},
      \underline{\diff^A\alpha}\right]_{_{RN}}$.
\end{enumerate}
\end{prop}
\begin{proof}
1) First notice that for all $X,Y\in \Gamma(A)$ we have $$\left[\underline{N},\underline{\alpha}\right]_{_{RN}}(X,Y)=\alpha(NX,Y)-\alpha(NY,X)=2\alpha_{_N}(X,Y).$$
 Since $\underline{N}$ and $\underline{\alpha}$ are both derivations, by Lemma \ref{lem:multiDer1} $\left[\underline{N},\underline{\alpha}\right]_{_{RN}}$ is a derivation and hence it is the unique extension of $2\alpha_{_N}$ by derivation.

2) It is a direct consequence of Lemma \ref{lem:deformingOidsN} together with Lemma \ref{lem:commutatorAlgebroidContraction}.

3) Using item 2 and Lemma \ref{lem:deformingOidsN}
\begin{equation*}
      \begin{array}{rcl}
      &&\left[\underline{N}+\underline{\alpha},\left[\underline{N}+\underline{\alpha},
            l_2^{\left[.,.\right]}\right]_{_{RN}}\right]_{_{RN}}\\
      &=&\left[\underline{N}+\underline{\alpha},l_2^{\left[.,.\right]_N}+\underline{\diff^A\alpha} \right]_{_{RN}}\\
      &=&l_2^{\left[.,.\right]_{N,N}}+\left[\underline{N},\underline{\diff^A\alpha}\right]_{_{RN}}+\left[\underline{\alpha},
      l_2^{\left[.,.\right]_{N}}\right]_{_{RN}}
      +\left[\underline{\alpha},\underline{\diff^A\alpha}\right]_{_{RN}}.\\
      \end{array}
\end{equation*}
But, using Lemma \ref{lem:deformingOidsN} and the graded Jacobi identity we have
\begin{equation*}
\begin{array}{rcl}
\left[\underline{\alpha},l_2^{\left[.,.\right]_{N}}\right]_{_{RN}}&=&
\left[\underline{\alpha},\left[\underline{N},l_2^{\left[.,.\right]}\right]_{_{RN}}\right]_{_{RN}}\\
&=&\left[\left[\underline{\alpha},\underline{N}\right]_{_{RN}},l_2^{\left[.,.\right]}\right]_{_{RN}}+
\left[\underline{N},\left[\underline{\alpha},l_2^{\left[.,.\right]}\right]_{_{RN}}\right]_{_{RN}}\\
&=&\left[-2\underline{\alpha_{_N}},l_2^{\left[.,.\right]}\right]_{_{RN}}+\left[\underline{N},\underline{\diff^A\alpha}\right]_{_{RN}}
\end{array}
\end{equation*}
and by Lemma \ref{underlinescommut} $\left[\underline{\alpha},\underline{\diff^A\alpha}\right]_{_{RN}}=0.$ Hence, since $N$ is Nijenhuis we get
\begin{equation*}
\begin{array}{rcl}
[\underline{N}+\underline{\alpha},[\underline{N}+\underline{\alpha},l_2^{\left[.,.\right]}]_{_{RN}}]_{_{RN}}&=& \left[\underline{N^2}-2\underline{\alpha_{_N}},l_2^{\left[.,.\right]}\right]_{_{RN}}+2\left[\underline{N},\underline{\diff^A\alpha}\right]_{_{RN}}\\
&=&\left[\underline{N^2},l_2^{\left[.,.\right]}\right]_{_{RN}}-2\underline{\diff^A\alpha_{_{N}}}+2\left[\underline{N},\underline{\diff^A\alpha}\right]_{_{RN}}.
\end{array}
\end{equation*}
\end{proof}
Proposition \ref{omega_N} in fact, gives a Nijenhuis form on what so called $\Omega N$-structure. Let us first recall this notion:
\begin{defi}
Let $(A,\left[.,.\right], \rho)$ be a Lie algebroid, with the associated de Rham differential $\diff^A$, $N$ be  a $(1-1)$-tensor field on $A$ and $\alpha \in \Gamma(\wedge^2 A^*)$ be a $2$-form. Let $\alpha_{_N}:\Gamma(A)\times\Gamma(A)\to\Gamma(A)$ be a bilinear map, defined as
\begin{equation*}
\alpha_{_N}(X,Y)=\alpha(NX,Y).
\end{equation*}
 Then, the pair $(N,\alpha)$ is called an $\Omega N$-structure on the Lie algebroid $A$ if, $\alpha_{_N}$ is skew-symmetric (and therefore a $2$-form on $A$) and $\alpha$ and $\alpha_{_N}$ are $\diff^A$-closed.
\end{defi}
\begin{cor}\label{NijenhuisonOmegaNstructures}
Let $(A,\left[.,.\right], \rho)$ be a Lie algebroid, with the associated de Rham differential $\diff^A$ and with the associated multiplicative GLA structure $l_2^{\left[.,.\right]}$ on the graded vector space $\Gamma(\wedge A)$. Let $(N,\alpha)$ be an $\Omega N$-structure on the Lie algebroid $A$, then $\underline{N}+\underline{\alpha}$ is a Nijenhuis vector valued form, with respect to $l_2^{\left[.,.\right]}$, with square $\underline{N^2}+\underline{\alpha_{_{N}}}$.
\end{cor}
\begin{proof}
It follows from item 3 in Proposition \ref{omega_N}, using the condition $\left[\underline{\alpha_{_{N}}},l_2^{\left[.,.\right]}\right]_{_{RN}}=0$, that
\begin{equation*}
\left[\underline{N}+\underline{\alpha},\left[\underline{N}+\underline{\alpha},
      l_2^{\left[.,.\right]}\right]_{_{RN}}\right]_{_{RN}}=
      \left[\underline{N^2}+\underline{\alpha_{_{N}}},l_2^{\left[.,.\right]}\right]_{_{RN}}.
\end{equation*}
Hence, it is enough to check that
\begin{equation*}
\left[\underline{N}+\underline{\alpha},\underline{N^2}+\underline{\alpha_{_{N}}}\right]_{_{RN}}=0.
\end{equation*}
But,
\begin{equation*}
\left[\underline{N}+\underline{\alpha},\underline{N^2}+\underline{\alpha_{_{N}}}\right]_{_{RN}}=\left[\underline{N},\underline{\alpha_{_N}}\right]_{_{RN}}+
\left[\underline{\alpha},\underline{N^2}\right]_{_{RN}}=2\underline{(\alpha_{_{N}})_{_N}}-2\underline{\alpha_{_{N^2}}}=0.
\end{equation*}
\end{proof}

We would like of course to include Poisson-Nijenhuis structures among our examples of Nijenhuis
structures on $L_\infty$-algebras.
Let us first, fix and recall some notations and notions.

Let $(A,\mu=\left[.,.\right], \rho)$ be a Lie algebroid, $\pi\in \Gamma(\wedge^2A)$ be a bi-vector and $N:\Gamma(A)\to \Gamma(A)$ be a $(1-1)$-tensor field. Let $\pi^{\#}: \Gamma(A^*)\to \Gamma(A)$ be the induced linear map given by $\langle \beta,\pi^{\#}\alpha\rangle=\pi(\alpha, \beta)$,  $N^*: \Gamma(A^*)\to \Gamma(A^*)$ be the induced linear map given by $\langle N^*\alpha,X\rangle=\langle \alpha,NX\rangle$, and $\pi_{_{N}}$ be the induced bi-vector defined by $\pi_{_{N}}(\alpha,\beta)=\langle \beta,N\pi^{\#}\alpha\rangle=\langle N^*\beta,\pi^{\#}\alpha\rangle$, for all $\alpha,\beta \in  \Gamma(A^*)$.
  A bracket $\{\cdot,\cdot\}_{_{\pi}}^{\mu}$ can be defined on $\Gamma(A^*)$, the space of $1$-forms on the Lie algebroid $(A,\mu=\left[.,.\right], \rho)$, as
\begin{equation*}
\{\alpha,\beta\RP=\mathcal{L}^A_{_{\pi^\#\alpha}}\beta -\mathcal{L}^A_{_{\pi^\#\beta}}\alpha-\diff^A(\pi(\alpha,\beta)),
\end{equation*}
for all $\alpha,\beta\in \Gamma(A^*)$. It is well-known that if $\pi$ is a Poisson bi-vector on the Lie algebroid $(A,\mu=\left[.,.\right], \rho)$, that is $\left[\pi,\pi\right]_{_{SN}}=0$, then $(\Gamma(A^*), \{.,.\RP)$ is a Lie algebra and if this is the case, then $\pi^{\#}$ is a Lie algebra morphism form the Lie algebra $(\Gamma(A^*), \{.,.\RP)$ to the Lie algebra $(\Gamma(A),\mu)$.
For every Poisson structure $\pi$ on a Lie algebroid $A $,
the triple $(\Gamma(\wedge A)[2], \left[.,.\right]_{_{SN}},\left[\pi,.\right]_{_{SN}}) $ is a DGLA,
so that the pair $(l_1^{\pi,\left[.,.\right]},l_2^{\left[.,.\right]})$ given by
\begin{equation}\label{l_1andl_2}
\begin{array}{rcl}
l_1^{\pi,\left[.,.\right]}(P)=\left[\pi,P\right]_{_{SN}}&\mbox{and}&l_2^{\left[.,.\right]}(P,Q) := (-1)^{(p-1)} \left[P,Q\right]_{_{SN}},
\end{array}
 \end{equation}
where $P\in \Gamma(\wedge^{p} A), Q\in \Gamma(\wedge^{q} A)$ is an $L_{\infty}$-structure on the graded vector space $\Gamma(\wedge A)[2]$, which is clearly multiplicative.
We call this $L_\infty$-structure the \emph{$L_\infty$-structure associated to the Poisson
structure $\pi$ and the Lie algebroid $A$}.
We recall the notion of Poisson-Nijenhuis structure:
\begin{defi}
Let $(A,\mu=\left[.,.\right], \rho)$ be a Lie algebroid, $\pi\in \Gamma(\wedge^2A)$ be a bi-vector and $N:\Gamma(A)\to \Gamma(A)$ be a $(1-1)$-tensor field. Then the pair $(\pi, N)$ is called a Poisson-Nijenhuis structure on the Lie algebroid $(A,\mu=\left[.,.\right], \rho)$ if
\begin{enumerate}
\item $N$ is a Nijenhuis $(1-1)$-tensor with respect to the Lie bracket $\mu$,
\item $\pi$ is a Poisson bi-vector,
\item $N\pi^{\#}=\pi^{\#} N^*$,
\item $(\{\alpha,\beta\RP)_{_{N^*}}=\{\alpha, \beta\}_{_{\pi}}^{\mu^{N}}$,
\end{enumerate}
for all $\alpha,\beta\in \Gamma(A^*)$, where $(\{.,.\RP)_{_{N^*}}$ is the deformation of the Lie bracket $\{\cdot,\cdot\}_{_{\pi}}^{\mu}$ by $N^*$ and $\{.,.\}_{_{\pi}}^{\mu^{N}}$ is the induced bracket by the pair $(\pi, \mu^{N}=\left[.,.\right]_{N})$.
\end{defi}
\begin{rem}\label{drough}
 It implies directly, from definitions that
\begin{equation*}
\pi_{_{N}}^{\#}=N\pi^{\#}=\pi^{\#}N^*
\end{equation*}
and hence, $\underline{N}(\pi)=\iota_{_{N^*}}\pi=2\pi_{_{N}}.$
\end{rem}
\begin{rem}\label{papereivet}
 Recall from \cite{YKS} that if $(\pi, N)$ is a Poisson-Nijenhuis structure on a Lie algebroid $(A,\mu=\left[.,.\right], \rho)$, then, the triples $\left(A,\mu^{N}=\left[.,.\right]_{N}, \rho\circ N\right)$, $\left(A^*, \{.,.\}_{\pi}^{\mu}, \rho\circ\pi^{\#}\right)$,\\ $\left(A^*, \left(\{.,.\}_{\pi}^{\mu}\right)_{N^*},\rho\circ\\\pi^{\#}\circ N^*\right)$, $\left(A^*, \{.,.\}_{\pi}^{\mu^{N}}, \rho\circ N\circ\pi^{\#}\right)$,$ \left(A^*, \{.,.\}_{\pi_{_{N}}}^{\mu}, \rho\circ\pi_{_{N}}^{\#}\right)$ are Lie algebroids such that all the triples $\left(A^*, \left(\{.,.\}_{\pi}^{\mu}\right)_{N^*}, \rho\circ\pi^{\#}\circ N^*\right)$, $\left(A^*, \{.,.\}_{\pi}^{\mu^{N}}, \rho\circ N\circ\pi^{\#}\right)$, $\left(A^*, \{.,.\}_{\pi_{_{N}}}^{\mu}, \rho\circ\pi_{_{N}}^{\#}\right)$ are identically the same Lie algebroids. Moreover, identifying the graded vector spaces $\Gamma(\wedge A^{**})$ and $\Gamma(\wedge A)$, the de Rham differential  $ \diff^{A^{**}}_{(\{.,.\}_{\pi}^{\mu})}$ coincide with the linear map $\left[\pi,.\right]_{_{SN}}$. Hence, the relation
 $\diff^{A^{**}}_{(\{.,.\}_{\pi}^{\mu^{N}})}=\diff^{A^{**}}_{(\{.,.\}_{\pi_{_{N}}}^{\mu})}$, which itself is a consequence of discussion in above, implies that $\left[\pi,.\right]'_{_{SN}}-\left[\pi_{_{N}},.\right]_{_{SN}}=0$, where $\left[.,.\right]'_{_{SN}}$ is the Schouten-Nijenhuis bracket with respect to the Lie bracket $\left[.,.\right]_N$.
 \end{rem}
\begin{rem}
Notice that the pair $(l_1^{\left[.,.\right],\pi},l_2^{\left[.,.\right]})$ introduced in (\ref{l_1andl_2})
defines a Lie bialgebroid, since it is a multiplicative $L_\infty$-structure on the graded vector space $\Gamma(A)[2]$.
\end{rem}

%

\begin{lem}\label{lamejanjali}
Let $(\pi,N)$ be a Poisson-Nijenhuis structure on a Lie algebroid $(A,\left[.,.\right], \rho)$. Then,
\begin{equation*}
\left[\underline{N},l_1^{\left[.,.\right],\pi}\right]_{_{RN}}(P)= \left[\pi,{\underline{ N}} (P)\right]_{_{SN}}-\underline {N} \left[\pi, P \right]_{_{SN}}  = \left[-\pi_N , P\right]_{_{SN}},
\end{equation*}
for all $P \in \Gamma(\wedge A)$.
\end{lem}
\begin{proof}
The first equality follows directly from the definition of Richardson-Nijenhuis bracket and definition of $l_1^{\left[.,.\right],\pi}$. For the second equality, observe that for all $P \in \Gamma(\wedge A)$ we have
\begin{equation*}
\left[\pi,P\right]'_{_{SN}}=\left[\underline{N}(\pi),P\right]_{_{SN}}+\left[\pi,\underline{N}(P)\right]_{_{SN}}-\underline{N}\left[\pi,P\right]_{_{SN}},
\end{equation*}
for all $P\in \Gamma(\wedge^pA)$, where $\left[.,.\right]'_{_{SN}}$ stands for the Schouten-Nijenhuis bracket with respect to the Lie bracket $\left[.,.\right]_{N}$.
Hence, using Remarks \ref{drough}  and \ref{papereivet} we have
\begin{equation*}
\begin{array}{rcl}
\left[\pi,{\underline{ N}} (P)\right]_{_{SN}}-\underline {N} \left[\pi, P \right]_{_{SN}}&=&\left[\pi,P\right]'_{_{SN}}-\left[\underline{N}(\pi),P\right]_{_{SN}}\\
               &=&\left[\pi,P\right]'_{_{SN}}-2\left[\pi_{_{N}},P\right]_{_{SN}}\\
               &=&\left(\left[\pi,P\right]'_{_{SN}}-\left[\pi_{_{N}},P\right]_{_{SN}}\right)-\left[\pi_{_{N}},P\right]_{_{SN}}\\
               &=&-\left[\pi_{_{N}},P\right]_{_{SN}}.
\end{array}
\end{equation*}
\end{proof}

\begin{prop}\label{NijenhuisonPoissonNijenhuisalgeroid}
Let $(N,\pi)$ be a Poisson-Nijenhuis structure on a Lie algebroid $(A,\left[.,.\right],\rho)$,
then the derivation $\underline N $ is a weak-Nijenhuis tensor
for the $L_\infty$-structure associated to the Poisson
structure $\pi$.

In this case, the deformed structure $[ \underline N,l_1^{\left[.,.\right],\pi} + l_2^{\left[.,.\right]}]_{_{RN}}$ is the $L_\infty$-structure associated to
the Poisson structure $-\pi_N $ on the Lie algebroid $(A,\left[.,.\right]_N,\rho^N)$.
\end{prop}
\begin{proof}
Lemmas \ref{lamejanjali} and \ref{lem:deformingOidsN} imply that
\begin{equation*}
\left[\underline{N},l_1^{\left[.,.\right],\pi}+l_2^{\left[.,.\right]}\right]_{_{RN}}=-l_1^{\left[.,.\right],\pi_{_{N}}}+l_2^{\left[.,.\right]_{_{N}}}.
\end{equation*}
Hence,
\begin{equation}\label{avalinesh}
\begin{array}{rcl}
\left[\underline{N},\left[\underline{N},l_1^{\left[.,.\right],\pi}+l_2^{\left[.,.\right]}\right]_{_{RN}}\right]_{{RN}}
&=&l_1^{\left[.,.\right],\pi_{_{N,N}}}+l_2^{\left[.,.\right]_{_{N,N}}}\\
&=&l_1^{\left[.,.\right],\pi_{_{N^2}}}+l_2^{\left[.,.\right]_{_{N^2}}}\\
&=&\left[\underline{N^2},-l_1^{\left[.,.\right],\pi}+l_2^{\left[.,.\right]}\right]_{_{RN}}\\
&=&\left[\underline{N^2},l_1^{\left[.,.\right],\pi}+l_2^{\left[.,.\right]}\right]_{_{RN}}-
2\left[\underline{N^2},l_1^{\left[.,.\right],\pi}\right]_{_{RN}}.
\end{array}
\end{equation}
Denoting $\mu=l_1^{\left[.,.\right],\pi}+l_2^{\left[.,.\right]}$ and using the fact that $\pi_{_{N^2}}$ is a Poisson bi-vector and hence $(\Gamma(\wedge A)[2],l_1^{\left[.,.\right],\pi_{_{N^2}}}+l_2^{\left[.,.\right]})$ is a symmetric DGLA, we have
\begin{equation}\label{dovominesh}
\begin{array}{rcl}
\left[\mu,\left[\underline{N},\left[\underline{N},\mu\right]_{_{RN}}\right]_{_{RN}}\right]_{_{RN}}
&=&\left[\mu,\left[\underline{N^2},\mu\right]_{_{RN}}\right]_{_{RN}}
-2\left[\mu,\left[\underline{N^2},l_1^{\left[.,.\right],\pi}\right]_{_{RN}}\right]_{_{RN}}\\
&=&-2\left[\mu,\left[\underline{N^2},l_1^{\left[.,.\right],\pi}\right]_{_{RN}}\right]_{_{RN}}\\
&=&2\left[\mu,l_1^{\left[.,.\right],\pi_{_{N^2}}}\right]_{_{RN}}\\
&=&2\left[l_1^{\left[.,.\right],\pi},l_1^{\left[.,.\right],\pi_{_{N^2}}}\right]_{_{RN}}+
2\left[l_2^{\left[.,.\right]},l_1^{\left[.,.\right],\pi_{_{N^2}}}\right]_{_{RN}}\\
&=&2\left[l_1^{\left[.,.\right],\pi},l_1^{\left[.,.\right],\pi_{_{N^2}}}\right]_{_{RN}}.
\end{array}
\end{equation}
But
\begin{equation}\label{sevominesh}
\begin{array}{rcl}
\left[l_1^{\left[.,.\right],\pi},l_1^{\left[.,.\right],\pi_{_{N^2}}}\right]_{_{RN}}(P)&=&
l_1^{\left[.,.\right],\pi_{_{N^2}}}(l_1^{\left[.,.\right],\pi}(P))+
l_1^{\left[.,.\right],\pi}(l_1^{\left[.,.\right],\pi_{_{N^2}}}(P))\\
&=&\left[\pi_{_{N^2}},\left[\pi,P\right]_{_{SN}}\right]_{_{SN}}+\left[\pi,\left[\pi_{_{N^2}},P\right]_{_{SN}}\right]_{_{SN}}\\
&=&\left[\left[\pi,\pi_{_{N^2}}\right]_{_{SN}},P\right]_{_{SN}}\\
&=&0.
\end{array}
\end{equation}
Therefore, $\left[\mu,\left[\underline{N},\left[\underline{N},\mu\right]_{_{RN}}\right]_{_{RN}}\right]_{_{RN}}=0$ which means that $\underline{N}$ is weak Nijenhuis vector valued form with respect to the symmetric DGLA structure $\mu=l_1^{\left[.,.\right],\pi}+l_2^{\left[.,.\right]}$ on the graded vector space
$\Gamma(\wedge A)[2]$.

%
%
%
%
\end{proof}
There is a second manner to see Poisson-Nijenhuis structures on a Lie algebroid as a Nijenhuis
form.
\begin{prop}
Let $(\pi,N)$ be a Poisson-Nijenhuis structure on a Lie algebroid $(A,\left[.,.\right],\rho)$. Then
 $\underline N + \pi $ is a weak Nijenhuis vector valued form with curvature, with respect to the multiplicative DGLA-structure $l_1^{\left[.,.\right],\pi}+l_2^{\left[.,.\right]}$ on the graded vector space $\Gamma(\wedge A)[2]$,
with square ${\underline{N^2}} $.
\end{prop}
\begin{proof}
It follows from Lemma \ref{lem:deformingOidsN} that
\begin{equation*} \label{ELONE}
\left[\underline{N}+\pi,l_1^{\left[.,.\right],\pi}\right]_{_{RN}}=
-l_1^{\left[.,.\right],\pi_{_{N}}}+\left[\pi,\pi\right]_{_{SN}}=-l_1^{\left[.,.\right],\pi_{_{N}}}.
\end{equation*}
Lemma \ref{lamejanjali} imply that
\begin{equation*} \label{ELTWO}
\left[\underline{N}+\pi,l_2^{\left[.,.\right]}\right]_{_{RN}}=
l_2^{\left[.,.\right]_N}+l_2^{\left[.,.\right]}(\pi,.)=l_2^{\left[.,.\right]_N}-l_1^{\left[.,.\right],\pi}.
\end{equation*}
Hence
\begin{equation}\label{divane}
\begin{array}{rcl}
&&\left[\underline{N}+\pi,\left[\underline{N}+\pi,l_1^{\left[.,.\right],\pi}+l_2^{\left[.,.\right]}\right]_{_{RN}}\right]_{_{RN}}\\
&=&\left[\underline{N}+\pi,-l_1^{\left[.,.\right],\pi_{_{N}}}+l_2^{\left[.,.\right]_N}-l_1^{\left[.,.\right],\pi}\right]_{_{RN}}\\
&=&l_1^{\left[.,.\right],\pi_{_{N,N}}}+l_1^{\left[.,.\right],\pi_{_{N}}}+l_2^{\left[.,.\right]_{N,N}}-l_1^{\left[.,.\right],\pi_{_{N}}}(\pi)
-l_1^{\left[.,.\right],\pi}(\pi)+l_2^{\left[.,.\right]_N}(\pi,.).\\
\end{array}
\end{equation}
But $l_1^{\left[.,.\right],\pi}(\pi)=\left[\pi,\pi\right]_{_{SN}}=0$, $l_1^{\left[.,.\right],\pi_{_{N}}}(\pi)=\left[\pi_{_{N}},\pi\right]_{_{SN}}=0$ and
\begin{equation*}
l_1^{\left[.,.\right],\pi_{_{N}}}(P)+l_2^{\left[.,.\right]_N}(\pi,P)=\left[\pi_{_{N}},P\right]_{_{SN}}-\left[\pi,P\right]'_{_{SN}}=-\mbox{concomitant}=0,
\end{equation*}
for all $P\in \Gamma(\wedge A)[2]$, where $\left[.,.\right]'_{_{SN}}$ is the Schouten-Nijenhuis bracket associated to the Lie bracket $\left[.,.\right]_N$.
Hence, \ref{divane} can be rewritten as
\begin{equation}
\begin{array}{rcl}
&&\left[\underline{N}+\pi,\left[\underline{N}+\pi,l_1^{\left[.,.\right],\pi}+l_2^{\left[.,.\right]}\right]_{_{RN}}\right]_{_{RN}}\\
&=&\left[\underline{N}+\pi,-l_1^{\left[.,.\right],\pi_{_{N}}}+l_2^{\left[.,.\right]_N}-l_1^{\left[.,.\right],\pi}\right]_{_{RN}}\\
&=&l_1^{\left[.,.\right],\pi_{_{N,N}}}+l_2^{\left[.,.\right]_{N,N}}.\\
\end{array}
\end{equation}
Similar computations as in (\ref{avalinesh}), (\ref{dovominesh}) and (\ref{sevominesh}) show that $\left[\mu,\left[\underline{N},\left[\underline{N},\mu\right]_{_{RN}}\right]_{_{RN}}\right]_{_{RN}}=0$ which means that $\underline{N}$ is weak Nijenhuis vector valued form with respect to the symmetric DGLA structure $\mu=l_1^{\left[.,.\right],\pi}+l_2^{\left[.,.\right]}$ on the graded vector space
$\Gamma(\wedge A)[2]$.
\end{proof}
We have already defined weak Nijenhuis structures.
For the purpose of these last lines, we shall introduce a notion that
is stronger than weak Nijenhuis but weaker than Nijenhuis admitting a square:
\begin{defi}
Let $E$ be a graded vector space and $\mu$ be a symmetric vector valued form on $E$ of degree $1$. A vector valued form ${\mathcal N}$ of degree zero is called
co-boundary Nijenhuis with respect to $\mu$ if there exists a vector valued form ${\mathcal K}$ of degree $0$ such that
       \begin{equation}\label{coboundN}
       \left[{\mathcal N},\left[{\mathcal N},\mu\right]_{_{RN}}\right]_{_{RN}}=\left[{\mathcal K},\mu\right]_{_{RN}},
       \end{equation}
  Such a ${\mathcal K} $ is called a square of ${\mathcal N} $.
  If $\mathcal{N}$ contains an element of the underlying graded
vector space, that is, $\mathcal{N}$ has a component which is a vector valued zero form, then $\mathcal{N}$
is called Nijenhuis vector valued form with curvature.
\end{defi}
Of course, if ${\mathcal K} $ commutes with $ {\mathcal N}$, this definition gives back the definition of Nijenhuis with square ${\mathcal K} $.
\begin{prop}\label{Poissonifandonlyif}
Let $(A,\left[.,.\right],\rho)$ be a Lie algebroid , $\pi\in \Gamma(\wedge^2 A)$ be a bi-vector and $N:\Gamma(A)\to \Gamma(A)$ be a $(1-1)$-tensor field such that
\begin{equation}\label{costezan}
N\pi^{\#}=\pi^{\#}N^*.
\end{equation}
Then
 $\underline{N} + \pi $ is a co-boundary Nijenhuis vector valued form with curvature, with respect to the multiplicative GLA-structure $l_2^{\left[.,.\right]}$ on the graded vector space $\Gamma(\wedge A)[2]$,
with square ${\underline{N^2}} $, if and only if $(N,\pi)$ is a Poisson-Nijenhuis structure on the Lie algebroid $(A,\left[.,.\right],\rho)$.
The deformed structure $[ \underline N, l_2^{\left[.,.\right]}]_{_{RN}}$ is the $L_\infty$-structure (indeed a DGLA) associated to
the Poisson structure $\pi $ on the Lie algebroid $(A,\left[.,.\right]_N,\rho \circ N) $.
\end{prop}
\begin{proof}
Assume that $(N,\pi)$ is a Poisson-Nijenhuis structure on the Lie algebroid $(A,\left[.,.\right],\rho)$. Then
\begin{equation*}
\left[\underline{N}+\pi,l_2^{\left[.,.\right]}\right]_{_{RN}}=l_2^{\left[.,.\right]_{N}}-l_1^{\left[.,.\right],\pi}
\end{equation*}
hence, by Remark \ref{papereivet} we get
\begin{equation*}
\begin{array}{rcl}
\left[\underline{N}+\pi\left[\underline{N}+\pi,l_2^{\left[.,.\right]}\right]_{_{RN}}\right]_{_{RN}}
&=&l_2^{\left[.,.\right]_{N,N}}+l_1^{\left[.,.\right],\pi_{_{N}}}-l_1^{\left[.,.\right]_N,\pi}\\
&=&l_2^{\left[.,.\right]_{N,N}}\\
&=&\left[\underline{N^2},l_2^{\left[.,.\right]}\right]_{_{RN}},
\end{array}
\end{equation*}
which means that $\underline{N}+\pi$ is a co-boundary Nijenhuis with respect to the multiplicative GLA-structure $l_2^{\left[.,.\right]}$ on the graded vector space
$\Gamma(\wedge A)[2]$, with square ${\underline{N^2}} $.

Conversely, let $\underline{N}+\pi$ be a co-boundary Nijenhuis with respect to the multiplicative GLA-structure $l_2^{\left[.,.\right]}$ on the graded vector space $\Gamma(\wedge A)[2]$, with square ${\underline{N^2}} $. Then
\begin{equation*}
\begin{array}{rcl}
\left[\underline{N}+\pi,\left[\underline{N}+\pi,l_2^{\left[.,.\right]}\right]_{_{RN}}\right]_{_{RN}}&=&
l_2^{\left[.,.\right]_{N,N}}+
\left(\left[\underline{N},l_2^{\left[.,.\right]}\right]_{_{RN}}+l_2^{\left[.,.\right]_{N}}(\pi,.)\right)+l_2^{\left[.,.\right]}(\pi,\pi)\\
&=&\left[\underline{N^2},l_2^{\left[.,.\right]}\right]_{_{RN}}\\
&=&l_2^{\left[.,.\right]_{N,N}}
\end{array}
\end{equation*}
implies that
\begin{equation*}
l_2^{\left[.,.\right]_{N,N}}=l_2^{\left[.,.\right]_{N^2}}
\end{equation*}
or
\begin{equation}\label{Nijenhuisshode}
\left[.,.\right]_{N,N}=\left[.,.\right]_{N^2},
\end{equation}
\begin{equation}\label{Poissonshode}
l_2^{\left[.,.\right]_{N}}(\pi,\pi)=0
\end{equation}
and
\begin{equation}\label{PoissonshodeNijenhuisham}
\left(\left[\underline{N},l_2^{\left[.,.\right]}\right]_{_{RN}}+l_2^{\left[.,.\right]_{N}}(\pi,.)\right)(P)=0,
\end{equation}
for all $P\in \Gamma(\wedge A)$.
Equation (\ref{Poissonshode}) means that $\pi$ is a Poisson element, while Equation (\ref{PoissonshodeNijenhuisham}) can be rewritten as
\begin{equation*}
-\left[\pi,\underline{N}(P)\right]_{_{SN}}+\underline{N}\left[\pi,P\right]_{_{SN}}-\left[\pi,P\right]'_{_{SN}}=0
\end{equation*}
or using (\ref{costezan})
\begin{equation}\label{javananchera}
\underline{N}\left[\pi_{_{N}},P\right]_{_{SN}}-\left[\pi,P\right]'_{_{SN}}=0,
\end{equation}
where $\left[.,.\right]'_{_{SN}}$ is the Schouten-Nijenhuis bracket with respect to the deformed bracket $\left[.,.\right]_{N}$.
Now, Equations (\ref{Poissonshode}), (\ref{javananchera}), (\ref{Nijenhuisshode}) and (\ref{costezan}) imply that $(N,\pi)$ is a Poisson-Nijenhuis structure on the Lie algebroid $(A,\left[.,.\right], \rho)$.
\end{proof}
Last, we shall say a few words about the so-called $\Pi\Omega$-structures.
Recall that a \emph{$\Pi\Omega$-structure on a Lie algebroid $(A,\rho,\left[.,.\right])$}
is a pair $(\pi,\omega)$ where $\pi \in \Gamma(\wedge^2 A)$ is a Poisson element
and $\omega \in \Gamma(\wedge^2 A^*)$ is a $2$-form, with $\diff \alpha =0$.
Defining a $1-1$ tensor $N := \pi^{\#} \circ \omega^b$, it is known that $(\pi,N) $ is always a Poisson-Nijenhuis structure
while $(N,\omega) $ is an $\Omega N$-structure.

\begin{prop}\label{NijenhuisonPoissonOmega}
Let $(\pi,\omega)$ be a  $\Pi\Omega$-structure on a Lie algebroid $(A,\left[.,.\right],\rho)$. Then,
 ${\mathcal N}=\underline{\omega} + \pi $ is a co-boundary Nijenhuis form, with curvature, with respect to the multiplicative GLA-structure $l_2^{\left[.,.\right]}$ on the graded vector space $\Gamma(\wedge A)[2]$,
with square $\underline{N} $, where $N = \pi^{\#} \circ \omega^b$.
The deformed structure is $-l_1^{\left[.,.\right],\pi}$.
\end{prop}
\begin{proof}
First, observe that
\begin{equation*}
l_1^{\left[.,.\right],\pi}(P)=\left[\pi,P\right]_{_{SN}}=-l_2^{\left[.,.\right]}(\pi,P)=-\left[\pi,l_2^{\left[.,.\right]}\right]_{_{RN}}(P)
\end{equation*}
for all $P\in \Gamma(\wedge^2 A)$  which means that
\begin{equation}\label{chetori}
l_1^{\left[.,.\right],\pi}=-\left[\pi,l_2^{\left[.,.\right]}\right]_{_{RN}}.
\end{equation}
Hence
\begin{equation} \label{samandoon}
\left[{\mathcal N},l_2^{\left[.,.\right]}\right]_{_{RN}} = -l_1^{\left[.,.\right],\pi}+ \underline{ \diff \omega} = -l_1^{\left[.,.\right],\pi},
\end{equation}
which proves the last the last claim (and proves that ${\mathcal N}$ is weak-Nijenhuis vector valued form with respect to $l_2^{\left[.,.\right]}$, since
 $ l_1^{\left[.,.\right],\pi}$ is an $L_{\infty}$-structure on $\Gamma(\wedge A)[2]$).
Now (\ref{samandoon}) and (\ref{chetori}) imply that
\begin{equation*}
\begin{array}{rcl}
\left[\mathcal{N},\left[\mathcal{N},l_2^{\left[.,.\right]}\right]_{_{RN}}\right]_{_{RN}}
&=&-\left[\mathcal{N},l_1^{\left[.,.\right],\pi}\right]_{_{RN}}\\
&=&-\left[\underline{\omega},l_1^{\left[.,.\right],\pi}\right]_{_{RN}}-\left[\pi,\pi\right]_{_{SN}}\\
&=&\left[\underline{\omega},\left[\pi,l_2^{\left[.,.\right]}\right]_{_{RN}}\right]_{_{RN}}\\
&=&\left[\left[\underline{\omega},\pi\right]_{_{RN}},l_2^{\left[.,.\right]}\right]_{_{RN}}.
\end{array}
\end{equation*}
This shows that $\mathcal{N}$ is a co-boundary Nijenhuis vector valued form with respect to the GLA-structure $l_2^{\left[.,.\right]}$,
 on the graded vector space $\Gamma(\wedge A)[2]$, with square $\left[\underline{\omega},\pi\right]_{_{RN}}.$
Now,  a direct computation shows that  $[\pi,\underline{\omega}]_{_{RN}}= \underline{N}$.
\end{proof}
%




\begin{thebibliography}{00}

\bibitem{antunes} P. Antunes, Crochets de Poisson gradu\'es et applications: structures compatibles et g\'en\'eralisations des structures hyperk\"ahl\'eriennes, Th\`ese de doctorat de l'\'Ecole Polytechnique, March 2010.
\bibitem{CJP} P. Antunes, C. Laurent-Gengoux, J. M. Nunes da Costa, Hierarchies and compatibility on Courant algebroids, \emph{Pac. J. Math.} {\bf 261}(1)(2013) 1-32.
\bibitem{AntunesCosta} P. Antunes, J. M. Nunes da Costa, Nijenhuis and compatible tensors on Lie and Courant algebroids, \emph{J. Geom. Phys.} {\bf 65}(2013) 66-79.
\bibitem{AJC} P. Aschieri, L.  Cantini, B. Jurco, Nonabelian Bundle Gerbes, their Differential Geometry and Gauge Theory,
        \emph{Comm. Math. Phys.} {\bf 254}(2)(2005) 367-400.
\bibitem{Said} I. Ayadi, S. Benayadi, Symmetric Novikov superalgebra,  \emph{J. Math. Phys.} {\bf 51}(2010) 023501, 15 pages.
\bibitem{BaezSchreiber} J. Baez, U. Schreiber, Higher gauge theory.  \emph{Contemp. Math.} {\bf 431}(2007) 7-30.
\bibitem{Bartels} T. Bartels, $2$-Bundles and Higher Gauge Theory, Ph.D. thesis, University of California, Riverside, 2004. arXiv:math/0410328
\bibitem{BehrendXu} K. Behrend, P. Xu,  Differentiable stacks and gerbes,
         \emph{J. Symplectic Geom.} {\bf  9}(2011) 285-341.
\bibitem{BergerGostiaux} M. Berger, B. Gostiaux, G\'eom\'etrie diff\'erentielle : vari\'et\'es, courbes et surfaces,
        \emph{Presses Universitaires de France} (1987).
\bibitem{BL} L. Breen, C. Laurent-Gengoux,  Nonabelian differential gerbes on local coordinates, private communication.
\bibitem{BM} L. Breen, W. Messing, Differential geometry of gerbes, \emph{Adv. Math.} {\bf 198}(2005) 732-846.
\bibitem{Brylinski} J. L. Brylinsky, Loop spaces, characteristic classes and geometric quantization, Progress in Mathematics {\bf 107} Birkh\"auser (1993).
\bibitem{Carchedi} D. Carchedi, Sheaf Theory for \'{E}tale Geometric Stacks , arXiv:1011.6070
\bibitem{CGM2} J.F Cari\~{n}ena, J. Grabowski, G. Marmo, Contraction: Nijenhuis and Saletan tensors for general algebraic structures, \emph{J. Phys. A} {\bf 34}(2001)(18)  3769-3789.
\bibitem{CGM} J.F Cari\~{n}ena, J. Grabowski, G. Marmo, Courant algebroid and Lie bialgebroid contractions, \emph{J. Phys. A} {\bf 37}(19)(2004) 5189-5202.
\bibitem{Cattaneo-Felder} A.S. Cataneo, G. Felder, Relative formality theorem and quantization of coisotropic submanifolds, \emph{Adv. Math.} {\bf 208}(2007) 521-548
\bibitem{Cattaneo-Schaaetz} A.S. Cataneo, F. Sch\"atz, Equivalences of higher derived brackets, \emph{J. Pure and Appl. Algebra} {\bf 212}(2008) 2450-2460.
\bibitem{Chatterjee} D. Chatterjee, On the construction of abelian gerbes, Ph.D. thesis (Cambridge) 1998.
\bibitem{Clemente-Nunes} J. Clemente-Gallardo, J. M. Nunes da Costa, Dirac-Nijenhuis structures, \emph{J. Phys. A: Math. Gen.} {\bf 37}(2004) 7267-7296
\bibitem{Courant} T. J. Courant, Dirac manifolds, \emph{Trans. Amer. Math. Soc.} {\bf 319}(1990) 631-661.
\bibitem{CrainicMoerdijk}  M. Crainic, I.  Moerdijk,  Foliation Groupoids and Their Cyclic Homology,
        \emph{Adv. Math.} {\bf 157}(2000) 177-197.
\bibitem{Dedecker} P.  Dedecker,  Sur la cohomologie non Ab\'elienne. I. (French) \emph{Canad. J. Math.}   {\bf 12}(1960) 231-251.
\bibitem{FRZ} Y. Fr\'{e}gier, C. L. Roger, M. Zambon, Homotopy moment maps. arXiv:1304.2051 [math.DG]
\bibitem{GR} K. Gawedzki, N. Reis, WZW branes and gerbes, \emph{Rev. Math. Phys.} {\bf 14}(2002) 1281-1334.
\bibitem{Getzler} E. Getzler, Higher derived bracket. arXiv:1010.5859v2
\bibitem{GinotStienon} G. Ginot, M. Sti\'enon, Groupoid extensions, principal $2$-group bundles and characteristic classes,
        arXiv:0801.1238.
\bibitem{Giraud} J. Giraud, Cohomologie non ab\'elienne, Die Grundlehren der mathematischen Wissenschaften {\bf 179} Springer-Verlag 1971.
\bibitem{GrB} J. Grabowski, Courant-Nijenhuis tensors and generalized geometries, In \emph{groups, geometry and Physics}, Mangr. Real. Acad. Ci. Exact. Fis.-Quim. Nat. Zaragoza {\bf 29}, Acad. Cienc. Exact. Fis. Quim. Nat. Zaragoza, Zaragoza (2006) 101-112.
\bibitem{Gra} J. Grabowski, G. Marmo, The graded Jacobi algebras and (co)homology, \emph{J. Phys A} {\bf 36}(2003) 161-181.
\bibitem{GrabowskiUrbanski} J. Grabowski, P. Urba\'nski, Lie algebroids and Poisson-Nijenhuis structures, \emph{Rep. Math. Phys.} {\bf 40}(1997)(2) 195-208.
\bibitem{Grothendieck1} A. Grothendieck, General Theory of Fiber Spaces, Report 4, University of Kansas, Lawrence, Kansas 1955.
\bibitem{Grothendieck2} A. Grothendieck, Le groupe de Brauer, Sem. Bourbaki {\bf 290} 1964/1965.
\bibitem{Hitchin} N. Hitchin, Lectures on special Lagrangian submanifolds, Winter School on Mirror Symmetry, Vector bundles and Lagrangian submanifolds (Cambridge, MA 1999)\emph{ AMS/IP Stud. Adv. Math.} {\bf 23} 151-182.
\bibitem{KaSt} H. Kajiwara, J. Stasheff, Homotopy algebras inspired by classical open-closed string field theory,\emph{Comm. Math. Phys.} {\bf 263}(2006) 553-581.
\bibitem{KMS} I. Kolar, P. W. Michor, J. Slovak, Natural Operations in Differential Geometry. Springer-Verlag 1993.
\bibitem{Kont} M. Kontsevich, Deformation quantization of Poisson manifolds, \emph{Lett. Math. Phys.}, {\bf 66}(2003), 157-216.
\bibitem{YKSbialgebroid} Y. Kosmann-Schwarzbach, The Lie bialgebroid of a Poisson-Nijenhuis manifold, \emph{Lett. Math. Phys.} {\bf 38}(1996) 421-428.
\bibitem{Kosmann} Y. Kosmann-Schwarzbach, Modular vector fields and Batalin-Vilkovisky algebras, in Poisson Geomtry, J. Grabowski and P. Urbanski, eds.,\emph{ Banach Center Publ.} {\bf 51}(2000) 109-129.
\bibitem{YKSquasi} Y. Kosmann-Schwarzbach, Quasi twisted and all that $\cdots$ in Poisson geometry and Lie algebroid theory, pp. 363-389 in \emph{The breath of symplectic and Poisson geometry}, J. Marsden and T. Ratiu eds., \emph{Prog. Math.} {\bf 232}, Birkha\"user, Boston, MA, 2005.
\bibitem{YKSBrazil} Y. Kosmann-Schwarzbach, Nijenhuis structures on Courant algebroids, \emph{Bull. Braz. Math. Soc. (N.S.)} {\bf 42}(2011)(4) 625-649.
\bibitem{YKS}Y. Kosmann-Schwarzabch, F. Magri, Poisson-Nijenhuis structures, \emph{Ann. Inst. Henri Poincar\'{e}, S\'{e}rie A,} {\bf 53}(1990), 35-81.
\bibitem{KoRou} Y. Kosmann-Schwarzbach, V. Rubtsov, Compatible structures on Lie algebroids and Monge-Ampere operators, \emph{Acta Appl. Math.} {\bf 109}(1)(2010) 101-135.
\bibitem{Lada-Sta} T. Lada, J. Stasheff, Introduction to SH Lie Algebras for Physisists, \emph{Int. J. Theor. Phys.} {\bf 32}(7)(1993) 1087-1103.
\bibitem{LaurentHabilitation} C. Laurent-Gengoux, Des groupo\"ides de Lie \`a la g\'eom\'etrie de Poisson, Habilitation Dissertation, Universit\'e de Poitiers (2009).
\bibitem{LaurentStienonXu} C. Laurent-Gengoux, M. Sti\'enon, P. Xu,  Non-abelian differentiable gerbes, \emph{ Adv. Math.}  {\bf 220}(2009) 1357-1427.
\bibitem{LTX} C. Laurent-Gengoux, J.L. Tu, P. Xu, Chern-Weil map for principal bundles over groupoids, \emph{Math. Z} {\bf 255}(2007) 451-491.
\bibitem{LMS} P.A.B. Lecomte, P.W. Michor, H.Schicketanz, The multigraded Nijenhuis-Richardson algebra, its universal property and applications,\emph{J. Pure Appl. Algebra} {\bf 77}(1992)(1) 87-102.
\bibitem{Lichnerowicz} A. Lichnerowicz, Les vari\'{e}t\'{e}s de Poisson et leurs alg\'ebres de Lie associ\'{e}es,\emph{ J. Diff. Geom.}
          {\bf 12}(1977) 253-300.
\bibitem{LSZ} Z. J. Liu, Y. Sheng, T. Zhang, \emph{Deformations of Lie $2$-algebras}, arXiv: 1306.6225.
\bibitem{LWX} Z. J. Liu, A.  Weinstein, P. Xu, Manin triples for Lie for Lie bialgebroids, \emph{J. Diff. Geom.} {\bf 45}(1997) 547-574.
\bibitem{McK}  K. Mackenzie, General theory of Lie groupoids and Lie algebroids,
           \emph{ London Mathematical Society Lecture Note series}, {\bf 213}, \emph{Cambridge University Press, Cambridge} (2005).
\bibitem{McX} K. Mackenzie P. Xu, Lie bialgebroids and Poisson groupoids, \emph{Duke Math. J.} {\bf 73}(1994) 415-452.
\bibitem{Magri} F. Magri, A simple model of the integrable Hamiltonian equation, \emph{J. Math. Phys.} {\bf 19}(1978) 1156-1162.
\bibitem{MM} F. Magri, C. Morosi, A geometrical characterization of integrable Hamiltonian systems through the theory of Poisson-Nijenhuis manifolds, Quaderno {\bf S 19}(1984) Univ. of Milan.
\bibitem{MMR} F. Magri, C. Morosi, O. Ragnisco, Reduction techniques for infinite dimensional Hamiltonian systems: some ideas and applications, \emph{Comm. Math. Phys} {\bf 99}(1)(1985) 115-140.
\bibitem{Murray} M. K Murray, D. Stevenson, Bundle gerbes: stable isomorphisms and local theory, \emph{J. London Math. Soc.}  {\bf 62}(2)(2000) 925-937.
\bibitem{NewlanderNirenberg} A. Newlander, L.Nirenberg, Complex analytic coordinates in almost complex manifolds,\emph{ Ann. of Math.} Second serries. {\bf 65}(1957) 391-404.
\bibitem{Nijenhuisformingsets} A. Nijenhuis, $X_{n-1}$ forming sets of eigenvectors, \emph{Indag. Math} {\bf 13}(1951) 202-212.
\bibitem{Nijenhuis-Jacobi-type} A. Nijenhuis, Jacobi-type identities for bilinear differential concommitants of certain tensor field I., \emph{Indag. Math.} {\bf17} (1955) 390-403.
\bibitem{NR} A. Nijenhuis, R. Richardson, Deformation of Lie algebra structures, \emph{J. Math. Mech.} {\bf 17}(1967) 89-105.
\bibitem{NikolausWaldorf} T. Nikolaus, K. Waldorf, Four Equivalent Versions of Non-Abelian Gerbes, arXiv:1103.4815v2
\bibitem{CRoger} C. L. Rogers, Higher symplectic geometry, Ph.D. thesis (2011) arXiv:1106.4068.
\bibitem{CRogerspaper} C. L. Rogers, $L_{\infty}$-algebras from multisymplectic geometry, \emph{Lett. Math. Phys.} {\bf 100}(2012) 29-50.
\bibitem{RoytenbergPh.D.} D. Roytenberg, Courant algebroids, derived brakets and even symplectic supermanifolds, Ph.D. thesis, UC Berkeley, 1999.
\bibitem{Roytenberg} D. Roytenberg, Quasi-Lie bialgebroids and twisted Poisson manifolds, \emph{Lett. Math. Phys.} {\bf 61}(2002) 123-137.
\bibitem{RoytenbergWeinstein} D. Roytenberg, A. Weinstein, Courant algebroids and strongly homotopy algebras, \emph{Lett. Math. Phys.} {\bf 46}(1998) 81-93.
\bibitem{SchatzThesis} F. Sch\"atz, Coisotropic submanifolds and the BFV-complex, Ph.D. Thesis, University of Z\"urich, 2009.
\bibitem{Scouten-differential-operators} J.A. Schouten, On the differential operators of the first order in tensor calculus, in: \emph{Convegno Int. Geom. Diff.} (Italia 1953) Ed. Cremonese (Roma, 1954)  pp.1-7.
\bibitem{SchreiberWaldorf} U. Schreiber, K. Waldorf, Connections on non-abelian Gerbes and their Holonomy, arXiv:0808.1923.
\bibitem{Uchino} K. Uchino, Remarks on the definition of a Courant Algebroid, \emph{Lett. Math. Phys.} {\bf 60}(2002) 171-175.
\bibitem{Vaintrob} A. Vaintrob, Lie algebroids and homological vector fields, \emph{Uspekhi. Mat. Nauk.}, {\bf 52}(2)(1997) 161-162, translated in \emph{Russian math. Surv.} {\bf 52}(2)(1997) 428-429.
\bibitem{Voronov} Th. Voronov, Higher derived brackets and homotopy algebras, \emph{J. Pure and Appl. Algebra} {\bf 202}(2005)(1-3) 133-153.
\bibitem{FWagemann} F. Wagemann, On Lie algebra crossed modules. \emph{comm. algebra} {\bf 34} (2006) 1699-1722.
\bibitem{AschieriCantiniJurco} P. Aschieri, L. Cantini, B. Jur\v{c}o, Non abellian bundle gerbes, their differential geometry and gauge theory. \emph{Comm. Math. Phys.}, {\bf 254}(2005) 367-400.
\bibitem{Waldorf} K. Waldorf, Multiplicative Bundle Gerbes with Connection, \emph{Diff. Geom. Appl.} {\bf 28}(2010) 313-340.
\bibitem{Witten} E. Witten, D-branes and $K$-theory, \emph{J. High Enery Phys.} {\bf 12}(1998) 19-44.
\bibitem{Zwiebach} B. Zwiebach, Closed string field theory: quantum action and the Batalin-Vilkovisky master equation, \emph{Nuclear Phys. B} {\bf390}(1)(1993) 33-152.
\end{thebibliography}
\end{document}